\documentclass[12pt,a4paper]{amsart}
\makeatletter
\renewcommand\normalsize{%
    \@setfontsize\normalsize{11.7}{14pt plus .3pt minus .3pt}%
    \abovedisplayskip 10\p@ \@plus4\p@ \@minus4\p@
    \abovedisplayshortskip 6\p@ \@plus2\p@
    \belowdisplayshortskip 6\p@ \@plus2\p@
    \belowdisplayskip \abovedisplayskip}
\renewcommand\small{%
    \@setfontsize\small{9.5}{12\p@ plus .2\p@ minus .2\p@}%
    \abovedisplayskip 8.5\p@ \@plus4\p@ \@minus1\p@
    \belowdisplayskip \abovedisplayskip
    \abovedisplayshortskip \abovedisplayskip
    \belowdisplayshortskip \abovedisplayskip}
\renewcommand\footnotesize{
    \@setfontsize\footnotesize{8.5}{9.25\p@ plus .1pt minus .1pt}
    \abovedisplayskip 6\p@ \@plus4\p@ \@minus1\p@
    \belowdisplayskip \abovedisplayskip
    \abovedisplayshortskip \abovedisplayskip
    \belowdisplayshortskip \abovedisplayskip}
\ifdefined\pdfpagewidth
\else
\fi
\calclayout
\makeatother

\usepackage{amsmath, amssymb, amsthm, fancyhdr, verbatim, graphicx}
\usepackage{enumerate}
\usepackage[all]{xy}
\usepackage[usenames,dvipsnames]{xcolor}
\usepackage{mathrsfs}
\usepackage{tikz-cd}
\usetikzlibrary{decorations.pathmorphing}
\usepackage{framed, hyperref}
\usepackage[titletoc]{appendix}
\usepackage{bbm}
\usepackage{lipsum}
\usepackage{adjustbox}
\usepackage{enumitem}

\numberwithin{equation}{section}

\newcommand{\LL}{\mathfrak{L}}

\DeclareFontFamily{U}{matha}{\hyphenchar\font45}
\DeclareFontShape{U}{matha}{m}{n}{
      <5> <6> <7> <8> <9> <10> gen * matha
      <10.95> matha10 <12> <14.4> <17.28> <20.74> <24.88> matha12
      }{}
\DeclareSymbolFont{matha}{U}{matha}{m}{n}
\DeclareFontFamily{U}{mathx}{\hyphenchar\font45}
\DeclareFontShape{U}{mathx}{m}{n}{
      <5> <6> <7> <8> <9> <10>
      <10.95> <12> <14.4> <17.28> <20.74> <24.88>
      mathx10
      }{}
\DeclareSymbolFont{mathx}{U}{mathx}{m}{n}

\DeclareMathSymbol{\obot}{2}{matha}{"6B}

\makeatletter
\newcommand*\leftdash{\rotatebox[origin=c]{-45}{$\dabar@\dabar@\dabar@$}}
\newcommand*\rightdash{\rotatebox[origin=c]{45}{$\dabar@\dabar@\dabar@$}}
\makeatother

\RequirePackage{xspace}

\newcommand{\sA}{\ensuremath{\mathscr{A}}\xspace}

\newcommand{\sM}{\ensuremath{\mathscr{M}}\xspace}

\newcommand{\sS}{\ensuremath{\mathscr{S}}\xspace}

\newcommand{\sZ}{\ensuremath{\mathscr{Z}}\xspace}

\newcommand{\rH}{\ensuremath{\mathrm{H}}\xspace}

\newcommand{\rV}{\ensuremath{\mathrm{V}}\xspace}

\usepackage{color}

\newcommand{\F}{\mathbf{F}}

\newcommand{\G}{\mathbf{G}}

\newcommand{\wt}[1]{\widetilde{#1}}

\newcommand{\Q}{\mathbf{Q}}
\newcommand{\Z}{\mathbf{Z}}

\newcommand{\mf}[1]{\mathfrak{#1}}

\newcommand{\ul}[1]{\underline{#1}}
\newcommand{\ol}[1]{\overline{#1}}
\newcommand{\wh}[1]{\widehat{#1}}

\newcommand{\Cal}[1]{\mathcal{#1}}
\newcommand{\A}{\mathbf{A}}
\DeclareMathOperator{\et}{\text{\'{e}t}}
\newcommand{\mbf}[1]{\mathbf{#1}}

\newcommand{\co}{\colon}
\newcommand{\mrm}[1]{\mathrm{#1}}

\newcommand{\bbm}[1]{\mathbbm{#1}}
\newcommand{\PP}{\mathbf{P}}

\newcommand{\DD}{\mathbf{D}}

\newcommand{\TT}{\mathbb{T}}
\newcommand{\inj}{\hookrightarrow}
\newcommand{\surj}{\twoheadrightarrow}

\newcommand{\tw}[1]{\sm{\langle #1 \rangle}}

\newcommand{\di}{\diamond}

\newcommand{\td}{\triangledown}
\newcommand{\tu}{\vartriangle}

\DeclareMathOperator{\Tot}{Tot}
\DeclareMathOperator{\Perf}{Perf}

\DeclareMathOperator{\GL}{GL}

\DeclareMathOperator{\Frob}{Frob}

\DeclareMathOperator{\Tr}{Tr}

\DeclareMathOperator{\Hom}{Hom}
\DeclareMathOperator{\dhom}{hom}
\DeclareMathOperator{\ext}{ext}
\newcommand{\cHom}{\Cal{H}om}
\newcommand{\cExt}{\Cal{E}xt}
\newcommand{\cRHom}{\Cal{R}\Cal{H}om}

\DeclareMathOperator{\rank}{rank}

\DeclareMathOperator{\Nm}{Nm}
\DeclareMathOperator{\Spec}{Spec\,}

\DeclareMathOperator{\End}{End}

\DeclareMathOperator{\Div}{Div}

\DeclareMathOperator{\Bun}{Bun}
\DeclareMathOperator{\Ext}{Ext}
\DeclareMathOperator{\Pic}{Pic}

\DeclareMathOperator{\Id}{Id}

\DeclareMathOperator{\QCoh}{QCoh}

\DeclareMathOperator{\Sym}{Sym}

\DeclareMathOperator{\pt}{pt}
\DeclareMathOperator{\Fib}{Fib}
\DeclareMathOperator{\Cofib}{Cofib}

\DeclareMathOperator{\Fix}{Fix}

\DeclareMathOperator{\Sht}{Sht}

\DeclareMathOperator{\pr}{pr}

\DeclareMathOperator{\Hk}{Hk}

\DeclareMathOperator{\Ch}{CH}

\DeclareMathOperator{\Map}{Map}

\DeclareMathOperator{\RHom}{RHom}
\DeclareMathOperator{\FT}{FT}

\DeclareMathOperator{\Mod}{Mod}
\DeclareMathOperator{\cc}{\mathfrak{c}}
\DeclareMathOperator{\Conv}{Conv}
\DeclareMathOperator{\cone}{cone}

\newcommand{\sm}[1]{{\scriptstyle  #1}}

\newcommand{\bL}{\mathbf{L}}
\newcommand{\bT}{\mathbf{T}}

\DeclareMathOperator{\BM}{BM}
\newcommand{\mBM}{\mrm{H}^{\BM}}

\DeclareMathOperator{\can}{can}

\DeclareMathOperator{\ev}{ev}

\newcommand{\bt}{\boxtimes}

\newcommand{\arith}{\mrm{arith}}

\newcommand{\incl}{\hookrightarrow}
\newcommand{\isom}{\stackrel{\sim}{\to}}

\newcommand{\Ql}{\Q_{\ell}}
\newcommand{\Qll}[1]{\Q_{\ell, #1}}
\newcommand{\Qlbar}{\overline{\Q}_\ell}

\renewcommand{\j}[1]{\langle{#1}\rangle}
\newcommand\un{\underline}
\newcommand{\bu}{\bullet}
\newcommand{\ov}{\overline}

\newcommand\xr{\xrightarrow}
\newcommand\op{\oplus}
\newcommand\ot{\otimes}

\renewcommand\c{\circ}

\newcommand\upH{\textup{H}}
\newcommand{\hBM}[2]{\textup{H}^{\textup{BM}}_{#1}({#2})}  

\renewcommand\a\alpha
\renewcommand\b\beta
\newcommand\g\gamma
\renewcommand\d\delta
\newcommand\D\Delta

\newcommand{\io}{\iota}

\newcommand{\ph}{\varphi}
\newcommand{\s}{\sigma}
\renewcommand{\t}{\tau}

\newcommand{\y}{\eta}
\newcommand{\z}{\zeta}

\newcommand{\om}{\omega}

\newcommand\na{\natural}
\newcommand\sh{\sharp}

\newcommand\cA{\mathcal{A}}

\newcommand\cC{\mathcal{C}}
\newcommand\cD{\mathcal{D}}
\newcommand\cE{\mathcal{E}}
\newcommand\cF{\mathcal{F}}
\newcommand\cG{\mathcal{G}}

\newcommand\cI{\mathcal{I}}

\newcommand\cK{\mathcal{K}}

\newcommand\cL{\mathcal{L}}
\newcommand\cM{\mathcal{M}}
\newcommand\cN{\mathcal{N}}
\newcommand\cO{\mathcal{O}}
\newcommand\cP{\mathcal{P}}
\newcommand\cQ{\mathcal{Q}}

\newcommand\cT{\mathcal{T}}
\newcommand\cU{\mathcal{U}}
\newcommand\cV{\mathcal{V}}
\newcommand\cW{\mathcal{W}}

\newcommand\cY{\mathcal{Y}}
\newcommand\cZ{\mathcal{Z}}

\newcommand\frL{\mathfrak{L}}
\newcommand\frM{\mathfrak{M}}

\newcommand\frc{\mathfrak{c}}
\newcommand\frd{\mathfrak{d}}

\newcommand\frp{\mathfrak{p}}
\newcommand\frq{\mathfrak{q}}

\newcommand\frs{\mathfrak{s}}

\newcommand\Corr{\mathrm{Corr}}
\newcommand\CoCorr{\mathrm{CoCorr}}
\newcommand\univ{\mathrm{univ}}

\def\upp{\textup{upp}}
\def\low{\textup{low}}
\def\out{\textup{out}}

\newtheorem{thm}{Theorem}[subsection]
\newtheorem{lemma}[thm]{Lemma}
\newtheorem{prop}[thm]{Proposition}
\newtheorem{cor}[thm]{Corollary}

\newtheorem{defn-prop}[thm]{Definition-Proposition}

\theoremstyle{remark}
\newtheorem{remark}[thm]{Remark} 
\newtheorem{defn}[thm]{Definition}

\newtheorem{example}[thm]{Example}

\newtheorem{observation}[thm]{Observation}

\makeatletter
\def\th@remark{%
  \thm@headfont{\bfseries}%
  \normalfont 
  \thm@preskip \thm@preskip 
  \thm@postskip\thm@preskip
}
\def\imod#1{\allowbreak\mkern5mu({\operator@font mod}\,\,#1)}
\makeatother

\numberwithin{equation}{subsection}

\title[]{Modularity of higher theta series I: \\ cohomology of the generic fiber}

\author{Tony Feng}
\address{University of California Berkeley, Department of Mathematics, Berkeley, CA 94720, USA}
\email{fengt@berkeley.edu}

\author{Zhiwei Yun}
\address{Massachusetts Institute of Technology, Department of Mathematics, 77 Massachusetts Avenue, Cambridge, MA 02139, USA}
\email{zyun@mit.edu}

\author{Wei Zhang}

\address{Massachusetts Institute of Technology, Department of Mathematics, 77 Massachusetts Avenue, Cambridge, MA 02139, USA}
\email{weizhang@mit.edu}

\begin{document}

\begin{abstract}
In a previous paper we constructed \textit{higher} theta series for unitary groups over function fields, and conjectured their modularity properties. Here we prove the generic modularity of the $\ell$-adic realization of higher theta series in cohomology. The proof debuts a new type of Fourier transform, occurring on the Borel-Moore homology of moduli spaces for shtuka-type objects, that we call the \textit{arithmetic Fourier transform}. Another novelty in the argument is a \textit{sheaf-cycle correspondence} extending the classical sheaf-function correspondence, which facilitates the deployment of sheaf-theoretic methods to analyze algebraic cycles. Although the modularity property is a statement within classical algebraic geometry, the proof relies on derived algebraic geometry, especially a nascent theory of \textit{derived Fourier analysis} on derived vector bundles, which we develop.
\end{abstract}

\maketitle

\tableofcontents

\section{Introduction}
The modularity of theta series has a long and storied history, beginning with Poisson, who famously applied Poisson summation to prove the modularity of Jacobi's theta series.\footnote{The proof of the modularity appears in Jacobi's paper \cite{Jac1828}, where he credits it to Poisson; cf. \cite[p.15, footnote $\ddag$]{Ed01}.} From the modern perspective of automorphic forms, theta series can be constructed much more generally, for all reductive groups fitting into \emph{dual reductive pairs}, and essentially the same Poisson summation argument generalizes to prove their modularity. So, to make a long story short, the modularity of theta series can ultimately be seen as a relatively simple (by modern standards) consequence of Fourier duality. 

Kudla introduced an analogue of theta series in arithmetic geometry, called \emph{arithmetic theta series}. These objects are again constructed as Fourier series, but with coefficients being \emph{algebraic cycles} rather than numbers. They are also conjectured to be modular, but that turns out to be much more difficult to prove (and even to formulate, in the ideal generality). For example, the modularity of arithmetic theta series of divisors on (the integral models of) unitary Shimura varieties was proved only recently, in \cite{BHKRY1}, and has interesting applications such as to the proof of the ``Arithmetic Fundamental Lemma'' \cite{Zh21}. For another example, the modularity of arithmetic theta series on (the integral models of) orthogonal Shimura varieties was proved even more recently, in \cite{HM20} (for divisors) and \cite{HM22}, and has interesting applications such as to the study of exceptional jumps of Picard ranks of K3 surfaces over number fields \cite{SSTT}. 

As these examples illustrate, the modularity of \emph{arithmetic} theta series has so far only been accessible through the codimension one (i.e., divisor) case; the reason for this is immediately clear from the proof strategy, which will be recalled below. In particular, it has so far been inaccessible in situations with no arithmetic theta series of codimension one\footnote{However, Kudla proved in \cite[Theorem 1.1]{Kud21} that the Beilinson-Bloch Conjecture can be used to deduce the modularity of the generating series for compact orthogonal Shimura varieties   (on the generic fiber) even in signatures that do not admit special divisors. The argument relies on particular features of the Hodge diamond of orthogonal Shimura varieties, and does not apply for unitary Shimura varieties.}, such as unitary groups with signature $(p,q)$ where both $p,q>1$. Moreover, this modularity has important and far-ranging consequences for other problems in arithmetic geometry. 

In the papers \cite{FYZ} and \cite{FYZ2}, the authors investigated theta series and arithmetic theta series over function fields, and discovered in that context that the story extends further: for each $r \geq 0$, there are \emph{higher theta series} $\wt{Z}^r$ that specialize to classical theta functions when $r=0$, and arithmetic theta series when $r=1$. The adjective ``higher'' refers to the fact, established in \cite{FYZ}, that these higher arithmetic theta series are related to \emph{higher derivatives} of Siegel-Eisenstein series, a generalization of the Siegel-Weil and arithmetic Siegel-Weil formulas. 

The main objective of \cite{FYZ2} was the construction of higher theta series, and the precise formulation of a Modularity Conjecture asserting their modularity property. The construction is itself a substantial task, because certain so-called ``singular'' Fourier coefficients comprising the Fourier series are especially subtle and complicated (the authors did not know how to define the singular Fourier coefficients at the time of writing \cite{FYZ}). In particular, we emphasize that the singular Fourier coefficients appear to be \emph{more} complicated in the function field setting than their counterparts over number fields (see the discussion in \S \ref{ssec: comparison to number fields}), at least from the perspective of classical algebraic geometry. A new insight of \cite{FYZ2}, however, was that \emph{all} Fourier coefficients -- singular or not -- have a uniform and concise description in terms of \emph{derived algebraic geometry}; we refer to the Introduction of \emph{loc. cit.} for more discussion of these issues. As was anticipated there, this phenomenon extends to the more traditional number field context of arithmetic theta series: Madapusi has recently found an interpretation of the virtual fundamental classes of special cycles on Shimura varieties in terms of derived algebraic geometry \cite{M22}. 

In the present paper, we prove the modularity of the higher theta series for all $r$, after realization in $\ell$-adic cohomology and restriction to the ``generic fiber'' (whose technical meaning will be explained below in \S \ref{ssec: formulation}). Notably, the argument is completely uniform in all parameters, including both $r$ and the codimension of the cycles, unlike what one has in the number field setting. We also believe that it will apply with little change for symplectic-orthogonal dual pairs, as well as unitary groups with different ``signatures'', although for comprehensibility we have not written it in the maximum generality here. The proof employs an \emph{arithmetic} incarnation of Fourier duality, which can be seen as a natural generalization to algebraic cycles of the Fourier-analytic argument for the modularity of the classical theta series (i.e., the case $r=0$). In particular, our proof is completely different from existing proofs of modularity for arithmetic theta series.

\subsection{Formulation of the results}\label{ssec: formulation} We turn next to a precise formulation of our results. Let $X' \rightarrow X$ be an \'{e}tale double cover of smooth projective curves over a finite field $\F_q$ of characteristic $p>2$. Fix integers $n\ge m \geq 1$, and $r \geq 0$.

We recall the following definitions from \cite[\S 4.5]{FYZ2}:
\begin{itemize}
\item Let $\Bun_{GU^-(2m)}$ be the moduli stack of triples $(\cG, \frM, h)$ where $\cG$ is a vector bundle of rank $2m$ over  $X'$, $\frM$ is a line bundle over $X$, and $h$ is a skew-Hermitian isomorphism $h: \cG\isom \s^{*}\cG^{\vee}\ot \nu^{*}\frM=\s^{*}\cG^{*}\ot\nu^{*}(\om_{X}\ot \frM)$. 
\item Let $\Bun_{\wt P_{m}}$ be the moduli stack of quadruples $(\cG,\frM, h,\cE)$ where $(\cG,\frM,h)\in \Bun_{GU^-(2m)}$, and $\cE\subset \cG$ is a Lagrangian sub-bundle (of rank $m$). 
\item Let $\Sht_{GU(n)}^r$ be the moduli stack of rank $n$ similitude Hermitian shtukas. 
\end{itemize}

In \cite[\S 4]{FYZ2}, we constructed the \emph{higher theta series}
\[
\wt Z^{r}_{m}: \Bun_{\wt P_{m}}(k) \to \Ch_{r(n-m)}(\Sht^{r}_{GU(n)}).
\]
The value of $\wt Z^{r}_{m}$ on a tuple $(\cG,\frM, h,\cE)$ is defined as a Fourier series, with Fourier coefficients $[\cZ_{\cE}^r(a)]$ where the Fourier parameter $a$ is a Hermitian map $\cE\to \s^{*}\cE^{\vee}\ot \nu^{*}(\om_X\ot \frM)$. 

The map $\Bun_{\wt P_{m}}(k) \rightarrow \Bun_{GU^-(2m)}(k)$, given by forgetting the Lagrangian sub-bundle $\cE \subset \cG$, is surjective, and \cite[Modularity Conjecture 4.15]{FYZ2} predicts that $\wt{Z}^r_m$ descends through this map to induce a function $Z^{r}_{m}: \Bun_{GU^-(2m)}(k) \to \Ch_{r(n-m)}(\Sht^{r}_{GU(n)})$, as in the diagram below. 
\[
\begin{tikzcd}
\Bun_{\wt P_{m}}(k)  \ar[d, twoheadrightarrow] \ar[dr, "\wt Z^{r}_{m}"] \\
 \Bun_{GU^-(2m)}(k) \ar[r, dashed, "Z^{r}_{m}"'] &  \Ch_{r(n-m)}(\Sht^{r}_{GU(n)})
\end{tikzcd}
\]
In other words, the Modularity Conjecture says that the function $\wt Z^{r}_{m}$, which a priori depends on $(\cG, \frM, h, \cE)$, is actually independent of the Lagrangian sub-bundle $\cE \subset \cG$.

Next we proceed to describe the main theorem of this paper. For $\ell \neq p$, there is an $\ell$-adic realization map 
\[
\Ch_{r(n-m)}(\Sht^{r}_{GU(n)}) \rightarrow \mBM_{2r(n-m)}(\Sht^{r}_{GU(n)})
\]
where $\mBM(Y)$ denotes the $\ell$-adic Borel-Moore homology of a space $Y \xrightarrow{\pi} \Spec k$, i.e., $\mBM_{2i}(Y) := \rH^{-2i}(Y; \pi^! \Qll{\Spec k} (-i) )$. We denote by $|\wt Z^{r}_{m}|_{\ell}$ the composition of the $\ell$-adic realization map with the higher theta series, which is a function 
\[
|\wt Z^{r}_{m}|_{\ell} \co \Bun_{\wt P_{m}}(k)  \rightarrow \mBM_{2r(n-m)}(\Sht^{r}_{GU(n)}).
\]

We consider a modification of $|\wt Z^{r}_{m}|_{\ell}$ according to the following structures:
\begin{itemize}
\item The stack $\Sht^{r}_{GU(n)}$ is locally of finite type, and admits a presentation as an inductive limit of finite type open substacks $\Sht^{r, \leq \mu}_{GU(n)}$ where $\mu$ is a Harder-Narasimhan polygon for $GU(n)$. Hence we have restriction maps 
\[
\mBM_{2r(n-m)}(\Sht^{r}_{GU(n)}) \rightarrow \mBM_{2r(n-m)}(\Sht^{r, \leq \mu}_{GU(n)})
\]
for every $\mu$. 

\item The stack $\Sht^{r}_{GU(n)}$ admits a ``leg map'' $\Sht^{r}_{GU(n)} \rightarrow (X')^r$. Let $\y=\Spec F'\to X'$ be the generic point. Let $\y^{r}=\Spec(F'\ot_{k}\cdots\ot_k F')\to (X')^{r}$. Note that $\y^{r}$ contains the generic point of $(X')^{r}$ but it also contains many more points such as the generic point of the diagonal $X'$. 
Hence we have a restriction map 
\[
\mBM_{2r(n-m)}(\Sht^{r}_{GU(n)}) \rightarrow \mBM_{2r(n-m)}(\Sht^{r}_{GU(n)} \times_{(X')^r} \eta^{r}).
\]
\end{itemize}
Our main result is that the function $|\wt Z^{r}_{m}|_{\ell}$ is modular after composing with the restriction maps of the bullet points above: 

\begin{thm}\label{thm: main}
The composition 
\[
\Bun_{\wt P_{m}}(k)  \xrightarrow{|\wt Z^{r}_{m}|_{\ell} }  \mBM_{2r(n-m)}(\Sht^{r}_{GU(n)}) \rightarrow \varprojlim_\mu \mBM_{2r(n-m)}(\Sht^{r, \leq \mu}_{GU(n)} \times_{(X')^r} \eta^{r})
\]
descends through $\Bun_{\wt P_{m}}(k)  \surj \Bun_{GU^-(2m)}(k)$. In other words, its value on $(\cG, \frM, h, \cE) \in \Bun_{\wt P_{m}}(k) $ is independent of the Lagrangian sub-bundle $\cE \subset \cG$. 
\end{thm}

\begin{remark} Theorem \ref{thm: main} implies the modularity of $|\wt Z^{r}_{m}|_{\ell}$ restricted to the generic point of $(X')^{r}$ (hence also the geometric generic fiber of $(X')^{r}$), but it contains more information. For example, it also implies the modularity of $|\wt Z^{r}_{m}|_{\ell}$ restricted to the generic point of the diagonal $\D(X')\incl (X')^{r}$ (i.e., all legs coincide), whose geometry is quite different from the generic fiber over $(X')^{r}$. 
\end{remark}

\begin{remark} We emphasize that Theorem \ref{thm: main} can be formulated completely within classical algebraic geometry, while its proof will draw upon the theory of \emph{derived algebraic geometry}. 
\end{remark}

\subsection{Comparison to number fields}\label{ssec: comparison to number fields}

As hinted earlier, the case $r=0$ of Theorem \ref{thm: main} is classical, while the case $r=1$ of Theorem \ref{thm: main} is parallel to the modularity of arithmetic theta series on the generic fiber of Shimura varieties. Therefore it is natural to compare Theorem \ref{thm: main} to analogous results for arithmetic theta series on the generic fiber, which we refer to as ``generic modularity''. 

One analogous result to (the $r=1$ case of) Theorem \ref{thm: main} is the landmark work of Kudla-Millson \cite{KM90}, establishing modularity in Betti cohomology of Shimura varieties, which was a vast generalization of a theorem of Hirzebruch-Zagier \cite{HZ76}. (Note however that this amounts to modularity in the \emph{geometric} generic fiber, which is weaker than modularity in the generic fiber.) This will be discussed further below. For the analogous problem on orthogonal Shimura varieties, a similar result is the work of Borcherds \cite{Bor99} which established the generic modularity in the Chow group for the codimension 1 case. Using Borcherds' work, the third author's thesis \cite{Zh09} proved the generic modularity in arbitrary codimensions conditionally upon a convergence hypothesis. Finally, this convergence hypothesis was established by Bruinier-Westerholt-Raum \cite{BWR15}, completing the proof of the generic modularity in arbitrary codimensions (in the orthogonal case). Of course, these achievements built upon work of many other people, whom we have not mentioned. 

Naturally, our initial attempts to prove the Modularity Conjecture started by looking to the proofs of the above results for inspiration. However, we did not find a way to adapt any of their ideas to the function field case, for reasons that we will briefly explain. 

\subsubsection{The work of Kudla-Millson} As mentioned above, the modularity of arithmetic theta series in the Betti cohomology of the geometric generic fiber was obtained by Kudla-Millson \cite{KM90}. Roughly speaking, they imitate the proof of modularity for theta functions, but replacing functions by differential forms on the complex points of the relevant Shimura varieties, which are then uniformized by complex hermitian domains. Unfortunately for us, this argument relies fundamentally on features that do not exist in positive characteristic, such as:
\begin{itemize}
\item An ``analytic description'' of cohomology classes in terms of automorphic forms, coming from de Rham theory. 
\item The control of Betti cohomology of locally symmetric spaces provided by $(\mf{g}, K)$-cohomology. 
\end{itemize}
By contrast, we have no analogous ``uniformization'' of $\Sht_{GU(n)}^r$, we cannot represent their $\ell$-adic cohomology classes by concrete objects close to automorphic functions, and their cohomology groups are comparatively very complicated (e.g., infinite dimensional).

We remark that although the statement of our Theorem \ref{thm: main} is formally analogous in the $r=1$ case to the results of Kudla-Millson on geometric modularity in cohomology, the actual arguments seem to have nothing in common. In particular, the reason we restrict to the generic fiber has nothing to do with the previous paragraph; for us the point is that we need to add level structure along certain points on the curve (that we have no control over), and the level-structure cover is generically finite but may fail to be finite when these points coincide with the legs. If the cover were proper over the whole curve, then we would be able to execute our argument over the whole curve. Also, we prove modularity in absolute cohomology, i.e., without having to pass to the \emph{geometric} generic fiber. 

\subsubsection{The work of Borcherds, etc.} Except in low rank cases that can be analyzed explicitly, all other approaches to modularity of arithmetic theta series are based on the method of Borcherds \cite{Bor99}. A summary of this method can be found in (for example) \cite[\S 1.2]{BHKRY1}. Roughly speaking, it proceeds by using \emph{Borcherds products} to lift weakly holomorphic modular forms to meromorphic forms on the unitary Shimura variety. The divisor of each such form provides a relation in the Chow group of the Shimura variety. Applying this to the entire space of weakly holomorphic modular forms, of the correct weight and level, leads to a host of such relations, which comprise the content of modularity, by Borcherds' modularity criterion. 

Unfortunately for us, no analogue of Borcherds lifting exists in positive characteristic. 

Moreover, one might say that the strategy above relies implicitly on the fact that the zero-th Fourier coefficient of the generating series has a simple form: it is the negative of the first Chern class of the line bundle of modular forms $\mbf{\omega}$. It is for this reason that constructing modular forms, i.e., sections of $\mbf{\omega}$, produces the right relations. In general, one expects roughly that the \emph{singular} Fourier coefficients of arithmetic theta series should be a power of this Chern class times a cycle that ``looks like'' a non-singular Fourier coefficient (see \cite{Kud04} for more precise formulations). 

By contrast, in \cite{FYZ2} we proposed a construction of the singular Fourier coefficients for higher theta series over function fields, which turned out to be much more complicated. For example, the constant term $\cZ_{\cE}^r(0)$ of the higher arithmetic theta series has a decomposition into \emph{infinitely} many (if $m>1$) open-closed pieces, indexed by sub-bundles $\cK \subset \cE$. The piece labeled by the sub-bundle $\cK = 0$ is what we call the ``least degenerate stratum'', while the piece labeled by $\cK = \cE$ dominates the whole moduli space of shtukas. Correspondingly, the virtual fundamental class $[\cZ_{\cE}^r(0)]$ is an infinite sum\footnote{This is a well-defined cycle because $\Sht_{GU(n)}^r$ is of infinite type. On any quasi-compact open substack $\Sht_{GU(n)}^{r, \leq \mu}$, only finitely many of these summands are supported.} of the form
\begin{equation}\label{eq: zero coefficient}
[\cZ_{\cE}^r(0)] := \sum_{\cK \subset \subset \cE} \left(\left(\prod_{i=1}^r c_{\mrm{top}}(p_i^* \sigma^* \cK^* \otimes \ell_i)\right)\cap [\cZ_{\cE/\cK}^{r,\circ}(\ol{0})]  \right).
\end{equation}
The notation is explained in \cite[\S 4]{FYZ2}; we do not explain it here as we only want to refer to coarse aspects of its form.
\begin{itemize}
\item The summand indexed by $\cK = \cE$ contributes the top Chern class of a vector bundle, which is analogous to the Hodge bundle $\omega^{-1}$ in the number-field case.
\item  The summand indexed by $\cK = 0$ (which we think of as the ``least degenerate'' piece) contributes a virtual class defined by certain, somewhat complicated, derived intersections of cycles. 
\item The intermediate terms, indexed by non-zero proper sub-bundles $\cK \subset \cE$, contribute some mixture of the above extremes: they are a Chern class times the virtual fundamental class of the ``least degenerate'' piece from a lower-dimensional situation. 
\end{itemize}
From this perspective, what happens over number fields is that only one summand from \eqref{eq: zero coefficient} appears (namely, the one corresponding to the ``most degenerate stratum''), because the other pieces are precluded by considerations at the archimedean place. 

In summary, the vastly more complicated form of \eqref{eq: zero coefficient}, as compared to the number field case, makes it difficult to imagine proving modularity by explicitly constructing all the necessary relations. 

\subsection{New ingredients}\label{ssec: intro ingredients of proof} Having explained why the pre-existing approaches to modularity do not seem applicable in our setting, we now proceed to describe the novel ingredients featuring into our proof of Theorem \ref{thm: main}.

\subsubsection{Derived fundamental classes}\label{sssec: derived fundamental classes}
The elementary but complicated definitions of the virtual fundamental cycles $[\cZ_{\cE}^r(a)]$, such as in \eqref{eq: zero coefficient}, are too unwieldy for us to work with effectively. A key point is to find a more conceptual description, which is uniform in the Fourier parameter $a$. 

An insight of \cite{FYZ2} is that the virtual fundamental classes admit an alternative description: they are the ``na\"{i}ve'' fundamental classes from a more sophisticated perspective. Namely, recall that the higher theta series is defined as a Fourier series, with Fourier coefficients $[\cZ_{\cE}^r(a)]$ where $a$ is the Fourier parameter. Here $\cZ_{\cE}^r(a)$ is a certain space which is finite over $\Sht_{GU(n)}^r$, but often of the ``wrong'' dimension, so the associated cycle class $[\cZ_{\cE}^r(a)]$ must be constructed as a \emph{virtual} fundamental class. However, it was discovered in \cite{FYZ2} that ``repeating'' the definition of $\cZ_{\cE}^r(a)$ in the natural way within derived algebraic geometry produces a derived stack $\sZ_{\cE}^r(a)$ which is \emph{quasi-smooth} of the ``correct'' dimension. As explained in \cite{KhanI}, a quasi-smooth derived stack $\sS$ has an intrinsic notion of \emph{fundamental class} $[\sS]$, which can be interpreted as a ``virtual fundamental class'' on its classical truncation. The classical truncation of $\sZ_{\cE}^r(a)$ is $\cZ_{\cE}^r(a)$, and we calculated that $[\sZ_{\cE}^r(a)]$ coincides with the elementary but complicated construction of the virtual class $[\cZ_{\cE}^r(a)]$; thus for example the ``na\"{i}ve'' notion of fundamental class $[\sZ_{\cE}^r(0)]$ of the derived stack $\sZ_{\cE}^r(0)$ agrees with the unwieldy formula \eqref{eq: zero coefficient}.

\subsubsection{Arithmetic Fourier transform} As mentioned at the beginning, the modularity of classical theta series is based on the Fourier transform. In order to prove Theorem \ref{thm: main}, we introduce an \emph{arithmetic Fourier transform} over $\mBM(\Sht^r_{GU(n)})$, which specializes to a relative version of the usual Fourier transform over finite fields when $r=0$. 

The construction of the arithmetic Fourier transform is formally analogous to that of the usual Fourier transform. We consider a level structure cover $``\Sht_V^r" \rightarrow \Sht^r_{GU(n)}$, which is an $\F_q$-vector space in stacks over $\Sht^r_{GU(n)}$. (The main reason for working on the generic fiber is that the geometry of the level structure cover is relatively simple over the generic fiber.) Therefore, it has a dual cover $``\Sht_{\wh{V}}^r" \rightarrow \Sht^r_{GU(n)}$, and an evaluation map 
\[
\Sht_V^r \times_{\Sht_{GU(n)}^r} \Sht_{\wh{V}}^r \xrightarrow{\ev} \F_q.
\]
For a nontrivial additive character $\psi$ of $\F_q$, the arithmetic Fourier transform 
\[
\FT^{\arith} \co \mBM(\Sht_V^r) \rightarrow \mBM(\Sht_{\wh{V}}^r)
\]
is defined in terms of the diagram
\[
\begin{tikzcd}
& \Sht_V^r \times_{\Sht_{GU(n)}^r} \Sht_{\wh{V}}^r \ar[r, "\ev"]  \ar[dl, "\pr_1"'] \ar[dr, "\pr_2"] & \F_q\\
\Sht_V^r & & \Sht_{\wh{V}}^r
\end{tikzcd}
\]
if $d$ denotes the rank of the relative $\F_q$-vector space, by sending $\alpha \in  \mBM(\Sht_V^r)$ to $(-1)^d\pr_{2!} (\ev^* \psi \cdot  \pr_1^* \alpha) \in \mBM(\Sht_{\wh{V}}^r)$. 

Recall that the Modularity Conjecture can be phrased as independence of the higher theta series $\wt{Z}_m^r$ on the Lagrangian sub-bundle $\cE \subset \cG$. In other words, its content is that for two different choices of Lagrangian sub-bundles $\cE_1, \cE_2 \subset \cG$, one has 
\begin{equation}\label{eq: intro indep of lagrangian}
\wt{Z}_m^r(\cG, \frM, h, \cE_1) = \wt{Z}_m^r (\cG, \frM, h, \cE_2).
\end{equation}
In the case of the classical theta series $r=0$, one can prove this relation as follows. Assuming for simplicity that $\cE_1$ and $\cE_2$ are \emph{transverse} Lagrangian sub-bundles, one can factor the special ``cycles'' through a level cover $\Sht_V^0 \rightarrow \Sht_{GU(n)}^0$ (adding level structure along a subset depending on $\cE_1, \cE_2$), which is an $\F_{q}$-vector space over $\Sht_{GU(n)}^0$ equipped with a self-duality. One can show that the Fourier transform of the special ``cycle'' $\sum_{a_1} [\cZ^0_{\cE_1}(a_1)]$ for $\cE_1$, which is really just a function on the discrete set $|\Sht^0_{V}|$, is essentially equal to the special cycle $\sum_{a_2}  [\cZ^0_{\cE_2}(a_2)]$ for $\cE_2$. The theta functions for $\cE_1, \cE_2$ are obtained by pairing the respective special cycles with a Gaussian, so then the equality \eqref{eq: intro indep of lagrangian} for $r=0$ follows from the Plancherel formula (i.e., unitarity of the finite Fourier transform) and the Fourier self-duality of Gaussians. 

We can formulate a generalization of this statement for higher $r$: the special cycles $\sum_{a_1} [\cZ^r_{\cE_1}(a_1)]$ and $\sum_{a_2}  [\cZ^r_{\cE_2}(a_2)]$ factor through a certain self-dual level cover $\Sht_V^r \rightarrow \Sht_{GU(n)}^r$, and: 
\begin{equation}\label{eq: intro quote}
  \parbox{\dimexpr\linewidth-8em}{%
    \strut
    The \emph{arithmetic} Fourier transform $\FT^{\arith} $ of the special cycle $\sum_{a_1} [\cZ^r_{\cE_1}(a_1)]$ should be essentially equal to the special cycle $\sum_{a_2}  [\cZ^r_{\cE_2}(a_2)]$. 
    \strut
  }
\end{equation}
However, for $r>0$ it is much less clear how one would prove such a statement, and this requires another innovation that we describe next. 

\begin{remark}\label{rem: AFT for number fields} 
The arithmetic transform does not depend on features specific to the function field context, such as the possibility of ``multiple legs'' or the existence of categorifications, so it makes sense even for Shimura varieties as in the traditional context of arithmetic theta series. It is therefore enticing to wonder how much of our strategy can be ported over to number fields. One new puzzle that arises when trying to do this is that the Archimedean place must be incorporated somehow (even when working on the generic fiber). 
\end{remark}

\subsubsection{The sheaf-cycle correspondence}

Grothendieck's sheaf-function correspondence associates to an $\ell$-adic sheaf on a variety over a finite field $\F_q$, a function on its $\F_q$-points. This allows one to bring the tools of sheaf theory to bear on functions, and its utility is by now well-documented in myriad applications.

In fact, Grothendieck's formalism \cite{SGA5} can also be applied to produce higher dimensional cohomology classes from sheaves, although we are not aware of any instance until now where this observation has been used. In order to prove a statement like \eqref{eq: intro quote}, we extend Grothendieck's formalism to a framework that we call a \emph{sheaf-cycle} correspondence, in order to bring the tools of sheaf theory to bear on the analysis of algebraic cycles. 

To hint at what this entails, we recall that in the sheaf-function correspondence, one begins with an endomorphism of a sheaf, and then extracts a function by taking the trace of the endomorphism. In the sheaf-cycle correspondence, one begins with a \emph{derived} endomorphism of a sheaf (i.e., a higher $\Ext$ class) and then extracts, by a generalization of the ``trace'' operation, a cycle class in cohomology. (By working with motivic sheaves, one can refine the trace to produce a class in the Chow group, but that is not considered in the present paper).

As usual, in practice it is useful to consider generalizations with complexes instead of sheaves, and correspondences instead of maps. Thus, in its general form, the sheaf-cycle correspondence applies a ``trace'' to extract a cycle class from a \emph{cohomological correspondence}, which is a certain map of complexes.

To prove the precise statement underlying \eqref{eq: intro quote}, we realize the virtual fundamental classes $\sum_{a_1} [\cZ_{\cE_1}^r(a_1)]$ and $\sum_{a_2} [\cZ_{\cE_2}^r(a_2)]$ as arising by the sheaf-cycle correspondence from cohomological correspondences $\cc_U$ and $\cc_{U^\perp}$, respectively. Then, we prove that a \emph{sheaf-theoretic} Fourier transform essentially takes $\cc_U$ to $\cc_{U^\perp}$. Finally, from this sheaf-theoretic statement we extract \eqref{eq: intro quote} by taking the trace. The above strategy uses crucially the additional flexibility afforded by sheaves (as opposed to cycles). 

\begin{remark}
The classical sheaf-function correspondence includes a compatibility with pushforward and pullback operations; in particular, the pullback compatibility is obvious there. For the sheaf-cycle correspondence, there is a form of pullback compatibility but it is much subtler, and seems to require derived algebraic geometry even to formulate (a reflection of the fact that pullback of algebraic cycles is a subtle operation, which is most robustly understood through derived algebraic geometry). Perhaps this is a reason why the sheaf-cycle correspondence has taken relatively long to materialize into applications. 
\end{remark}

\subsubsection{Derived Fourier analysis} Via the sheaf-cycle correspondence, the duality \eqref{eq: intro quote} ultimately comes out of a new apparatus that we call \emph{derived Fourier analysis}. This involves a generalization of Deligne-Laumon's theory of $\ell$-adic Fourier transform, which takes place on vector bundles, to a context that we call ``derived vector bundles'', which are spaces built out of \emph{perfect complexes}, generalizing how vector bundles are built from locally free coherent sheaves.

An example of a derived vector bundle is the \emph{derived} fibered product of a morphism of classical vector bundles $E' \rightarrow E$ with the zero section of $E$. Derived vector bundles also include certain types of classical stacks as well. Derived vector bundles have duals, and this duality interchanges the ``classical stacky'' and ``derived'' directions of derived vector bundles. In particular, the dual of a classical stack can have nontrivial derived structure, and vice versa. 

We can give a brief hint as to the role of derived vector bundles. For $r=0$, the special cycle $\sum_{a_1} \cZ_{\cE_1}^0(a_1)$ is a counting function on the set of (Hermitian) vector bundles $\cF$ over $X'$, which sends $\cF$ to the number of maps $\# \Hom_{X'}(\cE_1, \cF)$. For $r>0$, we want to let $\cF$ vary in moduli, but the vector spaces $\Hom_{X'}(\cE_1, \cF)$ do not assemble into a vector bundle as $\cF$ varies, for example because their dimensions jump discontinuously with $\cF$. However, the ``derived vector spaces'' $\RHom_{X'}(\cE_1, \cF)$ do (informally speaking) assemble into a derived vector bundle, which is locally of the form described in the first sentence of the preceding paragraph. 

In particular, the previously discussed cohomological correspondences $\cc_U$ and $\cc_{U^\perp}$ live on derived vector bundles of the above sort, and are defined using the notion of relative fundamental class for a quasi-smooth map of derived schemes. We therefore develop the theory of $\ell$-adic Fourier transform on derived vector bundles in order to compute with them. It turns out that there are several new technical challenges in the derived setting, which would be interesting for further study.

\begin{remark}
The primordial forms of derived Fourier analysis were discovered through computations in \cite{FW}, and the theory we develop here will also be applied in \emph{loc. cit.} (along with other ingredients) in order to categorify the Rankin-Selberg unfolding method for automorphic periods. Some of our results on the derived Fourier transform were inspired by work-in-progress of Adeel Khan investigating a derived Fourier transform for homogeneous sheaves. 
\end{remark}

\subsection{Organization of the paper} We provide some commentary on the organization of the paper. 

In the next section of the paper, \S \ref{sec: transverse Lagrangians}, we explain a proof of modularity in the special case $r=0$, as a template/toy model for the general case. Very roughly speaking, this proof will be geometrized from functions to sheaves, and then Theorem \ref{thm: main} will be extracted from the sheaf-theoretical level by an appropriate trace operation (which depends on $r$). The proof for $r=0$ is in \S \ref{ssec: r=0 transverse lagrangians}, and then in \S \ref{ssec: outline} we give an overview of the strategy for the general case, which relies on a setup that we call the ``transverse Lagrangian ansatz''. As the implementation of the strategy is quite long and involved, we recommend referring back to this overview repeatedly for guidance. In particular, we defer a discussion of the organization of some individual sections of this paper to \S \ref{ssec: outline}.

Part I, consisting of \S \ref{sec: zoo} -- \S \ref{sec: descent}, is devoted to the formalism of cohomological correspondences and their interaction with algebraic cycles through the sheaf-cycle correspondence. 

The notion of a cohomological correspondence, and the operation of extracting an algebraic cycle as the trace of a cohomological correspondence, are explained in \S \ref{ssec: trace of cc}. The majority of Part I is devoted to constructing the functoriality operations for cohomological correspondences, and establishing their compatibility with the formation of the trace. This story is much subtler than its analogue for the classical sheaf-function correspondence; in particular, derived algebraic geometry already arises naturally and crucially in the basic formulations. 


Part II, consisting of \S \ref{sec: FT} -- \S \ref{sec: arithmetic FT}, develops Fourier analysis in two new contexts.

In \S \ref{sec: FT}, we generalize the Deligne-Laumon Fourier transform for $\ell$-adic sheaves to a derived setting. We recall the notion of \emph{derived vector bundles}, which are built out of a perfect complex of coherent sheaves in a manner generalizing how vector bundles are built from locally free coherent sheaves. Then we define the \emph{derived Fourier transform} and state its basic properties, with the proofs deferred to Appendix \ref{app: A}. Actually, we are only able to establish one of these properties under a technical assumption of ``global presentability'', which appears to be an artefact of the proof. This is good enough for our purposes but it would be more satisfactory to remove it, which seems an interesting problem. 
 
 Then \S \ref{sec: FT cc} studies the interaction of the derived Fourier transform with cohomological correspondences between derived vector bundles. Next \S \ref{sec: arithmetic FT} introduces the arithmetic Fourier transform, establishes its basic properties, and relates it to the derived Fourier transform through the sheaf-cycle correspondence. 

Part III, consisting of \S \ref{ssec: modularity for cc} and \S \ref{sec: AFT and modularity}, assembles the preceding ingredients to complete the proof of Theorem \ref{thm: main}. We postpone an overview of the contents of this Part to \S \ref{ssec: outline}. 


\subsection{Acknowledgments} We thank Adeel Khan, Will Sawin, and Yakov Varshavsky for helpful correspondence. We are grateful to Haoshuo Fu, Adeel Khan, Tasuki Kinjo, Steve Kudla, Chao Li, and Yifeng Liu for comments on a draft of this paper. Parts of this work were carried out while the authors were visiting SLMath. TF was supported by an NSF Postdoctoral Fellowship (DMS \#1902927), an NSF Standard Grant (DMS \#2302520), a Viterbi Postdoctoral Fellowship at the Simons-Laufer Mathematical Sciences Institute, and the Hausdorff Institute for Mathematics (funded by the Deutsche Forschungsgemeinschaft (DFG, German Research Foundation) under Germany's Excellence Strategy – EXC-2047/1 -- 390685813). ZY was supported by the Simons Foundation and the Packard Foundation. WZ was supported by the NSF grant DMS \#1901642, \#2401548, and the Simons Foundation.

\subsection{Notation}

Throughout the paper, let $k=\F_{q}$ be a finite field.

\subsubsection{Notation related to spaces}\label{sssec:notation-on-spaces} In a previous article \cite{FYZ2}, we took care to use calligraphic fonts like $\cZ, \cM$ for classical stacks and script fonts like $\sZ, \sM$ for derived stacks. Starting in this paper, we will always work with derived stacks over $k$ by default (although many of them happen to have the \emph{property} of being classical, i.e., the natural map from the classical truncation is an isomorphism), and we do not use script fonts for derived objects. Hence when we say ``Cartesian square'' we mean what might be called ``derived Cartesian square'' (sometimes we keep the adjective ``derived'' for emphasis), unless noted otherwise. In particular, the notation departs from that of \cite{FYZ2}. 

``Derived Artin stacks'' are defined in \cite{TV08}. In this paper, \emph{all derived Artin stacks are assumed to be locally of finite type over a field}. 

\subsubsection{Notation related to $\ell$-adic sheaves}
Let $Y$ be a derived Artin stack. Then tautologically the classical truncation $Y_{cl}$ of $Y$ is a \emph{higher Artin stack} in the sense of \cite[\S 5.4]{LZa} (this notion goes back to To\"{e}n). We let $D(Y_{cl}) := D(Y_{cl}; \ol{\Q}_{\ell})_a$ be the bounded derived category of constructible \'{e}tale sheaves on $Y_{cl}$ as constructed in \cite[\S 1]{LZb}. This is the bounded subcategory of a homotopy category of a certain stable $\infty$-category $\cD(Y_{cl}; \ol{\Q}_{\ell})_a$ constructed in \emph{loc. cit.}, but we shall only need the six functors and their properties at the level of homotopy categories; for our purposes we prefer the framework of \cite{LZa, LZb} because of their generality in handling higher Artin stacks.

Let $f$ be a map of higher Artin stacks. In \cite[\S 1.3]{LZb} one finds the construction of $f^*$ and $f_*$ for general $f$, the construction of $f_!$ and $f^!$ for locally of finite type $f$, and the construction of $- \otimes_Y -$ and $\cRHom$. One also finds there (\cite[Theorem 1.3.9 and Theorem 1.3.10]{LZb}) the K\"{u}nneth formula, the base change isomorphism, the projection formula, and other ``usual'' properties of the six functors when $f$ is locally of finite type. 

Since $Y_{cl} \rightarrow Y$ induces an isomorphism of \'{e}tale sites by definition, we may  set $D(Y) := D(Y_{cl})$. For a map $f \co Y_1 \rightarrow Y_2$ of derived Artin stacks, we define $f^*, f_*$ to be the corresponding functor on classical truncations; if $f$ is locally of finite type then we define $f^!$ and $f_!$ to be the corresponding functor on classical truncations. In this way we may bootstrap all of \cite{LZa, LZb} to the setting of derived Artin stacks. 

\begin{remark}
It may seem at first that consideration of derived structure is totally irrelevant to the categories of $\ell$-adic sheaves. However this is not the case, as derived structures will be used to construct certain natural transformations of functors between such categories, namely the ``Gysin natural transformations'' associated to quasi-smooth morphisms $f$. This is analogous to how the Chow group of a derived stack is the same as that of its classical truncation, but the derived structure is still useful to construct a virtual fundamental class within the Chow group. 
\end{remark}

For a separated morphism $f\co X\to Y$ of derived Artin stacks that is representable in derived schemes, we let $\can(f)$ be the natural transformation $f_! \rightarrow f_*$ of functors $D(X) \rightarrow D(Y)$. If $f$ is proper, then $f_! = f_*$ and $\can(f)$ is the identity transformation. 

For a quasi-smooth morphism $f \co X \to Y$ of derived Artin stacks, we let $d(f)$ be the \emph{virtual dimension} of $f$ (i.e., Euler characteristic of its tangent complex), which we view as a locally constant function on the source. If $f \co E \rightarrow S$ is a \emph{derived vector bundle} in the sense of \S \ref{sssec: derived vector bundles}, then we call $d(f)$ the \emph{virtual rank} of $E$, and denote it by $\rank(E)$. 

For a locally of finite type morphism $f \co X \rightarrow Y$ of derived Artin stacks, we denote by $\DD_{X/Y} := f^! \Qll{Y}$ the relative dualizing sheaf. If $f$ is smooth of relative dimension $d$, then $f^! = f^*[2d](d)$. We denote by $\DD_{X/Y} (-) \co D(X) \rightarrow D(X)^{\mrm{op}}$ the relative Verdier duality functor, which is represented by the object $\DD_{X/Y}$. For $Y = \Spec k$, we abbreviate $\DD_X := \DD_{X/Y}$.

For any $\ell$-adic complex $\cK \in D(Y)$, we denote by $\cK \tw{i} := \cK[2i](i)$ the indicated shift and Tate twist.

\subsubsection{Notation related to coherent sheaves}

We let $\Perf(Y)$ be the triangulated category of perfect complexes on $Y$, i.e., the full subcategory of the derived category of quasicoherent sheaves on $Y$ spanned by objects locally quasi-isomorphic to finite complexes of finite rank vector bundles. 

For an animated ring $R$, we let $\Mod_R$ be the stable $\infty$-category of $R$-modules. 

For a torsion coherent sheaf $Q$ on a curve $X'$ we let $D_Q$ be its scheme-theoretic support, viewed as a divisor on $X'$, and $|Q| \subset X'$ its set-theoretic support. 

In \cite{FYZ2} we distinguished between the notion of a $\GL(n)$-torsor $\cF$ and the associated vector bundle $\rV(\cF)$, because we wanted to consider maps of the associated vector bundles that are not isomorphisms (and so do not come from maps of torsors). However, this would be too much of a notational burden in the present paper, so we use the same notation for $\cF$ and its associated vector bundle, trusting that context will make the usage clear.

\section{The transverse Lagrangians ansatz}\label{sec: transverse Lagrangians}
In this section, we will explain a proof of modularity of the higher theta series in the special case $r=0$. Although this result goes back at least to Weil \cite{Weil64}, our particular presentation serves as a toy model for the more general argument. In particular, the argument motivates the introduction of certain auxiliary spaces. 

To give a more precise overview of this section: 
\begin{enumerate}
\item In \S \ref{ssec: modularity conjecture} we review the formulation of the Modularity Conjecture, which says that a certain construction of higher theta series $\wt{Z}^r_m$, which a priori depends on a choice of a Lagrangian sub-bundle in a Hermitian bundle, is in fact independent of that choice. 
\item In \S \ref{ssec: reduction to transverse lagrangians}, we reduce the Modularity Conjecture for $\wt{Z}^r_m$ to a slightly weaker independence statement, namely that the values of $\wt{Z}^r_m$ on two \emph{transverse} Lagrangians coincide. 
\item In \S \ref{ssec: r=0 transverse lagrangians}, we prove this independence statement in the case $r=0$. This involves finite Fourier analysis on various auxiliary vector spaces. 
\item In \S \ref{ssec: outline}, we outline the proof of the general case, indicating in particular the ansatz of spaces and maps that will be used to generalize the modularity argument from $r=0$ to arbitrary $r$.  
\end{enumerate}

\subsection{The modularity conjecture for higher theta series}\label{ssec: modularity conjecture}

Recall from \cite[\S4.5]{FYZ2} that $\Bun_{GU(n)}$ parametrizes triples $(\cF,\frL, h)$, where $\cF$ is a vector bundle on $X'$ of rank $n$, $\frL$ is a line bundle on $X$, and $h: \cF\isom \s^{*}\cF^{\vee}\ot\nu^{*}\frL$ is an $\frL$-twisted Hermitian structure (i.e., $\s^{*}h^{\vee}=h$).

Recall from \cite[\S 4.6]{FYZ2} that $\Bun_{GU^{-}(2m)}$ parametrizes triples $(\cG,\frM, h)$, where $\cG$ is a vector bundle on $X'$ of rank $2m$, $\frM$ is a line bundle on $X$, and $h: \cG\isom \s^{*}\cG^{\vee}\ot\nu^{*}\frM$ is an $\frM$-twisted {\em skew-Hermitian} structure (i.e., $\s^{*}h^{\vee}=-h$). Alternatively, we can think of $h$ as an $\cO_{X'}$-bilinear perfect pairing
\begin{equation}\label{h to pairing}
(\cdot,\cdot)_{h}: \cG\times \s^{*}\cG\to \nu^{*}(\frM\ot \om_{X})
\end{equation}
satisfying $(\s^{*}\b,\s^{*}\a)_{h}=-\s^{*}(\a,\b)_{h}$  for local sections $\a$ and $\b$ of $\cG$ respectively.

Let $\Bun_{\wt P_{m}}$ be the moduli stack of quadruples $(\cG,\frM,h,\cE)$ where $(\cG,\frM, h)\in \Bun_{GU^{-}(2m)}$, and $\cE \subset \cG$ is a Lagrangian sub-bundle (i.e., $\cE$ has rank $m$ and the composition $\cE\subset \cG\xr{h}\s^{*}\cG^{\vee}\ot\nu^{*}\frM\to \s^{*}\cE^{\vee}\ot\nu^{*}\frM$ is zero). In \cite[\S 4.6]{FYZ2}, we defined for each $r \geq 0$ and $m \leq n$ a \emph{higher theta series} 
\[
\wt Z^{r}_{m}  \co \Bun_{\wt P_{m}}(k) \to \Ch_{r(n-m)}(\Sht^{r}_{GU(n)}).
\]

We briefly recall the definition of $\wt Z^{r}_{m}$. Let $\frL=\om_{X}\ot \frM$. Let $\Sht^{r}_{U(n),\frL}$ be the moduli stack of rank $n$ Hermitian shtukas $\cF_{\bullet}=((x_{i}), (\cF_{i}), (f_{i}), \ph: \cF_{r}\isom {}^{\t}\cF_{0})$ on $X'$ with $r$ legs and similitude line bundle $\frL$. For a vector bundle $\cE$ on $X'$ of rank $m$, we have the special cycle $\cZ^{r}_{\cE, \frL}$ parametrizing a point $\cF_{\bu}$ of $\Sht^{r}_{U(n),\frL}$, and maps $t_{i}: \cE\to \cF_{i}$ ($0\le i\le r$) compatible with the shtuka structure on $\cF_{\bullet}$. For details we refer to \cite[\S2.3]{FYZ2}. For a Hermitian map $a: \cE\to \s^{*}\cE^{\vee}\ot\nu^{*}\frL$, let $\cZ^{r}_{\cE, \frL}(a)$ be the open-closed substack of $\cZ^{r}_{\cE, \frL}$ consisting of $(\cF_{\bullet}, t_{\bullet})$ such that the Hermitian form on $\cF_{\bullet}$ induces the Hermitian map $a$ on $\cE$ via $t_{\bullet}$. Let $\z: \cZ^{r}_{\cE,\frL}(a)\to \Sht^{r}_{U(n),\frL}\subset \Sht^{r}_{GU(n)}$ be the  map forgetting $t_{\bu}$, which is known to be finite \cite[Proposition 7.5]{FYZ} and unramified.

In \cite[Definition 4.8]{FYZ2} we have defined a {\em virtual fundamental class} $[\cZ^{r}_{\cE,\frL}(a)]\in \Ch_{r(n-m)}(\cZ^{r}_{\cE,\frL}(a))$. Pushing forward along $\z$, we get Chow classes
\begin{equation*}
\z_{*}[\cZ^{r}_{\cE,\frL}(a)]\in \Ch_{r(n-m)}(\Sht^{r}_{U(n),\frL}).
\end{equation*}
The value of $\wt Z^{r}_{m}$ on $(\cG,\frM,h,\cE)$ (recall $\frM=\om_{X}^{-1}\ot\frL$), which we henceforth abbreviate as $(\cG,\cE)$,  is defined as
\begin{equation}
\wt Z^{r}_{m}(\cG,\cE)=\chi(\det\cE)q^{n(\deg \cE-\deg\frL-\deg\om_{X})/2}\sum_{a\in \cA_{\cE,\frL}(k)}\psi(\j{e_{\cG,\cE},a})\z_{*}[\cZ^{r}_{\cE,\frL}(a)].
\end{equation}
Here
\begin{itemize}
\item $\chi: \Pic_{X'}(k)\to \Qlbar^{\times}$ is a character whose restriction to $\Pic_{X}(k)$ is $\y^{n}$, where $\y: \Pic_{X}(k)\to \{\pm1\}$ is the character corresponding to the double cover $X'/X$. 
\item $\psi:\F_{q}\to \Qlbar^{\times}$ is a nontrivial character.
\item the summation of $a$ runs over the set $\cA_{\cE,\frL}(k)$ of all Hermitian maps $a: \cE\to \s^{*}\cE^{\vee}\ot\nu^{*}\frL$, including the singular ones.
\item Let $\cE'=\cG/\cE$. The pairing $(\cdot,\cdot)_{h}$ in \eqref{h to pairing} induces a perfect pairing $\cE\times \s^{*}\cE'\to \nu^{*}\frL$. This identifies $\cE'$ with $\s^{*}\cE^{*}\ot \nu^{*}\frL$. We thus have a short exact sequence
\begin{equation*}
\xymatrix{0\ar[r] & \cE\ar[r] & \cG\ar[r] & \s^{*}\cE^{*}\ot\nu^{*}\frL\ar[r] & 0}
\end{equation*}
giving an extension class $e_{\cG,\cE}\in \Ext^{1}(\s^{*}\cE^{*}\ot\nu^{*}\frL, \cE)$. 
\item The pairing $\j{-,-}$ is the Serre duality pairing between $\Ext^{1}(\s^{*}\cE^{*}\ot\nu^{*}\frL, \cE)$ and $\Hom(\cE, \s^{*}\cE^{\vee}\ot\nu^{*}\frL)$.
\end{itemize}

As explained after \cite[Conjecture 4.15]{FYZ2}, the modularity of $\wt Z^{r}_{m}$ can be formulated as the assertion that $\wt Z^{r}_{m}$ is actually independent of the choice of Lagrangian sub-bundle $\cE$.

\subsection{Reduction to the case of transverse Lagrangians}\label{ssec: reduction to transverse lagrangians}

We first argue that it suffices to show that whenever $\cE_1, \cE_2 \subset \cG$ are two \emph{transverse} Lagrangians in $\cG \in \Bun_{GU^{-}(2m)}(k)$, meaning that their intersection in the vector bundle $\cG$ is the $0$-section, then we have 
\[ 
\wt Z^{r}_{m}(\cG, \cE_1)  = \wt Z^{r}_{m}(\cG,\cE_2).
\]
Note that the condition that $\cE_1, \cE_2$ are transverse is equivalent to their intersection being zero on the (geometric) generic fiber of $X'$. 

\begin{lemma}\label{lem: transverse Lagrangian}
Let $F'/F$ be a non-trivial quadratic extension of fields of characteristic not equal to $2$, with non-trivial automorphism $\sigma$, and $V$ be a finite-dimensional (skew) $F'/F$-Hermitian space. Let $L_1, L_2$ be two Lagrangian subspaces of $V$. Then there exists a Lagrangian subspace $L \subset V$ such that 
\[
L_1 \cap L = L_2 \cap L = 0.
\] 
\end{lemma}

\begin{proof}The arguments for Hermitian and skew-Hermitian spaces are identical, so we give the Hermitian case only. Let $I := L_1 \cap L_2$, an isotropic subspace of $V$. Then (using that the characteristic of $F$ is not $2$) we may find an orthogonal decomposition $V \cong (I \oplus I^*) \oplus V'$ as Hermitian spaces, such that:
\begin{itemize}
\item $I^*$ is Lagrangian in $I \oplus I^*$ and the Hermitian form on $V$ induces a polarization $I \xrightarrow{\sim} I^*$. 
\item $L_1 = I \oplus L_1'$ and $L_2 = I \oplus L_2'$, with each $L_i'$ being Lagrangian in $V'$.  
\end{itemize}
We will take $L$ to be of the form $I^* \oplus L'$, where $L'$ is a Lagrangian in $V'$ transverse to both $L_1'$ and $L_2'$. To see that such $L'$ exists, note that the Lagrangian polarization $V' \cong L_1' \oplus L_2'$ induces an identification of $L_1'$ with $(L_2')^*$. Choosing any basis of $L_1'$ induces a dual basis for $L_2'$, and a corresponding decomposition of $V'$ into a direct sum of 2-dimensional Hermitian spaces. This reduces to the case $\dim_{F'}(V') = 2$. In this case, we may arrange $l_1 \in L_1'$ and $l_2 \in L_2'$ whose non-zero pairing under the Hermitian form lies in $(F')^{\sigma = -1}$. Then $l_1+l_2$ generates a Lagrangian subspace of $V'$ which is transverse to both $L_1'$ and $L_2'$. 

With this choice of $L$, it is clear that $L$ is transverse to both $L_1$ and $L_2$. 
\end{proof}

\begin{cor}
Suppose that for any $\cG \in \Bun_{GU^{-}(2m)}(k)$ and any two transverse Lagrangian sub-bundles $\cE_1, \cE_2 \subset \cG$, we have 
\[
\wt Z^{r}_{m} (\cG, \cE_1) = \wt Z^{r}_{m} (\cG, \cE_2).
\]
Then $\wt Z^{r}_{m}$ is modular. 
\end{cor}

\begin{proof} The meaning of modularity is that $\wt Z^{r}_{m} (\cG, \cE_1) = \wt Z^{r}_{m} (\cG, \cE_2)$ for any two (not necessarily transverse) Lagrangian sub-bundles $\cE_1, \cE_2 \subset \cG$.
By Lemma \ref{lem: transverse Lagrangian}, applied at the generic point of $X'$, we can find a Lagrangian subspace of the generic fiber of $\cG$ which is transverse to the generic fibers of both $\cE_1$ and $\cE_2$. This generic Lagrangian extends uniquely to a Lagrangian sub-bundle of $\cG$ by saturation. Thus we can link any two Lagrangian sub-bundles by a Lagrangian sub-bundle which is transverse to both. 
\end{proof}

Hence, in order to establish Theorem \ref{thm: main}, we are reduced to proving: 

\begin{thm}\label{thm: transverse modularity} For any $\cG \in \Bun_{GU^{-}(2m)}(k)$ and any two transverse Lagrangian sub-bundles $\cE_1, \cE_2 \subset \cG$, we have 
\[
|\wt Z^{r}_{m}(\cG, \cE_1)|_{\ell} = |\wt Z^{r}_{m}(\cG, \cE_2)|_{\ell} \in \varprojlim_{\mu} \mBM_{2(n-m)r} (\Sht_{GU(n)}^{r, \leq \mu} \times_{(X')^r} \eta^{r}).
\]
\end{thm}

The formulation in \cite{FYZ2} involves a similitude line bundle $\LL$ on $X'$. For sanity of notation, \emph{we will present the proof only in the case where $\LL$ is trivial}, so that it may be omitted entirely. The argument can be adapted to include $\LL$ in a completely straightforward manner. Accordingly, we assume henceforth that the similitude line bundle $\mf{M}$ for $\cG$ is $\om_{X}^{-1}$, i.e., the Hermitian form on $\cG$ is an isomorphism $\cG\isom \s^{*}\cG^{*}$.

\subsection{The case $r=0$}\label{ssec: r=0 transverse lagrangians}

The proof of Theorem \ref{thm: transverse modularity} will be long and complex. Some of the complications are caused by technical issues that are not present for $r=0$. Therefore, we will illustrate the argument for $r=0$, which can serve as a simplified model for the general case. 

\subsubsection{}
By definition \cite[Definition 4.13]{FYZ2}, the higher theta series for $r=0$ is a function on $\Bun_{U(n)}(\F_q)$, whose value on $\cF \in \Bun_{U(n)}(\F_q)$ is given by 
\begin{align*}
\wt Z^{0}_{m} (\cG,\cE_1)_{\cF} & := \chi(\det \cE_1) q^{n(\deg \cE_1 - \deg \omega_X)/2} \sum_{a \in \cA_{\cE_1}(k)} \psi(\langle e_{\cG,\cE_{1}}, a \rangle ) \#\cZ^0_{\cE_1}(a)_{\cF} \\
&=  \chi(\det \cE_1) q^{n(\deg \cE_1  - \deg \omega_X)/2} \sum_{t \in \Hom(\cE_1, \cF)} \psi(\langle e_{\cG,\cE_{1}}, a(t) \rangle )
\end{align*}
where $a(t) \in \Hom(\cE_{1}, \s^{*}\cE^{\vee}_{1})$ is the composition 
\[
\cE_1 \xrightarrow{t} \cF \xrightarrow{h_{\cF}} \s^{*}\cF^{\vee}  \xrightarrow{\sigma^* t^{\vee}} \s^{*}\cE_1^{\vee}.
\]
We find it more psychologically convenient to rewrite the index of summation as $\Hom(\cF^*, \cE_1^*) = \Hom(\cE_1, \cF)$. 

Similarly, the value of $\wt Z^{0}_{m} (\cG, \cE_2)$ at $\cF$ is 
\[
\wt Z^{0}_{m} (\cG, \cE_2)_{\cF}   := \chi(\det \cE_2) q^{n(\deg \cE_2  - \deg \omega_X)/2} \sum_{t \in \Hom(\cE_2, \cF )} \psi(\langle e_{\cG,\cE_{2}}, a(t) \rangle ).
\]

We will show that 
\begin{equation}
\wt Z^{0}_{m} (\cG,\cE_1)_{\cF}=\wt Z^{0}_{m} (\cG,\cE_2)_{\cF}, \quad \forall \cF\in \Bun_{U(n)}(\F_{q}),
\end{equation}
whenever $\cE_{1}$ and $\cE_{2}$ are transverse Lagrangians in $\cG$.

\subsubsection{Preliminaries} We begin with some preliminaries that are not specific to $r=0$. Recall that we assume the similitude line bundle $\frL$ is trivial, so that for $(\cG,h)\in \Bun_{GU^{-}(2m)}$, $h$ gives a perfect pairing
\begin{equation}\label{h to pairing O}
(\cdot,\cdot)_{h}: \cG\times \s^{*}\cG\to \cO_{X'}.
\end{equation}
Let $\cE_{1},\cE_{2}$ be Lagrangian sub-bundles of $\cG$. Each inclusion $\cE_i \inj \cG$ induces a short exact sequence
\begin{equation}\label{eq: Lagrangian SES}
0\to \cE_i \rightarrow (\cG  \xrightarrow{\sim} \sigma^* \cG^*) \rightarrow \sigma^* \cE_i^*\to 0.
\end{equation}
Here we are using the form $(\cdot,\cdot)_{h}$ to induce a perfect pairing $\cE_{i}\times \s^{*}(\cG/\cE_{i})\to \cO_{X'}$, which in turn induces an isomorphism $\cG/\cE_{i}\cong \s^{*}\cE_{i}^{*}$.

Another way to formulate the transversality of $\cE_1, \cE_2$ is as follows. If $\cE_1 \cap \cE_2$ vanishes, then the composition 
\[
b_{12}: \cE_1 \rightarrow \cG \rightarrow \sigma^* \cE_2^*
\]
has full rank generically, and therefore has torsion cokernel. Conversely, if the composite map $\cE_1 \rightarrow \sigma^* \cE_2^*$ has torsion cokernel, then $\cE_1 \cap \cE_2$ vanishes. We can think of $b_{12}$ as given by the pairing
\begin{equation}
\cE_{1}\times \s^{*}\cE_{2}\to \cO_{X'}
\end{equation}
obtained by restricting $(\cdot,\cdot)_{h}$ from \eqref{h to pairing O}. Similarly we have $b_{21}: \cE_{2}\to \s^{*}\cE_{1}^{*}$.

Since we are assuming that $\cE_1, \cE_2$ are transverse, we may define torsion sheaves $Q_{1}$ and $Q_{2}$ to fit into the short exact sequences
\begin{eqnarray}
\label{b12} 0\to \cE_{1}\xr{b_{12}} \s^{*}\cE^{*}_{2}\to Q_{2}\to 0,\\
\label{b21} 0\to \cE_{2}\xr{b_{21}} \s^{*}\cE^{*}_{1}\to Q_{1}\to 0.
\end{eqnarray}

Let $\un F'$ be the Zariski constant sheaf on $X'$ with stalks $F'$. For a torsion sheaf $\cT$, $\cT^{*}$ is defined to be $\cHom(\cT, \un F'/\cO_{X'})$. Taking the linear dual (composed with $\s^{*}$) of \eqref{b12}, we get a short exact sequence
\begin{equation}\label{dual b12}
0\to \cE_{2}\xr{\s^{*}b_{12}^{\vee}}\s^{*}\cE^{*}_{1}\to \s^{*}Q^{*}_{2}\to 0.
\end{equation}
Since $(\cdot,\cdot)_{h}$ is skew-Hermitian, and $b_{12}$ and $b_{21}$ can be interpreted as the restrictions of $(\cdot,\cdot)_{h}$ to $\cE_{1}\times \s^{*}\cE_{2}$ and $\cE_{2}\times\s^{*}\cE_{1}$ respectively, we have
\begin{equation}
\s^{*}b_{12}^{\vee}=-b_{21}.
\end{equation}
Comparing \eqref{b21} and \eqref{dual b12}, we get an isomorphism
\begin{equation}
\b_{12}: Q_{1}\isom \s^{*}Q_{2}^{*}
\end{equation}
compatible with the quotient maps $\s^{*}\cE^{*}_{1}\to Q_{1}$ and $\s^{*}\cE^{*}_{1}\to \s^{*}Q^{*}_{2}$ in \eqref{b21} and \eqref{dual b12}. 

Now switching the roles of $\cE_{1}$ and $\cE_{2}$, the same considerations give an isomorphism
\begin{equation}
\b_{21}: Q_{2}\isom \s^{*}Q_{1}^{*}
\end{equation}
compatible with the quotient maps $\s^{*}\cE^{*}_{2}\to Q_{2}$ in \eqref{b12} and $\s^{*}\cE^{*}_{2}\to \s^{*}Q^{*}_{1}$ obtained by dualizing \eqref{b21}.

\begin{lemma}\label{l:beta on Q skew}
The maps $\b_{12}$ and $\b_{21}$ satisfy
\begin{equation}
\s^{*}\b_{12}^{\vee}=-\b_{21}.
\end{equation}
\end{lemma}
\begin{proof}
Let $\cG^{\sh}=\s^{*}\cE_{2}^{*}\op \s^{*}\cE_{1}^{*}$; then $\cG$ is naturally a subsheaf of $\cG^{\sh}$ of the same rank. The form $(\cdot,\cdot)_{h}$ extends to a rational skew-Hermitian pairing 
\begin{equation}\label{rat pairing h}
(\cdot,\cdot)_{h}: \cG^{\sh}\times\s^{*}\cG^{\sh}\to \un F'.
\end{equation}
This restricts to pairings:
\begin{eqnarray}
\s^{*}\cE_{2}^{*}\times \cE^{*}_{1}\to \un F'\\
\s^{*}\cE_{1}^{*}\times \cE^{*}_{2}\to \un F'
\end{eqnarray}
which induce pairings
\begin{eqnarray}
\g_{21}: Q_{2}\times \s^{*}Q_{1}\to \un F'/\cO_{X'},\\
\g_{12}: Q_{1}\times \s^{*}Q_{2}\to \un F'/\cO_{X'}.
\end{eqnarray}
Unwinding the definitions, we see that $\b_{12}$ is induced from $\g_{12}$ and $\b_{21}$ is induced from $\g_{21}$. Since $(\cdot,\cdot)_{h}$ is skew-Hermitian, we see that for local sections $s_{1}$ of $Q_{1}$ and $s_{2}$ of $\s^{*}Q_{2}$, we have 
\begin{equation}\label{g12}
\g_{12}(s_{1}, s_{2})=-\s^{*}\g_{21}(\s^{*}s_{2},\s^{*}s_{1})
\end{equation}
This implies the desired equality for $\b_{12}$ and $\b_{21}$.
\end{proof}

\subsubsection{Self-duality of $Q$}\label{sssec: self-duality of Q} 
Recall from the proof of Lemma \ref{l:beta on Q skew} that $\cG^{\sh}=\s^{*}\cE_{2}^{*}\op \s^{*}\cE_{1}^{*}$ contains $\cG$ as a subsheaf of the same rank. Introduce the torsion sheaf $Q$:
\begin{equation}
Q:=\cG^{\sh}/\cG=(\s^{*}\cE_{2}^{*}\op \s^{*}\cE_{1}^{*})/\cG.
\end{equation}

Consider the commutative diagram of coherent sheaves on $X'$,
\begin{equation}\label{eq: Q diagram}
\begin{tikzcd}
\cE_1 \ar[rr, "b_{12}"] \ar[dr, hook] & & \sigma^* \cE_2^* \ar[rr, "\pi_2"] \ar[dr, hook, "i_{2}"] & &  Q_{2}  \ar[dr, "\io_{2}"] \\
& \cG \ar[ur] \ar[rr] \ar[dr]  & &  \cG^{\sh}\ar[rr] && Q\\
\cE_2 \ar[rr, "b_{21}"] \ar[ur, hook] & & \sigma^* \cE_1^* \ar[rr, "\pi_1"] \ar[ur, hook, "i_{1}"'] & & Q_{1}\ar[ur, "\io_{1}"'] 
\end{tikzcd}
\end{equation}

To explain the maps in this diagram:
\begin{itemize}
\item The diagonal arrows in and out of $\cG$ are as in \eqref{eq: Lagrangian SES}.
\item The horizontal sequences are short exact by definition.
\item The maps $i_{1}$ and $i_{2}$ are the inclusions as a summand by the definition of $\cG^{\sh}$, and $\io_{j}$ is induced from $i_j$. 
\end{itemize}

\begin{lemma}\label{lem: Q diagram} Maintaining our assumption that the Lagrangians $\cE_1, \cE_2 \subset \cG$ are transverse, then both $\io_{1}:Q_{1}\to Q$ and $\io_{2}:Q_{2}\to Q$ are isomorphisms.

%
%
\end{lemma}

\begin{proof} The maps connecting the first and second rows of \eqref{eq: Q diagram} give a map of short exact sequences. By definition $\cG^{\sh}$ is the pushout of $\cG$ along $b_{12}$ (and also along $b_{21}$), so the cokernels of the hook arrows $\cE_{1}\incl \cG$ and $\s^{*}\cE^{*}_{2}\incl \cG^{\sh}$ are both identified with $\s^{*}\cE^{*}_{1}$, hence we conclude that $\io_{2}$ is an isomorphism. The same argument applied to the second and third rows of \eqref{eq: Q diagram}  shows that $\io_{1}$ is an isomorphism.
\end{proof}

Therefore, composing the isomorphisms in Lemma \ref{lem: Q diagram} with either $\b_{12}$ or $\b_{21}$, we get isomorphisms
\begin{eqnarray}
\label{eq: h12} 
h_{12}: Q\xr{\io_{1}^{-1}}Q_{1}\xr{\b_{12}}\s^{*}Q_{2}^{*}\xr{\s^{*}(\io_{2}^{-1})^{\vee}}\s^{*}Q^{*}
\\
\label{eq: h21} 
h_{21}: Q\xr{\io_{2}^{-1}}Q_{2}\xr{\b_{21}}\s^{*}Q_{1}^{*}\xr{\s^{*}(\io_{1}^{-1})^{\vee}}\s^{*}Q^{*}.
\end{eqnarray}

\begin{lemma}\label{l:h12} Both $h_{12}$ and $h_{21}$ give Hermitian structures to $Q$, and 
\begin{equation}
h_{12}=-h_{21}.
\end{equation}
\end{lemma}
\begin{proof} 

Let
\begin{equation}
c_{12}: Q\times \s^{*}Q\to \un F'/\cO_{X'}
\end{equation}
be the pairing induced by $h_{12}$. Similarly define $c_{21}$. Then for local sections $s$ and $s'$ of $Q$, we have 
\begin{eqnarray}
c_{12}(s,s')=\g_{12}(\io_{1}^{-1}(s), \s^{*}(\io_{2}^{-1}(s'))),\\
c_{21}(s,s')=\g_{21}(\io_{2}^{-1}(s), \s^{*}(\io_{1}^{-1}(s'))).
\end{eqnarray}
The equality $h_{12}=-h_{21}$ is equivalent to $c_{12}=-c_{21}$, which is equivalent to
\begin{equation}\label{c1221}
\g_{12}(\io_{1}^{-1}(s), \s^{*}(\io_{2}^{-1}(s')))+\g_{21}(\io_{2}^{-1}(s), \s^{*}(\io_{1}^{-1}(s'))) = 0.
\end{equation}

Note $\cG^{\sh}\surj Q_{1}\op Q_{2}$. The rational pairing \eqref{rat pairing h} induces a pairing 
\begin{equation}
\g: (Q_{1}\op Q_{2})\times (\s^{*}Q_{1}\op \s^{*}Q_{2})\to \un F'/\cO_{X'}
\end{equation}
whose restriction to $Q_{i}\times \s^{*}Q_{i}$ is zero, and whose restriction to $Q_{1}\times \s^{*}Q_{2}$ (resp. $Q_{2}\times\s^{*}Q_{1}$) is $\g_{12}$ (resp. $\g_{21}$). Now the image $\ov\cG$ of $\cG$ in $Q_{1}\op Q_{2}$ is isotropic under the above pairing, and both projections $\ov \cG\to Q_{i}$ are isomorphisms. Therefore $\ov\cG$ is the graph of a unique isomorphism $\ph: Q_{1}\isom Q_{2}$. The fact that $\ov\cG$ is isotropic implies that for local sections $s_{1}, s'_{1}$ of $Q_{1}$, we have
$\g(s_{1}+\ph(s_{1}), \s^{*}s'_{1}+\s^{*}(\ph(s'_{1})))=0$, i.e.,
\begin{equation}\label{g1221}
\g_{12}(s_{1}, \s^{*}\ph(s'_{1}))+\g_{21}(\ph(s_{1}), \s^{*}s'_{1})=0.
\end{equation}
Note that $Q=(Q_{1}\op Q_{2})/\ov\cG$, hence 
\begin{equation}
\io^{-1}_{2}\c\io_{1}=-\ph: Q_{1}\isom Q_{2}.
\end{equation}
We can rewrite \eqref{g1221} as
\begin{equation}
\g_{12}(s_{1}, \s^{*}(\io_{2}^{-1}\io_{1}(s'_{1})))+\g_{21}(\io_{2}^{-1}\io_{1}(s_{1}), \s^{*}s'_{1})=0,
\end{equation}
which confirms \eqref{c1221} (letting $s=\io_{1}(s_{1}), s'=\io_{1}(s'_{1})$). This finishes the proof.
\end{proof}


\subsubsection{Self-duality of $\Hom(\cF^*, \s^{*}Q)$}\label{sssec: r=0 step 3}
Let
\begin{equation}
V:=\Hom(\cF^*,\s^{*}Q).
\end{equation}
The torsion sheaf $\cHom(\cF^{*}, \s^{*}Q)=\cF\ot \s^{*}Q$ carries a $\om_{X'}\ot \un F'/\cO_{X'}$-valued Hermitian form that is the tensor product of $h_{\cF}$ on $\cF$ and $\s^{*}c_{12}$ on $\s^{*}Q$:
\begin{equation}
h_{\cF}\ot c_{12}: \cHom(\cF^{*}, \s^{*}Q)\times \s^{*}\cHom(\cF^{*}, \s^{*}Q)\to \om_{X'}\ot \un F'/\cO_{X'}=\om_{F'}/\om_{X'}.
\end{equation}
Here, $\om_{F'}=\om_{X'}\ot_{\cO_{X'}}\un F'$ is the sheaf of rational differentials on $X'$. Taking global sections and applying the residue map, this gives a perfect symmetric $\F_{q}$-bilinear pairing
\begin{equation}
\tw{-, -}_{12} \co V\times V \rightarrow \F_q,
\end{equation}
which induces a quadratic form 
\[
\frq_{12} \co V \rightarrow \F_q .
\]

Concretely,  for $s: \cF^{*}\to \s^{*}Q$, we have the composition
\begin{equation}\label{expand quad s}
\xymatrix{\cF\ot\om_{X'}^{-1}\ar[r]^{h_{\cF}} & \s^{*}\cF^{*}\ar[r]^{\s^{*}s} & Q \ar[r]^{h_{12}}& \s^{*}Q^{*}\ar[r]^{s^{\vee}} & \cF[1]}
\end{equation}
giving an element in $\Ext^{1}(\cF\ot\om_{X'}^{-1},\cF)$. Then $\frq_{12}(s)$ is the image of this element under the trace map $\Ext^{1}(\cF\ot\om_{X'}^{-1},\cF)\to H^{1}(X',\om_{X'})=\F_{q}$.

Similarly, using $h_{21}$ on $\s^{*}Q$ instead of $h_{12}$, we obtain a perfect symmetric $\F_{q}$-bilinear pairing $\tw{-, -}_{21}$ on $\Hom(\cF^*, \s^{*}Q)$ and a quadratic form $\frq_{21}$. By Lemma \ref{l:h12}, we have 
\begin{equation}\label{q1221}
\frq_{12}=-\frq_{21}.
\end{equation}

\begin{lemma}\label{l:ts}
For $t\in \Hom(\cE_{1}, \cF)$, let $s\in V=\Hom(\cF^{*}, \s^{*}Q)$ be the composition 
$$\cF^{*}\xr{t^{\vee}}\cE_{1}^{*}\surj \s^{*}Q_{1}\xr{\s^{*}\io_{1}}\s^{*}Q.$$ 
Then we have
\begin{equation}
\j{e_{\cG,\cE_{1}}, a(t)}=\frq_{21}(s).
\end{equation}
\end{lemma}
\begin{proof} Let $\pi_{i}: \s^{*}\cE_{i}^{*}\surj Q_{i}$ be the projection. 
Let $e_{1}\in \Ext^{1}(Q_{2},\cE_{1})$ be the class of the top row of \eqref{eq: Q diagram}. Since the map $\cG\to \s^{*}\cE^{*}_{1}$ has a section over the subsheaf $\cE_{2}\incl \s^{*}\cE^{*}_{1}$, $e_{\cG,\cE_{1}}$ is the image of $e_{1}$ under the map
\begin{equation}
\Ext^{1}(Q_{2}, \cE_{1})\cong \Ext^{1}(Q_{1}, \cE_{1})\to \Ext^{1}(\s^{*}\cE^{*}_{1}, \cE_{1})
\end{equation}
induced by the projection $\pi_1 \co \s^{*}\cE^{*}_{1}\surj Q_{1}$ and the isomorphism $-\io_{2}^{-1}\io_{1}:Q_{1}\isom Q_{2}$. In other words, $e_{\cG,\cE_{1}}$ is the composition
\begin{equation}\label{eGE1}
\xymatrix{\s^{*}\cE^{*}_{1}\ar[r]^{\pi_1} & Q_{1}\ar[r]^-{\io_{1}}_-{\sim} & Q\ar[r]^{-\io_{2}^{-1}}_{\sim} & Q_{2}\ar[r]^-{e_{1}} & \cE_{1}[1]}
\end{equation}
Recall we have a commutative diagram of short exact sequences
\begin{equation}\label{map of extn Q}
\xymatrix{ 0\ar[r] & \cE_{1}\ar@{=}[d] \ar[r]^{b_{12}} & \s^{*}\cE_{2}^{*}\ar@{=}[d]\ar[r]^{\pi_2} & Q_{2}\ar[r]\ar[d]^{\b_{21}} & 0\\
0\ar[r] & \cE_{1} \ar[r]^{-\s^{*}b^{\vee}_{21}} & \s^{*}\cE_{2}^{*}\ar[r] & \s^{*}Q^{*}_{1}\ar[r] & 0
}
\end{equation}
This implies that $e_{1}: Q_{2}\to \cE_{1}[1]$ can be written as a composition
\begin{equation}
Q_{2}\xr{\b_{21}}\s^{*}Q_{1}^{*}\xr{-\s^{*}\pi^{\vee}_{1}}\cE_{1}[1].
\end{equation}
Using this and \eqref{eGE1}, we can rewrite $e_{\cG,\cE_{1}}$ as the composition
\begin{equation}
\xymatrix{\s^{*}\cE^{*}_{1}\ar[r]^{\pi_1} & Q_{1}\ar[r]^-{\io_{1}}_-{\sim} & Q\ar[r]^{-\io_{2}^{-1}}_{\sim} & Q_{2} \ar[r]^{\b_{21}} & \s^{*}Q^{*}_{1}\ar[r]^-{-\s^{*}\pi^{\vee}_{1}} & \cE_{1}[1]}
\end{equation}
Using the definition of $h_{21}$ in \eqref{eq: h21}, we see this is the composition
\begin{equation}
\xymatrix{\s^{*}\cE^{*}_{1}\ar[r]^{\io_{1}\pi_{1}}& Q\ar[r]^{h_{21}} & \s^{*}Q^{*} \ar[rr]^-{\s^{*}(\io_{1}\pi_{1})^{\vee}} && \cE_{1}[1]}
\end{equation}
Therefore $\j{e_{\cG,\cE_{1}}, a(t)}$ is the trace of 
\begin{equation}\label{eGE pair a}
\xymatrix{ \cF\ot\om_{X'}^{-1}\ar[r]^{h_{\cF}} & \s^{*}\cF^{*}\ar[r]^{\s^{*}t^{\vee}} & \s^{*}\cE^{*}_{1}\ar[r]^-{\io_{1}\pi_{1}} & Q\ar[r]^{h_{21}} & \s^{*}Q^{*} \ar[rr]^-{\s^{*}(\io_{1}\pi_{1})^{\vee}} && \cE_{1}[1]\ar[r]^-{t} & \cF[1]}
\end{equation}
Using the definition of $s$, we can rewrite \eqref{eGE pair a} as
\begin{equation}
\xymatrix{ \cF\ot\om_{X'}^{-1}\ar[r]^{h_{\cF}} & \s^{*}\cF^{*}\ar[r]^{\s^{*}s}& Q\ar[r]^{h_{21}} & \s^{*}Q^{*}\ar[r]^{s^{\vee}} & \cF[1]}
\end{equation}
which is $\frq_{21}(s)$ by comparing with the analogous concrete formula for $\frq_{21}$ obtained from \eqref{expand quad s} by replacing $h_{12}$ with $h_{21}$. This proves the lemma.
\end{proof}

\begin{remark} Combining Lemma \ref{l:ts} with \eqref{q1221}, we see that
\begin{equation}\label{eGE1 q12}
\j{e_{\cG,\cE_{1}}, a(t)}=-\frq_{12}(s).
\end{equation}
Switching the roles of $\cE_{1}$ and $\cE_{2}$ in Lemma \ref{l:ts}, we have the following formula. Let $t\in \Hom(\cE_{2}, \cF)$, let $s\in \Hom(\cF^{*}, \s^{*}Q)$ be the composition 
$$\cF^{*}\xr{t^{\vee}}\cE_{2}^{*}\surj \s^{*}Q_{2}\xr{\s^{*}\io_{2}}\s^{*}Q.$$ 
Then we have
\begin{equation}\label{eGE2 q21}
\j{e_{\cG,\cE_{2}}, a(t)}=\frq_{12}(s).
\end{equation}
\end{remark}

\subsubsection{Rewriting the theta series}\label{sssec: r=0 step 4}

Denote by  
\[
f \co \Hom(\cF^*, \cE_1^*) \rightarrow \Hom(\cF^*, \s^{*}Q)=V
\]
the map induced by the quotient $\s^{*}(\io_{1}\c\pi_{1}): \cE_1^* \surj \s^{*}Q_{1}\isom \s^{*}Q$. Using \eqref{eGE1 q12}, we may rewrite $\wt Z^{0}_{m}(\cG,  \cE_1)_{\cF}$ as
\begin{eqnarray}\label{eq: r=0 Z_1 rewrite}
\wt Z^{0}_{m}(\cG, \cE_1)_{\cF} 
&=& \chi(\det \cE_1) q^{n(\deg \cE_1  - \deg \omega_X)/2} \sum_{t\in \Hom(\cE_{1},\cF)} \psi(\j{e_{\cG,\cE_{1}}, a(t)}) \nonumber \\
&= & \chi(\det \cE_1) q^{n(\deg \cE_1  - \deg \omega_X)/2} \sum_{s \in V} \psi(\frq_{21}(s)) (f_!\bbm{1}_{\Hom(\cF^{*}, \cE_{1}^{*}) })(s)
\end{eqnarray}
Here $ \bbm{1}_{\Hom(\cF^{*},\cE_{1}^{*}) }$ is the constant function with value $1$ on the set $\Hom(\cF^{*},\cE_{1}^{*})$, so that for $s \in V$, $(f_!\bbm{1}_{\Hom(\cF^{*},\cE_{1}^{*}) })(s)$ is the number of maps $\cE_{1}\rightarrow \cF$ lying in the fiber over $s$. For functions $\varphi_1, \varphi_2$ on $V$, we denote
\[
\langle \varphi_1, \varphi_2 \rangle_V := \sum_{v \in V} \varphi_1(v) \varphi_2(v)
\]
so that \eqref{eq: r=0 Z_1 rewrite} becomes 
\begin{equation}
\wt Z^{0}_{m}(\cG, \cE_1)_{\cF}  = \chi(\det \cE_1) q^{n(\deg \cE_1  - \deg \omega_X)/2}  \j{\frq_{21}^{*}\psi, f_!\bbm{1}_{\Hom(\cF^{*}, \cE_{1}^{*}) }}_{V}.
\end{equation}

\subsubsection{More dualities}\label{sssec: r=0 step 2}
Applying $\s^{*}$ to the bottom row of \eqref{eq: Q diagram}, and using $\io_{1}$ to identify $Q_{1}$ with $Q$, we get a short exact sequence 
\begin{equation}\label{eq: E_1 E_2 SES 2}
0 \rightarrow \s^{*}\cE_2 \rightarrow \cE_1^{*} \rightarrow \s^{*}Q \rightarrow 0
\end{equation}
which induces a 5-term exact sequence 
\begin{equation}\label{eq: 5term 2}
\begin{gathered}
\begin{adjustbox}{max width=\textwidth}
$\displaystyle
\begin{tikzcd}[ampersand replacement=\&]
 \Hom(\cF^*, \s^{*}\cE_2) \ar[r] \& \Hom(\cF^*, \cE_1^*) \ar[r] \&  V \ar[r] \& \Ext^1(\cF^*, \s^{*}\cE_2) \ar[r] \& \Ext^1(\cF^*, \cE_1^{*}) .
 \end{tikzcd}
$
\end{adjustbox}
\end{gathered}
\end{equation}
From the top row of \eqref{eq: Q diagram} we get another short exact sequence
\begin{equation}\label{eq: E_1 E_2 SES 1}
0 \rightarrow \s^{*}\cE_1 \rightarrow  \cE_2^{*} \rightarrow \s^{*}Q \rightarrow 0
\end{equation}
 which induces a 5-term exact sequence 
\begin{equation}\label{eq: 5term 1}
\begin{gathered}
\begin{adjustbox}{max width=\textwidth}
$\displaystyle
\begin{tikzcd}[ampersand replacement=\&]
 \Hom(\cF^*, \s^{*}\cE_1) \ar[r] \& \Hom(\cF^*,  \cE_2^*) \ar[r] \&  V \ar[r] \& \Ext^1(\cF^*, \s^{*}\cE_1) \ar[r] \& \Ext^1(\cF^*,  \cE_2^{*}) .
 \end{tikzcd}
$
\end{adjustbox}
\end{gathered}
 \end{equation}


Serre duality exhibits certain dualities between the terms of \eqref{eq: 5term 1} and \eqref{eq: 5term 2}, as indicated in the diagram below: 
\begin{equation}\label{eq: comparison of 5term}
\begin{gathered}
\begin{adjustbox}{max width=\textwidth}
$\displaystyle
\begin{tikzcd}[ampersand replacement=\&]
 \Hom(\cF^*, \s^{*}\cE_2) \ar[r]  \ar[drrrr, dotted] \&  \Hom(\cF^*,  \cE_1^*) \ar[r]  \ar[drr, dotted] \&  V \ar[r] \ar[d, dotted]  \& \Ext^1(\cF^*, \s^{*}\cE_2) \ar[r]  \ar[dll, dotted] \& \Ext^1(\cF^*,\cE_1^*)   \ar[dllll, dotted] \\
 \Hom(\cF^*, \s^{*}\cE_1) \ar[r]  \ar[urrrr, dotted] \& \Hom(\cF^*,   \cE_2^*) \ar[r] \ar[urr, dotted] \&  V \ar[r]  \ar[u, dotted] \& \Ext^1(\cF^*, \s^{*}\cE_1) \ar[r] \ar[ull, dotted] \& \Ext^1(\cF^*,  \cE_2^*) \ar[ullll, dotted] 
\end{tikzcd}
$
\end{adjustbox}
\end{gathered}
\end{equation}
where the dotted arrows connect spaces that are dual. For example, the bottom left term $\Hom(\cF^*, \s^{*}\cE_1) \cong \Hom(\s^{*}\cE_1^*, \cF)$ is dual (via Serre duality) to 
\[
\Ext^1(\cF, \s^{*}\cE_1^* \otimes \omega_{X'}) \stackrel{h_{\cF}}{\cong} \Ext^1 (\sigma^* \cF^{\vee}, \s^{*}\cE_1^* \otimes \omega_{X'}) \cong  \Ext^1(\cF^*,  \cE_1^*),
\]
which is the top right term. The self-duality of the middle term $V$ has been explained in \S\ref{sssec: r=0 step 3}.

\begin{lemma}\label{l: duality up to sign}
Under either the pairing $\j{-,-}_{12}$ or $\j{-,-}_{21}$ on $V$ and the Serre duality pairings, the sequence of maps in the first row of \eqref{eq: comparison of 5term} is dual to the sequence of maps in the second row, up to sign.
\end{lemma}
\begin{proof}
This follows from applying $\RHom(\cF^{*}, -)$ to the diagram \eqref{map of extn Q}.
\end{proof}

\subsubsection{Fourier transform over finite fields}\label{sssec: finite FT} Next we consider the finite Fourier transform on $\Hom(\cF^*, \s^*Q)$. First we will have to set up our normalizations. 

Temporarily in this section, we let  $V$ denote a vector space of dimension $r$ over the finite field $k=\F_{q}$, and $\wh{V}$ the dual vector space over $k$. Recall $\psi$ is a nontrivial additive character of $k$. We define the Fourier transform of a function $\varphi$ on $V$ by the formula 
\[
\FT_V(\varphi)(\wh{v}) := (-1)^r \sum_{v \in V} \varphi(v) \psi \langle v, \wh{v} \rangle, \quad \wh{v} \in \wh{V}.
\]
This definition is compatible under the sheaf-function correspondence with the sheaf-theoretic Fourier transform to be defined in \S \ref{sec: FT}. With this normalization, we have the following properties of the Fourier transform: 
\begin{itemize}
\item (Involutivity) $\FT_{\wh{V}} \circ \FT_V (\varphi) = q^r [-1]^* \varphi$, where $[-1]^* \varphi(v) = \varphi(-v)$. 
\item (Plancherel formula)
\begin{equation}\label{eq: plancherel for functions}
q^r \sum_{v \in V} \varphi_1(v)  \varphi_2(v)  =  \sum_{\wh{v} \in \wh{V}} \FT_V([-1]^*\varphi_1) (\wh{v}) \FT_V(\varphi_2)(\wh{v}). 
\end{equation}
\item (Gaussians) Suppose we have an isomorphism $h \co V \xrightarrow{\sim} \wh V$ satisfying $\wh h=h$. 
This induces a quadratic form $\frq \co V \rightarrow \F_q$ given by $\frq(v) = \langle v, h(v) \rangle $ and a quadratic form $\wh{\frq} \co \wh{V} \rightarrow \F_q$ given by $\wh{\frq}(\wh{v}) = \langle h^{-1}(\wh{v}), \wh{v} \rangle$. Then we have
\begin{equation}\label{eq: FT of Gaussian function}
 \FT_V(\frq^* \psi ) = (-1)^{r}G(V,\frq)(-\frac{1}{4}\wh{\frq})^* \psi,
\end{equation}
where $G(V,\frq)$ is the Gauss sum
\begin{equation}
G(V,\frq)=\sum_{v \in V} \psi(\frq(v)).
\end{equation}

\item Let $f \co V' \rightarrow V$ be a linear map between vector spaces of ranks $r'$ and $r$ respectively. This induces a morphism $\wh{f} \co \wh{V} \rightarrow \wh{V}'$ of dual spaces. Then we have 
\begin{equation}\label{FT pushpull}
\FT_V (f_! \varphi') = (-1)^{r-r'} \wh{f}^* \FT_{V'} (\varphi')
\end{equation}
for all functions $\varphi'$ on $V'$, and 
\begin{equation}\label{FT pullpush}
\FT_{V'} (f^* \varphi) = (-1)^{r'-r} q^{r'-r} \wh{f}_! \FT_V (\varphi)
\end{equation}
for all functions $\varphi$ on $V$. 
\end{itemize}

\begin{example}
Let $\bbm{1}_V$ be the constant function on $V$ with value $1$ and $\delta_V$ be the delta function on $V$ with value $1$ at the origin. If $\dim_k (V) = r$ then we have 
\[
\FT (\bbm{1}_V) = (-1)^r q^r \delta_{\wh{V}}
\]
and 
\[
\FT(\delta_V) = (-1)^r \bbm{1}_{\wh{V}}.
\]
\end{example}

\subsubsection{The Gauss sum}\label{sssec: Gauss sum}
For a non-degenerate quadratic space $(V, \frq_{V})$ over $\F_{q}$, the normalized Gauss sum (with respect to the fixed additive character $\psi$ on $\F_{q}$) is defined as
\begin{equation}
\g(V,\frq_{V})=q^{-\dim V/2}\sum_{v\in V}\psi(\frq_{V}(v)).
\end{equation}
(Here $q^{-\dim V/2}$ means the positive square root of $q^{- \dim V}$.) We have the following well-known (at least as far back as Weil \cite[\S II]{Weil64}) facts:
\begin{enumerate}
\item $\g(V,\frq_{V})$ is a fourth root of unity in $\Qlbar$. 
\item\label{gs add} The function $(V,\frq_{V})\mapsto \g(V,\frq_{V})$ is additive, hence induces a homomorphism from the Witt group $\g: \mathrm{Witt}(\F_{q})\to \mu_{4}(\Qlbar)$.
\item\label{gs H+} For the hyperbolic plane $H_{+}:=\F_{q}^{2}$ and $\frq_{+}(x,y)=xy$, we have $\g(H_{+},\frq_{+})=1$. From this one deduces that if the quadratic form $\frq_{V}$ is split (i.e., there exists a Lagrangian subspace), then $\g(V,\frq_{V})=1$.
\item\label{gs H-} For $H_{-}:=\F_{q^{2}}$ and $\frq_{-}(x)=\Nm_{\F_{q^{2}}/\F_{q}}(x)$, we have $\g(H_{-}, \frq_{-})=-1$. From this one deduces that if $V$ is a $\F_{q^{2}}$-vector space with a nondegenerate $\F_{q^{2}}/\F_{q}$-Hermitian form $(-,-)$, and $\frq_{V}(x)=(x,x)$, then $\g(V,\frq_{V})=(-1)^{\dim_{\F_{q^{2}}}V}$.
\item\label{gs change field} Let $k'$ be a finite extension of $\F_{q}$, and $(V',\frq_{V'})$ be a quadratic space over $k'$. On one hand, we can define the normalized Gauss sum $\g_{k'}(V',\frq_{V'})$ using the additive character $\psi_{k'}=\psi\c \Tr_{k'/\F_{q}}$. On the other hand, we can view $V'$ as a vector space over $\F_{q}$ equipped with the quadratic form $\un \frq_{V'}:=\Tr_{k'/\F_{q}}\c \frq_{V'}$. Then we have $\g_{k'}(V',\frq_{V'})=\g(V', \un \frq_{V'})$. 
\end{enumerate}

Now we are back to the convention that $V=\Hom(\cF^{*}, \s^*Q)$. We want to compute $\g(V,\frq_{12})$.  Define the divisor $D'_{Q}$ on $X'$ to be the $\sum_{x'} d_{x'}x'$ where $d_{x'}$ is the length of $Q$ at $x'$. Since $Q$ carries a Hermitian form, we have $d_{x'}=d_{\s(x')}$, therefore $D'_{Q}=\nu^{*}D_{Q}$ for a unique divisor $D_{Q}$ on $X$.

\begin{lemma}\label{lem: sign of Gauss sum} We have
\[
\g(\Hom(\cF^{*}, \s^{*}Q),\frq_{12})=\g(\Hom(\cF^{*}, \s^{*}Q),\frq_{21})=\y_{F'/F}(D_{Q})^{n}.
\]
\end{lemma}
\begin{proof}
Decompose $Q=\op Q_{x}$ where $Q_{x}$ is the summand supported over a place $x\in |X|$. Then $(V,\frq_{12}):=(\Hom(\cF^{*}, \s^{*}Q),\frq_{12})$ is the orthogonal direct sum of quadratic spaces $(V_{x},\frq_{x})$ where $V_{x}=\Hom(\cF^{*}, \s^{*}Q_{x})$ and $\frq_{12,x}=\frq_{12}|_{V_{x}}$. By observation \eqref{gs add} above,  it suffices to consider the case $Q=Q_{x}$ for some $x\in |X|$.

If $x$ splits into $x'$ and $x''=\s(x')$ in $X'$, then we can write $Q_{x}=Q_{x'}\op Q_{x''}$ according to the support, and $\s^*Q_x=\s^*Q_{x'}\op\s^*Q_{x''}$. Relabeling the two summands by their supports, we get $V=V_{x'}\op V_{x''}$, where $V_{x'}=\Hom(\cF^{*}, \s^*Q_{x''})$ and $V_{x''}=\Hom(\cF^{*}, \s^*Q_{x'})$, so that $V_{x'}$ and $V_{x''}$ are both isotropic. By observation \eqref{gs H+} above, $\g(V,\frq_{12})=1$ in this case. On the other hand $\y_{F'/F}(D_{Q})=1$, hence $\y_{F'/F}(D_{Q})^{n}=1$.

Now we consider the case $x$ is inert in $X'$, and let $x'$ be the unique place above $x$.  By observation \eqref{gs change field}, we may rename the base fields, and thereby assume $k(x)=\F_{q}$ and $k(x')=\F_{q^{2}}$.  Recall the quadratic form $\frq_{12}$ on $V$ comes from a $\F_{q^{2}}/\F_{q}$ Hermitian form $(-,-)$ on $V$ by taking $\frq_{12}(v)=(v,v)$. By observation \eqref{gs H-}, we have
\[
\g(V,\frq_{12})=(-1)^{\dim_{\F_{q^{2}}}V}.
\]
Note that $\dim_{\F_{q^{2}}}V=nd$ where $d$ is the multiplicity of $D_{Q}$. Therefore
\[
\g(V,\frq_{12})=(-1)^{nd}=\y_{F'/F}(D_{Q})^{n}.
\]
The same argument works for $\frq_{21}$ in place of $\frq_{12}$.
\end{proof}

\subsubsection{Conclusion of the argument}\label{sssec: conclusion r=0}
We return to the notation $V=\Hom(\cF^*, \s^{*}Q)$. At the end of \S \ref{sssec: r=0 step 4}, we expressed
\begin{equation}\label{eq: pairing}
\wt Z^{0}_{m}(\cG, \cE_1)_{\cF}  = \chi(\det \cE_1) q^{n (\deg \cE_1  - \deg \omega_X)/2} \langle \frq_{21}^* \psi, f_!\bbm{1}_{\Hom(\cF^{*},\cE_{1}^{*})}  \rangle_{V}.
\end{equation}

By the Plancherel formula \eqref{eq: plancherel for functions}, we have 
\[
\langle \frq_{21}^* \psi, f_!\bbm{1}_{\Hom(\cF^{*},\cE^{*}_{1})}  \rangle_{V} = \frac{1}{q^{\dim V}} \langle [-1]^* \FT(\frq_{21}^* \psi), \FT(f_!\bbm{1}_{\Hom(\cF^{*},\cE_{1}^{*})}) \rangle_{V}.
\]

So we turn to analyze $\FT(\frq_{21}^* \psi)$ and $ \FT(f_!\bbm{1}_{\Hom(\cF^*,  \cE_1^*)})$. Below, we abbreviate $\hom(A,B) := \dim_{\F_q} \Hom(A,B)$ and $\ext^1(A,B) := \dim_{\F_q} \Ext^1(A,B)$.
\begin{itemize}
\item By \eqref{eq: FT of Gaussian function}, we have
\[
[-1]_* \FT(\frq_{21}^* \psi) = (-1)^{\dim V} G(V, \frq_{21})\cdot (-\frac{1}{4}\wh{\frq})^* \psi.
\]
As was explained in \S \ref{sssec: Gauss sum}, especially Lemma \ref{lem: sign of Gauss sum}, we have 
\[
G(V, \frq_{21})=q^{\dim V/2}\g(V, \frq_{21})=\eta_{F'/F}(D_Q)^n q^{\dim V/2}.
\]
On the other hand, if we use $\frq_{21}$ to identify $V$ with $\wh V$, then $\wh \frq_{21}=\frq_{21}$. By \eqref{q1221}, we have $-\frq_{21}=\frq_{12}$. Therefore, under this identification, we have
\begin{equation}
[-1]_* \FT(\frq_{21}^* \psi) = (-1)^{\dim V} \eta_{F'/F}(D_Q)^n q^{\dim V/2}\cdot (\frac{1}{4}\frq_{12})^* \psi.
\end{equation}

\item To analyze the Fourier transform of $f_! \bbm{1}_{\Hom(\cF^{*},\cE_{1}^{*}) }$, we have a long exact sequence
\begin{equation}\label{eq: boundary g}
0\to \Hom(\cF^{*}, \s^{*}\cE_{2})\to \Hom(\cF^{*}, \cE_{1}^{*})\xr{f} \Hom(\cF^{*},\s^{*}Q)\xr{g} \Ext^{1}(\cF^{*}, \s^{*}\cE_{2})
\end{equation}
coming from applying $\RHom(\cF^{*}, -)$ to the short exact sequence $\s^{*}\cE_{2}\to \cE_{1}^{*}\to \s^{*}Q$.  From this we get
\begin{equation}\label{eq: descent for functions}
f_! \bbm{1}_{\Hom(\cF^{*},\cE_{1}^*) } = q^{\dhom(\cF^*,\s^{*}\cE_2)} g^* \delta_{\Ext^1(\cF^*,\s^{*} \cE_2)},
\end{equation}
By \eqref{FT pullpush}, the Fourier transform on $V=\Hom(\cF^*, \s^{*}Q)$ sends  $g^* \delta_{\Ext^1(\cF^*, \s^{*}\cE_2)}$ to
\begin{eqnarray}\label{eq: FT of pullback function}
\FT(g^* \delta_{\Ext^1(\cF^*, \s^{*}\cE_2)}) &=& (-1)^{\dim V} q^{\dim V- \ext^1(\cF^*, \s^{*}\cE_2)} \wh{g}_! \bbm{1}_{\Ext^1(\cF^*, \s^{*}\cE_2)^*}\\
&=& (-1)^{\dim V} q^{\dim V - \ext^1(\cF^*, \s^{*}\cE_2)} \wh{g}_! \bbm{1}_{\Hom(\cF^*, \cE^{*}_{2})}
\end{eqnarray}
\end{itemize}
where the last equality uses the duality between the fourth term of the top row and the second term of the bottom row in the diagram \eqref{eq: comparison of 5term}, and Lemma \ref{l: duality up to sign}.

Putting these equations together and collecting factors yields: 
\begin{align}\label{eq: r=0 1}
\wt{Z}^0_m (\cG, \cE_1)_{\cF}  &= \chi(\det \cE_1) q^{n (\deg \cE_1   - \deg \omega_X)/2} \langle \frq_{21}^* \psi,  f_!\bbm{1}_{\Hom(\cF^*, \cE_1^*)}  \rangle_{V}  \nonumber \\
&=  \chi(\det\cE_1)  \eta_{F'/F}(D_Q)^n  q^{?} \langle (\frac{1}{4}\frq_{12})^* \psi,  \wh{g}_! \bbm{1}_{\Hom(\cF^*, \cE^{*}_2)} \rangle_{V}
\end{align}
where the exponent of $q$ is 
\begin{equation}\label{eq: r=0 2}
? = \frac{n}{2} (\deg \cE_1   - \deg \omega_X)+ \chi(\cF \otimes \s^{*}\cE_2) + \frac{1}{2}\dim V.
\end{equation}
We want to show that this agrees with 
\[
\wt{Z}^0_m ( \cG,  \cE_2)_{\cF}   = \chi(\det \cE_2 )  q^{n (\deg \cE_2   - \deg \omega_X)/2} \langle \frq_{12}^* \psi,  \wh{g}_! \bbm{1}_{\Hom(\cF^*, \cE^{*}_2)} \rangle_{V}
\]
so we separately compare the sign and the exponent of $q$. 

First, we claim that
\begin{equation}
\langle (\frac{1}{4}\frq_{12})^* \psi,  \wh{g}_! \bbm{1}_{\Hom(\cF^*, \cE^{*}_2)} \rangle_V =\langle \frq_{12}^* \psi,  \wh{g}_! \bbm{1}_{\Hom(\cF^*, \cE^{*}_2)} \rangle_V.
\end{equation}
Indeed, the factor $1/4=(1/2)^{2}$ can be eliminated on the left side since $\wh{g}_! \bbm{1}_{\Hom(\cF^*, \cE^{*}_2)}$ is invariant under the scaling action of $\F_{q}^{\times}$ on $\Hom(\cF^*, Q^{*})$.

Next we compare the exponents of $q$. We expand \eqref{eq: r=0 2} as 
\begin{equation}\label{eq: r=0 3}
? = \frac{1}{2}(n \deg \cE_1  - n \deg \omega_X) + n \deg (\cE_2) + m \deg (\cF) -\frac{1}{2} mn \deg \omega_{X'} + \frac{1}{2} \dim V.
\end{equation}
From \eqref{eq: E_1 E_2 SES 2} we have $- \deg \cE_1 = \deg \cE_1^* = \deg \cE_2 + \deg Q$. Also $\dim V=n\deg Q$. Substituting these into \eqref{eq: r=0 3} and simplifying yields 
\begin{equation}\label{eq: r=0 4}
? = \frac{1}{2} (n \deg \cE_2 - n \deg \omega_X) + m \deg \cF - mn \frac{\deg \omega_{X'}}{2}.
\end{equation}
Now, $\cF \xrightarrow{\sim} \s^{*}\cF^{\vee}$ implies that $\deg \cF = \frac{n \deg \omega_{X'}}{2}$, so \eqref{eq: r=0 4} equals $\frac{1}{2} (n \deg \cE_2 - n \deg \omega_X)$, as desired.

Next we compare the signs. From the top row of \eqref{eq: Q diagram} we have 
\begin{equation}\label{eq: r=0 5}
\chi(\det \cE_1) \chi(\nu^{*}D_Q) = \chi(\det \sigma^* \cE_2^*).
\end{equation}
Since $\chi\c\nu^{*}=\y^{n}$, we have $\chi(\nu^{*}D_Q)  = \eta_{F'/F}(D_Q)^n$, and we may rewrite \eqref{eq: r=0 5} as 
\begin{equation}\label{eq: r=0 6}
\chi(\det \cE_1)  \eta_{F'/F}(D_Q)^n  = \chi(\det \s^{*}\cE^{*}_2)=\chi(\det\s^{*}\cE_{2})^{-1}.
\end{equation}
Finally, we note that 
\[
\chi(\det \cE_2) \chi(\det \sigma^* \cE_2) = \chi(\nu^{*}\Nm \det \cE_2) = \eta(\Nm \det \cE_2)^n = 1,
\]
so that \eqref{eq: r=0 6} agrees with $\chi(\det\cE_2)$, as desired. 

\qed

\subsection[Outline of the proof of transverse modularity]{Outline of the proof of Theorem \ref{thm: transverse modularity}}\label{ssec: outline}
We now give a brief summary of the proof of Theorem \ref{thm: transverse modularity}. It is a generalization of the $r=0$ case in \S \ref{ssec: r=0 transverse lagrangians}, but the steps are of course much more complicated. We also take the opportunity to indicate, in more detail than in the Introduction, the role of the individual sections of this paper in this strategy. 

Choose transverse Lagrangians $\cE_1, \cE_2 \inj \cG$. Let $Q$ be as above. We want to show that 
\[
\wt{Z}^r_m(\cG, \cE_1) = \wt{Z}^r_m(\cG, \cE_2).
\]
\begin{enumerate}
\item We will realize the cycle $[\cZ_{\cE_1}^r]$, whose summands comprise the Fourier coefficients of $\wt{Z}^r_m(\cG, \cE_1)$, as the \emph{trace} of a \emph{cohomological correspondence}\footnote{These notions will be reviewed later, in \S \ref{sec: sheaf-cycle}.} ``$\cc_U$'', where $U$ is a derived vector bundle over $\Bun_{U(n)}$ geometrizing the $\F_q$-vector space $\Hom(\cF^*, \cE_1^*)$ which appeared in  \S \ref{sssec: r=0 step 2}. 

Similarly, we will realize the cycle $[\cZ_{\cE_2}^r]$, whose summands comprise the Fourier coefficients of $\wt{Z}^r_m(\cG, \cE_2)$, as the trace of a cohomological correspondence ``$\cc_{U^\perp}$'', where $U^\perp$ is a derived vector bundle over $\Bun_{U(n)}$ geometrizing the vector space $\Hom(\cF^*, \cE_2^*)$ which appeared in  \S \ref{sssec: r=0 step 2}.

The spaces $U$ and $U^\perp$ are defined in \S \ref{ssec: UVW}, the cohomological correspondences $\cc_U$ and $\cc_{U^\perp}$ are defined in \S \ref{ssec: comparison of cc}, and the computation of the traces is performed in \S \ref{ssec: calculation of traces}, based on general results established in \S \ref{ssec: VFC as trace}. 

\item We will construct maps $f \co U \rightarrow V$ and $f^\perp \co  U^\perp \rightarrow \wh{V}$, where $V$ is a vector bundle over $\Bun_{U(n)}$ geometrizing the $\F_q$-vector space $\Hom(\cF^*, Q)$ which appeared in \S \ref{sssec: r=0 step 2}. This occurs in \S \ref{ssec: UVW}. 

\item We will construct certain vector bundles $W$ and $W^\perp$, geometrizing the $\F_q$-vector spaces $\Ext^1(\cF^*, \sigma^* \cE_2)$ and $\Ext^1(\cF^*, \sigma^* \cE_1)$, respectively, which appeared in \S \ref{sssec: r=0 step 2}. In addition we will construct maps $g \co V \rightarrow W$ and $g^\perp \co \wh{V} \rightarrow W^\perp$ such that $g$ geometrizes the ``$g$'' from \eqref{eq: boundary g}, and $g^\perp$ plays the analogous role with $\cE_1$ and $\cE_2$ interchanged. This all occurs in \S \ref{ssec: UVW} and \S \ref{ssec: UVW dual}.

\item We will prove that a certain pushforward cohomological correspondence ``$f_! \cc_U$'' agrees with a pullback cohomological correspondence ``$g^* \cc_W$'', geometrizing the identity \eqref{eq: descent for functions}. We note that while the identity \eqref{eq: descent for functions} is trivial, its geometrization is highly non-obvious, and occupies the entirety of \S \ref{sec: descent}. 

\item We will prove that the sheaf-theoretic Fourier transform of cohomological correspondences (introduced in \S \ref{ssec: FT of cc}) takes the cohomological correspondence ``$g^* \cc_W$'' to the cohomological correspondence ``$f^\perp_! \cc_{U^\perp}$'' (up to shift and twist), geometrizing \eqref{eq: FT of pullback function}. This occurs in \S \ref{ssec: comparison of cc}, based on general results established in \S \ref{sec: FT cc}. 

At this point the situation is summarized by the following diagram, in which the dotted arrows connect dual vector bundles. 
\[
\begin{tikzcd}
{[\cZ_{\cE_1}^r]} & \cc_U \ar[d] \ar[l, dashed, "\Tr^{\Sht}"'] & U \ar[d, "f"]  \ar[ddrrr, dotted] & & & U^\perp \ar[d, "f^\perp = \wh{g}"] \ar[ddlll, dotted]  & \cc_{U^\perp} \ar[d] \ar[r, dashed, "\Tr^{\Sht}"]  & {[\cZ_{\cE_2}^r]} \\  
& f_! \cc_U & V \ar[d, "g"]  \ar[rrr, dotted] & & & \ar[lll, dotted]  \wh{V}\ar[d, "g^\perp   =\wh{f}"]  & f_!^\perp \cc_{U^\perp} \\
 & & W \ar[rrruu, dotted] & & &  W^\perp \ar[uulll, dotted] 
\end{tikzcd}
\]

\item Forming traces of $f_! \cc_U$ and $f_!^\perp \cc_{U^\perp}$ produces Borel-Moore homology classes of special cycles, generalizing the $f_! \bbm{1}_{\Hom(\cF^*, \cE_1^*)}$ and $\wh{g}_! \bbm{1}_{\Hom(\cF^*, \cE_2^*)}$ in \S \ref{sssec: r=0 step 4}. Therefore, the pairing of $\Tr^{\Sht}(f_! \cc_U)$ with an appropriate Gaussian produces the higher theta series $\wt{Z}^r_m(\cG, \cE_1)$, while the pairing of $\Tr^{\Sht}(f_!^\perp \cc_{U^\perp})$ with an appropriate Gaussian produces the higher theta series $\wt{Z}^r_m(\cG, \cE_2)$. This is detailed in \S \ref{ssec: conclusion}.

\[
\begin{gathered}
\begin{adjustbox}{max width=\textwidth}
$\displaystyle
\begin{tikzcd}[ampersand replacement=\&]
 \& f_! \bbm{1}_{\Hom(\cF^*,  \cE_1^*)}  \ar[rr, dotted] \&    \&\ar[ll, dotted]  \wh{g}_! \bbm{1}_{\Hom(\cF^*,  \cE_2^*)} \\
f_! \cc_U \ar[r, dashed, "\Tr^{\Sht}"] \& {[\cZ_{\cE_1}^r ]} \ar[rr, dotted, "\text{arithmetic Fourier transform}"] \ar[u, rightsquigarrow , "r=0"]  \&  \& {[\cZ_{\cE_2}^r]} \ar[ll, dotted] \ar[u, rightsquigarrow , "r=0"] \&  
f_!^\perp \cc_{U^\perp} \ar[l, dashed, "\Tr^{\Sht}"]  \\
\wt{Z}_m^r(\cG,\cE_1) \ar[r, equals] \& \langle {[\cZ_{\cE_1}^r]}, \text{Gaussian} \rangle   \& \&  \langle {[\cZ_{\cE_2}^r]}, \text{Gaussian} \rangle  \ar[r, equals] \&  \wt{Z}_m^r(\cG, \cE_2) 
\end{tikzcd}
$
\end{adjustbox}
\end{gathered}
\]

\item We introduce in \S \ref{sec: arithmetic FT} an \emph{arithmetic Fourier transform} on Borel-Moore homology classes, generalizing the finite Fourier transform from \S \ref{sssec: finite FT}. The previous steps imply that the arithmetic Fourier transform sends the special cycle $[\cZ_{\cE_1}^r]$ for $\cE_1$ to the special cycle $[\cZ_{\cE_2}^r]$ for $\cE_2$ (up to signs and powers of $q$), generalizing the observation that the finite Fourier transform sends $f_! \bbm{1}_{\Hom(\cF^*, \cE_1^*)}$ to $\wh{g}_! \bbm{1}_{\Hom(\cF^*, \cE_2^*)}$ (up to signs and powers of $q$) from \S \ref{sssec: r=0 step 4}. Noting again that the Gaussian $\frq^* \psi$ is essentially self-dual, we conclude that
\[
\langle {[\cZ_{\cE_1}^r]}, \text{Gaussian} \rangle =  \langle {[\cZ_{\cE_2}^r]}, \text{Gaussian} \rangle 
\]
using a version of the Plancherel formula \eqref{eq: plancherel for functions} for the arithmetic Fourier transform. This is carried out in \S \ref{ssec: conclusion}. 
\end{enumerate}

\part{Generalities on cohomological correspondences}

\section{Base change transformations}\label{sec: zoo}
The purpose of this section is to establish situations in which we can push forward or pull back cohomological correspondences (to be recalled in the next section). This is important for functoriality in the sheaf-cycle correspondence. 

We will start with a commutative square of derived Artin stacks (recall from \S \ref{sssec:notation-on-spaces} that these are automatically assumed to be locally of finite type over a field)
\begin{equation}\label{eq: commutative square}
\begin{tikzcd}
A \ar[r, "g'"] \ar[d, "f'"'] & B \ar[d, "f"] \\
C \ar[r, "g"] & D 
\end{tikzcd}
\end{equation}
We will investigate various situations in which we have base change maps between various combinations of pull/push functors in this diagram, and the compatibilities between them. 

\begin{remark}
Recall that restriction induces an equivalence of \'{e}tale sites of a derived stack and its classical truncation. Therefore, one might ask why we discuss derived stacks at all in this section, since all categories and functors are determined by the underlying classical stack. The answer is that when discussing pullbacks (in \S \ref{ssec: pull-pull}) we will need to invoke the relative fundamental class of a quasi-smooth morphism, which is a derived notion. 
\end{remark}

\subsection{Pushable and pullable squares}
Let $\wt B=C\times_{D}B$ (derived fiber product) such that the square \eqref{eq: commutative square} decomposes into a derived Cartesian square and two triangles
\begin{equation}\label{eq: pushable square}
\xymatrix{ A\ar[dr]^{a}\ar@/_1pc/[ddr]_{f'}\ar@/^1pc/[drr]^{g'} \\
& \wt B\ar[d]^{\wt f}\ar[r]^{\wt g} & B\ar[d]^{f}\\
& C \ar[r]^{g}& D
}
\end{equation} 
so that $g'=\wt g\circ a$ and $f'=\wt f\circ a$. 

\begin{defn}\label{def: pushable pullable} The outer square in \eqref{eq: pushable square} is called 
\begin{itemize}
\item {\em pushable}, if $a$ is proper.
\item {\em pullable}, if $a$ is quasi-smooth. In this case, we call the relative dimension $d(a)$ the {\em defect} of the pullable square $(A,B,C,D)$.
\end{itemize}
\end{defn}

Note that the notion of the outer square of \eqref{eq: pushable square} being pushable or pullable is invariant under flipping the square about the diagonal connecting $A$ and $D$.

\begin{example}\label{ex: pushable}
Here are some examples of pushable squares:
\begin{enumerate}
\item A square whose reduced classical truncation is Cartesian.
\item A square where $f$ and $f'$ are proper.
\item A square where $g'$ is proper and $g$ is separated and representable in derived schemes. 
\end{enumerate}
The pushability of (2) is implied by the pushability of (3) using the flipping symmetry, but we highlight it to make contact with previous constructions in the literature (see Example \ref{ex: pull-push f proper} below). Some special cases of the above examples are observed in \cite[\S 1.1.6]{Var07}. 
\end{example}

\begin{example}\label{ex: pullable}
Here are some examples of pullable squares:
\begin{enumerate}
\item A square whose reduced classical truncation is Cartesian. 
\item A square where $f$ and $f'$ are smooth.
\item A square where $g'$ is quasi-smooth and $g$ is smooth.
\end{enumerate}
Indeed, in the last case, $\wt g:\wt B\to B$ is smooth (as the base change of the smooth morphism $g$). Then because $g'=\wt g\circ a$ is quasi-smooth, we learn that $a$ is quasi-smooth. 

The pullability of (2) is implied by that of (3) using the flipping symmetry, but we highlight it because it will be of special importance. Some special cases of the above examples are observed in \cite[\S A.2]{GV20}.
\end{example}

\begin{remark}
After releasing the first draft of this paper, we learned that the notion of pushability and its significance for pushing forward cohomological correspondences had already been identified in work of Lu-Zheng, \cite[Construction 2.10]{LZ22}. We are not aware that the notion of pullability, which is most useful in the context of derived algebraic geometry, has previously been identified. 
\end{remark}

\subsection{Push-pull}\label{ssec: push-pull}
Referring to the diagram \eqref{eq: pushable square}, there are always base change natural transformations
\begin{equation}\label{eq: push-pull star star}
g^* f_* \rightarrow f'_* (g')^*
\end{equation}
and
\begin{equation}\label{eq: push-pull shriek shriek}
f'_!  (g')^! \rightarrow g^! f_! . 
\end{equation}
by adjunction. If the outer square in \eqref{eq: pushable square} is pushable, i.e., the map $a$ is proper, then we have a natural transformation
\begin{equation}\label{eq: pull-push gen}
g^{*}f_{!}\xr{\di} \wt f_{!}\wt g^{*}\to \wt f_{!}a_{*}a^{*}\wt g^{*}=\wt f_{!}a_{!}a^{*}\wt g^{*}=f'_{!}(g')^{*}
\end{equation}
Here and below we always label the proper base change isomorphism by $\di$. The second map above is the unit map $\Id\to a_{*}a^{*}$, and the next step uses $a_{!}=a_{*}$ because $a$ is proper.

We often denote such natural transformations coming from a pushable square by $\td$:
\begin{equation}
g^{*}f_{!}\xr{\td}f'_{!}(g')^{*}.
\end{equation}

\begin{example}\label{ex: pull-push f proper} In special cases, the map \eqref{eq: pull-push gen} has been observed before with a slightly different description. In particular, if $f$ and $f'$ are proper then $f_* = f_!$ and $f'_* = f'_!$, and Varshavsky \cite[\S 1.1.6]{Var07} observes that the base change map 
\begin{equation}\label{eq: push-pull mixed f proper}
g^* f_* \rightarrow  f'_* (g')^*
\end{equation}
therefore gives another natural transformation $g^* f_! \rightarrow f'_! (g')^*$. Since we will use some of Varshavsky's results concerning his natural transformation, we spell out for completeness why this natural transformation coincides with \eqref{eq: pull-push gen}. Under these assumptions \eqref{eq: pull-push gen} is adjoint to the composite map 
\[
f_*  \xrightarrow{\text{unit}(\wt g)}  f_* \wt{g}_* \wt{g}^* = g_* \wt{f}_* \wt{g}^* \xrightarrow{\text{unit}(a)}  g_* \wt{f}_* a_* a^* \wt{g}^* = g_* f'_* (g')^*
\]
which since $g' = \wt g \circ a$ coincides with 
\[
f_* \xrightarrow{\text{unit}(g')} f_*  (g')_* (g')^* = g_* f'_* (g')^* 
\]
which is adjoint to \eqref{eq: push-pull mixed f proper}.
\end{example}

\subsubsection{Compositions}\label{sssec: push-pull composition}

Suppose we have a commutative diagram 
\begin{equation}\label{eq: commutative composition}
\begin{tikzcd}
A \ar[r, "g''"] \ar[d, "f'"'] & B \ar[d, "f"] \\
C \ar[r, "g'"] \ar[d, "h'"'] & D \ar[d, "h"]  \\
E \ar[r, "g"] & F 
\end{tikzcd}
\end{equation}

\begin{lemma}[2 out of 3 for pushable squares]\label{lem: 2 out of 3} Consider the commutative diagram \eqref{eq: commutative composition}.
\begin{enumerate}
\item If both the upper square and the lower square are pushable, then the outer square formed by $(A,B,E,F)$ is also pushable.
\item Suppose $g$ and $g'$ are separated. If the outer square is pushable, then the upper square is also pushable. 
\end{enumerate}
\end{lemma}
\begin{proof}
Introduce the base changes $\wt D=E\times_{F}D, \wt B_{1}=\wt D\times_{D}B\cong E\times_{F}B$ and $\wt B=C\times_{D}B\cong C\times_{\wt D}\wt B_{1}$. We have a commutative diagram
\begin{equation}\label{eq: comp two pushable}
\xymatrix{ A \ar[r]_{a}\ar[dr]_{f'}\ar@/^1pc/[rrr]^{g''} & \wt B \ar[r]_{b}\ar[d]^{\wt f} & \wt B_{1}\ar[d]^{\wt f_{1}}\ar[r]_{\wt g_{1}} & B\ar[d]^{f}\\
& C\ar[r]_{c}\ar@/^1pc/[rr]^(.3){g'} \ar[dr]_{h'}& \wt D\ar[d]^{\wt h}\ar[r]_{\wt g} & D\ar[d]^{h}\\
& & E \ar[r]^{g}& F
}
\end{equation}
where all squares are Cartesian. 

In situation (1), both $a$ and $c$ are proper, hence $b$ is proper (being a base change of $c$), therefore $b\circ a$ is proper, i.e., the outer square is pushable.

In situation (2), we know that $b\circ a$ is proper. Since $g$ is separated, $\wt{g}$ is separated. Since $g' = \wt{g} \circ c$ is separated, $c$ is separated. Hence $b$ is separated, and $b \circ a $ proper implies $a$ is proper. 
\end{proof}

Assume both the upper square and the lower square are pushable. Then the outer square formed by $(A,B,E,F)$ is also pushable by the above lemma. In this case, we have two maps $g^{*}(h\c f)_{!}\to (h'\c f')_{!}(g'')^{*}$: one given by \eqref{eq: pull-push gen} for the outer pushable square, the other as the composition
\[
g^* h_! f_!   \to h'_! (g')^* f_!\to h'_! f'_! (g'')^*
\]
where both arrows are \eqref{eq: pull-push gen} applied to the upper and the lower squares. The next result says that these two maps are the same.

\begin{prop}\label{p: comp push}
Assume that both the upper square and the lower square in \eqref{eq: commutative composition} are pushable. Then the following diagram is commutative
\begin{equation}\label{eq: comp pull-push}
\xymatrix{g^* h_! f_! \ar@{=}[d] \ar[r]^-{\td f_{!}} & h'_! (g')^* f_! \ar[r]^-{h'_{!}\td} & h'_! f'_! (g'')^*\ar@{=}[d]\\
g^{*}(h\c f)_{!}\ar[rr]^-{\td} &&   (h'\c f')_{!}(g'')^{*}
}
\end{equation} 
\end{prop}
\begin{proof}
Consider the diagram \eqref{eq: comp two pushable}, in which $a,b$ and $c$ are proper.

Consider the diagram of natural transformations
\begin{equation}\label{comp big diagram}
\xymatrix{g^{*}h_{!}f_{!} \ar[dr]^{\di}\ar[r]^{\di}& \wt h_{!}\wt g^{*}f_{!}\ar[d]^{\di}\ar[r]^{u} & \wt h_{!}c_{*}c^{*}\wt g^{*}f_{!}\ar[d]^{\di}\ar@{=}[r] & h'_{!}c^{*}\wt g^{*}f_{!}\ar[d]^{\di}\\
& \wt h_{!}\wt f_{1!}\wt g_{1}^{*}\ar[r]^{u}\ar[d]^{u}\ar @{} [dr] |{\star} & \wt h_{!}c_{*}c^{*}\wt f_{1!}\wt g_{1}^{*}\ar@{=}[r]\ar[d]^{\di} & h'_{!}c^{*}\wt f_{1!}\wt g_{1}^{*}\ar[d]^{\di}\\
& \wt h_{!} \wt f_{1!} b_{*}b^{*}\wt g_{1}^{*}\ar[d]^{u}\ar@{=}[r] & \wt h_{!} c_{!}\wt f_{!} b^{*}\wt g_{1}^{*}\ar[d]^{u}\ar@{=}[r] & h'_{!}\wt f_{!}b^{*}\wt g_{1}^{*}\ar[d]^{u}\\
&  \wt h_{!}\wt f_{1!}b_{*}a_{*}a^{*}b^{*}\wt g_{1}^{*} \ar@{=}[r]
&\wt h_{!} c_{!}\wt f_{!} a_{*}a^{*}b^{*}\wt g_{1}^{*}\ar@{=}[r] & h'_{!}\wt f_{!}a_{*}a^{*}b^{*}\wt g_{1}^{*}
}
\end{equation}
Here the arrows labelled by $\di$ are the proper base change isomorphisms, and the arrows labelled by $u$ are the unit maps $\Id\to a_{*}a^{*}, \Id\to b_{*}b^{*}$ and $\Id\to c_{*}c^{*}$. Starting from the upper left corner, going along the top and then down to reach the lower right corner is the top row of \eqref{eq: comp pull-push}. On the other hand, going along the diagonal first, then down and then right gives the bottom row of \eqref{eq: comp pull-push}. Therefore it suffices to show that each square and triangle in \eqref{comp big diagram} is commutative. These are almost all clear except possibly the square labelled $\star$. To show $\star$ commutes, we refer to the Cartesian square
\begin{equation}
\xymatrix{ \wt B\ar[r]^{b}\ar[d]^{\wt f} & \wt B_{1}\ar[d]^{\wt f_{1}}\\
C\ar[r]^{c} & \wt D}
\end{equation}
and would like to show that
\begin{equation}
\xymatrix{ \wt f_{1!} \ar[r]^-{u} \ar[d]^{u} & c_{*}c^{*}\wt f_{1!}\ar[d]^{\di}\\
\wt f_{1!}b_{*}b^{*}\ar@{=}[r] &  c_{!}\wt f_{!}b^{*}}
\end{equation}
is commutative. Note that $b$ and $c$ are proper. To check this, it suffices to check on the geometric stalks. Therefore we may reduce to the case where $\wt D$ is a geometric point, in which case both compositions are identified with the map $R\Gamma_{c}(\wt B_{1}, -)\to R\Gamma_{c}(\wt B_{1}, -)\ot R\Gamma(C)$ given by $x\mapsto x\ot 1$.
\end{proof}

The following variant of Proposition \ref{p: comp push} will also be needed later. 

\begin{prop}\label{p: comp push horizontal} Suppose we have two pushable squares
\begin{equation}
\xymatrix{A\ar[r]^{g'}\ar[d]^{f''} & B\ar[r]^{h'}\ar[d]^{f'} & C\ar[d]^{f}\\
D\ar[r]^{g} & E\ar[r]^{h} & F}
\end{equation}
Then:
\begin{enumerate}
\item The outer square is also pushable.
\item The following diagram is commutative
\begin{equation}
\xymatrix{ g^{*}h^{*}f_{!} \ar@{=}[d] \ar[r]^{g^{*}\td} & g^{*}f'_{!}(h')^{*} \ar[r]^{\td (h')^{*}}& f''_{!} (g')^{*}(h')^{*}\ar@{=}[d]\\
(h\c g)^{*}f_{!} \ar[rr]^{\td} & & f''_{!}(h'\c g')^{*}
}
\end{equation}
\end{enumerate}
\end{prop}
\begin{proof}
Part (1) is the same as Lemma \ref{lem: 2 out of 3}(1). Part (2) follows from a similar argument as in the proof of Proposition \ref{p: comp push}. 

\end{proof}

\subsection{Push-push}\label{ssec: push-push}

Suppose the square \eqref{eq: commutative square} is Cartesian. Then we have a base change map 
\begin{equation}\label{eq: push-push mixed cartesian}
f_! g'_* \xrightarrow{\diamond} g_* f'_!.
\end{equation}
To construct this, we construct the adjoint map $g^* f_! g'_* \rightarrow f'_!$. For this, we start with the proper base change isomorphism $g^* f_! \xrightarrow{\di} f'_! (g')^*$, so composing with the counit of the $((g')^*, g'_*)$ adjunction gives a sequence of natural transformations
\[
g^* f_! g'_* \xrightarrow{\di} f'_! (g')^* g'_* \rightarrow f'_!.
\]

\begin{lemma}\label{lem: push push compatibility 1} Consider the commutative square \eqref{eq: commutative square}. Suppose it is Cartesian and the maps $g,g'$ are separated and locally of finite type.
\begin{enumerate}
\item The following diagram commutes:
\[
\begin{tikzcd}
f_! g'_!  \ar[r, equals] \ar[d, "\can(g')"'] & g_! f'_! \ar[d, "\can(g)"] \\ 
f_! g'_*  \ar[r, "\eqref{eq: push-push mixed cartesian}"]  & g_* f'_! 
\end{tikzcd}
\]

\item Assume additionally that $f$ is separated, hence $f'$ is also separated. Then the following diagram commutes:
\[
\begin{tikzcd}
f_! g'_*  \ar[r, "\eqref{eq: push-push mixed cartesian}"]   \ar[d, "\can(f)"'] & g_* f'_! \ar[d, "\can(f')"] \\ 
f_* g'_*    \ar[r, equals]  & g_* f'_* 
\end{tikzcd}
\]
\end{enumerate}
\end{lemma}

\begin{proof} Since $g$ is separated and locally of finite type, we may compactify $g$ into the composition of an open embedding followed by a proper map. The statement is obvious when $g$ is proper, so we reduce to the case where $g$ (hence also $g'$) is an open embedding. It suffices to check that the adjoint diagram 
\begin{equation}\label{eq: push-push compatibility}
\begin{tikzcd}
g^* f_! g'_!  \ar[r]   \ar[d, "\can(g')"'] &  f'_! \ar[d, equals, "\Id"] \\ 
g^*  f_! g'_*    \ar[r, equals]  &   f'_!
\end{tikzcd}
\end{equation}
commutes. By definition, the adjoint of the base change map $\di$ is the composition 
\[
g^* f_! g'_* \xrightarrow{\di}  f'_! (g')^* g'_*  \rightarrow f'_!.
\]
But since $g'$ is an open embedding, $(g')^* g'_!$ is the identity functor. Hence, after applying the isomorphism $g^*f_! \xrightarrow{\di} f'_! (g')^*$ to both of the left terms \eqref{eq: push-push compatibility}, the upper right path through \eqref{eq: push-push compatibility} is the identity natural transformation $f'_! \xrightarrow{\Id} f'_!$. Similarly, since $(g')^* g'_*$ is also the identity functor, the lower left path through \eqref{eq: push-push compatibility} is also evidently the identity natural transformation $f'_! \xrightarrow{\Id} f'_!$. The map $\can(g') \co (g')^* g'_! \rightarrow (g')^* g'_*$ is the identity natural transformation, so the diagram commutes. 

(2) Similar to (1). 
\end{proof}

Now, suppose that \eqref{eq: commutative square} is pushable. Then natural transformation 
\eqref{eq: pull-push gen} induces by adjunction a map
\begin{equation}\label{eq: push-push gen}
f_{!}g'_{*}\to g_{*}f'_{!}.
\end{equation}
Unravelling the construction, it is the composition
\begin{equation}
f_{!}g'_{*}=f_{!}\wt g_{*}a_{*}\xr{\di}g_{*}\wt f_{!}a_{*}=g_{*}\wt f_{!}a_{!}\to g_{*}f'_{!}.
\end{equation}
Here, the map $\di$ is \eqref{eq: push-push mixed cartesian}. 

The natural transformation \eqref{eq: push-push gen} has similar compositional properties to \eqref{eq: pull-push gen}. 

\begin{prop}
Assume both the upper square and the lower square in \eqref{eq: commutative composition} are pushable. Then the following diagram commutes:
\begin{equation}\label{eq: comp push-push}
\xymatrix{h_! f_! g''_* \ar@{=}[d] \ar[r]^-{h_! \td} & h_! g'_* f'_! \ar[r]^-{\td f'_!} & g_* h'_! f'_! \ar@{=}[d]\\
(h\c f)_{!} g''_* \ar[rr]^-{\td} &&   g_* (h'\c f')_{!}
}
\end{equation} 
\end{prop}

\begin{proof}
Similar to the proof of Proposition \ref{p: comp push}.
\end{proof}

\subsection{Recollections on relative fundamental classes}\label{ssec: relative fundamental nt} We review some properties of relative fundamental classes.

Let $f \co Y \rightarrow Z$ be a quasi-smooth map of derived stacks, of relative dimension $d(f)$. Then one has a relative fundamental class (also called ``Gysin map'') 
\begin{equation}\label{eq: Gysin}
f^* \Qll{Z} \xrightarrow{[f]} f^! \Qll{Z} \tw{-d(f)}.
\end{equation}
In fact, \eqref{eq: Gysin} induces a natural transformation $f^* \rightarrow f^! \tw{-d(f)}$. To see this, we note that there is a base change natural transformation for $\cK_{1},\cK_{2}\in D(Z)$
\begin{equation}\label{eq: shriek bifunctor}
f^*(\cK_1) \otimes f^!(\cK_2) \xrightarrow{\di} f^!(\cK_1 \otimes \cK_2)
\end{equation}
arising from the Cartesian square
\[
\begin{tikzcd}
Y \ar[r, "\Gamma_f"] \ar[d, "f"] & Y \times Z \ar[d, "f \times \Id"] \\
Z \ar[r, "\Delta"] & Z \times Z
\end{tikzcd}
\]
where $\Gamma_f$ is the graph of $f$. Taking $\cK_2 = \Qll{Z}$ in \eqref{eq: shriek bifunctor}, the relative fundamental class induces
\[
\begin{gathered}
\begin{aligned}
f^*(\cK_1)
&= f^*(\cK_1) \otimes f^*(\cK_2) \\
&\xrightarrow{\eqref{eq: Gysin}} f^* (\cK_1) \otimes f^!(\cK_2) \tw{-d(f)}
\xrightarrow{\eqref{eq: shriek bifunctor}} f^!(\cK_1 \otimes \cK_2) \tw{-d(f)}\\
&= f^! (\cK_1) \tw{-d(f)}.
\end{aligned}
\end{gathered}
\]

The base change property of relative fundamental classes \cite[Theorem 3.13]{KhanI} says that for a derived Cartesian square of derived Artin stacks 
\[
\begin{tikzcd}
Y' \ar[r, "f'"] \ar[d, "s'"] & Z' \ar[d, "s"] \\
Y \ar[r, "f"] & Z
\end{tikzcd}
\]
with $f$ (hence also $f'$) quasi-smooth of relative dimension $d(f) = d(f')$, the base change homomorphism $(s')^* f^! \rightarrow (f')^! s^*$ fits into a commutative diagram 
\begin{equation}\label{eq: base change relative class} 
\begin{tikzcd}
(s')^* f^* \ar[r, equals] \ar[d, "{[f]}"]  & (f')^* s^*  \ar[d, "{[f']}"]  \\
 (s')^* f^! \tw{-d(f)}\ar[r, "\di"] &  (f')^! s^* \tw{-d(f)}
\end{tikzcd}
\end{equation}
and dually, the diagram below commutes
\begin{equation}\label{eq: base change relative class dual}
\begin{tikzcd}
 (s')^! f^*   \ar[d, "{[f]}"]  & (f')^* s^!  \ar[d, "{[f']}"]   \ar[l, "\di"] \\
(s')^! f^!\tw{-d(f)} \ar[r, equals] &  (f')^! s^! \tw{-d(f)}
\end{tikzcd}
\end{equation}

\subsection{Pull-pull}\label{ssec: pull-pull}


Suppose the square \eqref{eq: commutative square} is Cartesian. Then we have a base change transformation 
\begin{equation}\label{eq: pull-pull mixed cartesian}
 (f')^* g^! \xrightarrow{\diamond} (g')^! f^*.
\end{equation}
To construct this, we construct the adjoint map $g^! \rightarrow f'_* (g')^! f^*$. We start with the proper base change natural isomorphism, $g^! f_* \xrightarrow{\sim} f'_* (g')^!$. Then we compose this with the unit of the $(f^*, f_*)$-adjunction: 
\[
g^! \xrightarrow{\text{unit}(f)} g^! f_* f^* \xrightarrow{\sim}  f'_* (g')^! f^*.
\]

Now suppose only that \eqref{eq: commutative square} is pullable, i.e., the map $a$ in \eqref{eq: pushable square} is quasi-smooth. Then we have a natural transformation
\begin{equation}\label{eq: pull-pull gen}
(f')^{*}g^{!} =  a^{*} \wt f^{*}g^{!} \xr{a^* \di}a^{*} \wt g^{!} f^{*}\xr{[a]}a^{!} \wt g^{!} f^{*}\tw{-d(a)}=(g')^{!} f^{*}\tw{-d(a)}.
\end{equation}
Here the map $\di$ is induced from the proper base change isomorphism as in \eqref{eq: pull-pull mixed cartesian}. We have also used the natural transformation
\[
[a]: a^* \rightarrow a^! \tw{-d(a)}
\]
induced by the derived fundamental class of the quasi-smooth map $a$. 

\subsubsection{Gysin natural transformation}\label{sssec: Gysin NT}

We often denote such a natural transformation induced by a pullable square by $\tu$,
\begin{equation}
(f')^{*}g^{!}\xr{\tu}(g')^{!} f^{*}\tw{-\d}
\end{equation}
where $\d=d(a)$ is the defect of the pullable square.

\begin{remark} \label{rem: p-p} If \eqref{eq: commutative square} is pullable, the map \eqref{eq: pull-pull gen} then induces the following maps by adjunction
\begin{equation}\label{eq: p-p}
g'_{!}(f')^{*}\to  f^{*}g_{!}\tw{-\d}.
\end{equation}
This resembles the pull-push map \eqref{eq: pull-push gen} (with the square flipped), but the arrow is in the opposite direction and there is a shift and twist.
\end{remark}

\begin{example}\label{ex: pull-pull f smooth}
If $f$ is smooth and $f'$ is quasi-smooth, then the square \eqref{eq: commutative square} is pullable with defect $d(f')-d(f)$. In this case, the base change map \eqref{eq: pull-pull gen} can be alternatively described as the composition
\begin{equation}
(f')^{*}g^{!} \xr{[f']} (f')^{!}g^{!}\tw{-d(f')}=(g')^{!}f^{!}\tw{-d(f')}\cong (g')^{!}f^{*}\tw{d(f)-d(f')}
\end{equation}
where the last isomorphism is induced from the inverse of the isomorphism $[f]: f^{*}\isom f^{!}\tw{-d(f)}$ since $f$ is smooth. The agreement of the composition above and \eqref{eq: pull-pull gen} can be proved by a similar argument to that of Example \ref{ex: pull-push f proper}. 
\end{example}

\subsubsection{Compositions} 

The setup is the same as in \S \ref{sssec: push-pull composition}. 

\begin{lemma}\label{l:pull 2 out of 3} Consider the commutative diagram \eqref{eq: commutative composition}.

(1) If both the upper square and the lower square are pullable, then the outer square formed by $(A,B,E,F)$ is also pullable.

(2) If $\d_{\upp},\d_{\low}$ and $\d_{\out}$ denote the defects of the upper, lower, and outer squares respectively, then
\begin{equation}
\d_{\out}=\d_{\upp}+\d_{\low}.
\end{equation}
\end{lemma}
\begin{proof}
(1) Consider the commutative diagram \eqref{eq: comp two pushable}, in which all squares are Cartesian. Both $a$ and $c$ are quasi-smooth, hence $b$ is also quasi-smooth. It follows that $b\circ a$ is quasi-smooth, i.e., the outer square is also pullable.


(2) The defect equality follows from the fact that $\d_{\upp}=d(a)$, $\d_{\low}=d(c)=d(b)$ and $\d_{\out}=d(b\circ a)$.
\end{proof}


\begin{prop}\label{p: comp pull} Suppose in the commutative diagram \eqref{eq: commutative composition} both the upper and the lower squares are pullable. Let $\d_{\upp},\d_{\low}$ and $\d_{\out}$ denote the defects of the upper, lower, and outer squares respectively. Then the following diagram is commutative
\begin{equation}\label{eq: comp pull}
\xymatrix{(f')^{*}(h')^{*}g^{!}\ar[r]^-{f'^{*}\tu}\ar@{=}[d] & (f')^{*}(g')^{!}h^{*}\tw{-\d_{\low}}\ar[r]^-{\tu h^{*}} & (g'')^{!}f^{*}h^{*}\tw{-\d_{\low}-\d_{\upp}}\ar@{=}[d]
\\
(h'\c f')^{*}g^{!} \ar[rr]^-{\tu} && (g'')^{!}(h\c f)^{*}\tw{-\d_{\out}}}
\end{equation}
\end{prop}

\begin{proof}
Consider the diagram \eqref{eq: comp two pushable}. Let $\wt{g}' \co \wt{B} \rightarrow B$ be the pullback of $g'$. Consider the diagram of natural transformations
\begin{equation}\label{eq: comp pull big diagram}
\begin{gathered}
\begin{adjustbox}{max width=\textwidth}
$\displaystyle
\begin{tikzcd}[ampersand replacement=\&]
(f')^* (h')^* g^! \ar[r, equals] \ar[dd, equals] \& (f')^* c^* \wt{h}^* g^! \ar[d, "\di"] \\
\& (f')^* c^* \wt{g}^! h^* \ar[d, equals] \ar[r, "{[c]}"] \&  (f')^* c^! \wt{g}^! h^*\tw{-\d_{\low}} = a^* \wt{f}^*c^! \wt{g}^! h^* \tw{-\d_{\low}} \ar[d, "\di"] \\
a^* b^* \wt{f}_1^* \wt{h}^* g^! \ar[d, "\di"]  \&  a^* b^* \wt{f}_1^* \wt{g}^! h^* \ar[r, "{[b]}"] \ar[d, "\di"]  \&  a^* b^! \wt{f}_1^* \wt{g}^! h^* \tw{-\d_{\low}} \ar[d, "\di"] \ar[r, "{[a]}"] \& a^! b^!  \wt{f}_1^* \wt{g}^! h^* \tw{-\d_{\out}} \ar[d, "\di"] \\
a^* b^* \wt{g}_1^! f^* h^* \ar[r, equals] \& a^* b^* \wt{g}_1^! f^* h^* \ar[r, "{[b]}"]   \& a^* b^! \wt{g}_1^! f^* h^* \tw{-\d_{\low}}  \ar[r, "{[a]}"]  \& a^! b^!  \wt{g}_1^! f^* h^*  \tw{-\d_{\out}}
\end{tikzcd}
$
\end{adjustbox}
\end{gathered}
\end{equation}
Tracing along the top and right edges gives the top and right path of \eqref{eq: comp pull}, while tracing along the left and bottom edges of \eqref{eq: comp pull big diagram} gives the bottom and left path of \eqref{eq: comp pull}, using that $[b] \circ [a]  = [b \circ a]$. Therefore, it suffices to verify that all the rectangles of \eqref{eq: comp pull big diagram} commute. For the rectangle on the left, this is because the proper base change natural isomorphisms are compatible for compositions of Cartesian squares. For the middle rectangle, the commutativity is an instance of \eqref{eq: base change relative class}. For the two rectangles in the bottom row, the commutativity is obvious. 

\end{proof}

\begin{prop}\label{p: comp pull horizontal} Suppose we have two pullable squares
\begin{equation}
\xymatrix{A\ar[r]^{g'}\ar[d]^{f''} & B\ar[r]^{h'}\ar[d]^{f'} & C\ar[d]^{f}\\
D\ar[r]^{g} & E\ar[r]^{h} & F}
\end{equation}
\begin{enumerate}
\item The outer square is also pullable. Moreover, if we let $\d_{\ell}$ and $\d_{r}$ be the defects of the left and right squares, respectively, then the defect of the outer square $\d_{\out}=\d_{\ell}+\d_{r}$. 
\item The following diagram is commutative
\begin{equation}
\xymatrix{
 (f'')^{*}g^{!}h^{!} \ar@{=}[d] \ar[r]^{\tu h^{!}} &
 (g')^{!}(f')^{*}h^{!}\tw{-\d_{\ell}} \ar[r]^{(g')^{!}\tu} &
 (g')^{!}(h')^{!}f^{*}\tw{-\d_{\ell}-\d_{r}} \ar@{=}[d]\\
 (f'')^{*}(h\c g)^{!} \ar[rr]^{\tu} && (h'\c g')^{!}f^{*}\tw{-\d_{\out}}
}
\end{equation}
\end{enumerate}
\end{prop}
\begin{proof}
(1) is the same as Lemma \ref{l:pull 2 out of 3}. The proof of (2) is similar to the proof of  Proposition \ref{p: comp pull}. 
\end{proof}

\section{The sheaf-cycle correspondence}\label{sec: sheaf-cycle}

As explained in the Introduction, an important part of our strategy is the realization of our cycle classes of interest as ``traces'' (in a suitable sense) of \emph{cohomological correspondences}. Furthermore, we will need notions of pushforward/pullback of cohomological correspondences, and we will need to know in some situations that they interact well with the notion of pushforward/pullback of cycles. We will explain this yoga in this section; it constitutes a framework that we call the ``sheaf-cycle correspondence'', extending the classical sheaf-function correspondence.

In \S \ref{ssec: trace of cc} we recall the notion of a cohomological correspondence, and the definition of the trace of a cohomological correspondence. In \S \ref{sssec: fix vs sht} we explain a variant of this construction, incorporating a Frobenius twist, that applies over finite fields. The early parts of this section are similar to material already appearing in \cite{SGA5} and \cite[\S 1]{Var07}, but we need much more generality than is handled there. 

In \S \ref{ssec: pushforward functoriality for CC} and \S \ref{ssec: pullback functoriality for CC} we will review some situations in which cohomological correspondences can be pushed forward or pulled back; this uses the material of \S \ref{sec: zoo}. In \S \ref{ssec: trace cc functoriality} we explain some situations in which this pushforward or pullback functoriality is compatible with the formation of traces, with proofs given in \S \ref{ssec: trace of smooth pullback}. In \S \ref{ssec: VFC as trace} we apply the theory to compute the trace for certain types of cohomological correspondences that will come up later. Finally, in \S \ref{ssec: coco} we introduce a dual notion of \emph{cohomological co-correspondence}, which plays a role in later analysis of the interaction with Fourier transform.

\subsection{The trace of a cohomological correspondence}\label{ssec: trace of cc}

\subsubsection{Basic definitions}
Let $Y_{0}$ and $Y_{1}$ be  derived Artin stacks. A {\em correspondence between $Y_{0}$ and $Y_{1}$} is a diagram of derived Artin stacks
\[
\begin{tikzcd}
& C \ar[dl, "c_0"'] \ar[dr, "c_1"]  \\
Y_0 & & Y_1 \end{tikzcd}
\]

Let $\cK_0 \in D(Y_0)$, $\cK_1 \in D(Y_1)$. A \emph{cohomological correspondence from $\cK_0$ to $\cK_1$ supported on $C$} is a map 
\[
c_0^* \cK_0 \rightarrow c_1^! \cK_1
\]
in $D(C)$. Let 
\begin{equation}
\Corr_{C}(\cK_{0}, \cK_{1}) :=\Hom_{C}(c_0^* \cK_0 , c_1^! \cK_1)
\end{equation}
denote the vector space of cohomological correspondences from $\cK_0$ to $\cK_1$ supported on $C$.

\subsubsection{Fixed points of a self-correspondence}
Now suppose that we have a fixed isomorphism $Y_0 \xrightarrow{\sim} Y_1$, which we will sometimes use to identify $Y_0$ with $Y_1$; however, it will also be convenient to distinguish them at times. Let $\Delta  \co Y_0 \rightarrow Y_0 \times Y_1$ be the diagonal embedding. Define $\Fix(C)$ as the fibered product 
\begin{equation}\label{eq: derived fixed point of c}
\begin{tikzcd}
\Fix(C) \ar[r, "\Delta'"] \ar[d, "c'"] & C \ar[d, "c"] \\
Y_0 \ar[r, "\Delta"] & Y_0 \times Y_1
\end{tikzcd}
\end{equation}
where $c = (c_0, c_1)$. 

\subsubsection{Trace of a cohomological self-correspondence}
Let $\cK \in D(Y_0)$. Then, following \cite[\S 1.2]{Var07}\footnote{Technically, Varshavsky does not consider the variant with the shifted Tate twist $\tw{-i}$.} we will define a trace map
\[
\Tr \co \Corr_{C}(\cK,\cK\tw{-i}) 
\rightarrow \mBM_{{2i}}(\Fix(C)) :=\rH^{-2i}(\Fix(C), \DD_{\Fix(C)}(-i)).
\]

\begin{observation}\label{obs: preliminaries for trace} Recall that one has the following isomorphisms:
\begin{enumerate}
\item For any $\cK_0 \in D(Y_0), \cK_1 \in D(Y_1)$, an isomorphism \cite[\S 0.3]{Var07}
\[
\cRHom_C(c_0^* \cK_0, c_{1}^! \cK_1) \cong c^! \cRHom_{Y_0 \times Y_1}(\pr_0^* \cK_0, \pr_1^! \cK_1)  \in D(C).
\]
\item For any $\cK_0 \in D(Y_0), \cK_1 \in D(Y_1)$, an isomorphism \cite[\S 0.4]{Var07} 
\[
\cRHom_{Y_0 \times Y_1}(\pr_0^* \cK_0, \pr_1^! \cK_1) \cong \DD(\cK_0) \boxtimes \cK_1 \in D(Y_0 \times Y_1).
\]
\end{enumerate}
\end{observation}

\begin{defn}\label{defn: trace}
Let $\cc \in \Corr_{C}(\cK,\cK\tw{-i})$. We will define its \emph{trace} $\Tr(\cc) \in \mBM_{2i}(\Fix(C))$. By Observation \ref{obs: preliminaries for trace}(1), we have  
\begin{equation}\label{eq: defining tr 1}
\cRHom_C(c_0^* \cK, c_1^! \cK \tw{-i}) \cong c^! \cRHom_{Y_0 \times Y_1}(\pr_0^*  \cK, \pr_1^! \cK \tw{-i}).
\end{equation}
Then Observation \ref{obs: preliminaries for trace}(2) gives an isomorphism
\begin{equation}\label{eq: defining tr 2}
c^!  \cRHom_{Y_0 \times Y_1}(\pr_0^*  \cK, \pr_1^! \cK \tw{-i}) \cong c^!(  \DD(\cK) \boxtimes \cK\tw{-i}).
\end{equation}
The evaluation map $\DD(\cK) \otimes \cK \tw{-i} \rightarrow \DD_{Y_0} \tw{-i}$ induces by adjunction a map 
\[
\DD(\cK) \boxtimes \cK\tw{-i} \rightarrow \Delta_* \DD_{Y_0} \tw{-i}.
\]
Composing this with \eqref{eq: defining tr 1} and \eqref{eq: defining tr 2} gives a map 
\begin{equation}\label{eq: defining tr 2.5}
\cRHom_C(c_0^* \cK, c_1^! \cK \tw{-i}) \rightarrow c^!  \Delta_* \DD_{Y_0} \tw{-i}.
\end{equation}
Finally, using (the adjoint form of) proper base change, we have isomorphisms 
\begin{equation}\label{eq: defining tr 3}
 c^! \Delta_* \DD_{Y_0} \tw{-i} \cong \Delta'_* (c')^! \DD_{Y_0} \tw{-i} = \Delta'_* \DD_{\Fix(C)} \tw{-i}.
\end{equation}
We may regard $\cc$ as a global section of $\cRHom_C(c_0^* \cK, c_1^! \cK \tw{-i})$. Then $\Tr(\cc) \in \rH^0( C, \Delta'_* \DD_{\Fix(C)} \tw{-i}) \cong \mBM_{2i}(\Fix(C))$ is defined as its image under the composition of \eqref{eq: defining tr 2.5} and \eqref{eq: defining tr 3}.
\end{defn}

\begin{remark}
It may not be clear why Definition \ref{defn: trace} is similar to the usual notion of ``trace''. The linear algebraic notion of trace has a natural generalization to symmetric monoidal categories, and it is possible to view the construction above as a special case of this general ``categorical trace'', at least when $i=0$. For a development along these lines, see \cite{LZ22}.
\end{remark}

\subsubsection{Shift and twist}\label{sss: cc sh tw}
Let $\cc\in \Corr_{C}(\cK_{0},\cK_{1})$. For $m,n\in\Z$, the same map $\cc: c_{0}^{*}\cK_{0}\to c_{1}^{!}\cK_{1}$ induces a map
\begin{equation*}
c_{0}^{*}\cK_{0}[m](n)\to c_{1}^{!}\cK_{1}[m](n)
\end{equation*}
which we denote by $\TT_{[m](n)}\cc$. Then $\cc\mapsto \TT_{[m](n)}\cc$ defines an isomorphism
\begin{equation*}
\TT_{[m](n)}: \Corr_{C}(\cK_{0},\cK_{1})\isom\Corr_{C}(\cK_{0}[m](n),\cK_{1}[m](n)).
\end{equation*}
When $Y_{0}=Y_{1}$ and $\cK_{1}=\cK_{0}\tw{-i}$, we have
\begin{equation*}
\Tr(\TT_{[m](n)}\cc)=(-1)^{m} \Tr(\cc)\in \mBM_{2i}(\Fix(C)).
\end{equation*}

\subsection{Fix vs Sht}\label{sssec: fix vs sht}
We will be working with correspondences of objects over the \emph{finite} field $k$. Therefore, our objects will have a Frobenius endomorphism $\Frob$, which in terms of the functor points is the absolute Frobenius $\Frob_q$ on the test scheme. We will use this to twist the preceding construction by Frobenius. 

For a correspondence
\[
\begin{tikzcd}
 Y & \ar[l, "c_0"'] C \ar[r, "c_1"] & Y
\end{tikzcd}
\]
over $k$, we will let $\Sht(C)$ (or sometimes $\Sht_Y$) be the derived fibered product 
\begin{equation}\label{eq: fix vs sht 1}
\begin{tikzcd}
\Sht(C) \ar[r] \ar[d] & C \ar[d, "{(c_0, c_1)}"] \\
Y \ar[r, "{(\Id, \Frob)}"] & Y \times Y 
\end{tikzcd}
\end{equation}
This derived fibered product can also be presented with the derived Cartesian square 
\begin{equation}\label{eq: fix vs sht 2}
\begin{tikzcd}
\Sht(C) \ar[r] \ar[d] & C \ar[d, "{(\Frob \circ c_0, c_1)}"] \\
Y \ar[r, "\Delta"] & Y \times Y 
\end{tikzcd}
\end{equation}
which is the ``fixed point Cartesian square'' for the correspondence datum 
\begin{equation}\label{eq: fix vs sht 3}
\begin{tikzcd}
 Y & \ar[l, "\Frob \circ c_0"'] C^{(1)} \ar[r, "c_1"] & Y
\end{tikzcd}
\end{equation}
where $C^{(1)} := C$ but with the left map twisted by $\Frob$. In other words, we have a canonical identification
\begin{equation}
\Sht(C)=\Fix(C^{(1)}).
\end{equation}

Given a cohomological correspondence $\cc \co c_0^* \cK_0 \rightarrow c_1^! \cK_1$ on $C$, with $\cK_i \in D(Y)$,  plus the canonical Weil structure $\Frob^* \cK_0 \cong \cK_0$ (because $Y$ is defined over $k=\F_{q}$), we have a cohomological correspondence $\cc^{(1)} \co (\Frob \circ c_0)^* \cK_0 \rightarrow c_1^! \cK_1$. In this way we obtain a linear isomorphism
\begin{equation*}
\Corr_{C}(\cK_{0}, \cK_{1})\isom \Corr_{C^{(1)}}(\cK_{0}, \cK_{1})
\end{equation*}
sending $\cc$ to $\cc^{(1)}$. If $\cK_1 = \cK_0 \tw{-i}$, then we define 
\[
\Tr^{\Sht}(\cc):=\Tr(\cc^{(1)}) \in \mBM_{2i}( \Fix(C^{(1)}))=\mBM_{2i}( \Sht(C))
\]
using the notion of trace in \S \ref{ssec: trace of cc} for the cohomological correspondence $\cc^{(1)}$ supported on $C^{(1)}$.

\begin{lemma}\label{lem: tangent complex of ShtM}
The tangent complex $\bT_{\Sht(C)/\F_q}$ is the restriction of $\bT_{c_1}$. In particular, if $c_1$ is quasi-smooth, then $\Sht(C)$ is quasi-smooth (over $\Spec \F_q$), of virtual dimension equal to $d(c_1)$. 
\end{lemma}

\begin{proof}The tangent complex of $\Sht(C)$ over $\F_q$ is the fibered product 
\[
\begin{tikzcd}
\bT_{\Sht(C)/\F_q} \ar[r] \ar[d] & \bT_{C/\F_q} |_{\Sht(C)}\ar[d, "{d(\Frob \circ c_0, c_1)}"] \\
\bT_{Y/\F_q}|_{\Sht(C)} \ar[r, "d\Delta"] & \bT_{Y/\F_q}|_{\Sht(C)} \oplus \bT_{Y/\F_q}|_{\Sht(C)}
\end{tikzcd}
\]
Since $\Frob$ induces the zero map of tangent complexes, this fibered product simplifies to $\bT_{c_1}|_{\Sht(C)}$. Since by assumption $c_1$ is quasi-smooth, this shows that $\Sht(C) \rightarrow \Spec \F_q$ is quasi-smooth. 
\end{proof}

Assume that $c_1$ is quasi-smooth, so that $[c_1] \in \mBM_{2d(c_1)}(C/ Y)$ exists. Then $\Sht(C)$ is quasi-smooth (over $\Spec \F_q$), so the derived fundamental class $[\Sht(C)]\in \mBM_{2d(c_1)}(\Sht(C))$ exists. 

On the other hand, regarding $[c_1]$ as a map $c_1^* \Qll{Y} \rightarrow c_1^! \Qll{Y} \tw{-d(c_1)}$, we have a cohomological correspondence 
\[
c_0^* \Qll{Y} = \Qll{C}  = c_1^* \Qll{Y} \xrightarrow{[c_1]} c_1^! \Qll{Y} \tw{-d(c_1)},
\]
whose composition we call $\cc_{Y}$, and regard as an element of
\[
\Corr_C(\Qll{Y}, \Qll{Y}\tw{-d(c_1)}).
\]
We equip $\Qll{Y}$ with the natural Weil structure $\Frob^* \Qll{Y} = \Qll{Y}$. Then we have $\Tr_{C}^{\Sht}(\cc_Y)  \in \mBM_{2d(c_1)}(\Sht(C))$. It is natural to ask when it will be true that
\begin{equation}\label{eq: general derived local term}
 \Tr_{C}^{\Sht}(\cc_Y) = [\Sht(C)] \in \mBM_{2d(c_1)}(\Sht(C)).
\end{equation}
We will see below in \S \ref{ssec: VFC as trace} that \eqref{eq: general derived local term} follows from Proposition \ref{prop: trace commutes with smooth pullback} whenever $Y$ is smooth. We expect that proving \eqref{eq: general derived local term} in more generality will be important for integral modularity statements. 

\subsubsection{Shift and twist} We view $\Ql(n)$ as a Weil sheaf where $\Frob$ acts as multiplication by $q^{-n}$. If $\cK_0$ is another Weil sheaf over $\Ql$, this equips $\cK_0(n)  := \cK_0 \otimes_{\Ql} \Ql(n)$ with a canonical Weil structure.  

For $\cc\in \Corr_{C}(\cK_{0},\cK_{0}\tw{-i})$ and $m,n\in\Z$, the shifted and twisted cohomological correspondence $\TT_{[m](n)}\cc\in \Corr_{C}(\cK_{0}[m](n),\cK_{0}[m](n)\tw{-i})$ is defined in \S\ref{sss: cc sh tw}. 
When $Y_{0}=Y_{1}$ and $\cK_{1}=\cK_{0}\tw{-i}$, we have
\begin{equation}\label{Tr Sht sh tw}
\Tr^{\Sht}(\TT_{[m](n)}\cc)=(-1)^{m}q^{-n}\Tr^{\Sht}(\cc)\in \mBM_{2i}(\Sht(C)).
\end{equation}

\subsection{Pushforward functoriality for cohomological correspondences}\label{ssec: pushforward functoriality for CC}

Suppose we have a commutative diagram of correspondences
\begin{equation}\label{eq: pushforward cohomological correspondences diagram}
\begin{tikzcd}
& C \ar[dl, "c_0"'] \ar[dr, "c_1"]  \ar[d, "f"] \\
Y_0 \ar[d, "f_0"'] & D \ar[dl, "d_0"'] \ar[dr, "d_1"]  & Y_1  \ar[d, "f_1"] \\
Z_0 & & Z_1 
\end{tikzcd}
\end{equation}
We assume that all morphisms are separated and representable in derived schemes.

\begin{defn}\label{def: corr pushable} The diagram of correspondences \eqref{eq: pushforward cohomological correspondences diagram} is called {\em left pushable} if the square with vertices $(C, Y_{0}, D, Z_{0})$ is pushable in the sense of Definition \ref{def: pushable pullable}. In other words, letting $\wt D_{0}=D\times_{Z_{0}}Y_{0}$, the natural map $\wt c_{0}=(f, c_{0}): C\to \wt D_{0}$ is proper.
\end{defn}

When \eqref{eq: pushforward cohomological correspondences diagram} is left pushable, for any cohomological correspondence $c_0^* \cK_0 \xrightarrow{\cc} c_1^! \cK_1$, we shall construct a ``pushforward correspondence'' $f_{!}(\cc): d_0^* f_{0!} \cK_0  \rightarrow d_1^!  f_{1!} \cK_1$, yielding a linear map
\begin{equation}
f_{!}: \Corr_{C}(\cK_{0}, \cK_{1})\to \Corr_{D}(f_{0!}\cK_{0}, f_{1!}\cK_{1}).
\end{equation}

Given $c_0^* \cK_0 \xrightarrow{\cc} c_1^! \cK_1$, we get a map
\begin{equation}\label{eq: push cc middle}
f_! c_0^* \cK_0 \xrightarrow{f_!(\cc)} f_! c_1^! \cK_1.
\end{equation}
Recall from \eqref{eq: push-pull shriek shriek} that there is always a base change morphism 
\begin{equation}\label{eq: push cc right}
f_! c_1^!  \cK_1 \rightarrow  d_1^! f_{1!} \cK_1.
\end{equation}
Since the left square of \eqref{eq: pushforward cohomological correspondences diagram} is pushable, \eqref{eq: pull-push gen} gives a base change morphism 
\begin{equation}\label{eq: push cc left}
d_0^* f_{0!} \cK_0 \rightarrow f_! c_0^* \cK_0.
\end{equation}
Precomposing \eqref{eq: push cc middle} with \eqref{eq: push cc left} and post-composing it with \eqref{eq: push cc right},  we get natural maps
\[
d_0^* f_{0!} \cK_0 \xrightarrow{\eqref{eq: push cc left}}	 f_! c_0^* \cK_0 \xrightarrow{\eqref{eq: push cc middle}} f_! c_1^! \cK_1 \xrightarrow{\eqref{eq: push cc right}}   d_1^! f_{1!} \cK_1
\]
whose composition is a cohomological correspondence from $f_{0!}\cK_0$ to $f_{1!} \cK_1$ that we denote by $f_{!}(\cc)$.

\subsection{Pullback functoriality for cohomological correspondences} \label{ssec: pullback functoriality for CC}
Consider the diagram of correspondences in \eqref{eq: pushforward cohomological correspondences diagram}.

\begin{defn}\label{def: corr pullable} The diagram of correspondences \eqref{eq: pushforward cohomological correspondences diagram} is called {\em right pullable} if the square with vertices $(C, Y_{1}, D, Z_{1})$ is pullable in the sense of Definition \ref{def: pushable pullable}. In other words, letting $\wt D_{1}=D\times_{Z_{1}}Y_{1}$, the natural map $\wt c_{1}=(f, c_{1}): C\to \wt D_{1}$ is quasi-smooth.

In this case, we define the {\em defect} $\d_{f}$ of the map of correspondences $f:C\to D$ to be the defect of the square $(C,Y_{1}, D,Z_{1})$, i.e., the relative dimension of $\wt{c}_1 \co C\to D\times_{Z_{1}}Y_{1}$. 
\end{defn}

When \eqref{eq: pushforward cohomological correspondences diagram} is right pullable, for any cohomological correspondence $d_0^* \cK_0 \xrightarrow{\cc} d_1^! \cK_1$, we shall construct a ``pullback correspondence'' $f^{*}(\cc): c_0^* f_{0}^{*} \cK_0  \rightarrow c_1^!  f_{1}^{*} \cK_1\tw{-\d_f}$, yielding a linear map
\begin{equation}
f^*: \Corr_{D}(\cK_{0}, \cK_{1})\to \Corr_{C}(f_{0}^{*}\cK_{0}, f_{1}^{*}\cK_{1}\tw{-\d_f}).
\end{equation}

Given $d_0^* \cK_0 \xrightarrow{\cc} d_1^! \cK_1$, we get a map
\begin{equation}\label{eq: pull cc middle}
f^* d_0^* \cK_0 \xrightarrow{f^*(\cc)} f^* d_1^! \cK_1.
\end{equation}
We have an obvious identification 
\begin{equation}\label{eq: pull cc left}
f^* d_0^* \cK_0 = c_0^* f_0^* \cK_0.
\end{equation}
Since the right square of \eqref{eq: pushforward cohomological correspondences diagram} is pullable, we get a base change morphism (cf. \S \ref{ssec: pull-pull})
\begin{equation}\label{eq: pull cc right}
f^* d_1^!  \cK_1 \rightarrow  c_1^! f_1^* \cK_1 \tw{-\d_f}.
\end{equation}
Precomposing \eqref{eq: pull cc middle} with \eqref{eq: pull cc left} and post-composing it with \eqref{eq: pull cc right}, we get natural maps
\[
c_0^* f_0^*  \cK_0 \xrightarrow{\eqref{eq: pull cc left}}	 f^* d_0^*  \cK_0 \xrightarrow{\eqref{eq: pull cc middle}} f^* d_1^! \cK_1 \xrightarrow{\eqref{eq: pull cc right}}   c_1^!  f_{1}^* \cK_1 \tw{-\d_f}
\]
whose composition is a cohomological correspondence between $f_0^* \cK_0$ and $f_{1}^* \cK_1\tw{-\d_f}$ that we denote by $f^{*}(\cc)$.

\subsection{Functoriality for the trace}\label{ssec: trace cc functoriality} Now we examine some situations in which pushforwards or pullbacks of cohomological correspondences are compatible with the formation of trace. 

Assume that in \eqref{eq: pushforward cohomological correspondences diagram} we have fixed identifications
\begin{equation}\label{eq: map of self-corr}
Y_{0}\cong Y_{1}, Z_{0}\cong Z_{1} \mbox{ and } f_{0}\cong f_{1}
\end{equation}
so that $\Fix(C)$ and $\Fix(D)$ are defined. 

\subsubsection{Proper pushforward} Assume the maps $f_{0}$ and $f$ are proper. Then \eqref{eq: pushforward cohomological correspondences diagram} is left pushable (cf. Example \ref{ex: pushable}), and the pushforward $f_{!}$ on cohomological correspondences is defined.  Moreover, the induced map $\Fix(f) \co \Fix(C) \rightarrow \Fix(D)$ is proper, so that we have a map $\Fix(f)_! \co \mBM_{2i}(\Fix(C)) \rightarrow \mBM_{2i}(\Fix(D))$. 

In this situation, we have the following compatibility between $f_!$ and the formation of trace, which generalizes \cite[Proposition 1.2.5]{Var07} and \cite[Corollary 2.22]{LZ22}. 

\begin{prop}\label{prop: trace commutes with proper push}
With notation as in \eqref{eq: pushforward cohomological correspondences diagram} and \eqref{eq: map of self-corr}, assume that $f_0, f, f_1$ are all proper and let $\cc\in \Corr_{C}(\cK,\cK\tw{-i})$. 
Then we have
\[
\Tr(f_{!}  \cc ) =\Fix(f)_! \Tr(\cc) \in \mBM_{2i}(\Fix(D)).
\]
\end{prop}

\begin{proof}
We stress at the outset that the result is essentially due to Lu--Zheng: their \cite[Theorem 2.21]{LZ22} proves the result when $i=0$ and all stacks involved are schemes. Furthermore, their proof can be adapted to our more general setting without much difficulty. To see the details written out, one can repeat the discussion of \cite[\S 5]{FK} verbatim with $\cD_{\mathrm{mot}}(-; \Q)$ replaced by our derived category $\cD_{\et}(-;\Ql)$ of $\ell$-adic sheaves. The statement in question is then the analogue of \cite[Theorem 5.4.2]{FK}.\footnote{From a pedagogical perspective this is somewhat backwards, as the point of the sequel paper \cite{FK} is to handle the more difficult setting of motivic sheaves. However, the referee accurately pointed out that this Proposition was important enough that a reference for its full generality should be given.} 
\end{proof}

\begin{remark}
Varshavsky's proof of \cite[Proposition 1.2.5]{Var07} does not generalize immediately to stacks (since he requires the existence of compactifications for all total spaces), but the proof of Lu-Zheng applies as is (cf. \cite[Remark 2.24]{LZ22}).
\end{remark}

\subsubsection{Smooth pullback}

Assume that the map $f_{1}$ (hence also $f_{0}$) is smooth, and $f$ is quasi-smooth. Then \eqref{eq: pushforward cohomological correspondences diagram} is right pullable (cf. Example \ref{ex: pullable}) of defect $\d=d(f)-d(f_{1})$, hence the map $f^{*}$ on cohomological correspondences is defined.  

On the other hand, we consider the map between fixed points. 

\begin{lemma}\label{l: Fix(f) qsm} Under the above assumptions, the map $\Fix(f): \Fix(C)\to \Fix(D)$ is quasi-smooth of dimension $\d=d(f)-d(f_{1})$. 
\end{lemma}
\begin{proof}
Define $E=\Fix(D)\times_{D}C$. Then $\Fix(f)$ factors as
\begin{equation*}
\Fix(f): \Fix(C)\xr{\a} E\xr{\b}\Fix(D).
\end{equation*}
We get an exact triangle of quasi-coherent sheaves on $\Fix(C)$ by taking cotangent complexes:
\begin{equation}\label{exact tri for LFix}
\a^{*}\bL_{\b}\to \bL_{\Fix(f)}\to \bL_{\a}.
\end{equation}
From the derived Cartesian square
\begin{equation*}
\xymatrix{ E\ar[d]^{\g} \ar[r]^{\b} &   \Fix(D)\ar[d]\\
C \ar[r]^{f} & D 
}
\end{equation*}
we get that $\bL_{\b}\cong \g^{*}\bL_{f}$. Denote both $Y_{0}$ and $Y_{1}$ by $Y$, and $Z_{0}$ and $Z_{1}$ by $Z$. From the derived Cartesian square
\begin{equation*}
\xymatrix{\Fix(C)\ar[r]^{\a}\ar[d]^{j} & E\ar[d]\\
Y\ar[r]^{\D_{Y/Z}} &  Y\times_{Z} Y
}
\end{equation*}
we conclude that $\bL_{\a}\cong j^{*}\bL_{\D_{Y/Z}}\cong j^{*}\bL_{f_{1}}[1]$. Combining these observations with \eqref{exact tri for LFix}, we get an exact triangle
\begin{equation}\label{exact tri for LFix 2}
\a^{*}\g^{*}\bL_{f}\to \bL_{\Fix(f)}\to j^{*}\bL_{f_{1}}[1].
\end{equation}

Since $f_{0}=f_{1}: Y\to Z$ is smooth, we see that $j^* \bL_{f_{1}}[1]$ is a perfect complex with tor-amplitude in $[-1, \infty)$. Since $f$ is quasi-smooth, we see that $\a^* \gamma^* \bL_{f}$ is a perfect complex with tor-amplitude in $[-1, \infty)$. By \eqref{exact tri for LFix 2}, we conclude that $\bL_{\Fix(f)}$ is a perfect complex with tor-amplitude in $[-1, \infty)$, i.e., $\Fix(f)$ is quasi-smooth. The relative dimension calculation of $\Fix(f)$ also follows from the exact triangle \eqref{exact tri for LFix 2}.
\end{proof}

By Lemma \ref{l: Fix(f) qsm}, the pullback map on Borel-Moore homology
\begin{equation*}
\Fix(f)^{*} \co \mBM_{2i}(\Fix(D)) \rightarrow \mBM_{2i+2\d}(\Fix(C))
\end{equation*}
is defined. It is induced from the map 
\begin{equation}\label{Fix(f) pullback}
\Fix(f)^{*}\DD_{\Fix(D)}\xr{[\Fix(f)]}\Fix(f)^{!}\DD_{\Fix(D)}\tw{-\d}=\DD_{\Fix(C)}\tw{-\d}
\end{equation}
by taking global sections.

\begin{prop}\label{prop: trace commutes with smooth pullback}
With notation as in \eqref{eq: pushforward cohomological correspondences diagram} and \eqref{eq: map of self-corr}, assume that $f_1$ is smooth and  $f$ is quasi-smooth. Let  $\d=d(f)-d(f_{1})$. Let $\mf{d} \in \Corr_{D}(\cK,\cK\tw{-i})$. 
Then we have
\[
\Tr(f^{*} \mf{d}) = \Fix(f)^{*} \Tr(\mf{d}) \in \mBM_{2i+2\d}(\Fix(C)).
\]
\end{prop}

The proof of Proposition \ref{prop: trace commutes with smooth pullback} is long, and will be given over the course of the next subsection.

\subsection[Proof of trace commutes with smooth pullback]{Proof of Proposition \ref{prop: trace commutes with smooth pullback}}\label{ssec: trace of smooth pullback} We begin with some preliminary technical observations.

\subsubsection{Another construction of the trace of a cohomological correspondence} Let notation be as in \S \ref{ssec: trace of cc}. We fix an identification $Y_0 \cong Y_1 =: Y$.

The following alternative formulation of the trace will be useful. By Observation \ref{obs: preliminaries for trace}, we may interpret a cohomological correspondence $\cc \co c_0^* \cK  \rightarrow c_1^! \cK \tw{-i}$ as an element of $\rH^{-2i}(C, c^! (\DD_Y(\cK) \boxtimes \cK)(-i))$. We have a map 
\begin{equation}\label{eq: alternative trace 1}
\rH^*(C, c^! (\DD_Y(\cK) \boxtimes \cK)) \rightarrow \rH^*(\Fix(C), (\Delta')^* c^!(\DD_Y(\cK) \boxtimes \cK)).
\end{equation}
Then the base change transformation $(\Delta')^* c^! \xrightarrow{\di} (c')^! \Delta^*$ for the Cartesian square \eqref{eq: derived fixed point of c} induces a map 
\begin{equation}\label{eq: alternative trace 2}
\rH^*(\Fix(C), (\Delta')^* c^!(\DD_Y(\cK) \boxtimes \cK) ) \xrightarrow{\di} \rH^*(\Fix(C), (c')^! \Delta^* (\DD_Y(\cK) \boxtimes \cK) ) .
\end{equation}
Applying the trace map $\Delta^* (\DD_Y(\cK) \boxtimes \cK) \cong \DD_Y(\cK) \otimes \cK \rightarrow \DD_Y$, we get a map 
\begin{equation}\label{eq: alternative trace 3}
\begin{gathered}
\begin{aligned}
\rH^{-2i}(\Fix(C), (c')^! \Delta^* (\DD_Y(\cK) \boxtimes \cK)(-i) )
&\rightarrow \rH^{-2i}(\Fix(C), (c')^! \DD_Y(-i)) \\
&\cong  \rH^{-2i}(\Fix(C), \DD_{\Fix(C)}(-i)) \\
&=: \mBM_{2i}(\Fix(C)).
\end{aligned}
\end{gathered}
\end{equation}

\begin{lemma}\label{lem: trace as BC}
The map $\Tr \co \Corr_{C}(\cK,\cK\tw{-i})\cong \rH^{-2i}(C, c^! (\DD_Y(\cK) \boxtimes \cK)(-i))  \rightarrow \mBM_{2i}(\Fix(C))$ coincides with the composition of the maps \eqref{eq: alternative trace 1}, \eqref{eq: alternative trace 2}, and \eqref{eq: alternative trace 3}.
\end{lemma}

\begin{proof} In the proof we denote all pull-pull base change transformations \eqref{eq: pull-pull mixed cartesian} induced by Cartesian squares with the symbol $\di$.

Comparing the definitions, the Lemma amounts to the commutativity of the diagram 
\begin{equation}\label{eq: alt trace eq 1}
\begin{tikzcd}
c^! \ar[d, equals] \ar[r, "\mrm{unit}(\Delta')"] & \Delta'_* (\Delta')^* c^! \ar[r, "\di"] & \Delta'_* (c')^! \Delta^* \ar[d, equals] \\
c^! \ar[r, "\mrm{unit}(\Delta)"] & c^! \Delta_* \Delta^* \ar[r, "\di"] & \Delta'_*(c')^! \Delta^*
\end{tikzcd}
\end{equation}

By definition, the morphism $(\Delta')^* c^! \xrightarrow{\di} (c')^! \Delta^* $ is the composition 
\[
(\Delta')^* c^! \xrightarrow{\mrm{unit}(\Delta)} (\Delta')^* c^! \Delta_* \Delta^* \xrightarrow{\di} (\Delta')^* \Delta'_* (c')^! \Delta^* \xrightarrow{\mrm{counit}(\Delta')} (c')^! \Delta^*.
\]
Therefore, the top row of \eqref{eq: alt trace eq 1} fits into a commutative diagram 
\begin{equation}\label{eq: alternative trace 4}
\begin{tikzcd}[column sep = huge]
c^! \ar[d, equals] \ar[rr, "\mrm{unit}(\Delta')"] &&  \Delta'_* (\Delta')^* c^! \ar[r, "\di"] \ar[d, "\mrm{unit}(\Delta)"] & \Delta'_* (c')^! \Delta^*  \\
c^! \ar[rr, "\mrm{unit}(\Delta') \circ \mrm{unit}(\Delta)"] & & \Delta'_* (\Delta')^* c^! \Delta_* \Delta^* \ar[r, "\di"] & \Delta'_* (\Delta')^* \Delta'_*  (c')^! \Delta^*  \ar[u, "\mrm{counit}(\Delta')"]
\end{tikzcd}
\end{equation}
Similarly we have a commutative diagram 
\begin{equation}\label{eq: alternative trace 5}
\begin{tikzcd}
c^! \Delta_* \Delta^* \ar[r, "\di"] \ar[d, "\mrm{unit}(\Delta')"] & \Delta'_* (c')^! \Delta^* \ar[d, "\mrm{unit}(\Delta')"] \\
\Delta'_* (\Delta')^* c^! \Delta_* \Delta^* \ar[r, "\di"] & (\Delta')_* (\Delta')^* (\Delta'_*) (c')^! \Delta^*
\end{tikzcd}
\end{equation}
Putting together \eqref{eq: alternative trace 4} and \eqref{eq: alternative trace 5}, we get a commutative diagram 
\[
\begin{tikzcd}
c^! \ar[d, equals] \ar[rr, "\mrm{unit}(\Delta')"] & & \Delta'_* (\Delta')^* c^! \ar[r, "\di"] \ar[d, "\mrm{unit}(\Delta)"] & \Delta'_* (c')^! \Delta^*  \\
c^! \ar[rr, "\mrm{unit}(\Delta') \circ \mrm{unit}(\Delta)"] \ar[d, equals] &  & \Delta'_* (\Delta')^* c^! \Delta_* \Delta^* \ar[r, "\di"] & \Delta'_* (\Delta')^* \Delta'_* (c')^! \Delta^*  \ar[u, "\mrm{counit}(\Delta')"] \\
c^! \ar[rr, "\mrm{unit}(\Delta)"] & &  c^! \Delta_* \Delta^* \ar[r, "\di"] \ar[u, "\mrm{unit}(\Delta')"'] & (\Delta'_*)(c')^! \Delta^* \ar[u, "\mrm{unit}(\Delta')"] 
\end{tikzcd}
\]
whose upper and lower rows agree with those in \eqref{eq: alt trace eq 1}. It remains to verify that the vertical composition on the right column, $\mrm{counit}(\Delta') \circ \mrm{unit}(\Delta')$, is the identity map. But this is the triangle identity for the adjunction $((\Delta')^*,\Delta'_*)$.
\end{proof}

\subsubsection{Alternative description of pullback cohomological correspondence}\label{sss: alt pullback corr} 
Recall the setup of Proposition \ref{prop: trace commutes with smooth pullback}: suppose we are given commutative diagram of derived Artin stacks
\[
\begin{tikzcd}
Y_0 \ar[d, "f_0"]  & \ar[l, "{c_0}"'] C\ar[r, "{c_1}"] \ar[d, "f"] &  Y_1 \ar[d, "f_1"]  \\
Z_0 & \ar[l, "{d_0}"'] D \ar[r, "{d_1}"] &  Z_1 
\end{tikzcd}
\]
such that $f_{1}$ is smooth and $f$ is quasi-smooth with $d(f)-d(f_{1})=\d$. We suppose a fixed identification of $f_0 \co Y_0 \rightarrow Z_0$ with $f_1 \co Y_1 \rightarrow Z_1$. 

Since $f_{1}$ is smooth and $f$ is quasi-smooth, we may use Example \ref{ex: pull-pull f smooth} to give an alternative description of the pullback map $f^{*}$ on cohomological correspondences: given $\frd: d_{0}^{*}\cK_{0}\to d_{1}^{!}\cK_{1}$ (where $\cK_{i}\in D(Z_{i})$), $f^{*}\frd$ is the composition
\begin{equation*}
c_{0}^{*}f_{0}^{*}\cK_{0}=f^{*}d_{0}^{*}\cK_{0}\xr{f^{*}\frd}f^{*}d_{1}^{!}\cK_{1}\xr{[f]}f^{!}d_{1}^{!}\cK_{1}\tw{-d(f)}=c_{1}^{!}f_{1}^{!}\cK_{1}\tw{-d(f)}\stackrel{[f_{1}]}{\cong} c_{1}^{!}f_{1}^{*}\cK_{1}\tw{-\d}.
\end{equation*}
By Observation \ref{obs: preliminaries for trace}, we have 
\begin{eqnarray}\label{Corr D H0}
\Corr_{D}(\cK_{0}, \cK_{1})&=&\upH^{0}(D, d^! (\DD_{Z_0}(\cK_0) \boxtimes \cK_1)),\\
\label{Corr C H0}\Corr_{C}(f_{0}^{*}\cK_{0}, f_{1}^{*}\cK_{1}\tw{-\d})&=&\upH^{0}(C, c^! (\DD_{Y_0}(f_{0}^{*}\cK_0) \boxtimes f_{1}^{*}\cK_1\tw{-\d})).
\end{eqnarray}
Consider the commutative square
\begin{equation*}
\xymatrix{ C \ar[d]^{c}\ar[r]^{f} & D \ar[d]^{d}\\ 
Y_{0}\times Y_{1}\ar[r]^{f_{0}\times f_{1}} & Z_{0}\times Z_{1}}
\end{equation*}
Since $f_{0}\times f_{1}$ is smooth and $f$ is quasi-smooth, this square is pullable (cf. Example \ref{ex: pullable}). The base change map (cf. \S\ref{ssec: pull-pull}) in this situation is defined and reads
\begin{equation*}
f^{*}d^{!}\xr{\tu}c^{!}(f_{0}\times f_{1})^{*}\tw{-d(f)+d(f_{0})+d(f_{1})}.
\end{equation*}
Using the identifications \eqref{Corr D H0} and  \eqref{Corr C H0}, the map $f^{*}: \Corr_{D}(\cK_{0}, \cK_{1})\to \Corr_{C}(f^{*}_{0}\cK_{0}, f_{1}^{*}\cK_{1}\tw{-\d})$ is induced by the following map upon taking $\upH^{0}$:
\begin{equation}\label{alt pullback corr}
\begin{gathered}
\begin{aligned}
f^{*}d^{!}(\DD\cK_{0}\bt\cK_{1})
&\xr{\tu}c^{!}(f_{0}^{*}\DD\cK_{0}\bt f^{*}_{1}\cK_{1})\tw{-d(f)+d(f_{0})+d(f_{1})} \\
&\cong  c^{!}(\DD(f_{0}^{!}\cK_{0})\bt f^{*}_{1}\cK_{1})\tw{-d(f)+d(f_{0})+d(f_{1})}\\
&\stackrel{[f_{0}]}{\cong} c^{!}(\DD(f_{0}^{*}\cK_{0})\bt f^{*}_{1}\cK_{1})\tw{-\d}.
\end{aligned}
\end{gathered}
\end{equation}

\subsubsection{}

Consider the commutative diagram
\begin{equation}\label{eq: bc diagram 1}
\xymatrix{\Fix(C) \ar[r]^{i_{Y}} \ar[d]^{j_{Y}} & C \ar[r]^{f} \ar[d]^{c} & D \ar[d]^{d} \\
Y_0 \ar[r]^{\D_{Y}} & Y_0 \times Y_1 \ar[r]^{f_{0}\times f_{1}} & Z_0 \times Z_1 
}
\end{equation}
The left square is derived Cartesian and the right is pullable. Therefore by Lemma \ref{l:pull 2 out of 3}, the outer square is also pullable.  By Proposition \ref{p: comp pull horizontal}, the pull-pull base change map for the outer square is the composition of the base change maps for the two inner squares
\begin{equation}\label{eq: bc map 1}
i_Y^* f^* d^! \xr{i_{Y}^{*}\tu} i_{Y}^{*}c^{!}(f_{0}\times f_{1})^{*}\tw{-\d+d(f_{0})}\xr{\di(f_{0}\times f_{1})^{*}} j_Y^! \Delta_Y^* (f_0\times f_1)^*\tw{-\d+d(f_{0})}.
\end{equation}

On the other hand, consider the commutative diagram
\begin{equation}\label{eq: bc diagram 2}
\xymatrix{
\Fix(C) \ar[r]^{\Fix(f)} \ar[d]^{{j_{Y}}} &  \Fix(D) \ar[r]^{i_{Z}} \ar[d]^{j_{Z}} & D \ar[d]^{d} \\
Y_0 \ar[r]^{f_{0}} & Z_0  \ar[r]^{\D_{Z}} & Z_0 \times Z_1
}
\end{equation}
where the right square is derived Cartesian and the left square is pullable, since $f_{0}$ is smooth and $\Fix(f)$ is quasi-smooth by Lemma \ref{l: Fix(f) qsm}. Again by Proposition \ref{p: comp pull horizontal}, the pull-pull base change map for the outer square is the composition of the base change maps for the two inner squares
\begin{equation}\label{eq: bc map 2}
\Fix(f)^* i_Z^* d^! \xr{\Fix(f)^{*}\di} \Fix(f)^{*}j_{Z}^{!}\D_{Z}^{*}\xr{\tu\D_{Z}^{*}} j_Y^! f_0^* \Delta_Z^*\tw{-\d+d(f_{0})}\cong j_Y^! f_0^! \Delta_Z^*\tw{-\d}.
\end{equation}

Since the outer  squares of both diagrams \eqref{eq: bc diagram 1} and \eqref{eq: bc diagram 2} are the same, Proposition \ref{p: comp pull horizontal} implies that both \eqref{eq: bc map 1} and \eqref{eq: bc map 2} give the base change map for the same square.  We thus get a commutative diagram
\begin{equation}\label{bc map 12 agree}
\begin{gathered}
\begin{adjustbox}{max width=\textwidth}
$\displaystyle
\xymatrix{i_Y^* f^* d^! \ar@{=}[d]\ar[rr]^-{i_{Y}^{*}\tu} && i_{Y}^{*}c^{!}(f_{0}\times f_{1})^{*}\ar[rr]^-{\di(f_{0}\times f_{1})^{*}}\tw{-\d+d(f_{0})} &&  j_Y^! \Delta_Y^* (f_0\times f_1)^*\tw{-\d+d(f_{0})}\ar@{=}[d]\\
\Fix(f)^* i_Z^* d^! \ar[rr]^-{\Fix(f)^{*}\di} && \Fix(f)^{*}j_{Z}^{!}\D_{Z}^{*}\ar[rr]^-{\tu\D_{Z}^{*}} &&  j_Y^! f_0^* \Delta_Z^*\tw{-\d+d(f_{0})}
}
$
\end{adjustbox}
\end{gathered}
\end{equation}

\subsubsection{Completion of the proof}
Now we may complete the proof of Proposition \ref{prop: trace commutes with smooth pullback}. 

Consider applying \eqref{eq: bc map 1} to an element of the form $\DD(\cK_0) \boxtimes \cK_1 \in D(Z_0 \times Z_1)$: we get a map
\begin{eqnarray}\label{bc map 1 on KK}
i_Y^* f^* d^! (\DD(\cK_0) \boxtimes \cK_1) &\to& i_{Y}^{*}c^{!}(f_{0}^{*}\DD(\cK_{0})\bt f_{1}^{*}\cK_{1})\tw{-\d+d(f_{0})}\\
&\to&  j_Y^! \Delta_Y^* (f_0^*\DD(\cK_{0}) \boxtimes f_1^*\cK_{1})\tw{-\d+d(f_{0})}
\\
&\cong&    j_Y^! (\DD(f_{0}^{*}\cK_0) \ot f_{1}^{*}\cK_1\tw{-\d}).
\end{eqnarray}
Here we use $f_{0}^{*}\DD\tw{d(f_{0})}\cong f_{0}^{!}\DD \cong \DD f_{0}^{*}$. As we have observed in \S\ref{sss: alt pullback corr}, for a cohomological correspondence $\mf{d} \in \Corr_D(\cK_0, \cK_1)$, viewed as a global section of $d^! (\DD(\cK_0) \boxtimes \cK_1)$, $f^{*}\frd$ is the image of $\frd$ under the map on $\upH^{0}$ induced by the base change map \eqref{alt pullback corr}, which is the first step in \eqref{bc map 1 on KK}.

Suppose now that $\cK_1 = \cK_0 \tw{-i}$ (and recall that $f_0 = f_1$). Then we may compose \eqref{eq: bc map 1} with the map $\Tr_{Y_0} \co \DD_{Y_0}(f_0^* \cK_0) \otimes (f_1^* \cK_0 \tw{-i-\d}) \rightarrow \DD_{Y_0} \tw{-i-\d}$ to get a sequence of maps  
\begin{equation}\label{eq: trace of pull}
\begin{gathered}
\begin{aligned}
i_Y^* f^* d^{!}(\DD\cK_{0}\bt\cK_{0}\tw{-i})
&\xrightarrow{\eqref{bc map 1 on KK}}
j_Y^! \bigl(\DD_{Y_0}(f_0^* \cK_0) \otimes (f_1^* \cK_0 \tw{-i-\d})\bigr)\\
&\xrightarrow{\Tr_{Y_0}}  j_Y^! \DD_{Y_0} \tw{-i-\d}
\cong \DD_{\Fix(C)} \tw{-i-\d}.
\end{aligned}
\end{gathered}
\end{equation}
Then by Lemma \ref{lem: trace as BC} and the preceding paragraph, the image of $i_Y^* f^* \mf{d}$ under \eqref{eq: trace of pull} is $\Tr(f^* \mf{d}) \in \mBM_{2i+2\d}(\Fix(C))$.

Next consider applying \eqref{eq: bc map 2} to $\DD(\cK_0) \boxtimes \cK_1 \in D(Z_0 \times Z_1)$: we get a map
\begin{equation}\label{bc map 2 on KK}
\Fix(f)^* i_Z^* d^! (\DD(\cK_0) \boxtimes \cK_1) \to
\Fix(f)^{*}j_{Z}^{!}(\DD(\cK_{0})\otimes\cK_{1}) \to 
  j_Y^! f_0^!  (\DD(\cK_0) \ot \cK_1) \tw{-\d}.
\end{equation}

Suppose now $\cK_1 = \cK_0 \tw{-i}$. Then we may compose \eqref{eq: bc map 2} with the map $\Tr_{Z_0} \co \DD(\cK_0) \otimes \cK_0 \tw{-i} \rightarrow \DD_{Z_{0}}\tw{-i}$ to get a commutative diagram
\begin{equation}\label{eq: pull of trace}
\begin{gathered}
\begin{adjustbox}{max width=\textwidth}
$\displaystyle
\xymatrix{
\Fix(f)^* i_Z^* d^! (\DD(\cK_0) \boxtimes \cK_0\tw{-i}) \ar[r] & \Fix(f)^{*}j_{Z}^{!}(\DD(\cK_{0})\otimes\cK_{0}\tw{-i}) \ar[d]^{ \Fix(f)^{*}j_{Z}^{!}\Tr_{Z_0}} \ar[r] & j_Y^!  f_0^! (\DD(\cK_0) \otimes \cK_0 \tw{-i-\d})\ar[d]^{j_Y^!  f_0^!\Tr_{Z_0}}\\
&    \Fix(f)^{*}j_{Z}^{!}\DD_{Z_{0}}\tw{-i}\ar@{=}[d] \ar[r] &  j_Y^!  f_0^!  \DD_{Z_{0}} \tw{-i-\d}\ar@{=}[d] \\
& \Fix(f)^{*}\DD_{\Fix(D)}\tw{-i} \ar[r]^-{\eqref{Fix(f) pullback}} & \DD_{\Fix(C)} \tw{-i-\d}.
}
$
\end{adjustbox}
\end{gathered}
\end{equation}
Again we view the cohomological correspondence $\mf{d} \in \Corr_D(\cK_0, \cK_0 \tw{-i})$ as a global section of $ d^! (\DD(\cK_0) \boxtimes \cK_0\tw{-i}) $. By Lemma \ref{lem: trace as BC}, the image of this global section under $i_Z^* \Tr_{Z_{0}}$ is $\Tr( \mf{d}) \in \mBM_{2i}(\Fix(D))$. Since the bottom map induces $\Fix(f)^{*}$ on Borel-Moore homology, we get that the image of $\Fix(f)^* i_Z^* \mf{d}$ in the lower right corner of \eqref{eq: pull of trace} is $\Fix(f)^* \Tr( \mf{d})  \in \mBM_{2i+2\d}(\Fix(C))$. 

Finally, the commutativity of \eqref{bc map 12 agree} establishes that \eqref{eq: trace of pull} agrees with the upper-left to lower-right composition of \eqref{eq: pull of trace}, so then $\Tr(f^* \mf{d}) = \Fix(f)^* \Tr( \mf{d})  \in \mBM_{2i+2\d}(\Fix(C))$, as desired. \qed

\subsection{Derived fundamental class as a trace}\label{ssec: VFC as trace}

Let notation be as in \S \ref{ssec: trace of cc}, and set $Y := Y_0 = Y_1$. 

\begin{prop}\label{prop: smooth derived local terms} Assume that $Y$ is smooth over $k$ and $c_1$ is quasi-smooth. Then we have 
\[
\Tr(\cc_Y) = [\Fix(C)] \in \mBM_{2d(c_1)} (\Fix(C)).
\]
\end{prop}

\begin{remark} This result almost appears as \cite[Theorem 1.7]{O15}, but we need more generality than is treated there. We will see that the proof is almost a trivial application of Proposition \ref{prop: trace commutes with smooth pullback}. 
\end{remark}

\begin{proof}
Consider the map of correspondences 
\[
\begin{tikzcd}
Y  \ar[d, "\pi_0"] & \ar[l, "c_0"'] C \ar[r, "c_1"] \ar[d, "\pi"]  & Y \ar[d, "\pi_1"] \\
\pt &  \ar[l] \pt \ar[r] & \pt 
\end{tikzcd}
\]
On the bottom row we have the trivial correspondence $\cc_{\pt} = \Id \in \Corr_{\pt}(\Ql, \Ql)$. By assumption $\pi_1$ is smooth and $c_1$ is quasi-smooth, so $\pi$ is quasi-smooth, hence the diagram is right pullable, and it is immediate from the definitions that 
\[
\pi^* \cc_{\pt} = \cc_Y \in \Corr_C(\Qll{Y}, \Qll{Y}\tw{-d(c_1)}).
\]
Then Proposition \ref{prop: trace commutes with smooth pullback} implies that 
\[
\Tr(\pi^* \cc_{\pt}) = \Fix(\pi)^* \Tr(\cc_{\pt}) = \Fix(\pi)^*[\pt] = [\Fix(C)] \in \mBM_{2d(c_1)}(\Fix(C)).
\]
\end{proof}

\begin{cor}
In the setup of \S \ref{sssec: fix vs sht}, suppose $Y$ is smooth over $\F_q$ and $c_1$ is quasi-smooth. Then we have
\[
 \Tr_{C}^{\Sht}(\cc_Y) = [\Sht(C)] \in \mBM_{2d(c_1)}(\Sht(C)).
 \]
\end{cor}

\subsection{Cohomological co-correspondences}\label{ssec: coco}
Later when we study how cohomological correspondences interact with Fourier duality, we will naturally encounter variants of the above constructions in terms of what we call \emph{co-correspondences} and \emph{cohomological co-correspondences}. 

\subsubsection{Definitions}\label{sssec: cohomological co-correspondence} 
A \emph{co-correspondence} between derived Artin stacks $Y_{0}$ and $Y_{1}$ is a diagram\footnote{We deliberately denote the map $Y_{0}\to C'$ by $c'_{1}$. This is not a typo. This will be justified in the situation when we can complete the diagram into one of the form \eqref{eq: corr and co-corr diagram}, in which we follow the convention that the names of arrows at opposite sides of a square differ by a prime.}
\begin{equation}
\xymatrix{Y_0 \ar[dr]_{c'_{1}}& & Y_{1}\ar[dl]^{c'_{0}}\\
& C'
}
\end{equation}
We define a \emph{cohomological co-correspondence} from $\cK_0 \in D(Y_0)$ to $\cK_1 \in D(Y_1)$ to be an element of $\Hom_{C'} (c'_{1!} \cK_0, c'_{0*} \cK_1)$. Let
\begin{equation}
\CoCorr_{C'}(\cK_{0}, \cK_{1}):=\Hom_{C'} (c'_{1!} \cK_0, c'_{0*} \cK_1).
\end{equation}

To see the relation between cohomological correspondences and co-correspondences, suppose we have a Cartesian square

\begin{equation}\label{eq: corr and co-corr diagram}
\begin{tikzcd}
& C^{\flat} \ar[dl, "c_0"'] \ar[dr, "c_1"] \\
Y_0 \ar[dr, "c_1'"']  & & Y_1 \ar[dl, "c_0'"] \\
& C^{\sharp}
\end{tikzcd}
\end{equation}
Then for $\cK_0 \in D(Y_0)$ and $\cK_1 \in D(Y_1)$, there is a canonical isomorphism of vector spaces
\begin{equation}\label{eq: coh corr to coh co-corr}
\g_{C}: \Corr_{C^{\flat}}(\cK_{0}, \cK_{1})\isom \CoCorr_{C^{\sharp}}(\cK_{0}, \cK_{1})
\end{equation}
given by adjunctions and proper base change
\begin{equation*}
\begin{gathered}
\begin{aligned}
\Hom_{C^{\flat}}(c_0^* \cK_0 , c_1^! \cK_1)
&\cong \Hom_{Y_1}(c_{1!} c_0^* \cK_0 , \cK_1)\\
&\cong  \Hom_{Y_1}((c_0')^* c'_{1!}  \cK_0 , \cK_1)\\
&\cong  \Hom_{C^{\sharp}} (c'_{1!} \cK_0, c'_{0*} \cK_1).
\end{aligned}
\end{gathered}
\end{equation*}


\subsubsection{Pushforward and pullback of cohomological co-correspondences}\label{sssec: pull push coco}
We have dual notions of pushability and pullability for co-correspondences.  Consider a morphism of co-correspondences, i.e., a commutative diagram
\begin{equation}\label{map of coco}
\xymatrix{Y_0\ar[d]_{f_{0}}\ar[dr]_{c'_{1}} & & Y_{1}\ar[dl]^{c'_{0}}\ar[d]^{f_{1}}\\
Z_{0}\ar[dr]_{d'_{1}}  & C'\ar[d]^{f'} & Z_{1}\ar[dl]^{d'_{0}} \\
& D' 
}
\end{equation}

\begin{defn} The diagram of co-correspondences \eqref{map of coco} is called {\em left pullable}, if the square with vertices $(Y_{0}, C', Z_{0}, D')$ is pullable in the sense of Definition \ref{def: pushable pullable}. We define the defect $\d_{f'}$ to be the defect of the square $(Y_{0}, C', Z_{0}, D')$.

Similarly, \eqref{map of coco} is called {\em right pushable}, if the square with vertices $(Y_{1}, C', Z_{1}, D')$ is pushable in the sense of Definition \ref{def: pushable pullable}.
\end{defn}

When \eqref{map of coco} is left pullable, we have a pullback map of cohomological co-correspondences (where $\cK_{i}\in D(Z_{i})$)
\begin{equation}
(f')^{*}: \CoCorr_{D'}(\cK_{0}, \cK_{1})\to \CoCorr_{C'}(f_{0}^{*}\cK_{0}, f_{1}^{*}\cK_{1}\tw{-\d_{f'}})
\end{equation}
defined as follows. For $\cc': d'_{1!}\cK_{0}\to d'_{0*}\cK_{1}$, we define $(f')^{*}(\cc')$ as the composition
\begin{equation}
c'_{1!}f_{0}^{*}\cK_{0}\xr{\tu}(f')^{*}d'_{1!}\cK_{0}\tw{-\d_{f'}}\xr{(f')^{*}(\cc')}(f')^{*}d'_{0*}\cK_{1}\tw{-\d_{f'}}\to c'_{0*}f_{1}^{*}\cK_{1}\tw{-\d_{f'}}
\end{equation}
where the first map is the natural transformation \eqref{eq: p-p}. 

When \eqref{map of coco} is right pushable, we have a pushforward map of cohomological co-correspondences (where $\cK_{i}\in D(Y_{i})$)
\begin{equation}
f'_{!}: \CoCorr_{C'}(\cK_{0}, \cK_{1})\to \CoCorr_{D'}(f_{0!}\cK_{0}, f_{1!}\cK_{1})
\end{equation}
defined as follows. For $\cc': c'_{1!}\cK_{0}\to c'_{0*}\cK_{1}$, we define $f'_{!}(\cc')$ as the composition
\begin{equation}
d'_{1!}f_{0!}\cK_{0}\cong f'_{!}c'_{1!}\cK_{0}\xr{f'_{!}(\cc')}f'_{!}c'_{0*}\cK_{1}\to d'_{0*}f_{1!}\cK_{1}
\end{equation}
where the last map is the push-push natural transformation \eqref{eq: push-push gen}.

\subsubsection{Compatibility}

Suppose we are given a commutative diagram
\begin{equation}\label{eq: co-corr to corr pull/push}
\xymatrix{&  C^{\flat} \ar[dl]_{c_{0}}   \ar[dr]^{c_{1}} \ar@/^1pc/[ddd]_{f^{\flat}}\\ 
 Y_0 \ar[dr]_{c'_{1}} \ar[ddd]_{f_{0}} && Y_1 \ar[dl]^{c'_{0}}   \ar[ddd]^{f_{1}} \\ 
 & C^{\sharp} \ar@/_1pc/[ddd]^{f^{\sharp}}\\
&  D^{\flat} \ar[dl]_{d_{0}}   \ar[dr]^{d_{1}} \\ 
 Z_0 \ar[dr]_{d'_{1}} && Z_1 \ar[dl]^{d'_{0}}  \\ 
 & D^{\sharp} 
} 
\end{equation}
where the top and bottom diamonds are derived Cartesian. We view $f^{\flat}: C^{\flat}\to D^{\flat}$ as a map of correspondences, and $f^{\sh}: C^{\sh}\to D^{\sh}$ as a map of co-correspondences.

\begin{prop}\label{prop: comp coco} In the above situation,
\begin{enumerate}
\item If $f^{\sh}$ is left pullable, then $f^{\flat}$ is right pullable, and $\d_{f^{\sh}}=\d_{f^{\flat}}$. In this case, we have a commutative diagram for any $\cK_{0}\in D(Z_{0})$ and $\cK_{1}\in D(Z_{1})$
\begin{equation}\label{pull coco and co}
\xymatrix{\CoCorr_{D^{\sh}}(\cK_{0}, \cK_{1})\ar[r]^{\g_{D}^{-1}}\ar[d]^{(f^{\sh})^{*}} & \Corr_{D^{\flat}}(\cK_{0}, \cK_{1})\ar[d]^{(f^{\flat})^{*}}\\
\CoCorr_{C^{\sh}}(f_{0}^{*}\cK_{0}, f_{1}^{*}\cK_{1}\tw{-\d_{f^{\sh}}})\ar[r]^{\g_{C}^{-1}}
& \Corr_{C^{\flat}}(f_{0}^{*}\cK_{0}, f_{1}^{*}\cK_{1}\tw{-\d_{f^{\flat}}})
}
\end{equation}
\item If $f^{\sh}$ is right pushable, then $f^{\flat}$ is left pushable. In this case, we have a commutative diagram for any $\cK_{0}\in D(Y_{0})$ and $\cK_{1}\in D(Y_{1})$
\begin{equation}
\xymatrix{\CoCorr_{C^{\sh}}(\cK_{0}, \cK_{1})\ar[r]^{\g_{C}^{-1}}\ar[d]^{f^{\sh}_{!}} & \Corr_{C^{\flat}}(\cK_{0}, \cK_{1})\ar[d]^{f^{\flat}_{!}}\\
\CoCorr_{D^{\sh}}(f_{0!}\cK_{0}, f_{1!}\cK_{1})\ar[r]^{\g_{D}^{-1}}
& \Corr_{D^{\flat}}(f_{0!}\cK_{0}, f_{1!}\cK_{1})
}
\end{equation}
\end{enumerate}
\end{prop}
\begin{proof}
(1) Using that the top diamond in \eqref{eq: co-corr to corr pull/push} is derived Cartesian, the fact that the square $(Y_{0}, C^{\sh}, Z_{0}, D^{\sh})$ is pullable implies that the square $(C^{\flat}, Y_{1}, Z_{0}, D^{\sh})$ is pullable with the same defect. Using that the bottom diamond is derived Cartesian, we have $Y_{1}\times_{Z_{1}}D^{\flat}\isom Y_{1}\times_{D^{\sh}}Z_{0}$, from which we see that the square $(C^{\flat}, Y_{1}, D^{\flat}, Z_{1})$ is pullable with the same defect as $(C^{\flat}, Y_{1}, Z_{0}, D^{\sh})$.

We prove the commutativity of the diagram \eqref{pull coco and co}. Let $\cc': d'_{1!}\cK_{0}\to d'_{0*}\cK_{1}$ be an element in $\CoCorr_{D^{\sh}}(\cK_{0}, \cK_{1})$. Consider the following diagram
\begin{equation}
\begin{gathered}
\begin{adjustbox}{max width=\textwidth}
$\displaystyle
\xymatrix{c_{1!}c_{0}^{*}f_{0}^{*}\cK_{0}\ar[r]^-{\di}\ar@{=}[d] & (c'_{0})^{*}c'_{1!}f_{0}^{*}\cK_{0}\ar[r]^-{(c'_{0})^{*}\tu} &  (c'_{0})^{*}(f^{\sh})^{*}d'_{1!}\cK_{0}\tw{-\d_{f^{\sh}}} \ar[r]^-{(c'_{0})^{*}(f^{\sh})^{*}\cc'}\ar@{=}[d] & (c'_{0})^{*}(f^{\sh})^{*}d'_{0*}\cK_{1}\tw{-\d_{f^{\sh}}} \ar[d]\\
c_{1!}(f^{\flat})^{*}d_{0}^{*}\cK_{0} \ar[r]^{\tu d_{0}^{*}} & f_{1}^{*}d_{1!}d_{0}^{*}\cK_{0}\tw{-\d_{f^{\flat}}}\ar[r]^-{\di} & f_{1}^{*}(d'_{0})^{*}d'_{1!}\cK_{0}\tw{-\d_{f^{\flat}}} \ar[r]^-{f_{1}^{*}\cc''} & f_{1}^{*}\cK_{1}\tw{-\d_{f^{\flat}}}}
$
\end{adjustbox}
\end{gathered}
\end{equation}
The right vertical map comes from the identity $(c'_{0})^{*}(f^{\sh})^{*}=f^{*}_{1}(d'_{0})^{*}$ by adjunction. The map $\cc'':(d'_{0})^{*}d'_{1!}\cK_{0}\to \cK_{1}$ that appears on the bottom row is obtained from $\cc'$ by adjunction.

In the above diagram, the two ways of getting from $c_{1!}c_{0}^{*}f_{0}^{*}\cK_{0}$ to $f_{1}^{*}\cK_{1}\tw{-\d_{f^{\flat}}}$, by going right and down, and by going down and right, are obtained from $\g_{C}^{-1}(f^{\sh})^{*}(\cc')$ and $(f^{\flat})^{*}\g_{D}^{-1}(\cc')$ by adjunction respectively. Therefore it suffices to show that the above diagram is commutative. The equality in the middle divides the diagram into two parts. The right square is tautologically commutative. The top and bottom rows of the left part both induce base change maps for the square $(C^{\flat}, Y_{1}, Z_{0}, D^{\sh})$
\begin{equation}
c_{1!}(f_{0}\c c_{0})^{*}=c_{1!}(d_{0}\c f^{\flat})^{*}\to (d'_{0}\c f_{1})^{*}d'_{1!}\tw{-\d_{f^{\flat}}}=(f^{\sh}\c c'_{0})^{*}d'_{1!}\tw{-\d_{f^{\sh}}}
\end{equation}
by decomposing it in two different ways
\begin{equation}
\xymatrix{C^{\flat}\ar[d]_{c_{0}}\ar[r]^{c_{1}} & Y_{1}\ar[d]^{c'_{0}} & & C^{\flat}\ar[d]_{f^{\flat}}\ar[r]^{c_{1}}  & Y_{1}\ar[d]^{f_{1}}\\
Y_{0}\ar[d]_{f_{0}}\ar[r]^{c'_{1}}  & C^{\sh}\ar[d]^{f^{\sh}} & & D^{\flat}\ar[d]_{d_{0}}\ar[r]^{d_{1}}  & Z_{1}\ar[d]^{d'_{0}}\\
Z_{0}\ar[r]^{d'_{1}}  & D^{\sh} & & Z_{0}\ar[r]^{d'_{1}}  & D^{\sh}}
\end{equation}
The two base change maps agree by Proposition \ref{p: comp pull}, after adjunction.

(2) The proof is similar. We omit it as the statement is not used in the sequel.
\end{proof}

\section{Base change for cohomological correspondences}\label{sec: descent}
In this section we formulate and prove the \emph{Base Change Theorem} for cohomological correspondences (Theorem \ref{thm: descent of pushforward correspondence}). This will be used in Step (4) of the outline in \S \ref{ssec: outline}.

\subsection{The setup}\label{ssec: coh corr for descent}

\subsubsection{Correspondences}\label{sssec: descent assumptions} 

Suppose we are given a commutative diagram of derived Artin stacks, all of whose morphisms are separated and representable in derived schemes,
\begin{equation}\label{eq: two cubes}
\xymatrix{ U_{0}\ar[ddd]^{\pi_{0}}\ar[dr]^{f_{0}} && \ar[ll]_-{a_{0}}\ar[rr]^-{a_{1}}C_{U}\ar[ddd]^{\pi}\ar[dr]^{f} && U_{1}\ar[ddd]^{\pi_{1}}\ar[dr]^{f_{1}} \\
& V_{0}\ar[ddd]^{g_{0}} && \ar[ll]_-{b_{0}}\ar[rr]^-{b_{1}} C_{V}\ar[ddd]^{g} && V_{1}\ar[ddd]^{g_{1}}\\
\\
S_{0}\ar[dr]^{z_{0}} && \ar[ll]^-{h_{0}}\ar[rr]_-{h_{1}} C_{S}\ar[dr]^{z} && S_{1}\ar[dr]^{z_{1}}\\
& W_{0} && \ar[ll]_-{c_{0}}\ar[rr]^-{c_{1}} C_{W} && W_{1} 
}
\end{equation}
satisfying the following conditions:  
\begin{enumerate}
\item The middle vertical parallelogram
\begin{equation}
\xymatrix{C_{U}\ar[ddd]^{\pi}\ar[dr]^{f}\\
& C_{V}\ar[ddd]^{g}\\
\\
C_{S}\ar[dr]^{z}\\
& C_{W}
}
\end{equation}
is derived Cartesian.

\item The three squares in the following diagram are pushable
\begin{equation}\label{eq:  UVW derived cartesian}
\xymatrix{U_{0}\ar[ddd]^{\pi_{0}}\ar[dr]^{f_{0}} && \ar[ll]_{a_{0}}C_{U}\ar[dr]^{f} && \\
& V_{0}\ar[ddd]^{g_{0}} &&  \ar[ll]_{b_{0}}C_{V}&& \\
\\
S_{0}\ar[dr]^{z_{0}} && \ar[ll]_(.6){h_{0}}C_{S}\ar[dr]^{z} &&\\
& W_{0} && \ar[ll]_{c_{0}}C_{W} && 
}
\end{equation}

\item The three squares in the following diagram are pullable
\begin{equation}
\xymatrix{&& \ar[rr]^-{a_{1}}C_{U}\ar[ddd]^{\pi} && U_{1}\ar[ddd]^{\pi_{1}}\ar[dr]^{f_{1}} \\
&&& \ar[rr]^-{b_{1}} C_{V}\ar[ddd]^{g} && V_{1}\ar[ddd]^{g_{1}}\\
\\
&& \ar[rr]^(.4){h_{1}} C_{S} && S_{1}\ar[dr]^{z_{1}}\\
&&& \ar[rr]^-{c_{1}} C_{W} && W_{1} 
}
\end{equation}
Moreover, the right square $(U_{1}, V_{1}, S_{1}, W_{1})$ above has defect zero.
\end{enumerate}

We view $C_{S}$ as a correspondence between $S_{0}$ and $S_{1}$, and similarly for $C_{U}, C_{V}$ and $C_{W}$. 

\subsubsection{Push and pull of cohomological correspondences}
Let $\cK_{i}\in D(S_{i})$ for $i=0,1$. Let $\frs\in \Corr_{C_{S}}(\cK_{0}, \cK_{1})$.

Consider the back face of the diagram \eqref{eq: two cubes}, viewed as a map of correspondences $\pi: C_{U}\to C_{S}$:
\begin{equation}
\xymatrix{U_0 \ar[d]^{\pi_{0}} &  C_{U}\ar[l]_{a_{0}}\ar[r]^{a_{1}}  \ar[d]^{\pi}  &  U_{1} \ar[d]^{\pi_{1}} \\
S_{0} & \ar[l]_{h_{0}} C_S\ar[r]^{h_{1}} & S_{1}
}
\end{equation}
By assumption, this map of correspondences is right pullable. Therefore, by \S\ref{ssec: pullback functoriality for CC}, the map
\begin{eqnarray}
\pi^{*}: \Corr_{C_{S}}(\cK_{0}, \cK_{1})\to \Corr_{C_{U}}(\pi_{0}^{*}\cK_{0}, \pi_{1}^{*}\cK_{1}\tw{-\d_{\pi}})
\end{eqnarray}
is defined (where the defect $\d_{\pi}$ is defined in Definition \ref{def: corr pullable}).

Consider the top face of the diagram \eqref{eq: two cubes}, viewed as a map of correspondences $f: C_{U}\to C_{V}$:
\begin{equation}
\xymatrix{U_0 \ar[d]^{f_{0}} &  C_{U}\ar[l]_{a_{0}}\ar[r]^{a_{1}}  \ar[d]^{f}  &  U_{1} \ar[d]^{f_{1}} \\
V_{0} & \ar[l]_{b_{0}} C_V\ar[r]^{b_{1}} & V_{1}
}
\end{equation}
By assumption, this map of correspondences is left pushable. Therefore,  by \S\ref{ssec: pushforward functoriality for CC}, the map
\begin{eqnarray}
f_{!}: \Corr_{C_{U}}(\pi^{*}_{0}\cK_{0}, \pi_{1}^{*}\cK_{1}\tw{-\d_{\pi}})\to \Corr_{C_{V}}(f_{0!}\pi_{0}^{*}\cK_{0}, f_{1!}\pi_{1}^{*}\cK_{1}\tw{-\d_{\pi}})
\end{eqnarray}
is defined. 

The composition of the two maps give an element
\begin{equation}
f_{!}\pi^{*}(\frs)\in \Corr_{C_{V}}(f_{0!}\pi_{0}^{*}\cK_{0}, f_{1!}\pi_{1}^{*}\cK_{1}\tw{-\d_{\pi}}).
\end{equation}

Similarly, the bottom face of the diagram \eqref{eq: two cubes}
\begin{equation}
\xymatrix{S_0 \ar[d]^{z_{0}} &  C_{S}\ar[l]_{h_{0}}\ar[r]^{h_{1}}  \ar[d]^{z}  &  S_{1} \ar[d]^{z_{1}} \\
W_{0} & \ar[l]_{c_{0}} C_W\ar[r]^{c_{1}} & W_{1}
}
\end{equation}
is left pushable and the front face
\begin{equation}
\xymatrix{V_0 \ar[d]^{g_{0}} &  C_{V} \ar[l]_{b_{0}}\ar[r]^{b_{1}}  \ar[d]^{g}  &  V_{1}\ar[d]^{g_{1}} \\
W_{0} & \ar[l]_{c_{0}} C_W\ar[r]^{c_{1}} & W_{1}
}
\end{equation}
is right pullable. Therefore the cohomological correspondence
\begin{equation}
g^{*}z_{!}(\frs)\in \Corr_{C_{V}}(g_{0}^{*}z_{0!}\cK_{0}, g_{1}^{*}z_{1!}\cK_{1}\tw{-\d_{g}})
\end{equation}
is defined. 

\subsubsection{Matching source and target}
We will now formulate a Base Change Theorem for cohomological correspondences, which is the main result of this section. The Base Change Theorem asserts that, referring to the diagram \eqref{eq: two cubes}, for any $\frs\in \Corr_{C_{S}}(\cK_{0}, \cK_{1})$, the cohomological correspondences $f_{!}\pi^{*}(\frs)$ and $g^{*}z_{!}(\frs)$ on $C_{V}$ ``agree''. To make sense of this statement, we need to relate the source and targets of the respective cohomological correspondences. 

By assumption, the square $(U_{0},V_{0}, S_{0}, W_{0})$ in \eqref{eq: two cubes} is pushable. We get a base change map
\begin{equation}\label{eq: match source}
g^{*}_{0}z_{0!}\xr{\td} f_{0!}\pi_{0}^{*}: D(S_{0})\to D(V_{0}).
\end{equation}
By assumption, the square $(U_{1},V_{1}, S_{1}, W_{1})$ in \eqref{eq: two cubes} is pullable with defect zero. We get a pull-pull map
\begin{equation}
\pi_{1}^{*}z_{1}^{!}\xr{\tu} f_{1}^{!}g^{*}_{1}: D(W_{1})\to D(U_{1}).
\end{equation}
By adjunction (cf. Remark \ref{rem: p-p}), it gives a base change map
\begin{equation}\label{eq: match target}
 f_{1!}\pi_{1}^{*}\to g^{*}_{1}z_{1!}: D(S_{1})\to D(V_{1}).
\end{equation}

\begin{lemma}\label{l:match defects}
We have an equality of defects (Definition \ref{def: corr pullable}) $\d_{\pi}=\d_{g}$.
\end{lemma}
\begin{proof}
By Lemma \ref{l:pull 2 out of 3}, both defects are equal to the defect of the pullable square $(C_{U}, V_{1}, C_{S}, W_{1})$, using that both $(C_{U}, C_{V}, C_{S}, C_{W})$ and $(U_{1}, V_{1}, S_{1}, W_{1})$ have defect zero.
\end{proof}

\begin{example}
A special case (which will be our case of interest) is when both $(U_{0},V_{0}, S_{0}, W_{0})$ and $(U_{1},V_{1}, S_{1}, W_{1})$ are derived Cartesian. In this case, the sources and targets of $f_{!}\pi^{*}(\frs)$ and $g^{*}z_{!}(\frs)$ are matched by the proper base change isomorphisms
\begin{equation}\label{eq: match source target Cart}
f_{0!}\pi_{0}^{*}\cK_{0}\xr{\di} g^{*}_{0}z_{0!}\cK_{0}, \quad f_{1!}\pi_{1}^{*}\cK_{1}\tw{-\d_{\pi}}\xr{\di} g^{*}_{1}z_{1!}\cK_{1}\tw{-\d_{g}}.
\end{equation}
\end{example}

\begin{thm}[Base Change Theorem for cohomological correspondences]\label{thm: descent of pushforward correspondence}Under the assumptions in \S\ref{sssec: descent assumptions}, for any $\frs\in \Corr_{C_{S}}(\cK_{0}, \cK_{1})$, the following diagram is commutative
\begin{equation}
\xymatrix{g^{*}_{0}z_{0!}\cK_{0}\ar[rr]^-{g^{*}z_{!}(\frs)}\ar[d]^{\eqref{eq: match source}} & & g^{*}_{1}z_{1!}\cK_{1}\tw{-\d_{g}} \\
f_{0!}\pi_{0}^{*}\cK_{0} \ar[rr]^-{f_{!}\pi^{*}(\frs)} && f_{1!}\pi_{1}^{*}\cK_{1}\tw{-\d_{\pi}}\ar[u]_{\eqref{eq: match target}}
}
\end{equation}
Here we use Lemma \ref{l:match defects} to match the twists. 

In particular, when both $(U_{0},V_{0}, S_{0}, W_{0})$ and $(U_{1},V_{1}, S_{1}, W_{1})$ are derived Cartesian, we have an equality of cohomological correspondences on $C_{V}$
\begin{equation}
f_{!}\pi^{*}(\frs)=g^{*}z_{!}(\frs)
\end{equation}
under the isomorphisms \eqref{eq: match source target Cart}.
\end{thm}

\begin{remark} In applications to the modularity theorem, we are interested in the special case where both $(U_{0}, V_{0}, S_{0}, W_{0})$ and $(U_{1}, V_{1}, S_{1}, W_{1})$ are derived Cartesian, $h_{1}$ is quasi-smooth, $\cK_{0}=\Qll{S_{0}}$, $\cK_{1}=\Qll{S_{1}}\tw{-d(h_{1})}$, and the map $\frs: h_{0}^{*}\cK_{0}=\Qll{C_S}\to h_{1}^{!}\cK_{1}=h_{1}^{!}\Qll{S_{1}}\tw{-d(h_{1})}$ is given by the relative fundamental class $[h_{1}]$.
\end{remark}

\subsection[Proof of descent of pushforward correspondence]{Proof of Theorem \ref{thm: descent of pushforward correspondence}}
Unravelling the constructions of $g^{*}z_{!}(\frs)$ and $f_{!}\pi^{*}(\frs)$, we see that they appear as the top and bottom rows of the following diagram
\begin{equation}\label{eq: proof des}
\begin{gathered}
\begin{adjustbox}{max width=\textwidth}
$\displaystyle
\xymatrix{b_{0}^{*}g_{0}^{*}z_{0!}\cK_{0}\ar[d]^{b_{0}^{*}\td}\ar@{=}[r] & g^{*}c_{0}^{*}z_{0!}\cK_{0}\ar[r]^{g^{*}\td} & g^{*}z_{!}h_{0}^{*}\cK_{0} \ar[r]^{g^{*}z_{!}\frs}\ar[d]^{\di h_{0}^{*}} &  g^{*}z_{!}h_{1}^{!}\cK_{1}\ar[r] & g^{*}c_{1}^{!}z_{1!}\cK_{1}\ar[r]^{\tu z_{1!}} &  b_{1}^{!}g_{1}^{*}z_{1!}\cK_{1}\tw{-\d_{g}}\\
b_{0}^{*}f_{0!}\pi_{0}^{*}\cK_{0}\ar[r]^{\td\pi_{0}^{*}} & f_{!}a_{0}^{*}\pi_{0}^{*}\cK_{0}\ar@{=}[r] & f_{!}\pi^{*}h_{0}^{*}\cK_{0}\ar[r]^{f_{!}\pi^{*}\frs} & f_{!}\pi^{*}h_{1}^{!}\cK_{1} \ar[r]^{f_{!}\tu}\ar[u]_{\di h_{1}^{!}}
 & f_{!}a_{1}^{!}\pi_{1}^{*}\cK_{1}\tw{-\d_{\pi}}\ar[r] & b_{1}^{!} f_{1!}\pi_{1}^{*}\cK_{1}\tw{-\d_{\pi}}\ar[u]^{b_{1}^{!}\tu}
}
$
\end{adjustbox}
\end{gathered}
\end{equation}
Here, the arrows marked by $\td$ are the base change maps obtained from a pushable square as in \S\ref{ssec: push-pull}; the arrows marked by $\tu$ are the base change maps obtained from a pullable square as in \S\ref{ssec: pull-pull}, and the arrows marked by $\di$ are the proper base change isomorphisms, which are special cases of both $\td$ and $\tu$. The unmarked arrows are the tautological base change maps from commutative squares.

To prove the theorem, we need to check that all three rectangles in \eqref{eq: proof des} commute.  The middle square is clearly commutative. Below we check separately that the left and right sleeves commute.

\subsubsection{Left sleeve}
We need to show that the diagram of natural transformations
\begin{equation}\label{eq: left rec des}
\xymatrix{b_{0}^{*}g_{0}^{*}z_{0!}\ar[d]^{b_{0}^{*}\td}\ar@{=}[r] & g^{*}c_{0}^{*}z_{0!}\ar[r]^{g^{*}\td} & g^{*}z_{!}h_{0}^{*}\ar[d]^{\td h_{0}^{*}}\\
b_{0}^{*}f_{0!}\pi_{0}^{*}\ar[r]^{\td\pi_{0}^{*}} & f_{!}a_{0}^{*}\pi_{0}^{*}\ar@{=}[r] & f_{!}\pi^{*}h_{0}^{*}
}
\end{equation}
is commutative. Here we have replaced the $\di$ with the $\td$ in the right vertical arrow because the proper base change isomorphism $g^{*}z_{!}\isom f_{!}\pi^{*}$ is a special case of the $\td$ map from a pushable square. 

Let 
\begin{eqnarray}
\a_{0}:=\pi_{0}\c a_{0}=h_{0}\c \pi: C_{U}\to S_{0},\\
\b_{0}:=g_{0}\c b_{0}=c_{0}\c g: C_{V}\to W_{0}.
\end{eqnarray}

The left cube in the diagram \eqref{eq: two cubes} provides two decompositions of the commutative square
\begin{equation}\label{eq: left diag push}
\xymatrix{C_{U}\ar[r]^{f} \ar[d]^{\a_{0}} & C_{V} \ar[d]^{\b_{0}} \\
S_0 \ar[r]^{z_{0}} & W_0
}
\end{equation}

The first is
\begin{equation}
\begin{tikzcd}
C_{U} \ar[r, "f"]   \ar[d, "a_0"] & C_{V} \ar[d, "b_0"]  \\
U_0 \ar[r, "f_0"] \ar[d, "\pi_0"] & V_0 \ar[d, "g_0"] \\
S_0 \ar[r, "z_0"]  & W_0 
\end{tikzcd}
\end{equation}
in which both the upper and lower squares are pushable. After flipping the two squares about the diagonal, Proposition \ref{p: comp push horizontal} shows that the natural transformation $\td$ for the square \eqref{eq: left diag push} agrees with the composition
\begin{equation}
\xymatrix{\b^{*}_{0}z_{0!} \ar@{=}[r] & b_{0}^{*}g_{0}^{*}z_{0!} \ar[r]^{b_{0}^{*}\td} & b_{0}^{*}f_{0!}\pi_{0}^{*} \ar[r]^{\td\pi_{0}^{*}} & f_{!}a_{0}^{*}\pi_{0}^{*}\ar@{=}[r] &  f_{!}\a_{0}^{*}}.
\end{equation}
This is the lower composition of the diagram \eqref{eq: left rec des}.

The second decomposition of \eqref{eq: left diag push} is
\begin{equation}
\begin{tikzcd}
C_{U} \ar[r, "f"] \ar[d, "\pi"] & C_{V}\ar[d, "g"] \\
C_{S} \ar[r, "z"] \ar[d, "h_0"] & C_{W} \ar[d, "c_0"] \\
S_0 \ar[r, "z_0"] &  W_0  
\end{tikzcd}
\end{equation}
in which both the upper and lower squares are pushable. After flipping the two squares about the diagonal, Proposition \ref{p: comp push horizontal} again shows that the natural transformation $\td$ for the square \eqref{eq: left diag push} agrees with the composition
\begin{equation}
\xymatrix{\b^{*}_{0}z_{0!} \ar@{=}[r] & g^{*}c_{0}^{*}z_{0!} \ar[r]^{g^{*}\td} & g^{*}z_{!}h_{0}^{*} \ar[r]^{\td h_{0}^{*}} & f_{!}\pi^{*}h_{0}^{*}\ar@{=}[r] & f_{!}\a^{*}_{0}.
}
\end{equation}
This is the upper composition of the diagram \eqref{eq: left rec des}. Therefore both compositions in \eqref{eq: left rec des} compute the same base change map $\td: \b^{*}_{0}z_{0!}\to f_{!}\a^{*}_{0}$. This proves that \eqref{eq: left rec des} is commutative.

\subsubsection{Right sleeve}
We need to show that the diagram of natural transformations
\begin{equation}\label{eq: right rec des}
\xymatrix{g^{*}z_{!}h_{1}^{!}\ar[r] & g^{*}c_{1}^{!}z_{1!}\ar[r]^-{\tu z_{1!}} &  b_{1}^{!}g_{1}^{*}z_{1!}\tw{-\d_{g}}\\
f_{!}\pi^{*}h_{1}^{!} \ar[r]^-{f_{!}\tu}\ar[u]_{\tu h_{1}^{!}}
 & f_{!}a_{1}^{!}\pi_{1}^{*}\tw{-\d_{\pi}}\ar[r] & b_{1}^{!} f_{1!}\pi_{1}^{*}\tw{-\d_{\pi}}\ar[u]^{b_{1}^{!}\tu}
}
\end{equation}
is commutative. Here we have replaced the $\di$ with the $\tu$ in the left vertical arrow because the proper base change isomorphism $ f_{!}\pi^{*}\isom g^{*}z_{!}$ is a special case of the $\tu$ map from a pullable square. 

Let 
\begin{eqnarray}
\a_{1}:=f_{1}\c a_{1}=b_{1}\c f: C_{U}\to V_{1},\\
\b_{1}:=z_{1}\c h_{1}=c_{1}\c z: C_{S}\to W_{1}.
\end{eqnarray}

The right cube in the diagram \eqref{eq: two cubes} provides two decompositions of the commutative square
\begin{equation}\label{eq: right diag pull}
\xymatrix{C_{U}\ar[r]^{\pi} \ar[d]^{\a_{1}} & C_{S} \ar[d]^{\b_{1}} \\
V_{1}\ar[r]^{g_{1}} & W_{1}
}
\end{equation}

The first is
\begin{equation}
\begin{tikzcd}
C_{U} \ar[r, "\pi"]   \ar[d, "a_{1}"] & C_{S} \ar[d, "h_{1}"]  \\
U_{1} \ar[r, "\pi_{1}"] \ar[d, "f_{1}"] & S_{1} \ar[d, "z_{1}"] \\
V_{1}\ar[r, "g_{1}"]  & W_{1}
\end{tikzcd}
\end{equation}
in which both the upper and lower squares are pullable. After flipping the two squares about the diagonal, Proposition \ref{p: comp pull horizontal} shows that the natural transformation $\tu$ for the square \eqref{eq: right diag pull} agrees with the composition
\begin{equation}
\xymatrix{\pi^{*}\b_{1}^{!}\ar@{=}[r] & \pi^{*}h_{1}^{!}z_{1}^{!} \ar[r]^{\tu z_{1}^{!}} & a_{1}^{!}\pi_{1}^{*}z_{1}^{!}\tw{-\d_{\pi}}\ar[r]^{a_{1}^{!}\tu} & a_{1}^{!}f_{1}^{!}g_{1}^{*}\tw{-\d_{\pi}}\ar@{=}[r] & \a_{1}^{!}g_{1}^{*}\tw{-\d_{\pi}}}
\end{equation}

The second decomposition of \eqref{eq: right diag pull} is 
\begin{equation}
\begin{tikzcd}
C_{U} \ar[r, "\pi"]   \ar[d, "f"] & C_{S} \ar[d, "z"]  \\
C_{V} \ar[r, "g"] \ar[d, "b_{1}"] & C_{W} \ar[d, "c_{1}"] \\
V_{1}\ar[r, "g_{1}"]  & W_{1}
\end{tikzcd}
\end{equation}
in which both the upper and lower squares are pullable. After flipping the two squares about the diagonal, Proposition \ref{p: comp pull horizontal} again shows that the natural transformation $\tu$ for the square \eqref{eq: right diag pull} agrees with the composition
\begin{equation}
\xymatrix{\pi^{*}\b_{1}^{!}\ar@{=}[r] & \pi^{*}z^{!}c_{1}^{!} \ar[r]^{\tu c_{1}^{!}} & f^{!}g^{*}c_{1}^{!}\ar[r]^{f^{!}\tu} & f^{!}b_{1}^{!}g_{1}^{*}\tw{-\d_{g}}\ar@{=}[r] & \a_{1}^{!}g_{1}^{*}\tw{-\d_{g}}}
\end{equation}
Combining the two expressions of $\tu: \pi^{*}\b_{1}^{!}\to \a_{1}^{!}g_{1}^{*}\tw{-\d_{g}}=\a_{1}^{!}g_{1}^{*}\tw{-\d_{\pi}}$, we get a commutative diagram
\begin{equation}\label{eq: two pull-pull}
\xymatrix{\pi^{*}z^{!}c_{1}^{!}\ar@{=}[d]\ar[r]^-{\tu c_{1}^{!}} & f^{!}g^{*}c_{1}^{!}\ar[r]^-{f^{!}\tu} & f^{!}b_{1}^{!}g_{1}^{*}\tw{-\d_{g}}\ar@{=}[d]\\
\pi^{*}h_{1}^{!}z_{1}^{!}\ar[r]^{\tu z_{1}^{!}} & a^{!}_{1}\pi_{1}^{*}z_{1}^{!}\tw{-\d_{\pi}}\ar[r]^-{a_{1}^{!}\tu} & a_{1}^{!}f_{1}^{!}g_{1}^{*}\tw{-\d_{\pi}}
}
\end{equation}

Compare \eqref{eq: two pull-pull} with \eqref{eq: right rec des}. Starting in both diagrams from the lower left corner, going upward then turning right to arrive at the upper right corner, we see the two maps 
\begin{eqnarray}
\pi^{*}h_{1}^{!}z_{1}^{!}\to f^{!}b_{1}^{!}g_{1}^{*}\tw{-\d_{g}},\\
f_{!}\pi^{*}h_{1}^{!}\to b_{1}^{!}g_{1}^{*}z_{1!}\tw{-\d_{g}},
\end{eqnarray}
are related by adjunctions $(f_{!},f^{!})$ and $(z_{1!}, z_{1}^{!})$. Similarly, starting in both diagrams from the lower left corner, going right and then turning upward to arrive at the upper right corner, we get two maps of the above shape that are again related by adjunctions. Since \eqref{eq: two pull-pull} is commutative, we conclude that \eqref{eq: right rec des} is also commutative. This completes the proof of Theorem \ref{thm: descent of pushforward correspondence}. \qed

\part{Generalities on Fourier transform}

\section{Derived Fourier analysis}\label{sec: FT}
In this section we introduce a package of results that constitute what we call ``derived Fourier analysis'', because it occurs on generalizations of vector bundles that we call \emph{derived vector bundles}. These are the ``total spaces'' of perfect complexes, generalizing how vector bundles are the ``total spaces'' of locally free coherent sheaves. Then we generalize the Deligne-Laumon Fourier transform for $\ell$-adic sheaves from vector bundles to derived vector bundles; we call this the \emph{derived Fourier transform}.  This theory of ``derived Fourier analysis'' is needed to lift the function-theoretic Fourier analysis in the proof for $r=0$ to the level of sheaves. We establish several properties of the derived Fourier transform generalizing familiar ones, deferring most of the proofs to Appendix \ref{app: A}.

\subsection{Derived Fourier transform}\label{ssec: Derived FT}

Let $S$ be a derived Artin stack, and $E \rightarrow S$ a vector bundle. 

For $\wh{E}$ the linear dual of $E$, we have the tautological evaluation pairing $\ev \co E \times_S \wh{E} \rightarrow \A^1$. 

Let $\psi \co \F_q \rightarrow \Q_{\ell}^\times$ be a nontrivial additive character and $\cL_{\psi}$ be the corresponding Artin-Schreier sheaf on $\A^1$.

The \emph{Deligne-Laumon Fourier transform} is the functor 
\[
\FT_E^{\psi} \co D(E) \rightarrow D(\wh E)
\]
given by (following \cite[D\'{e}finition 1.2.1.1]{Lau87} in our normalizations)
\[
\cK \mapsto \pr_{1!} (\pr_0^* \cK \otimes \ev^* \cL_{\psi})[r]
\]
where $r := \rank(E)$ is the rank of $E$, and the maps are as in the diagram 
\[
\begin{tikzcd}
& E \times_S \wh{E} \ar[dl, "\pr_0"'] \ar[dr, "\pr_1"]  \ar[r,"\ev"] & \A^1 \\
E & & \wh{E}
\end{tikzcd}
\]

We will extend the Deligne-Laumon Fourier transform to certain ``derived linear spaces'' that we call \emph{derived vector bundles}, generalizing vector bundles. 

\subsubsection{Derived vector bundles}\label{sssec: derived vector bundles}
Let $S$ be a derived Artin stack. There is a functor $\Tot_S$ from the category $\Perf(S)$ of perfect complexes on $S$ to the category of derived stacks over $S$, which extends the usual construction of a vector bundle from a locally free coherent sheaf. As far as we know, the construction is due to To\"en and is documented in \cite[p.200-201]{To14} (and essentially goes back at least to \cite{To06}); \emph{however, be warned that To\"en's convention differs from ours: what he calls $\Tot_S(\cE)$ is what we would call $\Tot_S(\cE^*)$.} (The construction of $\Tot_S$ also appears in \cite{KhanI}, who agrees with To\"en's convention and therefore disagrees with ours.)

The elegant definition from \cite[p. 201]{To14} explains that as a functor from (derived) $S$-schemes to anima, $\Tot_S(\cE)$ sends $u \co T \rightarrow S$ to $\Map_{\QCoh(T)}(\cO, u^*\cE)$, after taking into account the reversal of convention as mentioned above. Here $\Map_{\QCoh(T)}$ invokes the enrichment of $\QCoh(T)$ over anima. 

\begin{defn}
For $\cE \in \Perf(S)$, we will call $E := \Tot_S(\cE)$ the \emph{derived vector bundle} associated to $\cE$. The \emph{virtual rank} of $E$, still denoted $\rank(E)$, is the locally constant function on $S$ given by $s\mapsto \chi(\cE_{s})$, the Euler characteristic of the fiber of $\cE$ at a geometric point $s$. \textbf{In general, for perfect complexes denoted with calligraphic letters such as $\cE, \cE'$, etc., the corresponding roman letters such as $E, E'$, etc. denote their associated total spaces.} 
\end{defn}

The map $0 \rightarrow \cE$ equips $E$ with a \emph{zero-section}
\[
z_{E} \co S \rightarrow E.
\]
We caution that $z_E$ is \emph{not} necessarily a closed embedding -- it is a closed embedding exactly when $E$ comes from a perfect complex $\cE$ with tor-amplitude in $[0, \infty)$. 

Also, the map $\cE \rightarrow 0$ equips $E$ with a \emph{projection}
\[
\pi_{E} \co E \rightarrow S. 
\]
We caution that $\pi_E$ is \emph{not} necessarily representable in derived schemes -- it is representable exactly when $E$ comes from a perfect complex $\cE$ with tor-amplitude in $[0, \infty)$.

\begin{remark}
If the perfect complex $\cE$ has tor-amplitude in $(-\infty,0]$, then the morphism $E \rightarrow S$ is represented by stacks which are \emph{classical} in the sense of being isomorphic to their classical truncations. If $\cE$ has tor-amplitude in $[0, \infty)$ then the morphism $E \rightarrow S$ is represented by derived schemes. Therefore, duality of derived vector bundles interchanges ``stackiness'' with ``derivedness''. 
\end{remark}

\subsubsection{Geometry of derived vector bundles}\label{sssec:geometry-of-dvb}

We record some observations about the geometric interpretation of derived vector bundles.

\begin{example}\label{ex:dvb-schematic}
Suppose $\cE$ has tor-amplitude in $[0, \infty)$ (sometimes referred to as $\cE$ being \emph{co-connective}). Then $\cE^*$ has tor-amplitude in $(-\infty, 0]$, and in particular is an animated $\cO_S$-module. The forgetful functor from animated $\cO_S$-algebras to animated $\cO_S$-modules has a left adjoint, the \emph{derived symmetric algebra} functor $\Sym_S^{\bu}$ (see for example \cite[\S 25.2.2]{SAG}). Then $\Tot_S(\cE)$ is the relative spectrum of $\Sym_S^{\bu}(\cE^*)$. 
 \end{example}
 
\begin{remark}
We will spell out a possibly more concrete way to think about the geometric object $\Tot_S(\cE)$. We start by explicating $\Tot_S(\cE)$ in the special case where $\cE \in \Perf(S)$ has tor-amplitude in $(-\infty, 0]$ (sometimes referred to as $\cE$ being \emph{connective}). Then $\Tot_S(\cE)$ represents the functor given by (derived) global sections, which at the level of 0-cells assigns to a derived affine test scheme $T \rightarrow S$ the sections $R\Gamma(T, \cE|_T)$ viewed as an animated abelian group.

Next we explain how to generalize the construction of the preceding paragraph to general $\cE \in \Perf(S)$. It is immediate from the definition that the construction $\cE \mapsto \Tot_S(\cE)$, defined so far on connective $\cE$, preserves limits (if they exist as connective complexes), hence sends exact triangles $\cE_1 \rightarrow \cE_2 \rightarrow \cE_3$ to derived Cartesian squares 
\[
\begin{tikzcd}
\Tot_S(\cE_1) \ar[d] \ar[r] & \Tot_S(\cE_2) \ar[d] \\
\Tot_S(0) \ar[r] & \Tot_S(\cE_3)
\end{tikzcd}
\]
Note that $\Tot_S(0) \cong S$; we call the map $\Tot_S(0) \rightarrow \Tot_S(\cE)$ the \emph{zero-section}. We extend the construction of $\Tot_S(-)$ to all $\cE \in \Perf(S)$ by this condition. Concretely, $\Tot_S(\cE[-1])$ is the derived self-intersection\footnote{Although we call this an ``intersection'', the zero-section will typically not be a closed embedding.} of the zero-section of $\Tot_S(\cE)$, 
\[
\begin{tikzcd}
\Tot_S(\cE[-1]) \ar[r] \ar[d] & \Tot_S(0)  \ar[d] \\
\Tot_S(0) \ar[r] & \Tot_S(\cE)
\end{tikzcd}
\]
and for a general $\cE \in \Perf(S)$ there exists some $d$ such that $\cE[d]$ is connective; then $\Tot_S(\cE)$ is obtained from $\Tot_S(\cE[d])$ by iterating the procedure of forming derived self-intersection of the zero-section $d$ times. 
\end{remark}

\begin{lemma}\label{lem:closed-embedding}
Let $S$ be a derived Artin stack and $f \co \cE \rightarrow \cF$ be a map in $\Perf(S)$. Denote the associated map of derived vector bundles by $f \co E \rightarrow F$. 

\begin{enumerate}
\item The map $f$ is smooth if $\Fib(f)$ has tor-amplitude in $(-\infty, 0]$. 

\item The map $f$ is quasi-smooth if $\Fib(f)$ has tor-amplitude in $(-\infty, 1]$. 

\item The map $f$ is a closed embedding if $\Fib(f)$ has tor-amplitude in $[1, \infty)$. 

\item The map $f$ is representable in derived affine schemes if $\Fib(f)$ has tor-amplitude in $[0, \infty)$. 
\end{enumerate}
\end{lemma}

\begin{proof}
Note that if $\cE$ is a perfect complex on $S$, then the relative tangent complex $\bT_{E/S}$ of the associated derived vector bundle is the pullback of $\cE$ to $E$. It follows that the relative tangent complex $\bT_f$ is the pullback to $E$ of $\Fib(f)$. From this, (1) and (2) follow immediately. 

For (3) and (4), we may work locally, and thus assume that $\cF$ is represented by a complex of finite locally free $\cO_S$-modules
\[
\ldots \rightarrow \cF^{-1} \rightarrow \cF^0 \rightarrow \cF^1 \rightarrow \ldots.
\]
By (1), the map $F^{\geq 0} \rightarrow F$ is a smooth cover. We can check the properties in question after smooth base change on $F$, so we can replace $F$ by $F^{\geq 0}$ to reduce to the case where $\cF$ has tor-amplitude in $[0,\infty)$. If $\Fib(f)$ has tor-amplitude in $[0,\infty)$, then $\cE$ also has tor-amplitude in $[0, \infty)$. From the explicit description of derived vector bundles in this co-connective case (Example \ref{ex:dvb-schematic}) we immediately see (4). 

For (3), the property of being a closed embedding can be checked after classical truncation and smooth base change, so we may assume that the base $S = \Spec R$ is a classical affine scheme. Since $\Fib(f)$ has tor-amplitude in $[1,\infty)$, its dual has tor-amplitude in $(-\infty,-1]$. Dualizing the exact triangle $\Fib(f) \rightarrow \cE \rightarrow \cF$ gives an exact triangle $\cF^* \rightarrow \cE^* \rightarrow \Fib(f)^*$, which shows that $H^0(\cF^*) \twoheadrightarrow H^0(\cE^*)$ is surjective. Since, after the preceding reduction, $\cE$ and $\cF$ have tor-amplitude in $[0,\infty)$, the complexes $\cE^*$ and $\cF^*$ are animated $R$-modules; for an animated $R$-module $M$, one has $\pi_0\Sym_R^{\bu}(M)\cong \Sym_R^{\bu}(\pi_0 M)$. It follows that $\Sym_R^{\bu}(\cF^*)\rightarrow \Sym_R^{\bu}(\cE^*)$ is surjective on $\pi_0$. By Example \ref{ex:dvb-schematic}, this says exactly that the map $E\rightarrow F$ is a closed embedding on classical truncations.
\end{proof}

\subsubsection{The $\ell$-adic Fourier transform for derived vector bundles}

Let $S$ be a derived Artin stack, $\cE \in \Perf(S)$. For $\cE^* \in \Perf(S)$ the linear dual of $\cE$, we have a tautological evaluation pairing $\cE \times \cE^* \rightarrow \cO_S$. Setting $E := \Tot_S(\cE)$ and $\wh{E} := \Tot_S(\cE^*)$, this induces on total spaces a map 
\[
\ev \co E \times_S \wh{E} \rightarrow \A^1.
\]
The Fourier transform 
\[
\FT^\psi_E \co D(E) \rightarrow D(\wh{E})
\]
is defined as
\[
\cK \mapsto \pr_{1!} (\pr_0^* (\cK) \otimes \ev^* \cL_{\psi})[r]
\]
where $r := \rank(E)$ is the virtual rank of $E \rightarrow S$ (which is constant on each connected component of $S$, so the shift $[r]$ makes sense over each connected component of $S$), and the maps are as in the diagram 
\[
\begin{tikzcd}
& E \times_S \wh{E} \ar[dl, "\pr_0"'] \ar[dr, "\pr_1"] \ar[r, "\ev"] & \A^1 \\
E & & \wh{E}
\end{tikzcd}
\]

This extends the Deligne-Laumon Fourier transform for $\ell$-adic sheaves for vector bundles, which corresponds to the case where $\cE$ is a locally free coherent sheaf. When the additive character $\psi$ is understood, we will simply omit it from the notation. 


\subsection{Properties of the derived Fourier transform}\label{ssec: properties of FT}
We now tabulate some basic properties of the derived Fourier transform. Below we let $r$ be the virtual rank of $E \rightarrow S$. The nontrivial proofs are all found in Appendix \ref{app: A}. Below we use the adjective ``canonical'' to describe an isomorphism whose construction does not depend on any auxiliary choices; we emphasize this because the proofs require considering, at intermediate stages, natural isomorphisms that a priori depend on auxiliary choices (but are ultimately seen to be independent of such choices a posteriori).

\subsubsection{Fourier transform of Gaussians}\label{sssec: FT of Gaussian}
Suppose $h_E \co E \xrightarrow{\sim} \wh{E}$ is a symmetric isomorphism. This induces:
\begin{itemize}
\item a quadratic form $q \co E \rightarrow \A^1$ given by $q(e) := \langle e, h_E(e) \rangle $, and 
\item a quadratic form $\wh{q} \co \wh{E} \rightarrow \A^1$ given by $\wh q(\wh{e}) := \langle h_E^{-1}(\wh{e}), \wh{e} \rangle$.
\end{itemize}
Then one has a canonical isomorphism
\[
[2]^* \FT_E(q^* \cL_{\psi}) \cong (-\wh{q})^* \cL_{\psi} \otimes (\pi_{\wh{E}}^* \pi_{E!} q^* \cL_{\psi} [r] ).
\]
The proof is the same as for classical vector bundles, which is found in \cite[Proposition 1.2.3.3]{Lau87}.

\subsubsection{Base change}\label{sssec: FT base change} Let $h \co \wt{S} \rightarrow S$. For a derived vector bundle $E \rightarrow S$, let $\wt{E} \rightarrow E$ be its base change along $h$. So we have derived Cartesian squares 
\[
\begin{tikzcd}
\wt{E} \ar[r, "h^E"] \ar[d] & E \ar[d] \\
\wt{S} \ar[r, "h"] & S
\end{tikzcd} 
\hspace{1cm}
\begin{tikzcd}
\wh{\wt{E}} \ar[r, "h^{\wh{E}}"] \ar[d] & \wh{E} \ar[d] \\
\wt{S} \ar[r, "h"] & S
\end{tikzcd}
\]

Then there are canonical natural isomorphisms of functors $D(E) \rightarrow D(\wh{\wt{E}})$
\begin{equation}\label{eq: FT bc h^*}
\FT_{\wt{E}} \circ (h^E)^* \cong (h^{\wh{E}} )^* \circ \FT_E
\end{equation}
\begin{equation}\label{eq: FT bc h^!} 
\FT_{\wt{E}} \circ (h^E)^! \cong (h^{\wh{E}} )^! \circ \FT_E
\end{equation}
and canonical natural isomorphisms of functors $D(\wt{E}) \rightarrow D(\wh{E})$
\begin{equation}\label{eq: FT bc h_!}
\FT_{E} \circ (h^E)_!  \cong (h^{\wh{E}} )_! \circ \FT_{\wt{E}}
\end{equation}
\begin{equation}\label{eq: FT bc h_*}
\FT_{E} \circ (h^E)_*  \cong (h^{\wh{E}} )_* \circ \FT_{\wt{E}}.
\end{equation}
The isomorphisms \eqref{eq: FT bc h^*} and \eqref{eq: FT bc h_!} follow directly from proper base change. The other natural isomorphisms will be constructed in \S \ref{sssec: FT bc proof}.

\subsubsection{Involutivity}  There is a canonical natural isomorphism $\FT_{\wh{E}} \circ \FT_{E} \cong [-1]^* (-r)$ of functors $D(E) \rightarrow D(E)$, where $[-1]$ is multiplication by $-1$ on $E$. 

\begin{remark}
The construction of this natural isomorphism appears to be significantly more involved than in the situation of the Deligne-Laumon Fourier transform, and occupies much of Appendix \ref{app: A}.
\end{remark}

\subsubsection{Functoriality}\label{sssec: FT functoriality} Let $f \co E' \rightarrow E$ be a linear map of derived vector bundles having virtual ranks $r', r$ respectively. This induces a morphism $\wh{f} \co \wh{E} \rightarrow \wh{E}'$ of dual derived bundles. Then we have canonical natural isomorphisms of functors $D(E') \rightarrow D(\wh{E})$: 
\begin{enumerate}
\item\label{item: functoriality 1} $\wh{f}^* \circ \FT_{E'} \cong \FT_E \circ f_! [r'-r]$,
\item\label{item: functoriality 1.25} $\wh{f}^! \circ \FT_{E'} \cong  \FT_E \circ f_*[r-r'](r-r')$,
\end{enumerate}
and canonical natural isomorphisms of functors $D(E) \rightarrow D(\wh{E}')$: 
\begin{enumerate}[resume]
\item\label{item: functoriality 1.5} $\FT_{E'} \circ f^* \cong \wh{f}_! \circ \FT_E [r-r'] (r-r')$,
\item\label{item: functoriality 2} $\FT_{E'} \circ f^! \cong \wh{f}_* \circ \FT_E[r'-r]$. 
\end{enumerate}
We record for convenience that in the Fourier dual coordinates, the previous two isomorphisms become natural isomorphisms of functors $D(\wh{E}') \rightarrow D(E)$: 
\begin{enumerate}[resume]
\item $\FT_{\wh{E}}\circ \wh{f}^*  \cong f_!\circ \FT_{\wh{E}'} [r'-r](r'-r)$,
\item\label{item: functoriality 3} $\FT_{\wh{E}}\circ \wh{f}^! \cong f_*\circ \FT_{\wh{E}'} [r-r']$. 
\end{enumerate}
Moreover, we shall see in \S \ref{sssec: more functoriality} that by construction, the natural isomorphisms above intertwine the adjunction $(f_!, f^!)$ with the adjunction $(\wh{f}^*, \wh{f}_*)$ (up to shift and twist), and the adjunction $(f^*, f_*)$ with the adjunction $(\wh{f}_!, \wh{f}^!)$ (up to shift and twist). In particular, $\FT_E$ sends the unit and counit 
\[
\Id \rightarrow f^! f_! \hspace{1cm}  f_! f^! \rightarrow \Id
\]
to the unit and the counit\footnote{In some normalizations of the involutivity isomorphism, a sign would appear here. We have set up our formalism so that no sign issues intervene here.}
\[
\FT_E \rightarrow \wh{f}_* \wh{f}^* \FT_E \hspace{1cm}  \wh{f}^* \wh{f}_* \FT_E \rightarrow \FT_E 
\]
under the above identifications, and similarly for the other adjunction.

\begin{defn}[The delta-sheaf] For a derived vector bundle $E \rightarrow S$, recall that $z_{E} \co S \rightarrow E$ is the zero-section, which may not be a closed embedding. We define $\delta_{E} := z_{E!} \Qll{S}$. 
\end{defn}

\begin{example}[Delta-constant duality]\label{ex: delta-constant duality}
Suppose $E\rightarrow S$ is a derived vector bundle of rank $r$. Then the natural transformation \eqref{item: functoriality 1} gives an isomorphism
\[
\FT_E(\delta_E) \cong \Qll{\wh{E}}[r]
\]
and natural transformation \eqref{item: functoriality 1.5} gives an isomorphism
\[
\FT_{\wh{E}} \Qll{\wh{E}} \cong  \delta_E [-r](-r).
\]
\end{example}

\subsubsection{Verdier duality}\label{sssec: Verdier duality} Letting $\DD_E$ (resp. $\DD_{\wh{E}}$) denote the Verdier duality functor on $E$ (resp. $\wh{E}$), there is a canonical isomorphism naturally in $\cK\in D(E)$
\[
\DD_{\wh{E}} ( \FT_E^{\psi} (\cK))\cong \FT_E^{\psi^{-1}}( \DD_E(\cK)) (r).
\]

\subsubsection{Convolutions}\label{sssec:convolutions} For $\cK_0, \cK_1 \in D(E)$, we write 
\[
\cK_0 \star \cK_1 = +_! (\pr_0^* \cK_0 \otimes \pr_1^* \cK_1)
\]
where maps are as in the diagram 
\[
\begin{tikzcd}
& E \times_S E \ar[r, "+"] \ar[dl, "\pr_0"'] \ar[dr, "\pr_1"] & E \\
E & & E
\end{tikzcd}
\]
For $r = \rank(E)$, there is a canonical natural isomorphism 
\[
\FT_E(\cK_0 \star \cK_1) \cong \FT_{E}(\cK_0) \otimes \FT_{E}(\cK_1)[-r]
\]
which is constructed formally from the functorialities of \S \ref{sssec: FT functoriality}.

\subsubsection{Plancherel formula}\label{sssec: Plancherel} There is a canonical natural isomorphism 
\begin{equation}\label{eq:plancherel-eq}
 \pi_{\wh{E}!} (\FT_E(\cK_0) \otimes \FT_E(\cK_1)) \cong \pi_{E!} (\cK_0 \otimes [-1]^* \cK_1) (-r).
\end{equation}

To see this, we write the LHS as $ \pi_{\wh{E}!} \Delta^* \FT_{E \times_S E}( \pr_0^* (\cK_0) \otimes \pr_1^* (\cK_1))$, then use the dual-side Cartesian square involving $-\Delta:\wh{E}\to \wh{E}\times_S\wh{E}$ and $+_{\wh E}$ to rewrite $(z_{\wh E})^*(+_{\wh{E}})_!$ by proper base change.

\subsection{Proper base change} We will need the compatibility of the Fourier transform with proper base change, at least under a ``global presentation'' hypothesis.

\subsubsection{Globally presented derived vector bundles}
We introduce the following definition for technical reasons: 

\begin{defn}
A perfect complex $\cE \in \Perf(S)$ is \emph{globally presented} if it is quasi-isomorphic to a bounded complex of vector bundles on $S$,
\[
\cE \cong (\ldots \rightarrow \cE^{-1} \rightarrow \cE^0 \rightarrow \cE^1 \rightarrow \ldots ).
\]
By definition of $\Perf(S)$ such a presentation exists Zariski-locally on $S$, but we are asking for its existence globally.

We say that the associated derived vector bundle $E = \Tot_S(\cE)$ is \emph{globally presented} if $\cE$ is globally presented.

We say that a map $f \co E' \rightarrow E$ of derived vector bundles lying over $h \co S' \rightarrow S$ is \emph{globally presented} if there exist global presentations $\cE^{\bu}$ for $E$ and $(\cE^{\bu})'$ for $E'$ and $f$ is induced by a map of complexes $(\cE^{\bu})' \rightarrow h^* \cE^{\bu}$. 

More generally, we say that a diagram of derived vector bundles is \emph{globally presented} if there exist global presentations for all derived vector bundles such that all maps between derived vector bundles are induced by maps of these presentations. 

\end{defn}

\begin{example}
Any diagram of classical vector bundles is globally presented. 
\end{example}

The role of the notion of global presentation is the following. Certain proofs (deferred to Appendix \ref{app: A}) towards the results already mentioned in this section rely on the notion of global presentation at intermediate stages. Furthermore, in the statement of Proposition \ref{prop: ft pbc} below, we impose a global presentation assumption. Various results in later sections depend on Proposition \ref{prop: ft pbc} and will therefore also require a global presentation assumption. We expect that this assumption is not actually necessary, but it provides a ``shortcut'' for the proof, ultimately because a globally presented map can be factored into the composition of a closed embedding and a smooth map (Lemma \ref{lem: factorize smooth/proper}).

\subsubsection{Compatibility with proper base change}\label{sssec: pbc}

Consider a Cartesian square of globally presented derived vector bundles, along with the dual Cartesian square
\[
\begin{tikzcd}
& B \ar[dl, "g'"'] \ar[dr, "f'"] \\
A  \ar[dr, "f"'] & & D  \ar[dl, "g"] \\
& C
\end{tikzcd}\hspace{1cm} \text{and} \hspace{1cm}
\begin{tikzcd}
& \wh{C} \ar[dl, "\wh{f}"'] \ar[dr, "\wh{g}"] \\
\wh{A}  \ar[dr, "\wh{g}'"'] & & \wh{D}  \ar[dl, "\wh{f}'"] \\
& \wh{B}
\end{tikzcd}
\]
Then proper base change gives natural isomorphisms 
\begin{equation}\label{eq: app BC}
g^* f_! \cong (f')_! (g')^*  \hspace{1cm} \text{and} \hspace{1cm} \wh{g}_! \wh{f}^* \cong (\wh{f}')^* \wh{g}'_!
\end{equation}
Let $d=d(f),\delta := d(g)$. According to \S \ref{sssec: FT functoriality}, there are natural isomorphisms 
\begin{equation}\label{eq: app FT}
\wh{g}_! \wh{f}^*\FT_A \cong  \FT_D g^* f_! [d+\delta](\delta)  \hspace{1cm} \text{and} \hspace{1cm}  (\wh{f}')^* \wh{g}'_! \FT_A \cong \FT_D f'_! (g')^* [d+\delta](\delta).
\end{equation}

\begin{prop}\label{prop: ft pbc} 
Assume that $f$ and $g$ are globally presented (in particular, $A,C,D$ are globally presented). Then the diagram 
\begin{equation}\label{eq: FT BC compatibility}
\begin{tikzcd}
\wh{g}_! \wh{f}^*\FT_A  \ar[r, "\sim"]   \ar[d, "\sim"]  &   \FT_D g^* f_! [d+\delta](\delta)  \ar[d, "\sim"]  \\
(\wh{f}')^* \wh{g}'_! \FT_A \ar[r, "\sim"]   &  \FT_D f'_! (g')^* [d+\delta](\delta) 
\end{tikzcd}
\end{equation}
commutes, where the identifications are as in \eqref{eq: app BC} and \eqref{eq: app FT}.
\end{prop}

This innocuous-looking statement turns out to be rather involved to prove, so the proof will be deferred to \S \ref{app: A} (see Proposition \ref{prop: FT hexagon}). As discussed above, we believe that the technical assumption that $f,g$ are globally presented is an artefact of the proof.

\subsection{Fourier transform of the Gysin map}\label{ssec: FT of Gysin}

Let $f \co E' \rightarrow E$ be a quasi-smooth morphism of derived vector bundles over $S$, which by Lemma \ref{lem:closed-embedding} is equivalent to $\Fib(\cE' \rightarrow \cE)$ having tor-amplitude in $(-\infty, 1]$. Then $f$ has a relative fundamental class $[f]$, which induces a Gysin natural transformation $f^* \rightarrow f^! \tw{-d(f)}$, as explained in \S \ref{ssec: relative fundamental nt}. Dualizing and using Lemma \ref{lem:closed-embedding} again, this is equivalent to the condition that the dual map $\wh{f} \co \wh{E} \rightarrow \wh{E'}$ is separated, and therefore has a ``forget supports'' natural transformation $\can(\wh{f}) \co \wh{f}_! \rightarrow \wh{f}_*$.

\begin{example}\label{ex: factorize bundle map} If $E$ and $E'$ are classical vector bundles over $S$, then the map $f$ is automatically LCI (and therefore quasi-smooth). Indeed, the graph of $f$ provides a factorization 
\[
E' \inj E \times_S E' \xrightarrow{\pr_1} E
\]
which is a composition of a regular embedding and a smooth morphism. 
\end{example}

We need the following identification of the Fourier transform of the Gysin map. After some inquiries, we found that this statement was unknown to experts even in the case where $f$ is a map of classical vector bundles, hence automatically LCI by Example \ref{ex: factorize bundle map}.

\begin{prop}\label{prop: fourier of Gysin}
Let $f \co E' \rightarrow E$ be a globally presented quasi-smooth map of derived vector bundles and let $\wh{f} \co \wh{E} \rightarrow \wh{E}'$ be the dual map to $f \co E' \rightarrow E$. Then the diagram of functors $D(E) \rightarrow D(\wh{E}')$
\begin{equation}\label{eq: fourier of gysin}
\begin{tikzcd}
\wh{f}_! \FT_{E} \ar[d, "\sim"] \ar[r, "\can(\wh{f})"] & \wh{f}_* \FT_{E} \ar[d, "\sim"] \\
\FT_{E'} f^* [d(f)](d(f)) \ar[r, "{[f]}"] & \FT_{E'} f^! [-d(f)]
\end{tikzcd}
\end{equation}
commutes. 
\end{prop}

\begin{remark}
The significance of Proposition \ref{prop: fourier of Gysin} is to describe the derived (relative) fundamental class $[f]$ in terms of classical geometry \emph{in the Fourier dual space}. 
\end{remark}

The proof of Proposition \ref{prop: fourier of Gysin} is rather lengthy. We will begin with several reductions. 

\subsubsection{Reduction to a classical base} 	The first observation is that for the purpose of proving Proposition \ref{prop: fourier of Gysin}, we may assume that the base $S$ is classical. To justify this, we note that derived base change along the classical truncation $S_{\mathrm{cl}} \inj S$ induces an equivalence of \'etale topoi, compatible with all functors such as $\wh{f}_!$, $\wh{f}_*$, $f^*$, $f^!$, and $\FT_E, \FT_{E'}$. It is also clearly compatible with the natural transformation $\can(\wh f)$ since this does not depend on the derived structure. For the natural transformation $[f]$, which \emph{does} depend on the derived structure, the compatibility has more content, but follows from compatibility of formation of relative fundamental classes with respect to change of base \cite[Theorem 3.13]{KhanI}. Thus, for the rest of the proof of Proposition \ref{prop: fourier of Gysin} we may and do assume that our derived vector bundles are over a classical higher Artin stack. 

\subsubsection{Reduction to smooth derived vector bundles} 
We begin by reducing Proposition \ref{prop: fourier of Gysin} to the case where $E$ and $E'$ are both smooth (but still potentially stacky) over $S$. 

\begin{lemma}\label{lem: pullback from connective}
We can find globally presented $\cF, \cF'$ with tor-amplitude in $(-\infty, 0]$ and a globally presented quasi-smooth map $g \co F' \rightarrow F$ fitting into a derived Cartesian square 
\begin{equation}\label{eq: pullback from connective}
\begin{tikzcd}
E' \ar[d, "f"] \ar[r, "h'"] & F' \ar[d, "g"] \\
E \ar[r, "h"] & F
\end{tikzcd}
\end{equation}
\end{lemma}

\begin{proof}We have by assumption that $f$ is represented by a map of complexes of vector bundles
\begin{equation}\label{eq: presentation 1}
\begin{tikzcd}
\ldots \ar[d] \ar[r]  & (\cE')^{m-1} \ar[r]  \ar[d, "f_{m-1}"] &  (\cE')^m \ar[r]  \ar[d, "f_m"] & 0\\
\ldots \ar[r]  & \cE^{m-1} \ar[r] &  \cE^m \ar[r] & 0
\end{tikzcd}
\end{equation}
We induct on the statement: as long as $m \geq 1$ and $f$ is quasi-smooth, any such diagram is up to homotopy equivalence pulled back from one in which both rows are complexes of vector bundles which vanish in degrees at least $m$ (in both rows). 

To prove this, we will replace \eqref{eq: presentation 1} via homotopy equivalences by a map of complexes for which $f_m = \Id$, for then 
\[
\begin{tikzcd}
\cE' \ar[d, "f"] \ar[r] & (\cE')^{< m}  \ar[d, "f_{< m}"]\\
\cE \ar[r] & \cE^{< m}
\end{tikzcd}
\]
is a pullback square, where $(\ldots)^{<m}$ refers to the naive truncation, and we may take $g \co \cF' \rightarrow \cF$ to be $f_{ <m} \co (\cE')^{< m}   \rightarrow \cE^{< m}$. This $g$ is globally presented by the truncated complexes, and the pullback square identifies $\Fib(g)$ with $\Fib(f)$, so $g$ is quasi-smooth. Indeed, the assumption that $f$ is quasi-smooth implies that it induces an isomorphism on $H^{\geq 2}$ and a surjection on $H^1$, so by the assumption that $m \geq 1$ the map $\cE^{m-1} \oplus (\cE')^m \xrightarrow{d+ f_m} \cE^m$ is surjective. (Note that this makes sense only because we have reduced to a classical base.) We may replace $(\cE')^{m-1} \xrightarrow{d'} (\cE')^m$ by $\cE^{m-1}\oplus (\cE')^{m-1} \xrightarrow{\Id\oplus d'} \cE^{m-1} \oplus (\cE')^m$, and replace the diagram \eqref{eq: presentation 1} by
\begin{equation*}
\begin{tikzcd}
\ldots \ar[d] \ar[r] & \cE^{m-1}\op (\cE')^{m-1} \ar[r, "\Id\oplus d'"] \ar[d, "\Id+f_{m-1}"] & \cE^{m-1}\oplus (\cE')^m  \ar[d, "d+f_m"] \ar[r] & 0  \\
\ldots \ar[r] &  \cE^{m-1} \ar[r, "d"] & \cE^m \ar[r] & 0 
\end{tikzcd}
\end{equation*}
Hence we reduce to the case where $f_{m}$ is surjective. Therefore, naive pullback along $f_m$ computes the derived pullback. Then we may replace the bottom row in \eqref{eq: presentation 1} by its pullback along the surjection $f_{m}\co (\cE')^m \to \cE^m$. The map in degree $m$ is now the identity map, as desired. 
\end{proof}


\begin{lemma}\label{lem: reduce to connective}
If Proposition \ref{prop: fourier of Gysin} holds for the map $g \co F' \rightarrow F$ in the right column of \eqref{eq: pullback from connective}, then it holds for the map $f \co E' \rightarrow E$. 
\end{lemma}

\begin{proof}
In the construction of Lemma \ref{lem: pullback from connective}, the maps $h$ and $h'$ are obtained by successively removing terms in strictly positive degrees, and induce equivalences on \'etale topoi. It is consequently enough to verify the desired natural transformation on objects of the form $h^*\cL$, with $\cL\in D(F)$.

By the base change property for relative fundamental classes \cite[Theorem 3.13]{KhanI}, we have $(h')^* [g] = [f]$, meaning that the following diagram commutes:
\begin{equation}\label{eq:connective-fundamental-class-base-change}
\begin{tikzcd}
(h')^* g^* \Qll{F} \ar[r, "{(h')^*[g]}"]  \ar[d, equals] & (h')^* g^! \Qll{F} \tw{-d(g)}  \ar[d, "\di"] \\
f^* h^* \Qll{F} \ar[r, "{[f]}"]  & f^! h^* \Qll{F} \tw{-d(f)} 
\end{tikzcd}
\end{equation}
The commutative diagram \eqref{eq:connective-fundamental-class-base-change} is the special case $\cL=\Qll{F}$ of the Lemma; tensoring it with $h^*\cL$ and using the construction of the Gysin natural transformation from the relative fundamental class gives the corresponding base-change identity for every $h^*\cL\in D(E)$.
We are granted that $\FT_{F'}([g]) = \can(\wh g)$. Then applying $\FT_{E'}$ to \eqref{eq:connective-fundamental-class-base-change}, using \S \ref{sssec: FT functoriality}, we have that 
\[
\FT_{E'}([f]) = \FT_{E'}((h')^* [g]) = \wh{h}'_! \FT_{F'}([g])= \wh{h}'_!  \can(\wh g) = \can(\wh f)
\]
where the last equality used Lemma \ref{lem: push push compatibility 1} (note that both maps $\wh g$ and $\wh f$ are separated and locally of finite type). 
\end{proof}

Hence we have reduced the proof of Proposition \ref{prop: fourier of Gysin} to the case where $E'$ and $E$ are smooth over $S$, and in the rest of the argument we will assume this to be the case.

\subsubsection{Equivalence of formulations} Recall from Example \ref{ex: delta-constant duality} that for $r:= \rank(E)$, we have $\FT(\Qll{E}) \cong \delta_{\wh{E}}[-r](-r)$. Therefore, a special case of Proposition \ref{prop: fourier of Gysin} is the following Lemma.  

\begin{lemma}\label{lem: FT of delta}
The relative fundamental class $f^* \Qll{E} \xrightarrow{[f]} f^! \Qll{E}\tw{-d(f)}$ is sent by $\FT_{E'}$ to 
\[
\wh{f}_! \delta_{\wh{E}} [-r'](-r') \xrightarrow{\can(\wh{f})} \wh{f}_* \delta_{\wh{E}}[-r'](-r').
\]
\end{lemma}

However, the converse is also true, at least under the given assumptions.

\begin{lemma} Let $E', E$ be derived vector bundles smooth over $S$, and $f \co E' \rightarrow E$ be quasi-smooth and globally presented. If Lemma \ref{lem: FT of delta} holds for $f$, then Proposition \ref{prop: fourier of Gysin} holds for $f$. 
\end{lemma}

\begin{proof}
Indeed, let $\cK \in D(E)$. Then $f^* \cK \rightarrow f^! \cK \tw{-d(f)}$ is the composition 
\begin{equation}\label{eq: long composition}
f^* \cK = f^* \cK \otimes f^* \Qll{E} \xrightarrow{\Id \otimes [f]} f^* \cK \otimes f^! \Qll{E} \tw{-d(f)} \xrightarrow{\di} f^!(\cK \otimes \Qll{E}) \tw{-d(f)} = f^! \cK \tw{-d(f)}
\end{equation}
where we recall that the second arrow is the base change map for the Cartesian square 
\[
\begin{tikzcd}
E' \ar[r, "{(\Id,f)}"] \ar[d, "f"] & E' \times E \ar[d, "f \times \Id"] \\
E \ar[r, "\Delta"] & E \times E
\end{tikzcd}
\]
Note that a global presentation for $f$ induces a global presentation for this diagram. 

Let us factor \eqref{eq: long composition} into two halves: 
\begin{equation}\label{eq: long composition first}
f^* \cK = f^* \cK \otimes f^* \Qll{E}  \xrightarrow{\Id \otimes [f]} f^* \cK \otimes f^! \Qll{E} \tw{-d(f)} 
\end{equation}
and 
\begin{equation}\label{eq: long composition second}
f^* \cK \otimes f^! \Qll{E} \tw{-d(f)} \xrightarrow{\di} f^!(\cK \otimes \Qll{E}) \tw{-d(f)} = f^! \cK \tw{-d(f)}.
\end{equation}
Abbreviate $\wh{\cK} := \FT_E(\cK)$. By hypothesis, the first half \eqref{eq: long composition first} is sent (up to shift and twist by $r-r'$) by $\FT_{E'}$ to 
\[
\wh{f}_! \wh{\cK} = \wh{f}_! \wh{\cK} \star \wh{f}_! \delta_{\wh{E}}   \xrightarrow{\Id \star \can(\wh{f})} \wh{f}_! \wh{\cK} \star\wh{f}_* \delta_{\wh{E}}. 
\]
Applying Proposition \ref{prop: ft pbc}, we see that the second half \eqref{eq: long composition second} is sent (up to shift and twist by $r-r'$) by $\FT_{E'}$ to 
\[
 \wh{f}_! \wh{\cK} \star\wh{f}_* \delta_{\wh{E}}  \xrightarrow{\di} \wh{f}_*(\wh{\cK} \star \delta_{\wh{E}}) = \wh{f}_* \wh{\cK} 
 \]
 where the arrow comes from base change for the Cartesian square 
 \[
 \begin{tikzcd}
 \wh{E} \times \wh{E} \ar[r, "+"] \ar[d, "{\Id\times \wh{f}}"]  & \wh{E} \ar[d, "\wh{f}"] \\
 \wh{E} \times \wh{E}'  \ar[r, "\wh{f} + \Id "] & \wh{E}'
 \end{tikzcd}
 \]
 which has a global presentation induced by that of $f$. 
 
To complete the proof, we need to show that the above composition agrees with the ``forget supports'' map for $\wh{f}$. This follows from the compatibility statement in Lemma \ref{lem: push push compatibility 1} that the following diagram commutes:
\[
\begin{tikzcd}
\wh{f}_! \wh{\cK} \star \wh{f}_! \delta_{\wh{E}}  \ar[r, "\sim"] \ar[d, equals] & \wh{f}_! \wh{\cK} \star \wh{f}_* \delta_{\wh{E}}   \ar[d] \\
\wh{f}_! (\wh{\cK} \star \delta_{\wh{E}}) \ar[r, "\can(\wh{f})"] & \wh{f}_* (\wh{\cK} \star  \delta_{\wh{E}})
\end{tikzcd}
\]
where the top horizontal isomorphism is justified by the chain of identifications (using that $z_{\wh{E}}$ and $z_{\wh{E}'}$ are closed embeddings because $\cE, \cE'$ are connective by assumption)
\begin{equation}\label{eq: deltas}
\wh{f}_! \delta_{\wh{E}}  = \wh{f}_! z_{\wh{E}!} \Qll{S} = z_{\wh{E}'!} \Qll{S} = z_{\wh{E}'*} \Qll{S} = \wh{f}_* z_{\wh{E}*} \Qll{S} = \wh{f}_* \delta_{\wh{E}}.
\end{equation}
Indeed, tracing the diagram along the left and bottom gives $\can(\wh{f})$, while tracing along the top and right gives (up to shift and twist by $r-r'$) the natural transformation $\FT_{E'}([f])$. This concludes the proof.
\end{proof}

\subsubsection[Proposition on Fourier transform of Gysin maps]{Proposition \ref{prop: fourier of Gysin}} We now complete the proof of Proposition \ref{prop: fourier of Gysin} by establishing Lemma \ref{lem: FT of delta}. By Lemma \ref{lem: pullback from connective} and Lemma \ref{lem: reduce to connective}, we may assume that $E', E$ are smooth over $S$. Then the map $f$ is LCI, hence $f^! \Qll{E}\tw{-d(f)}$ is isomorphic to $f^* \Qll{E} \cong \Qll{E'}$. Note that the space $\Hom_{E'}(f^* \Qll{E}, f^! \Qll{E}\tw{-d(f)})$ then identifies with $\rH^0(E'; \Q_\ell)  \cong \rH^0(S; \Q_{\ell})$. 

Similarly we see $\wh{f}_! \delta_{\wh{E}} \cong \delta_{\wh{E}'} \cong \wh{f}_* \delta_{\wh{E}}$ and $\Hom_{\wh{E}'}(\wh{f}_! \delta_{\wh{E}}, \wh{f}_* \delta_{\wh{E}}) \cong \rH^0(\wh{E}'; \Q_\ell) \cong \rH^0(S; \Q_{\ell})$. Of course, this also follows formally from the previous paragraph, since $\FT_{E'}$ is an equivalence of categories. 

By examining each connected component at a time, we reduce to the case where $S$ is connected. Then $\rH^0(S; \Q_\ell) \cong \Q_\ell$, so the two maps in question differ by some scalar. We want to verify that this scalar is $1$; it suffices to check this after pulling back to a geometric point of $S$: Fourier transform is compatible with pullback on the base, the relative fundamental class is compatible with base change by \cite[Theorem 3.13]{KhanI}, and the natural transformation $\can$ and the delta-sheaf identifications are functorial under base change. Thus we reduce to the case where $S$ is a point. 

For $f^!$, the six-functor formalism is developed (cf. \cite{Ver67}) so that for a smooth morphism $f$, the shriek pullback $f^!$ is \emph{equal} to $f^*\tw{d(f)}$. Hence if $f \co E' \rightarrow E$ is surjective, then the Gysin map is \emph{the} identity map, and since $\wh{f}$ is a closed embedding the forget supports map $\wh{f}_! \rightarrow \wh{f}_*$ is \emph{the} identity as well. So Proposition \ref{prop: fourier of Gysin} is evident in this case. 

By factoring $f \co E' \rightarrow E$ as the composition of a linear closed embedding and a linear smooth map (cf. Lemma \ref{lem: factorize smooth/proper}), we may assume that $f$ is a closed embedding (since the statement of Proposition \ref{prop: fourier of Gysin} is compatible with compositions). Since $S$ is a point, $f$ necessarily has a splitting $\pi \co E \rightarrow E'$, which is necessarily smooth. Then since $\pi \circ f = \Id$, we have 
\[
\Qll{E'} = f^* \pi^* \Qll{E'} \xrightarrow{[\pi]} f^* \pi^! \Qll{E'}\tw{-d(\pi)} \xrightarrow{[f]} f^! \pi^! \Qll{E'} \tw{-d(\pi)-d(f)} = \Qll{E'}.
\]
Since $\pi$ is smooth, $[\pi]$ is the identity map. The composition $[\pi \circ f ] = [\Id]$ is also the identity map, so we deduce that with respect to the identifications $f^* \pi^! \Qll{E'}\tw{-d(\pi)}  = \Qll{E'}$ and $f^! \pi^! \Qll{E'} \tw{-d(\pi)-d(f)} = \Qll{E'}$ induced by the equalities in the equation above, $[f]$ is \emph{the} identity map. 

Similarly, we have
\[
\delta_{\wh{E'}} = \wh{f}_! z_{\wh{E}!} \Qll{S} \xrightarrow{\can(z_{\wh{E}})} \wh{f}_! z_{\wh{E}*}  \Qll{S} \xrightarrow{\can(\wh{f})} \wh{f}_* z_{\wh{E}*} \Qll{S} = \delta_{\wh{E}'}.
\]
Since $z_{\wh{E}}$ is a closed embedding by the hypothesis that $\cE$ is connective, $\can(z_{\wh{E}})$ is the identity map. The composition in the above diagram is $\can(\wh{f} \circ z_{\wh{E}})$, which is also the identity map because $\wh{f} \circ z_{\wh{E}} = z_{\wh{E}'}$ is also a closed embedding. Therefore, $\can(\wh{f})$ is also \emph{the} identity map.

\qed

\section{Fourier analysis of cohomological correspondences}\label{sec: FT cc}
In this section we study how derived Fourier transform interacts with cohomological correspondences. This provides the main content towards Step (5) of the outline \S \ref{ssec: outline}. \emph{In this section, all maps between derived vector bundles (over the same base) are assumed to be linear, i.e., induced by maps of perfect complexes.}

\subsection{Fourier transform of cohomological correspondences}\label{ssec: FT of cc}
In \S \ref{ssec: coco} we defined the notion of cohomological co-correspondence. These arise naturally as the Fourier transforms of cohomological correspondences, as we now explain. 
 
\subsubsection{Over the same base}\label{sssec: FT of cc}
Suppose we have a Cartesian square of derived vector bundles over a base $S$, 
 \begin{equation}\label{eq: Cart square of dvb}
 \begin{tikzcd}
&  C^{\flat} \ar[dl, "p_0"']   \ar[dr, "p_1"] \\ 
 E_0 \ar[dr, "q_0"']&& E_1 \ar[dl, "q_1"]  \\ 
 & C^{\sharp} 
 \end{tikzcd}
 \end{equation}
Dualizing, we obtain a Cartesian square of derived vector bundles over $S$:
 \[
\begin{tikzcd}
& \wh{C^{\sharp}}  \ar[dl, "\wh{q}_0"']   \ar[dr, "\wh{q}_1"] \\ 
\wh{E}_0 \ar[dr, "\wh{p}_0"'] && \wh{E}_1 \ar[dl, "\wh{p}_1"]  \\ 
 & \wh{C^{\flat} }
\end{tikzcd}
\]
We may apply \S \ref{sssec: FT functoriality} to see that the cohomological correspondence $p_0^* \cK_0 \rightarrow p_1^!  \cK_1$ is taken by $\FT_{C^{\flat}}$ to a cohomological co-correspondence on $\wh{C^{\flat}}$
 \[
 \wh{p}_{0!} \FT_{E_{0}}(\cK_0)[-d(p_0)](-d(p_0)) \rightarrow \wh{p}_{1*}  \FT_{E_{1}}(\cK_1) [d(p_1)].
 \]
This may be converted as in \eqref{eq: coh corr to coh co-corr} to a cohomological correspondence on $\wh{C^{\sharp}}$
 \[
\wh{q}_{0}^* \FT_{E_{0}}(\cK_0)\rightarrow \wh{q}_{1}^!  \FT_{E_{1}}(\cK_1) [d(p_0)+d(p_1)](d(p_0))
 \]
 which we call $\FT_{C^\flat}(\cc)$. The construction $\cc\mapsto \FT_{C^\flat}(\cc)$ defines an isomorphism of vector spaces
\begin{equation}\label{FT corr same base}
\FT_{C^{\flat}}: \Corr_{C^{\flat}}(\cK_{0}, \cK_{1})\isom \Corr_{\wh{C^{\sharp}}}(\FT_{E_{0}}(\cK_0), \FT_{E_{1}}(\cK_1)[d(p_0)+d(p_1)](d(p_0)))
\end{equation} 
thanks to the near-involutivity of the derived Fourier transform.

\subsubsection{Varying the base}\label{sssec: FT cc vary base} In \S\ref{sssec: FT of cc} we explained that cohomological correspondences on certain correspondences of derived vector bundles over a base $S$ could be Fourier transformed to a dual correspondence. 

We will define the Fourier transform of a cohomological correspondence in a slightly more general situation, where the base of the derived vector bundles is also permitted to change. Suppose we are given a map of correspondences 
\begin{equation}
\xymatrix{& C^{\flat}\ar[d] \ar[dl]_{p_{0}}\ar[dr]^{p_{1}}\\
E_{0}\ar[d] & C_{S}\ar[dl]_{h_{0}}\ar[dr]^{h_{1}} & E_{1}\ar[d] \\
S_{0} && S_{1}}
\end{equation}
where $E_0, C^{\flat}$ and $E_1$ are derived vector bundles on $S_0, C_{S}$ and $S_1$ respectively.

Let $\wt{E}_0$ and  $\wt{E}_1$ be the pullbacks of $E_{0}$ and $E_{1}$ to $C_S$ via $h_{i}$. We can canonically extend the correspondence $E_{0}\xleftarrow{p_{0}}C^{\flat}\xr{p_{1}} E_{1}$ into a commutative diagram
 \[
 \begin{tikzcd}
& &  C^{\flat}=\wt C^{\flat} \ar[dl, "\wt{p}_0"]   \ar[dr, "\wt{p}_1"']  \ar[ddll, bend right, "p_0"'] \ar[ddrr, bend left, "p_1"] \\ 
& \wt{E}_0 \ar[dr, "\wt p'_1"'] \ar[dl, "h_0^E"] && \wt{E}_1 \ar[dl, "\wt p'_0"]  \ar[dr, "h_1^E"'] \\ 
 E_0 & & \wt C^{\sharp}  & & E_1
 \end{tikzcd}
 \]
Here $\wt C^{\sh}$ is defined to be the pushout of the correspondence of vector bundles $\wt E_{0}\xleftarrow{\wt p_{0}}C^{\flat}\xr{\wt p_{1}} \wt E_{1}$, so that $\wt C^{\sh}$ is also a derived vector bundle over $C_{S}$, and the inner diamond is derived Cartesian. When we view $\wt C^{\flat}$ as a correspondence between $E_{0}$ and $E_{1}$, we denote it by $C^{\flat}$; when we view it as a correspondence between $\wt E_{0}$ and $\wt E_{1}$, we denote it by $\wt C^{\flat}$.  
 
Taking dual derived vector bundles we get a diagram
\[
 \begin{tikzcd}
& &  \wh {C^{\sh}}=\wh{\wt C^{\sharp}} \ar[dl, "\wh{\wt p'_1}"]   \ar[dr, "\wh{\wt p'_0}"']  \ar[ddll, bend right, "\wh{p'_1}"'] \ar[ddrr, bend left, "\wh{p'_0}"]\\ 
& \wh{\wt{E}}_0 \ar[dr, "\wh{\wt{p}}_0"'] \ar[dl, "h_0^{\wh E}"] && \wh{\wt{E}}_1 \ar[dl, "\wh{\wt{p}}_1"]  \ar[dr, "h_1^{\wh E}"'] \\ 
 \wh{E}_0 & & \wh{\wt C^{\flat}}  & & \wh{E}_1
 \end{tikzcd}
\] 
Again, when we view $\wh{\wt C^{\sh}}$ as a correspondence between $\wh E_{0}$  and $\wh E_{1}$, we denote it by $\wh C^{\sh}$.

For $\cK_{i}\in D(E_{i})$, $i=0,1$, we define an isomorphism of vector spaces
\begin{equation}\label{eq: FT of corr}
\FT_{C^{\flat}}: \Corr_{C^{\flat}}(\cK_{0}, \cK_{1})\isom \Corr_{\wh{C}^{\sharp}}(\FT_{E_{0}}(\cK_{0}), \FT_{E_{1}}(\cK_{1})\sm{[d(\wt p_0) + d(\wt p_1)](d(\wt p_0))})
\end{equation}
as the composition of isomorphisms
\begin{eqnarray*}
\Corr_{C^{\flat}}(\cK_{0}, \cK_{1})&=&\Corr_{\wt C^{\flat}}(( h_0^E)^* \cK_0,(h_1^E)^! \cK_1 )\\
&\xr{\FT_{\wt C^{\flat}}}& \Corr_{\wh{\wt C^{\sh}}}(\FT_{\wt E_{0}}(( h_0^E)^* \cK_0),\FT_{\wt E_{1}}((h_1^E)^! \cK_1 )\sm{[d(\wt p_0) + d(\wt p_1)](d(\wt p_0))})\\
&\cong & \Corr_{\wh{\wt C^{\sh}}}((h_0^{\wh E})^* \FT_{E_0}(\cK_0), (h_1^{\wh E})^! \FT_{E_1}(\cK_1)\sm{[d(\wt p_0) + d(\wt p_1)](d(\wt p_0))})\\
&=& \Corr_{\wh{C^{\sharp}}}(\FT_{E_{0}}(\cK_{0}), \FT_{E_{1}}(\cK_{1})\sm{[d(\wt p_0) + d(\wt p_1)](d(\wt p_0))}).
\end{eqnarray*}
Here we have used \S\ref{sssec: FT base change} in the second-to-last isomorphism.

\subsection{Functoriality}\label{sssec: FT functoriality of cc}
We state and prove functorial properties of the Fourier transform of cohomological correspondences constructed in \S\ref{sssec: FT of cc} and more generally in \S\ref{sssec: FT cc vary base}.

\subsubsection{Functoriality over the same base}Suppose we have a commutative diagram 
\begin{equation}\label{eq: FT cc functoriality diagram}
 \begin{tikzcd}
&  C^{\flat} \ar[dl, "p_0"']   \ar[dr, "p_1"] \ar[ddd, "f^{\flat}"', bend left] \\ 
 E_0 \ar[dr, "p_1'"'] \ar[ddd, "f_0"] && E_1 \ar[dl, "p_{0}'"]   \ar[ddd, "f_1"] \\ 
 & C^{\sharp} \ar[ddd, "f^{\sharp}", bend right]\\
&  D^{\flat} \ar[dl, "q_0"']   \ar[dr, "q_1"] \\ 
 F_0 \ar[dr, "q_1'"'] && F_1 \ar[dl, "q_0'"]  \\ 
 & D^{\sharp} 
 \end{tikzcd}
\end{equation}
of derived vector bundles over $S$, where the top and bottom diamonds are derived Cartesian and all maps are linear.

The dual diagram to \eqref{eq: FT cc functoriality diagram} is
\begin{equation}\label{eq: FT cc functoriality diagram dual}
\begin{tikzcd}
 & \wh{D^{\sharp}}  \ar[dl, "\wh{q}'_1"']   \ar[dr, "\wh{q}'_0"]   \ar[ddd, "\wh{f}^{\sharp}"', bend left] \\ 
\wh{F}_0 \ar[dr, "\wh{q}_0"'] \ar[ddd, "\wh{f}_0"]  && \wh{F}_1 \ar[dl, "\wh{q}_1"]   \ar[ddd, "\wh{f}_1"]\\ 
 & \wh{D^{\flat}}  \ar[ddd, "\wh{f}^{\flat}", bend right]  \\
 & \wh{C^{\sharp}}   \ar[dl, "\wh{p}_1'"']   \ar[dr, "\wh{p}_0'"] \\ 
\wh{E}_0 \ar[dr, "\wh{p}_0"']&& \wh{E}_1 \ar[dl, "\wh{p}_1"]  \\ 
 & \wh{{C}^{\flat} }
 \end{tikzcd}
\end{equation}

\begin{lemma}\label{lem: left pushable dualize to right pullable} In the above setup, the following are equivalent:
\begin{enumerate}
\item The map of correspondences $f^{\flat}: C^{\flat}\to D^{\flat}$ is left pushable;
\item The map of co-correspondences $f^{\sh}: C^{\sh}\to D^{\sh}$ is right pushable;
\item The map of co-correspondences $\wh f^{\flat}: \wh {D^{\flat}}\to \wh{C^{\flat}}$ is left pullable;
\item The map of correspondences $\wh f^{\sh}: \wh{D^{\sh}}\to \wh{C^{\sh}}$ is right pullable.
\end{enumerate}
Moreover, when (3) and (4) hold, we have $\d_{\wh f^{\flat}}=\d_{\wh f^{\sh}}$.
\end{lemma}
\begin{proof} Using that the bottom diamond in \eqref{eq: FT cc functoriality diagram} is derived Cartesian, we see that the following diagram is derived Cartesian
\begin{equation}
\xymatrix{ C^{\flat}\ar[r]^-{c^{\flat}}\ar[d]_{p_{1}} & E_{0}\times_{F_{0}}D^{\flat}\ar[d]^{p'_{1}\times_{q'_{1}} q_{1}}\\
E_{1}\ar[r]^-{e_{1}} &  C^{\sh}\times_{D^{\sh}}F_{1}
}
\end{equation}
Now (1) $\iff$ $c^{\flat}$ is a closed embedding of derived vector bundles $\iff$ $e_{1}$ is a closed embedding of derived vector bundles (by the above derived Cartesian diagram) $\iff$ (2). This proves (1) $\iff$ (2). The same argument shows (3) $\iff$ (4) and that $\d_{\wh f^{\flat}}=\d_{\wh f^{\sh}}$.

It remains to show that (1) $\iff$ (3).  Let $\cE_{i}, \cF_{i}, \cC^{\flat},\cdots$ be perfect complexes over $S$ whose total spaces are $E_{i}, F_{i}, C^{\flat}, \cdots$. We christen the maps between these perfect complexes with the same names of the induced maps between their total spaces. Let $P$ be the derived fiber of the linear map $(p_{0}, f^{\flat}): C^{\flat}\to E_{0}\times_{F_{0}} D^{\flat}$ over the zero section $0_{S}$; let $Q$ be the derived fiber of $(\wh q_{0}, \wh f_{0}):\wh F_{0}\to \wh{D^{\flat}}\times_{\wh {C^{\flat}}}\wh E_{0}$ over the zero section $0_{S}$.

Note that $P$ is the total space of the perfect complex $\cP$ over $S$ that is obtained by taking the total complex of the double complex
\begin{equation}
\cC^{\flat}\xr{(p_{0}, f^{\flat})} \cE_{0}\op \cD^{\flat}\xr{f_{0}-q_{0}} \cF_{0}.
\end{equation}
Here the displayed complex means the totalization with $\cC^{\flat}$ in horizontal degree $0$, $\cE_{0}\op \cD^{\flat}$ in horizontal degree $1$, and $\cF_{0}$ in horizontal degree $2$, while retaining the internal cohomological gradings. By the preceding closed-embedding criterion, (1) is equivalent to
\begin{enumerate}
\item[(1')] $\cP$ is locally represented by a complex of vector bundles in degrees $\ge1$. 
\end{enumerate}
 
Similarly, $Q$ is the total space of the perfect complex $\cQ$ over $S$ that is obtained by taking the total complex of the double complex
\begin{equation}
\wh\cF_{0}\xr{(\wh q_{0}, \wh f_{0})} \wh{\cD^{\flat}}\op \wh\cE_{0}\xr{\wh f^{\flat}-\wh p_{0}} \wh{\cC^{\flat}}.
\end{equation}
Here $\wh\cF_{0}$ denotes the derived dual $\cF_{0}^{*}$, and the same horizontal-degree convention is used. Since pullability is quasi-smoothness of the corresponding comparison map, Lemma \ref{lem:closed-embedding} gives that (3) is equivalent to
\begin{enumerate}
\item[(3')] $\cQ$ is locally represented by a complex of vector bundles in degrees $\le1$. 
\end{enumerate}

With these totalization conventions, dualizing the three-term complex defining $\cP$ reverses the horizontal degree and identifies $\cQ$ with $\cP^{*}[-2]$. Thus $\cP$ has tor-amplitude in degrees $\ge 1$ if and only if $\cQ$ has tor-amplitude in degrees $\le 1$, proving (1') $\iff$ (3') and hence (1) $\iff$ (3).

\end{proof}

Assume that $f^{\flat}: C^{\flat}\to D^{\flat}$ is left pushable. Then for $\cK_{i}\in D(E_i)$ ($i=0,1$), the pushforward map
\begin{equation}
(f^{\flat})_{!}: \Corr_{C^{\flat}}(\cK_{0}, \cK_{1})\to \Corr_{D^{\flat}}(f_{0!}\cK_{0}, f_{1!}\cK_{1})
\end{equation}
is defined. By  Lemma \ref{lem: left pushable dualize to right pullable}, the map of correspondences $\wh{f^{\sh}}: \wh{D^{\sh}}\to \wh{C^{\sh}}$ is right pullable. Hence for $\cL_{i}\in D(\wh E_{i})$ ($i=0,1$), the pullback map
\begin{equation}
(\wh{f^{\sh}})^{*}: \Corr_{\wh{C^{\sh}}}(\cL_{0}, \cL_{1})\to \Corr_{\wh{D^{\sh}}}(\wh{f}_{0}^{*}\cL_{0}, \wh{f}_{1}^{*}\cL_1\tw{-\d_{\wh{f^{\sh}}}})
\end{equation}
is defined. 

On the other hand, applying Fourier transform to $f_{0!}\cK_{0}$ and $f_{1!}\cK_{1}$, we have by \S\ref{sssec: FT functoriality}
\begin{equation}
\FT_{F_{0}}(f_{0!}\cK_{0})\cong \wh f_{0}^{*}\FT_{E_{0}}(\cK_{0})[-d(f_{0})], \quad \FT_{F_{1}}(f_{1!}\cK_{1})\cong \wh f_{1}^{*}\FT_{E_{1}}(\cK_{1})[-d(f_{1})].
\end{equation}

\begin{lemma}\label{l:FT corr dim}
We have
\begin{enumerate}
\item $\d_{\wh{f^{\sh}}}=d(p_{0})-d(q_{0})$.
\item $d(p_{1})-d(q_{1})+d(f_{1})=d(p_{0})-d(q_{0})+d(f_{0})$.
\end{enumerate}
\end{lemma}
\begin{proof}
(1) Let $\cP$ be the perfect complex computing the derived fiber of $C^{\flat}\to E_{0}\times_{F_{0}}D^{\flat}$, as in the proof of Lemma \ref{lem: left pushable dualize to right pullable}. That proof identifies the fiber complex defining the defect of $\wh f^{\sh}$ with $\cP^{*}[-2]$, whose Euler characteristic is the same as that of $\cP$. Since $\cP$ is the total complex of $\cC^{\flat}\to \cE_{0}\op \cD^{\flat}\to \cF_{0}$, this Euler characteristic is
\[
\rank(C^{\flat})-\rank(E_{0})-\rank(D^{\flat})+\rank(F_{0})=d(p_{0})-d(q_{0}).
\]

(2) Both sides are equal to $\rank(C^{\flat})-\rank(D^{\flat})$. Indeed,
\[
\begin{aligned}
d(p_{i})-d(q_{i})+d(f_{i})
&=\bigl(\rank(C^{\flat})-\rank(E_{i})\bigr)-\bigl(\rank(D^{\flat})-\rank(F_{i})\bigr)\\
&\quad+\bigl(\rank(E_{i})-\rank(F_{i})\bigr)
=\rank(C^{\flat})-\rank(D^{\flat})
\end{aligned}
\]
for $i=0,1$.
\end{proof}

\begin{prop}\label{prop: cc push/pull functoriality tilde}Assume that diagram \eqref{eq: FT cc functoriality diagram} is globally presented. 

(1) Suppose the map of correspondences $f^{\flat}: C^{\flat}\to D^{\flat}$ is left pushable. Let $\cK_i \in D(E_i)$ for $i=0,1$. Then the following diagram commutes
\begin{equation}
\begin{gathered}
\begin{adjustbox}{max width=\textwidth}
$\displaystyle
\xymatrix{\Corr_{C^{\flat}}(\cK_{0}, \cK_{1}) \ar[rr]^-{\FT_{C^{\flat}}}\ar[d]^{(f^{\flat})_{!}}& &  \Corr_{\wh{C^{\sh}}}(\FT_{E_{0}}(\cK_0), \FT_{E_{1}}(\cK_1)\sm{[d(p_0)+d(p_1)](d(p_0))})\ar[d]^{(\wh{f}^{\sh})^{*}}\\
\Corr_{D^{\flat}}(f_{0!}\cK_{0}, f_{1!}\cK_{1})\ar[rr]^-{\TT_{[d(f_0)](0)}\FT_{D^{\flat}}} &  & \Corr_{\wh{D^{\sh}}}(\wh{f}_{0}^{*}\FT_{E_{0}}(\cK_0), \wh{f}_{1}^{*}\FT_{E_{1}}(\cK_1)\sm{[d(q_0)+d(q_1)+d(f_{0})-d(f_{1})](d(q_0))})
}
$
\end{adjustbox}
\end{gathered}
\end{equation}
Here we use Lemma \ref{l:FT corr dim} to match the differences of the twists that appear in the right vertical map with $\tw{-\d_{\wh{f^{\sh}}}}$, which is the correct twist for $(\wh{f^{\sh}})^{*}$.

(2) Suppose the map of correspondences $f^{\flat}: C^{\flat}\to D^{\flat}$ is right pullable. Let $\cK_i \in D(F_i)$ for $i=0,1$. Then the following diagram commutes
\[
\begin{gathered}
\begin{adjustbox}{max width=\textwidth}
$\displaystyle
\begin{tikzcd}[column sep = huge, ampersand replacement=\&]
\Corr_{D^{\flat}}(\cK_{0}, \cK_{1}) \ar[rr, "\FT_{D^{\flat}}"] \ar[d, "(f^{\flat})^*"] \& \& {\Corr_{\wh{D^{\sh}}}(\FT_{F_{0}}(\cK_0), \FT_{F_{1}}(\cK_1)\sm{[d(q_0)+d(q_1)](d(q_0))})} \ar[d, "(\wh{f}^{\sh})_{!}"] \\
\Corr_{C^{\flat}}(f_{0}^*\cK_{0}, f_{1}^*\cK_{1}\tw{-\d_{f^{\flat}}})   \ar[rr, "\TT_{[d(f_0)](d(f_0))} \FT_{C^{\flat}}"]  \& \& {\Corr_{\wh{C^{\sh}}}(\wh{f}_{0!}\FT_{F_{0}}(\cK_0), \wh{f}_{1!}\FT_{F_{1}}(\cK_1)\sm{[d(q_0)+d(q_1)](d(q_0))})} 
\end{tikzcd}
$
\end{adjustbox}
\end{gathered}
\]
In (2), the lower-left term includes the defect twist from the right-pullback operation, and here $\d_{f^{\flat}}=d(p_{1})-d(q_{1})$. With this identity and Lemma \ref{l:FT corr dim}, the shift and twist introduced by $\TT_{[d(f_0)](d(f_0))}\FT_{C^{\flat}}$ match the displayed lower-right term.
\end{prop}

\begin{proof} (1) Below, to shorten notation, we write $\wh\cK_{i}:=\FT_{E_{i}}(\cK_{i})$. By definition, $\FT_{C^{\flat}}$ is the composition of two isomorphisms
\begin{equation}
\begin{gathered}
\begin{aligned}
\FT_{C^{\flat}}: \Corr_{C^{\flat}}(\cK_{0}, \cK_{1})
&\xr{\FT'_{C^{\flat}}}\CoCorr_{\wh{ C^{\flat}}}(\wh\cK_{0}, \wh\cK_{1}\sm{[d(p_{0})+d(p_{1})](d(p_{0}))})\\
&\xr{\g_{C}^{-1}}\Corr_{\wh {C^{\sh}}}(\wh\cK_{0}, \wh\cK_{1}\sm{[d(p_{0})+d(p_{1})](d(p_{0}))}).
\end{aligned}
\end{gathered}
\end{equation}
Here $\g_C^{-1}$ denotes the inverse of the comparison isomorphism \eqref{eq: coh corr to coh co-corr} for the dual Cartesian diamond.
Similar remarks apply to $\FT_{D^{\flat}}$. Therefore it suffices to prove the commutativity of the following two diagrams separately
\begin{equation}\label{FT push pull coco}
\begin{gathered}
\begin{adjustbox}{max width=\textwidth}
$\displaystyle
\xymatrix{\Corr_{C^{\flat}}(\cK_{0}, \cK_{1}) \ar[rr]^-{\FT'_{C^{\flat}}}\ar[d]^{(f^{\flat})_{!}}& & \CoCorr_{\wh{C^{\flat}}}(\wh\cK_0, \wh\cK_1\sm{[d(p_0)+d(p_1)](d(p_0))})\ar[d]^{(\wh{f}^{\flat})^{*}}\\
\Corr_{D^{\flat}}(f_{0!}\cK_{0}, f_{1!}\cK_{1})\ar[rr]^-{\TT_{[d(f_0)]} \FT'_{D^{\flat}}} & &  \CoCorr_{\wh{D^{\flat}}}(\wh{f}_{0}^{*}\wh\cK_0, \wh{f}_{1}^{*}\wh\cK_1\sm{[d(q_0)+d(q_1)+d(f_{0})-d(f_{1})](d(q_0))})
}
$
\end{adjustbox}
\end{gathered}
\end{equation}
and
\begin{equation}\label{coco to co diagram}
\xymatrix{\CoCorr_{\wh{C^{\flat}}}(\wh\cK_0, \wh\cK_1) \ar[d]^{(\wh f^{\flat})^{*}}\ar[r]^{\g_{C}^{-1}}& 
\Corr_{\wh{C^{\sh}}}(\wh\cK_0, \wh\cK_1)\ar[d]^{(\wh f^{\sh})^{*}}\\
\CoCorr_{\wh{D^{\flat}}}(\wh f_{0}^{*}\wh\cK_0, \wh f_{1}^{*}\wh\cK_1) \ar[r]^{\g_{D}^{-1}}& 
\Corr_{\wh{D^{\sh}}}(\wh f_{0}^{*}\wh\cK_0, \wh f_{1}^{*}\wh\cK_1)
}
\end{equation}
Here the pullback map of co-correspondences $(\wh f^{\flat})^{*}$ is defined in \S\ref{sssec: pull push coco}. In \eqref{coco to co diagram}, the lower terms are displayed with the common defect twists suppressed: the second arguments carry the twists attached to $(\wh f^\flat)^*$ and $(\wh f^\sh)^*$, and Lemma \ref{l:FT corr dim} identifies the corresponding defects. The commutativity of the diagram \eqref{coco to co diagram} is proved in Proposition \ref{prop: comp coco}.

It remains to show that \eqref{FT push pull coco} is commutative. Since the statement does not involve $C^{\sh}$ and $D^{\sh}$, we will omit the superscript $\flat$ from the notations and denote $C^{\flat}, D^{\flat}$ by $C$ and $D$. For $\cc\in \Corr_{C}(\cK_{0}, \cK_{1})$, denoting  $\wh\cK_{i}=\FT_{E_{i}}(\cK_{i})$, we have to show the commutativity of the outer square of the diagram
\begin{equation}\label{FT bc cc}
\xymatrix{\FT(q^{*}_{0}f_{0!}\cK_{0}) \ar@{-}[d]_{\wr}\ar[r]^{\td} & \FT(f_{!}p_{0}^{*}\cK_{0})\ar@{-}[d]_{\wr} \ar[rr]^-{\FT(f_{!}\cc)}&& \FT(f_{!}p_{1}^{!}\cK_{1})\ar@{-}[d]_{\wr}\ar[r] & \FT(q_{1}^{!}f_{1!}\cK_{1}) \ar@{-}[d]_{\wr}\\
\wh q_{0!}\wh f_{0}^{*}\wh\cK_{0}[?](?)\ar[r]^{\tu} & \wh f^{*}\wh p_{0!}\wh\cK_{0}[?](?) \ar[rr]^-{\wh f^{*}\FT(\cc)} && \wh f^{*}\wh p_{1*}\wh \cK_{1}[?](?)\ar[r] & \wh q_{1*}\wh f_{1}^{*}\wh\cK_{1}[?](?)
}
\end{equation}
Here the arrows marked by $\td$ and $\tu$ are the base change maps attached to the pushable square $(C, E_{0}, D, F_{0})$ and the pullable square $(\wh F_{0}, \wh D, \wh E_{0}, \wh C)$. The unmarked arrows are induced by the natural transformation  $f_{!}p_{1}^{!}\to q_{1}^{!}f_{1!}$ attached to the square $(C, E_{1}, D, F_{1})$ and the natural transformation $\wh f^{*}\wh p_{1*}\to \wh q_{1*}\wh f_{1}^{*}$ attached to the dual square. The vertical isomorphisms are from \S\ref{sssec: FT functoriality}. We have omitted the shifts and twists in the bottom row.

The middle square above is commutative by the naturality of the isomorphisms in \S\ref{sssec: FT functoriality}. The right square is commutative: write $f_{!}p_{1}^{!}\to q_{1}^{!}f_{1!}$ as the composition
\begin{equation}
f_{!}p_{1}^{!}\xr{\textup{unit}} f_{!}p_{1}^{!}f_{1}^{!}f_{1!}=f_{!}f^{!}q_{1}^{!}f_{1!}\xr{\textup{counit}}q_{1}^{!}f_{1!}
\end{equation}
Therefore it suffices to note that $\FT$ transforms the unit map $\Id\to f^!_1 f_{1!}$ (resp. counit map $f_{!}f^{!}\to \Id$) to 
the unit map $\Id\to \wh f_{1*}\wh f_1^{*}$ (resp. the counit map $\wh f^{*}\wh f_{*}\to \Id$), as explained in \S \ref{sssec: FT functoriality}.

It remains to show that the left square in \eqref{FT bc cc} commutes.  Let $C^{\na}=E_{0}\times_{F_{0}}D$, with the induced map $c: C\to C^{\na}$. The square $(C, E_{0}, D, F_{0})$ can be decomposed as the composition of two squares of derived vector bundles over $S$
\begin{equation}\label{decomp CDEF}
\xymatrix{E_{0}\ar@{=}[d]  & C \ar[d]^{c}\ar[l]_{p_{0}}\\
E_{0}\ar[d]_{f_{0}} & C^{\na}\ar[d]^{f^{\na}}\ar[l]_{p_{0}^{\na}}\\
F_{0} & D\ar[l]_{q_{0}}\\
}
\end{equation}
where $c$ is proper by assumption, and the bottom square is derived Cartesian by definition. Note that the assumptions guarantee that \eqref{decomp CDEF} is globally presented. Using the compatibility of the base change maps with composition of squares proved in Proposition \ref{p: comp push} and Proposition \ref{p: comp pull}, we reduce to showing the commutativity of the left square in \eqref{FT bc cc} separately for the two squares in \eqref{decomp CDEF}, i.e., for two special cases:
\begin{enumerate}
\item The square $(C, E_{0}, D, F_{0})$ is derived Cartesian (hence so is its dual square).
\item The map $f_{0}=\Id: E_{0}\to F_{0}=E_{0}$ is the identity map, and $f^{\na}$ is proper (i.e., a closed embedding of derived vector bundles; dually $\wh{f^{\na}}$ is smooth).
\end{enumerate}
The first case follows from Proposition \ref{prop: ft pbc}  (proved in Proposition \ref{prop: FT hexagon}), which applies because \eqref{decomp CDEF} is globally presented. In the second case, we reduce to showing that the outer square in the following diagram is commutative 
\begin{equation}
\begin{adjustbox}{max width=\textwidth}
\xymatrix{\FT q_{0}^{*}\ar@{-}[d]_{\wr} \ar[r]^-{\textup{unit}}& \FT (f^{\na})_{*}(f^{\na})^{*}q_{0}^{*}\ar@{-}[d]_{\wr}\ar@{=}[r] & \FT (f^{\na})_{*}(p_{0}^{\na})^{*}\ar@{-}[d]_{\wr} & \FT (f^{\na})_{!}(p_{0}^{\na})^{*}\ar[l]_{\cong}\ar@{-}[d]_{\wr}\\
\wh q_{0!}\FT\sm{[?](?)}\ar[r]^-{\textup{unit}} & (\wh{f^{\na}})^{!}(\wh{f^{\na}})_{!}\wh q_{0!}\FT\sm{[?](?)} \ar@{=}[r]& (\wh{f^{\na}})^{!}\wh{p_{0}^{\na}}_{!}\FT\sm{[?](?)} & (\wh{f^{\na}})^{*}\wh{p_{0}^{\na}}_{!}\FT\sm{[?](?)}\tw{d(\wh{f^{\na}})}\ar[l]_-{[\wh{f^{\na}}]}
}
\end{adjustbox}
\end{equation}
Here $?=-d(q_{0})$. (The global presentability assumption is used here to produce the first and second vertical maps, since $q_0$ is not assumed to be smooth or a closed embedding.) Now the left and middle squares commute because $\FT$ takes the unit to the unit. The right square commutes by Proposition \ref{prop: fourier of Gysin}.

(2) An analogous argument can be applied. Alternatively, (2) follows from (1) using the near-involutivity of $\FT$. 
\end{proof}

\subsubsection{Functoriality over varying bases}\label{sss: push varying base} We next extend the preceding discussion to the situation where the base space may vary. Consider a diagram of maps of correspondences
\begin{equation}\label{linear corr varying base}
\xymatrix{ & C^{\flat}\ar[dl]_{p_{0}}\ar[dr]^{p_{1}}\ar[d]^{f^{\flat}}\\
E_{0}\ar[d]_{f_{0}} & D^{\flat}\ar[dl]_{q_{0}}\ar[dr]^{q_{1}}\ar[d]^{} & E_{1}\ar[d]^{f_{1}}\\
F_{0}\ar[d] & C_{S}\ar[dl]_{h_{0}}\ar[dr]^{h_{1}} & F_{1}\ar[d]\\
S_{0} && S_{1}
}
\end{equation}
where $E_{i}$ and $F_{i}$ are derived vector bundles over $S_{i}$ (for $i=0,1$),  and $C^{\flat}$ and $D^{\flat}$ are derived vector bundles over $C_{S}$. All maps between derived vector bundles are assumed to be linear.

Let $\wt{E}_i \rightarrow C_S$, $\wt{F}_i \rightarrow C_S$ and $\wt f_{i}: \wt{E}_i \rightarrow \wt{F}_{i}$ be the base changes of $E_{i}, F_{i}$ and $f_{i}$ along $h_i \co C_S \rightarrow S_i$. Using the discussion in \S\ref{sssec: FT cc vary base}, we can canonically extend the upper part of the diagram \eqref{linear corr varying base} into a  commutative diagram 
\begin{equation}\label{eq: FT cc functoriality diagram varying base}
 \begin{tikzcd}
& &  C^{\flat}=\wt C^{\flat} \ar[dl, "\wt{p}_0"] \ar[dd, phantom, "\di"]  \ar[dr, "\wt{p}_1"']  \ar[ddll, bend right, "p_0"']   \ar[ddrr, bend left, "p_1"] \ar[ddd, "f^{\flat}=\wt{f}^{\flat}"', bend left] \\ 
& \ar[dl, "h_0^E"]  \wt{E}_0 \ar[dr, "p'_{1}"'] \ar[ddd, "\wt{f}_0"] && \wt{E}_1 \ar[dl, "p'_0"]   \ar[ddd, "\wt{f}_1"]  \ar[dr, "h_1^E"'] \\ 
E_0 \ar[ddr, phantom, "\di"]\ar[ddd, "f_0"]  &  & \wt C^{\sharp}  \ar[ddd, "\wt{f}^{\sh}", bend right] & & E_1 \ar[ddd, "f_1"]\ar[ddl, phantom, "\di"] \\
& &  D^{\flat}=\wt{D}^{\flat} \ar[dd, phantom, "\di"]\ar[dl, "\wt q_0"]   \ar[dr, "\wt q_{1}"']  \ar[ddll, bend right, "q_0"']  \ar[ddrr, bend left, "q_1"] \\ 
&\ar[dl, "h_0^{F}"]   \wt{F}_0 \ar[dr, "\wt q'_1"'] && \wt{F}_1 \ar[dl, "\wt q'_0"] \ar[dr, "h_1^{F}"']    \\ 
F_0  & &  \wt{D}^{\sharp} & & F_1
 \end{tikzcd}
\end{equation}
where the squares labeled by $\di$ are derived Cartesian.

Since the leftmost parallelogram is derived Cartesian, the square $(C^{\flat}, E_{0}, D^{\flat}, F_{0})$ is pushable if and only if the square $(\wt C^{\flat}, \wt E_{0}, \wt D^{\flat}, \wt F_{0})$ is pushable. In other words, the morphism $f^{\flat}: C^{\flat}\to D^{\flat}$ of correspondences is left pushable if and only if the morphism $\wt f^{\flat}: \wt C^{\flat}\to \wt D^{\flat}$ of correspondences is left pushable. When any of these equivalent conditions holds, we have a pushforward map
\begin{equation}
f^{\flat}_{!}: \Corr_{C^{\flat}}(\cK_{0}, \cK_{1})\to \Corr_{D^{\flat}}(f_{0!}\cK_{0}, f_{1!}\cK_{1}).
\end{equation}

\subsubsection{}\label{sss: pull varying base} The dual diagram to \eqref{eq: FT cc functoriality diagram varying base} is: 
\begin{equation}\label{eq: FT cc functoriality diagram varying base dual}
\begin{tikzcd}
 & & \wh{D^{\sharp}}=\wh{\wt D^{\sh}}  \ar[dl, "\wh{\wt q'_1}"]   \ar[dr, "\wh{\wt q'_0}"']    \ar[ddd, "\wh{f}^{\sh}=\wh{\wt{f}^{\sh}}"', bend left] \ar[ddll, bend right, "\wh{q'_1}"']   \ar[ddrr, bend left, "\wh{q'_0}"]\\ 
&\ar[dl, "h_0^{\wh F}"]  \wh{\wt{F}}_0 \ar[dr, "\wh{\wt{q}_0}"'] \ar[ddd, "\wh{\wt{f}}_0"]  && \wh{\wt{F}}_1 \ar[dl, "\wh{\wt{q}}_1"]   \ar[ddd, "\wh{\wt{f}}_1"]  \ar[dr, "h_1^{\wh F}"']  \\ 
\wh{F}_0  \ar[ddd, "\wh{f}_0"] &  & \wh{\wt{D}^{\flat}}  \ar[ddd, "\wh{\wt f}^{\flat}", bend right]    & & \wh{F}_1  \ar[ddd, "\wh{f}_1"] \\
 & & \wh{C^{\sh}}=\wh{\wt C^{\sharp}}   \ar[dl, "\wh{\wt p'_1}"]   \ar[dr, "\wh{\wt p'_0}"'] \ar[ddll, bend right, "\wh{p'_1}"']   \ar[ddrr, bend left, "\wh{p'_0}"]\\ 
& \wh{\wt{E}}_0 \ar[dr, "\wh{\wt{p}}_0"']  \ar[dl, "h_0^{\wh E}"]  && \wh{\wt{E}}_1 \ar[dr, "h_1^{\wh E}"']   \ar[dl, "\wh{\wt{p}}_1"]  \\ 
\wh{E}_0 &  & \wh{\wt{C}^{\flat} }  && \wh{E}_1
 \end{tikzcd}
\end{equation}
Since the rightmost parallelogram is derived Cartesian, the square $(\wh{D^{\sh}}, \wh F_{1}, \wh{C^{\sh}}, \wh E_{1})$ is pullable if and only if the square $(\wh{\wt D^{\sh}}, \wh{\wt F_{1}}, \wh{\wt C^{\sh}}, \wh{\wt E_{1}})$ is pullable. In other words, the morphism $\wh f^{\sh}: \wh{D^{\sh}}\to \wh{C^{\sh}}$ of correspondences is right pullable if and only if the morphism $\wh{\wt f^{\sh}}: \wh{\wt D^{\sh}}\to \wh{\wt C^{\sh}}$ of correspondences is right pullable. When any of these equivalent conditions holds, we have a pullback map
\begin{equation}
(\wh{f^{\sh}})^{*}: \Corr_{\wh{C^{\sh}}}(\cL_{0}, \cL_{1})\to \Corr_{\wh{D^{\sh}}}(\wh{f}_{0}^{*}\cL_{0}, \wh{f}_{1}^{*}\cL_{1}\tw{-\d_{\wh{ f^{\sh}}}}).
\end{equation}

Moreover, by Lemma \ref{lem: left pushable dualize to right pullable}, $f^{\flat}$ is left pushable if and only if $\wh{f^{\sh}}$ is right pullable.

\begin{prop}\label{prop: cc push/pull functoriality}  Assume the diagram \eqref{eq: FT cc functoriality diagram varying base dual} is globally presented.

(1) Suppose the map of correspondences $f^{\flat}: C^{\flat}\to D^{\flat}$ is left pushable. Let $\cK_i \in D(E_i)$ for $i=0,1$. Then the following diagram commutes:
\begin{equation}\label{FT cc push pull diagram}
\begin{gathered}
\begin{adjustbox}{max width=\textwidth}
$\displaystyle
\xymatrix{\Corr_{C^{\flat}}(\cK_{0}, \cK_{1}) \ar[rr]^-{\FT_{C^{\flat}}}\ar[d]^{(f^{\flat})_{!}}& &  \Corr_{\wh{C^{\sh}}}(\FT_{E_{0}}(\cK_0), \FT_{E_{1}}(\cK_1)\sm{[d(\wt p_0)+d(\wt p_1)](d(\wt p_0))})\ar[d]^{(\wh{f}^{\sh})^{*}}\\
\Corr_{D^{\flat}}(f_{0!}\cK_{0}, f_{1!}\cK_{1})\ar[rr]^-{\TT_{[d(f_0)]}\FT_{D^{\flat}}} & & \Corr_{\wh{D^{\sh}}}(\wh{f}_{0}^{*}\FT_{E_{0}}(\cK_0), \wh{f}_{1}^{*}\FT_{E_{1}}(\cK_1)\sm{[d(\wt q_0)+d(\wt q_1)+d(f_{0})-d(f_{1})](d(\wt q_0))})
}
$
\end{adjustbox}
\end{gathered}
\end{equation}
Here we use Lemma \ref{l:FT corr dim} to match the differences of the twists that appear in the right vertical map with $\tw{-\d_{\wh{f^{\sh}}}}$, which is the correct twist for $(\wh{f^{\sh}})^{*}$. 

(2) Suppose the map of correspondences $f^{\flat}: C^{\flat}\to D^{\flat}$ is right pullable. Let $\cK_i \in D(F_i)$ for $i=0,1$. Then the following diagram commutes
\[
\begin{gathered}
\begin{adjustbox}{max width=\textwidth}
$\displaystyle
\begin{tikzcd}[column sep = huge, ampersand replacement=\&]
\Corr_{D^{\flat}}(\cK_{0}, \cK_{1}) \ar[rr, "\FT_{D^{\flat}}"] \ar[d, "(f^{\flat})^*"] \& \& {\Corr_{\wh{D^{\sh}}}(\FT_{F_{0}}(\cK_0), \FT_{F_{1}}(\cK_1)\sm{[d(\wt{q}_0)+d(\wt{q}_1)](d(\wt{q}_0))})} \ar[d, "(\wh{f}^{\sh})_{!}"] \\
\Corr_{C^{\flat}}(f_{0}^*\cK_{0}, f_{1}^*\cK_{1}\tw{-\d_{f^{\flat}}})   \ar[rr, "\TT_{[d(f_0)](d(f_0))} \FT_{C^{\flat}}"] \& \&  {\Corr_{\wh{C^{\sh}}}(\wh{f}_{0!}\FT_{F_{0}}(\cK_0), \wh{f}_{1!}\FT_{F_{1}}(\cK_1)\sm{[d(\wt{q}_0)+d(\wt{q}_1)](d(\wt{q}_0))})} 
\end{tikzcd}
$
\end{adjustbox}
\end{gathered}
\]
As in the same-base case, the lower-left term includes the defect twist for the right pullback; in tilded notation $\d_{f^{\flat}}=d(\wt p_{1})-d(\wt q_{1})$.
\end{prop}

\begin{proof} (1) Let $\wt\cK_{0}= (h^{E}_{0})^{*}\cK_{0}$ and $\wt\cK_{1}=(h^{E}_{1})^{!}\cK_{1}$. Using the definition of $\FT_{C^{\flat}}$ and $\FT_{D^{\flat}}$ from \S\ref{sssec: FT cc vary base}, we decompose the square \eqref{FT cc push pull diagram} into three squares
\begin{equation}\label{FT cc to tilde}
\xymatrix{\Corr_{C^{\flat}}(\cK_{0}, \cK_{1}) \ar[r]^{\sim} \ar[d]^{(f^{\flat})_{!}}& \Corr_{\wt C^{\flat}}(\wt\cK_0, \wt \cK_1)\ar[d]^{(\wt{f}^{\flat})_{!}}\\
\Corr_{D^{\flat}}(f_{0!}\cK_{0}, f_{1!}\cK_{1})\ar[r]^{\sim} & \Corr_{\wt D^{\flat}}(\wt f_{0!}\wt\cK_0, \wt f_{1!}\wt\cK_1)
}
\end{equation}
\begin{equation}\label{FT cc comp tilde}
\begin{gathered}
\begin{adjustbox}{max width=\textwidth}
$\displaystyle
\xymatrix{\Corr_{\wt C^{\flat}}(\wt \cK_{0}, \wt \cK_{1}) \ar[rr]^-{\FT_{\wt C^{\flat}}}\ar[d]^{(\wt f^{\flat})_{!}}& & \Corr_{\wh{\wt C^{\sh}}}(\FT(\wt\cK_0), \FT(\wt\cK_1)\sm{[d(\wt p_0)+d(\wt p_1)](d(\wt p_0))})\ar[d]^{(\wh{\wt f^{\sh}})^{*}}\\
\Corr_{\wt D^{\flat}}(\wt f_{0!}\wt\cK_{0}, \wt f_{1!}\wt \cK_{1})\ar[rr]^-{\TT_{[d(f_0)]} \FT_{\wt D^{\flat}}} & & \Corr_{\wh{\wt D^{\sh}}}((\wh{\wt f_{0}})^{*}\FT(\wt\cK_0), (\wh{\wt f_{1}})^{*}\FT(\wt\cK_1)\sm{[d(\wt q_0)+d(\wt q_1)+d(f_{0})-d(f_{1})](d(\wt q_0))})
}
$
\end{adjustbox}
\end{gathered}
\end{equation}
and
\begin{equation}\label{FT cc from tilde}
\xymatrix{\Corr_{\wh{\wt C^{\sh}}}(\FT(\wt\cK_0), \FT(\wt\cK_1)) \ar[d]^{(\wh{\wt f^{\sh}})^{*}} \ar[r]^{\sim} & \Corr_{\wh{C^{\sh}}}(\FT(\cK_0), \FT(\cK_1))\ar[d]^{(\wh f^{\sh})^{*}}\\
\Corr_{\wh{\wt D^{\sh}}}((\wh{\wt f_{0}})^{*}\FT(\wt\cK_0), (\wh{\wt f_{1}})^{*}\FT(\wt\cK_1)) \ar[r]^{\sim} & \Corr_{\wh{D^{\sh}}}(\wh f_{0}^{*}\FT(\cK_0), \wh f_{1}^{*}\FT(\cK_1))}
\end{equation}
As in \eqref{coco to co diagram}, the defect twists in the second arguments of the lower terms in \eqref{FT cc from tilde} are suppressed from the notation.
The commutativity of \eqref{FT cc comp tilde} is proved in Proposition \ref{prop: cc push/pull functoriality tilde}. 

Let us prove the commutativity of \eqref{FT cc to tilde}. After  unraveling definitions, the non-obvious part is to show the commutativity of the following diagram of functors $D(E_{0})\to D(D^{\flat})$
\begin{equation}
\xymatrix{\wt q_{0}^{*}(h^{F}_{0})^{*}f_{0!}\ar@{=}[d]\ar[r]^{\wt q_{0}^{*}\di} & \wt q_{0}^{*}\wt f_{0!}(h^{E}_{0})^{*} \ar[r]^{\td(h^{E}_{0})^{*}}& f^{\flat}_{!}\wt p_{0}^{*}(h^{E}_{0})^{*}\ar@{=}[d]\\
q_{0}^{*}f_{0!} \ar[rr]^{\td} & & f^{\flat}_{!}p_{0}^{*} }
\end{equation}
This follows from Proposition \ref{p: comp push horizontal} applied to the two squares
\begin{equation}
\xymatrix{ C^{\flat}\ar[r]^{\wt p_{0}} \ar[d]^{f^{\flat}}& \wt E_{0}\ar[r]^{h^{E}_{0}}\ar[d]^{\wt f_{0}} & E_{0}\ar[d]^{f_{0}}\\
D^{\flat} \ar[r]^{\wt q_{0}}& \wt F_{0}\ar[r]^{h^{F}_{0}} & F_{0}
}
\end{equation}

Let us prove the commutativity of \eqref{FT cc from tilde}. After  unraveling definitions, the non-obvious part is to show the commutativity of the following diagram of functors $D(\wh E_{1})\to D(\wh D^{\sh})$
\begin{equation}
\xymatrix{(\wh f^{\sh})^{*}(\wh{\wt p'_{0}})^{!}(h^{\wh E}_{1})^{!} \ar@{=}[d]\ar[r]^{\tu(h^{\wh E}_{1})^{!}} & (\wh{\wt q'_{0}})^{!}(\wh{\wt f_{1}})^{*}(h^{\wh E}_{1})^{!} \ar[r]^{(\wh{\wt q'_{0}})^{!}\di}& (\wh{\wt q'_{0}})^{!}(h^{\wh F}_{1})^{!}(\wh f_{1})^{*}\ar@{=}[d]\\
(\wh f^{\sh})^{*} (\wh{p'_{0}})^{!}\ar[rr]^{\tu} & & (\wh {q'_{0}})^{!}(\wh f_{1})^{*} }
\end{equation}
This follows from Proposition \ref{p: comp pull horizontal} applied to the two squares
\begin{equation}
\xymatrix{ \wh D^{\sh}\ar[r]^{\wh{\wt q'_{0}}} \ar[d]^{\wh f^{\sh}}& \wh{\wt F_{1}}\ar[r]^{h^{\wh F}_{1}}\ar[d]^{\wh{\wt f_{1}}} & \wh F_{1}\ar[d]^{\wh f_{1}}\\
\wh C^{\sh} \ar[r]^{\wh{\wt p'_{0}}}& \wh{\wt E_{1}}\ar[r]^{h^{\wh E}_{1}} & \wh E_{1}
}
\end{equation}

(2) The proof is similar. Alternatively, it follows from (1) using the near-involutivity of $\FT$. 

\end{proof}

\section{Arithmetic Fourier transform}\label{sec: arithmetic FT} In this section we introduce an ``arithmetic'' variant of the Fourier transform, which will be used to do Fourier analysis on the Borel-Moore homology of moduli spaces for shtuka-type objects.
\[
\begin{tikzcd}
\Sht_V^r \ar[r] \ar[d] & \Hk^\flat_V \ar[d, "{(b_0, b_r)}"] \\
V \ar[r, "{(\Id, \Frob)}"] & V \times V
\end{tikzcd}
\]
When specialized to $r=0$, the arithmetic Fourier transform recovers the finite Fourier transforms used in \S \ref{ssec: r=0 transverse lagrangians} to prove modularity for $r=0$.

\subsection{The general setup}\label{sssec: AFT setup}
Let $T$ be a derived Artin stack and $Y \rightarrow T$ be a finite rank $\F_q$-vector space over $T$, meaning that $Y$ is a finite-dimensional $\F_q$-vector space object in stacks over $T$. In particular, $Y \rightarrow T$ is representable (in schemes) and finite \'{e}tale. Let $\wh{Y} \rightarrow T$ be the dual $\F_q$-vector space, i.e., at the level of \'{e}tale sheaves over $T$ we have 
\[
\wh{Y} = \cRHom_T(Y, \ul{\F_q}	). 
\]
Note that $\wh{\wh{Y}} \cong Y$.

Then we have an ``evaluation'' map 
\[
\begin{tikzcd}
\ev \co Y \times_T \wh{Y} \rightarrow \ul{\F_q},
\end{tikzcd}
\]
where $\ul{\F_q}$ is the set $\F_q$ viewed as a discrete scheme. Consider the diagram 
\[
\begin{tikzcd}
& Y \times_T \wh{Y}   \ar[dl, "\pr_0"'] \ar[dr, "\pr_1"]  \ar[r, "\ev"] & \ul{\F_q}\\
Y \ar[dr, "\pi"']  & & \wh{Y} \ar[dl, "\wh{\pi}"] \\
& T
\end{tikzcd}
\]

\begin{defn}\label{def: arithmetic FT}
Let $\psi$ be a nontrivial additive character of $\ul{\F}_q$. Let $d$ be the rank of $Y$ as an $\F_q$-vector space over $T$. We define the \emph{arithmetic Fourier transform (with respect to $\psi$)} to be the map 
\[
\FT^{\arith, \psi}_Y \co \mBM_*(Y) \rightarrow \mBM_{*} (\wh{Y})
\]
given by
\[
\alpha \mapsto   (-1)^d \pr_{1!}(\pr_0^*(\alpha) \cdot \ev^* \psi).
\]
Here, we used that $\mBM_*(-)$ is a module over $\rH^0(-; \Ql)$, or in other words, we can multiply Borel-Moore homology classes by locally constant functions.  More generally, $\mBM_*(-)$ is a module over $\rH^{2*}(-; \Q_\ell(*))$.

 Similarly, we define the \emph{arithmetic Fourier transform (on cohomology)} to be the map
\[
\FT^{\arith, \psi} \co \rH^*(Y; \Qll{Y}) \rightarrow \rH^* (\wh{Y}; \Qll{\wh Y})
\]
given by 
\[
\alpha \mapsto   (-1)^d\pr_{1!}(\pr_0^*(\alpha) \cdot \ev^* \psi).
\]
When $\psi$ is understood, we will suppress it from the notation, writing $\FT^{\arith} = \FT^{\arith, \psi}$. 
\end{defn}

\subsection{Basic properties} We establish some basic properties of the arithmetic Fourier transform, parallel to those of the usual finite Fourier transform (\S \ref{sssec: finite FT}).

\begin{lemma}[Plancherel property]\label{lem: Plancherel for cycles}
Let $\alpha_1 \in \mBM_{2i}(Y)$ and
\[
\beta_2 \in \rH^{2j}(\wh{Y}; \Qll{\wh{Y}}(j)).
\]
Then
\[
\pi_!(\alpha_1 \cdot \FT^{\arith}(\beta_2)) = \wh{\pi}_!(\FT^{\arith}(\alpha_1) \cdot \beta_2) \in \mBM_{2i-2j}(T).
\]
\end{lemma}

\begin{proof}
We have 
\begin{align*}
\pi_!(\alpha_1 \cdot\FT^{\arith}(\beta_2))  &=(-1)^d \pi_! (\alpha_1 \cdot \pr_{0!} (\pr_1^* \beta_2 \cdot \ev^* \psi)) \\
&=(-1)^d \pi_! \pr_{0!} (\pr_0^* \alpha_1 \cdot \pr_1^* \beta_2 \cdot \ev^* \psi) \\
&=(-1)^d \wh{\pi}_! \pr_{1!} (\pr_0^* \alpha_1 \cdot \pr_1^* \beta_2 \cdot \ev^* \psi) \\
&= (-1)^d\wh{\pi}_! ( \pr_{1!} (\pr_0^* \alpha_1  \cdot \ev^* \psi) \cdot \beta_2 ) \\
&= \wh{\pi}_! (\FT^{\arith}(\alpha_1) \cdot \beta_2).
\end{align*}
\end{proof}

Next we examine compatibility of the arithmetic Fourier transform with base change. Let $\varphi \co T' \rightarrow T$ be a proper map. Let $Y' \rightarrow T'$ be the pullback of $Y \rightarrow T$ and $\wh{Y}' \rightarrow T'$ be the pullback of $\wh{Y} \rightarrow T$. So we have Cartesian squares 
\[
\begin{tikzcd}
Y' \ar[r, "\varphi"] \ar[d] & Y \ar[d] \\
T' \ar[r, "\varphi"] & T
\end{tikzcd}
\hspace{1cm}
\begin{tikzcd}
\wh{Y}' \ar[r] \ar[d] & \wh{Y} \ar[d] \\
T' \ar[r, "\varphi"] & T
\end{tikzcd}
\]
The properness of $\varphi$ ensures the existence of maps $\varphi_! \co \mBM_*(Y') \rightarrow \mBM_{*}(Y)$ and $\varphi_! \co \mBM_*(\wh{Y}') \rightarrow \mBM_{*}(\wh{Y})$.

\begin{lemma}\label{lem: arithmetic FT proper pushforward} Let $\varphi \co T' \rightarrow T$ be proper and maintain the above notation. Then we have $\FT^{\arith} \circ \varphi_! = \varphi_!  \circ \FT^{\arith}$ as maps $\mBM_{*}(Y') \rightarrow  \mBM_{*}(\wh{Y}) $. 
\end{lemma}

\begin{proof}
The pushforward $\varphi_!$ satisfies base change against: smooth pullback, proper pushforward, and tensoring with $\rH^0(-, \Q_\ell)$, and $\FT^{\arith}$ is a composition of such operations.
\end{proof}

Now suppose that $\varphi \co T' \rightarrow T$ is quasi-smooth. Let $Y' \rightarrow T'$ be the pullback of $Y \rightarrow T$ and $\wh{Y}' \rightarrow T'$ be the pullback of $\wh{Y} \rightarrow T$. The quasi-smoothness of $\varphi$ ensures the existence of maps $\varphi^* \co \mBM_*(Y) \rightarrow \mBM_{*+2d(\varphi)}(Y')$ and $\varphi^* \co \mBM_*(\wh{Y}) \rightarrow \mBM_{*+2d(\varphi)}(\wh{Y}')$.

\begin{lemma}\label{lem: arithmetic FT qs pullback} Let $\varphi \co T' \rightarrow T$ be quasi-smooth and maintain the above notation. Then we have $\FT^{\arith} \circ \varphi^* = \varphi^*  \circ \FT^{\arith}$ as maps $\mBM_{*}(Y) \rightarrow  \mBM_{*+2d(\varphi)}(\wh{Y}') $. 
\end{lemma}

\begin{proof}
The pullback $\varphi^*$ satisfies base change against: smooth pullback, proper pushforward, and tensoring with $\rH^*(-, \Q_\ell)$, and $\FT^{\arith}$ is a composition of such operations.
\end{proof}

\begin{lemma}[Involutivity]\label{lem: involutivity}
We have $\FT^{\arith}_{\wh{Y}} \circ  \FT^{\arith}_{Y} = q^d [-1]^*$, where $[-1]$ is multiplication by $-1$ on $Y$ using its $\F_q$-vector space structure over $T$. 
\end{lemma}

\begin{proof}
First suppose that $Y \rightarrow T$ is a \emph{split} \'{e}tale $\F_q$-vector space over $T$, i.e., there exists a (finite-dimensional) $\F_q$-vector space $Y_0$ such that $Y = Y_0 \times T$. In this case, $\mBM_*(Y)  = \mBM_*(T) \otimes_{\Ql} \Ql[Y_0]$ and the arithmetic Fourier transform simplifies to the identity on $\mBM_*(T)$ tensored with the usual Fourier transform on $\Ql[Y_0]$. Therefore, the identity follows as for the usual finite Fourier transform (\S \ref{sssec: finite FT}). 

Now, since $Y \rightarrow T$ is finite \'{e}tale, it is split by some finite \'{e}tale pullback $\varphi \co T' \rightarrow T$. Letting $Y' \rightarrow T'$ be the pullback of $Y \rightarrow T$ along $\varphi$, the previous paragraph shows that 
\begin{equation}\label{eq: invol after pullback}
\FT^{\arith}_{\wh{Y}'} \circ  \FT^{\arith}_{Y'} = q^d [-1]^*.
\end{equation}
According to Lemma \ref{lem: arithmetic FT proper pushforward} and Lemma \ref{lem: arithmetic FT qs pullback}, $\FT^{\arith}_{\wh{Y}'}  \circ \varphi^* = \varphi^* \circ \FT^{\arith}_{\wh{Y}}$ and $\FT^{\arith}_{Y'} \circ \varphi^*  = \varphi^* \circ \FT^{\arith}_{Y}$, so composing \eqref{eq: invol after pullback} with $\varphi^*$ shows that 
\[
\varphi^*\FT^{\arith}_{\wh{Y}} \circ  \FT^{\arith}_{Y} = \varphi^* q^d [-1]^*.
\]
Now apply $\varphi_!$. Since $\varphi_! \circ \varphi^*$ is multiplication by $\deg \varphi$ (which is invertible), we get the desired equation. 
\end{proof}

Let $\varphi \co Y' \rightarrow Y$ be an $\F_q$-linear map of $\F_q$-vector spaces over $T$. In particular, $\varphi$ is finite \'{e}tale, so that we have maps 
\begin{align*}
\varphi^* \co \mBM_*(Y) & \rightarrow \mBM_{*}(Y'), \\
\varphi_! \co \mBM_*(Y') & \rightarrow \mBM_*(Y).
\end{align*}
Let $\wh{\varphi} \co \wh{Y} \rightarrow \wh{Y}'$ be the dual map. As $\wh{\varphi}$ is also finite \'{e}tale, we have maps 
\begin{align*}
\wh{\varphi}^* \co \mBM_*(\wh{Y}') & \rightarrow \mBM_*(\wh{Y}), \\
\wh{\varphi}_! \co \mBM_*(\wh{Y}) & \rightarrow \mBM_{*}(\wh{Y}').   
\end{align*}

\begin{prop}[Functoriality for linear maps]\label{prop: arithmetic FT functoriality} Keep the above notation. Let $d, d'$ be the ranks of $Y,Y'$ as $\F_q$-vector spaces over $T$. Then we have the following identities. 
\begin{enumerate}
\item $\FT^{\arith}_Y \circ \varphi_! =  (-1)^{d-d'} \wh{\varphi}^* \circ \FT^{\arith}_{Y'}$ as maps $\mBM_*(Y') \rightarrow \mBM_*(\wh{Y})$. 
\item $\FT^{\arith}_{Y'} \circ \varphi^* = (-1)^{d'-d} q^{d'-d} \wh{\varphi}_! \circ \FT^{\arith}_Y$ as maps $\mBM_*(Y) \rightarrow \mBM_*(\wh{Y}')$. 
\end{enumerate}
\end{prop}

\begin{proof}First suppose that the entire linear map $\varphi:Y'\to Y$ is split over $T$: there are finite-dimensional $\F_q$-vector spaces $Y'_0,Y_0$ and an $\F_q$-linear map $\varphi_0:Y'_0\to Y_0$ such that $Y'=Y'_0\times T$, $Y=Y_0\times T$, and $\varphi=\varphi_0\times\Id_T$. Then $\mBM_*(Y)=\mBM_*(T)\otimes_{\Ql}\Ql[Y_0]$ and $\mBM_*(Y')=\mBM_*(T)\otimes_{\Ql}\Ql[Y'_0]$, the arithmetic Fourier transforms are the identity on $\mBM_*(T)$ tensored with the finite Fourier transforms, and $\varphi_!,\varphi^*,\wh\varphi_!$, and $\wh\varphi^*$ are the identity on $\mBM_*(T)$ tensored with the usual finite pushforward and pullback attached to $\varphi_0$ and its dual. The two asserted identities therefore follow from \eqref{FT pushpull} and \eqref{FT pullpush}.

In the general case, after working componentwise if necessary, choose a finite \'{e}tale surjective map $s:T_0\to T$ which splits the finite \'{e}tale sheaves $Y$, $Y'$ and the morphism $\varphi:Y'\to Y$ in the sense of the previous paragraph. Then the proof reduces to that in the split case, as in the proof of Lemma \ref{lem: involutivity}.
\end{proof}

\subsection{Compatibility with sheaf-theoretic Fourier transform}\label{sssec: AFT and FT}

Let $S$ be a derived stack and $p: E\to S$ be a vector bundle (note that we deliberately do not allow more general derived vector bundles here). Suppose $c=(c_{0},c_{1}): C\to S\times S$ is a correspondence and we are given an isomorphism of vector bundles over $C$
\begin{equation}
\io: c_{0}^{*}E\simeq c_{1}^{*}E.
\end{equation}
Let $C_{E}$ be the total space of $c_{0}^{*}E$ and of $c_{1}^{*}E$, identified via $\io$. Let $e_{i}: C_{E}\simeq c_{i}^{*}E\to E$ be the projection, for $i=0,1$. Then we get a correspondence $e=(e_{0},e_{1}): C_{E}\to E\times E$ that fits into a commutative diagram
\begin{equation}\label{corr vb}
\xymatrix{E\ar[d]^{p} &  C_{E}\ar[d]^{p_{C}}\ar[l]_-{e_{0}}\ar[r]^{e_{1}} & E\ar[d]^{p}\\
S & C\ar[l]_-{c_{0}}\ar[r]^{c_{1}} & S }
\end{equation}
such that both squares are Cartesian. 

The above data induces a correspondence $\wh e: C_{\wh E}\to \wh E\times \wh E$ by passing to the dual vector bundles. Let $\wh p:\wh E\to S$ and $p_{\wh E}: C_{\wh E}\to \wh E$ be the projections. 

Let $\cK\in D^{b}_{c}(E)$ and $\frc: e_{0}^{*}\cK\to e_{1}^{!}\cK{\tw{-n}}$ be a cohomological correspondence.  Applying Fourier transforms, using that Fourier transform commutes with arbitrary $*$ and $!$ base change (\S \ref{sssec: FT base change}), $\frc$ induces a cohomological correspondence of $\FT_{E}(\cK)$ as the composition:
\begin{equation}
\FT_{C_{E}}(\frc): \wh e_{0}^{*}\FT_{E}(\cK)\simeq \FT_{C_{E}}(e_{0}^{*}\cK)\xr{\FT_{C_{E}}(\frc)}
\FT_{C_{E}}(e_{1}^{!}\cK\tw{-n})\simeq \wh e_{1}^{!}\FT_{E}(\cK){\tw{-n}}.
\end{equation}

Consider the map $c^{(1)}=(\Frob_{S}\c c_{0}, c_{1}): C \rightarrow S\times S$. This makes $C$ into a self-correspondence of $S$ in a different way, which we denote by $C^{(1)}$. Similarly, we define $C^{(1)}_{E}$ (a self-correspondence of $E$) and $C^{(1)}_{\wh E}$ (a self-correspondence of $\wh E$). Recall notation
\begin{equation}
\Sht(C)=\Fix(C^{(1)}), \quad 
\Sht(C_{E})=\Fix(C^{(1)}_{E}), \quad \Sht(C_{\wh E})=\Fix(C^{(1)}_{\wh E}).
\end{equation}

\begin{lemma}
The projections 
\begin{equation}
\pi: \Sht(C_{E})\to \Sht(C), \quad  \wh \pi: \Sht(C_{\wh E})\to \Sht(C)
\end{equation}
are relative $\F_{q}$-vector spaces over $\Sht(C)$ that are dual to each other. 
\end{lemma}

\begin{proof}
Evident. 
\end{proof}

Let $\cK\in D(E)$. Then $\cK$ is equipped with a canonical Weil structure $\Frob_{E}^{*}\cK\simeq \cK$. A cohomological correspondence $\frc: e_{0}^{*}\cK\to e_{1}^{!}\cK{\tw{-n}}$ induces a cohomological correspondence $\frc^{(1)}$ supported on $C_{E}^{(1)}$:
\begin{equation}
\frc^{(1)}: e_{0}^{*}\Frob_{E}^{*}\cK\simeq e_{0}^{*}\cK\xr{\frc}e_{1}^{!}\cK\tw{-n}.
\end{equation}
Taking trace we get
\begin{equation}
\Tr^{\Sht}(\frc):=\Tr(\frc^{(1)})\in \hBM{2n}{\Sht(C_E)}.
\end{equation}
Similarly, we have the cohomological correspondence $\FT_{C_{E}}(\frc)^{(1)}$ of $\FT_{E}(\cK)$ supported on $C_{\wh E}^{(1)}$, and its trace
\begin{equation}
\Tr^{\Sht}(\FT_{C_{E}}(\frc)):=\Tr(\FT_{C_{E}}(\frc)^{(1)})\in \hBM{2n}{\Sht(C_{\wh E})}.
\end{equation}

\begin{thm}\label{thm: trace compatible with FT} In the above situation, we have
\begin{equation}
\Tr^{\Sht}(\FT_{C_{E}}(\frc))=\FT^{\arith}_{\Sht(C_{E})}(\Tr^{\Sht}(\frc))\in \hBM{2n}{\Sht(C_{\wh E})}.
\end{equation}
\end{thm}

\begin{proof}It is easy to see that 
\begin{equation}
\FT_{C_{E}}(\frc^{(1)})=\FT_{C_{E}}(\frc)^{(1)}.
\end{equation}
Therefore we need to show
\begin{equation}
\Tr(\FT_{C_{E}}(\frc^{(1)}))=\FT^{\arith}_{\Sht(C_{E})}(\Tr(\frc^{(1)})).
\end{equation}
Consider the following diagram of correspondences (correspondences are written horizontally and morphisms between correspondences are written vertically)
\begin{equation}
\xymatrix{ E & C^{(1)}_{E}\ar[l]_-{\Frob\c e_{0}}\ar[r]^-{e_{1}} & E\\
E\times_{S}\wh E\ar[u]^{\pr_{1}}\ar[d]^{\pr_{2}} & C^{(1)}_{E\times_{S}\wh E}\ar[u]^{\pr_{1}}\ar[d]^{\pr_{2}}\ar[l]_-{\Frob\c d_{0}}\ar[r]^-{d_{1}} & E\times_{S}\wh E\ar[u]^{\pr_{1}}\ar[d]^{\pr_{2}}\\
\wh E & C^{(1)}_{\wh E}\ar[l]_-{\Frob\c\wh e_{0}}\ar[r]^-{\wh e_{1}} & \wh E
}
\end{equation} 
Note that all squares are \'{e}tale topologically Cartesian. The normalized Fourier transform $\FT(\frc^{(1)})$ is obtained by applying the following unshifted correspondence operation and then the normalizing shift by $\rank(E)$:
\begin{equation}
\FT(\frc^{(1)})=\TT_{[\rank(E)]}\pr_{2!}(t_{\psi}(\pr_{1}^{*}(\frc^{(1)}))).
\end{equation}
After applying $\Tr$, this shift contributes the sign $(-1)^{\rank(E)}$, which is exactly the sign in the arithmetic Fourier transform of Definition \ref{def: arithmetic FT}. Thus the unshifted calculation below proves the stated normalized trace identity.
Here,  for a cohomological correspondence $\frd: d_{0}^{*}\cF\to d_{1}^{!}\cF{\tw{-n}}$,  $t_{\psi}(\frd)$ means the cohomological correspondence 
\begin{equation}
t_{\psi}(\frd): d_{0}^{*}(\cF\ot \ev^{*}\cL_{\psi})\simeq d_{0}^{*}\cF\ot \wt\ev^{*}\cL_{\psi}\xr{\frd\ot\Id}d_{1}^{!}\cF{\tw{-n}}\ot \wt\ev^{*}\cL_{\psi}\simeq d_{1}^{!}(\cF\ot \ev^{*}\cL_{\psi}){\tw{-n}}
\end{equation}
where $\ev: E\times_{S}\wh E\to \G_{a}$ and $\wt\ev: C_{E\times_{S}\wh E}=C_{E}\times_{C}C_{\wh E}\to \G_{a}$ are the tautological evaluation pairings. 

On the other hand, let 
\[
\pr^{\Sht}_{1}: \Sht(C_{E\times_{S}\wh E})=\Sht(C_{E})\times_{\Sht(C)}\Sht(C_{\wh E})\to \Sht(C_{E})
\]
\[
\pr^{\Sht}_{2}: \Sht(C_{E\times_{S}\wh E})=\Sht(C_{E})\times_{\Sht(C)}\Sht(C_{\wh E})\to \Sht(C_{\wh E})
\]
be the projections to the two factors (the maps are finite \'etale). Then the unshifted part of $\FT^{\arith}$ is the composition $\pr_{2!}^{\Sht}\c m_{\psi}\c\pr^{\Sht *}_{1}$, and the prefactor $(-1)^{\rank(E)}$ is the sign discussed above, where $m_{\psi}$ means multiplying the Borel-Moore classes on $\Sht(C_{E})\times_{\Sht(C)}\Sht(C_{\wh E})$ by the function $\ev^{\Sht*}\psi$, where $\ev^{\Sht}: \Sht(C_{E})\times_{\Sht(C)}\Sht(C_{\wh E})\to \F_{q}$ is the evaluation pairing.

By Proposition \ref{prop: trace commutes with smooth pullback}, $\Tr(\pr_{1}^{*}(\frc^{(1)}))=\pr_{1}^{\Sht *} \Tr(\frc^{(1)})$. It is also clear that taking trace intertwines $t_{\psi}$ and $m_{\psi}$. Hence
\begin{equation}
\Tr(t_{\psi}\pr_{1}^{*}(\frc^{(1)}))=m_{\psi}\pr_{1}^{\Sht*}(\Tr(\frc^{(1)})).
\end{equation}
Let $\cF=\pr_{1}^{*}\cK\ot\ev^{*}\cL_{\psi}\in D(E\times_{S}\wh E)$, and $\frd=t_{\psi}\pr_{1}^{*}(\frc): d_{0}^{*}\cF\to d_{1}^{!}\cF{\tw{-n}}$ be the cohomological correspondence of $\cF$ supported on $C_{E\times_{S}\wh E}$. It remains to show
\begin{equation}
\Tr(\pr_{2!}(\frd^{(1)}))=\pr^{\Sht}_{2!}\Tr(\frd^{(1)}).
\end{equation}
This follows by applying Lemma \ref{l:push corr along vector bundle} below to the map of correspondences $\pr_{2}: C_{E\times_{S}\wh E}\to C_{\wh E}$ (viewing $E\times_{S}\wh E$ as a vector bundle over $\wh E$).
\end{proof}

\begin{lemma}\label{l:push corr along vector bundle}
Consider the situation \eqref{corr vb}. Let $\cK\in D(E)$ and $\frc$ be a cohomological correspondence $\frc: e_{0}^{*}\cK\to e_{1}^{!}\cK\tw{-n}$. Let $p^{\Sht}: \Sht(C_{E})\to \Sht(C)$ be the induced map on fixed points of $C_{E}^{(1)}$ and $C^{(1)}$, which is a relative $\F_{q}$-vector space (in particular a finite morphism). Then
\begin{equation}
\Tr(p_{C!}(\frc^{(1)}))=p^{\Sht}_{!}\Tr(\frc^{(1)})\in \hBM{2n}{\Sht(C)}.
\end{equation}
\end{lemma}
\begin{proof}
Consider the projective bundle $\ov E=\PP(E\op \cO)\to S$ that contains $E$ as an open substack. Similarly let $C_{\ov E}=\PP(c_{0}^{*}E\op \cO)\cong \PP(c_{1}^{*}E\op \cO)\to C$  be the pullback projective bundle over $C$, via either $c_{0}$ or $c_{1}$. Then $C_{\ov E}$ is a self-correspondence of $\ov E$ with a proper map to $C$:
\begin{equation}
\xymatrix{\ov E\ar[d]^{\ov p} &  C_{\ov E}\ar[d]^{\ov p_{C}}\ar[l]_-{\ov e_{0}}\ar[r]^{\ov e_{1}} & \ov E\ar[d]^{\ov p}\\
S & C\ar[l]_-{c_{0}}\ar[r]^{c_{1}} & S
}
\end{equation}
Let $E_{\infty}=\ov E-E$ be the divisor at infinity, which is isomorphic to $\PP(E)$. Similarly define $C_{E_{\infty}}=C_{\ov E}-C_{E}$, which is a self-correspondence of $E_{\infty}$. We define $C^{(1)}_{\ov E}$ and $C^{(1)}_{E_{\infty}}$ by composing the first projections by the Frobenius. 
For the fixed points, we have an open and closed decomposition
\begin{equation}
\Sht(C_{\ov E})=\Sht(C_{E})\coprod \Sht(C_{E_{\infty}}).
\end{equation}
Fiberwise over $\Sht(C)$, this decomposition takes the form $\PP^{N}(\F_{q})\simeq \F_{q}^{N}\coprod \PP^{N-1}(\F_{q})$, where $N$ is the rank of $E$.

Let $j:E\incl \ov E$ and $j_{C}: C_{E}\incl C_{\ov E}$ be the open inclusions. The map of correspondences 
\[
\begin{tikzcd}
E \ar[d, hook, "j"] & C_E \ar[l, "e_0"'] \ar[r, "e_1"] \ar[d, hook, "j_C"] & E \ar[d, hook, "j"] \\
\ol{E} & C_{\ol{E}}  \ar[l, "\ol{e}_0"'] \ar[r, "\ol{e}_1"] & \ol{E}
\end{tikzcd}
\]
has both squares Cartesian, so it is left pushable. Therefore we have the pushforward cohomological correspondence
\begin{equation}
\ov\frc:=j_{C!}\frc \in \Corr_{C_{\ov E}}(j_! \cK, j_! \cK\tw{-n}).
\end{equation}
Here $j_{C!}$ is the left pushforward of cohomological correspondences along the displayed open immersion of correspondences. Since this open immersion and the Frobenius twist are obtained by base change from the same maps, $j_{C!}(\frc^{(1)})=(j_{C!}\frc)^{(1)}=\ov\frc^{(1)}$, and $p_{C!}(\frc^{(1)})=\ov p_{C!}j_{C!}(\frc^{(1)})$.
Since $\ov p$ and $\ov p_{C}$ are proper, we have by Proposition \ref{prop: trace commutes with proper push},
\begin{equation}
\Tr(p_{C!}(\frc^{(1)}))=\Tr(\ov p_{C!}j_{C!}(\frc^{(1)}))=\Tr(\ov p_{C!}(\ov\frc^{(1)}))=\ov p^{\Sht}_{!}\Tr(\ov\frc^{(1)})
\end{equation}
where $\ov p^{\Sht}: \Sht(C_{\ov E})\to \Sht(C)$ is the obvious map. On the other hand, since $j_{C}$ is an open embedding, we have by Proposition \ref{prop: trace commutes with smooth pullback} that $\Tr(\ov\frc^{(1)})|_{\Sht(C_{E})}=\Tr(\frc^{(1)})\in \hBM{2n}{\Sht(C_{E})}$, therefore
\begin{equation}
p^{\Sht}_{!}(\Tr(\frc^{(1)}))=p^{\Sht}_{!}(\Tr(\ov\frc^{(1)})|_{\Sht(C_{E})}).
\end{equation}
It remains to show that
\begin{equation}
\ov p^{\Sht}_{!}\Tr(\ov\frc^{(1)})=p^{\Sht}_{!}(\Tr(\ov\frc^{(1)})|_{\Sht(C_{E})}),
\end{equation}
which would follow from the vanishing
\begin{equation}\label{tr van along infty}
\Tr(\ov\frc^{(1)})|_{\Sht(C_{E_{\infty}})}=0\in \hBM{2n}{\Sht(C_{E_{\infty}})}.
\end{equation}

The divisor $E_{\infty}$ is invariant under $C_{\ov E}$, since $\ov e_{1}^{-1}(E_{\infty})=C_{E_{\infty}}$ and this closed substack maps to $E_{\infty}$ by $\ov e_{0}$. Hence $E_{\infty}$ is also invariant under the Frobenius-twisted correspondence $C_{\ov E}^{(1)}$. We claim that $C_{\ov E}^{(1)}$ is contracting near $E_{\infty}$ in the sense of \cite[Definition 2.1.1(b)]{Var07}. Let $\cI_{\infty}$ be the defining ideal of $E_{\infty}$ in $\ov E$. Because $C_{E_{\infty}}=\ov e_{0}^{-1}(E_{\infty})=\ov e_{1}^{-1}(E_{\infty})$, both $\ov e_{0}^{*}\cI_{\infty}$ and $\ov e_{1}^{*}\cI_{\infty}$ are the defining ideal of $C_{E_{\infty}}$ inside $C_{\ov E}$. For the twisted left leg $\Frob\circ\ov e_{0}$, therefore,
\[
(\Frob\circ\ov e_{0})^{*}\cI_{\infty}=\Frob^{*}(\ov e_{0}^{*}\cI_{\infty})=(\ov e_{0}^{*}\cI_{\infty})^{q}=(\ov e_{1}^{*}\cI_{\infty})^{q}\subset (\ov e_{1}^{*}\cI_{\infty})^{2},
\]
which is precisely Varshavsky's contracting ideal condition.

Now we can apply the $\ell$-adic version of \cite[Theorem 7.5.1]{FK} to the correspondence $C^{(1)}_{\ov E}$, the closed substack $Z=E_{\infty}\subset \ov E$, and each connected component $\b$ of $\Sht(C_{E_{\infty}})$ as an open-closed substack of the fixed point locus of $C^{(1)}_{\ov E}$. Since $C^{(1)}_{\ov E}$ is contracting near $E_{\infty}$, we conclude that the $\b$-part of the trace satisfies
\begin{equation}
\Tr(\ov\frc^{(1)})|_{\b}=\Tr(\ov\frc^{(1)}|_{E_{\infty}})|_{\b}\in \hBM{2n}{\b}.
\end{equation}
Here $\ov\frc^{(1)}|_{E_{\infty}}$ is the restriction of the cohomological correspondence $\ov\frc^{(1)}$ as defined in \cite[1.5.6(a)]{Var07}, which is a cohomological correspondence of $(j_{!}\cK)|_{E_{\infty}}=0$. Therefore the right side above is zero, hence $\Tr_{\b}(\ov\frc^{(1)})=0$ for any connected component of $\Sht(C_{E_{\infty}})$. This proves \eqref{tr van along infty} and finishes the proof of the lemma.
\end{proof}

\part{Modularity}
\section{Modularity for cohomological correspondences}\label{ssec: modularity for cc}
In this section we carry out Steps (2), (3), and (5) of the proof outline from \S \ref{ssec: outline}. 

First, in \S \ref{ssec: UVW} we set up the derived vector bundles $U,V$, and $W$ and their Hecke correspondences; in the terminology of the outline in \S \ref{ssec: outline}, $U$ belongs to Step (1), $V$ to Step (2), and $W$ to Step (3). Switching the role of $\cE_{1}$ and $\cE_{2}$, in \S \ref{ssec: UVW dual} we define the corresponding dual-side derived vector bundles $U^\perp,\wh V,W^\perp$, belonging respectively to the same three steps. The main result of this section is Theorem \ref{thm: FT of pushforward correspondence}, showing that the Fourier transform on $\Hk_V^\flat$ of the cohomological correspondence $f_! \cc_U$ agrees (up to shift and twist) with the parallel cohomological construction obtained by interchanging $\cE_1$ and $\cE_2$. The proof of Theorem \ref{thm: FT of pushforward correspondence} uses the general results proved in \S\ref{sec: descent} and \S\ref{sec: FT cc}. We regard Theorem \ref{thm: FT of pushforward correspondence} as an incarnation of modularity at the level of cohomological correspondences -- the modularity of the higher theta series will be extracted from it by taking a trace in the sense of the sheaf-cycle correspondence -- which explains the title of the section.

\subsection{The stacks $U, V, W$ and their Hecke correspondences}\label{ssec: UVW} 
We begin by defining various spaces of interest. 

Suppose we are given a short exact sequence of coherent sheaves on $X'$:
\begin{equation}\label{eq: descent ses 1}
0 \rightarrow  \s^{*}\cE_2 \rightarrow \cE_1^* \rightarrow \wt{Q}_1 \rightarrow 0
\end{equation}
where $\cE_1, \cE_2$ are vector bundles and $\wt{Q}_1$ is torsion.

\subsubsection{Definition of $U, V$ and $W$}

Below we denote $S$ for a Harder-Narasimhan truncation $\Bun_{U(n)}^{\leq \mu}$. (The reason for this truncation is to guarantee global presentability; see \S \ref{ssec: safety verification}.) We define several derived vector bundles over $S$. Let $\cF_{\univ}$ be the universal Hermitian bundle over $X'\times S$. Let $R$ be any animated $\F_{q}$-algebra.
\begin{itemize}

\item Define $\cU:=\ul{\RHom(\cF^*_{\univ}, \cE_1^*)}$ to be the perfect complex on $S$ whose pullback to an $R$-point $\cF \in S(R)$ is naturally in $R$ isomorphic to  $\RHom_{X_R'}(\cF^*, \cE_1^*\ot R)$ regarded in $\Mod_R$. Let $U := \Tot_S(\cU)$ be the associated derived vector bundle over $S$.

\item Define $\cV:=\ul{\RHom(\cF^*_{\univ},\wt{Q}_1)}$ to be the perfect complex on $S$ whose pullback to an $R$-point $\cF \in S(R)$ is naturally in $R$ isomorphic to $\RHom_{X_R'}(\cF^*, \wt{Q}_1\ot R)$ regarded in $\Mod_R$. Let $V := \Tot_S(\cV)$ be the associated derived vector bundle over $S$. Since $\Ext^1_{X_R'}(\cF^*,\wt{Q}_1\ot R) =0$ as $\wt{Q}_1$ is torsion, $\cV$ is in fact a locally free coherent sheaf on $S$, so that $V \rightarrow S$ is actually a vector bundle in the classical sense. This fact is not important for this section, although it plays a role in later sections.

\item Define $\cW:=\ul{\RHom(\cF^*_{\univ}, \sigma^* \cE_2[1])}$ to be the perfect complex on $S$ whose pullback to an $R$-point $\cF \in S(R)$ is naturally in $R$ isomorphic to $\RHom_{X_R'}(\cF^*, \sigma^* \cE_2[1]\ot R)$ regarded in $\Mod_R$. Let $W := \Tot_S(\cW)$ be the associated derived vector bundle over $S$. 
\end{itemize} 
From \eqref{eq: descent ses 1} we get an exact triangle of sheaves on $X'$,
\begin{equation}\label{eq: descent tri}
\cE_1^* \rightarrow \wt{Q}_1 \rightarrow \sigma^* \cE_2[1] 
\end{equation}
which induces the exact triangle in $\Perf(S)$,
\begin{equation}\label{eq: exact triangle on S E_1}
\cU\to \cV\to \cW.
\end{equation}
Forming total spaces, this is equivalent to the derived Cartesian square of derived vector bundles over $S$: 
\begin{equation}\label{eq: derived cartesian E_1 S}
\begin{tikzcd}
U \ar[r] \ar[d]  & V \ar[d] \\
S \ar[r, "z"] & W
\end{tikzcd}
\end{equation}
where $z: S\to W$ denotes the zero section.

\begin{remark}\label{rem: derived hitchin stack}
By \cite[\S 5.7]{FYZ2}, $U$ is isomorphic to the derived pullback of the derived Hitchin stack ``$\sM_{\GL(m)', U(n)} \rightarrow \Bun_{\GL(m)'} \times \Bun_{U(n)}$'' from \cite[\S 5]{FYZ2} 
along the map $\{ \cE_1 \} \times S   \rightarrow \Bun_{\GL(m)'} \times \Bun_{U(n)}$. 
\end{remark}

\subsubsection{Hecke stacks}
Recall the Hecke stack $\Hk_{U(n)}^{r}$ from \cite[\S 5.4]{FYZ2}. The $R$-points of $\Hk_{U(n)}^r$ are diagrams 
\begin{equation}\label{eq: flat hecke diagram}
\begin{tikzcd}
 & \cF_{1/2}^{\flat} \ar[dl] \ar[dr]   & & \ldots  & & \cF_{r-1/2}^{\flat} \ar[dl] \ar[dr] \\ 
\cF_0 \ar[rr, dashed]  &  & \cF_{1} \ar[r, dashed] & \ldots \ar[r, dashed] & \cF_{r-1} \ar[rr, dashed] & & \cF_r 
\end{tikzcd}   
\end{equation}
where each $\cF_i \in \Bun_{U(n)}(R)$, and each $\cF_{1/2}^{\flat}$ is a rank $n$ vector bundle on $X'_R$, satisfying some conditions (for example, the maps are injective). We shall abbreviate such diagrams as $(\cF_{\star}) \in \Hk_{U(n)}^r(R)$. There are maps 
\begin{equation}\label{eq: h_i untruncated}
h_i \co \Hk_{U(n)}^{r} \rightarrow \Bun_{U(n)}, \quad i =0, \ldots, r
\end{equation}
projecting to the datum of $\cF_i$, as well as $\pr\co \Hk_{U(n)}^{r} \rightarrow (X')^r$ projecting to the ``legs'' ($r$-tuple of points on $X'$ at which the dashed maps have poles). 

We define the open substack $\Hk_S^r \subset \Hk_{U(n)}^r$ as 
\[
\Hk_S^r := h_0^{-1}(S) \cap h_1^{-1} (S) \cap \ldots \cap h_r^{-1}(S).
\]
Therefore the maps \eqref{eq: h_i untruncated} restrict to give 
\[
h_i \co \Hk_{S}^r \rightarrow S, \quad i = 0, \ldots, r.
\]
We will at some points find it convenient to distinguish the different copies of $S$, so we will also sometimes use $S_i$ to denote a copy of $S$ and write $h_i  \co \Hk_{S}^r \rightarrow S_i$.

\subsubsection{Hecke stacks for $U, V$ and $W$}\label{sss: Hk UVW}  Given a diagram $(\cF_{\star}) \in \Hk_{S}^r(R)$, define $\cF_{\bu}^{\flat}$ to be the homotopy limit of the diagram\footnote{For $r=0$, this is the empty diagram and the homotopy limit is interpreted to be the $0$ object.}
\[
\begin{tikzcd}
& \cF_{1/2}^\flat \ar[dr] &  & \ar[dl] \ldots  \ar[dr] & & \cF_{r-1/2}^\flat \ar[dl]  \\
 & & \cF_1 & \ldots & \cF_{r-1} & & 
\end{tikzcd}
\]
This is a perfect complex on $X'_R$ (as the limit of a finite diagram of locally free sheaves). Explicitly, it admits a compact description as the complex (placed in degrees $0,1$)
\[
(\cF_{1/2}^{\flat} \oplus \ldots \oplus \cF_{r-1/2}^{\flat}) \rightarrow (\cF_1 \oplus \ldots \oplus \cF_{r-1}),
\]
where the map sends $(s_{1/2},\cdots, s_{r-1/2})$ to $(s_{1/2}-s_{3/2}, \cdots, s_{r-3/2}-s_{r-1/2})$. This description makes clear that $\cF_{\bu}^{\flat}$ has tor-amplitude in $(-\infty, 1]$. Note that even for classical $R$, $\cF_{\bu}^{\flat}$ may have non-zero cohomology in both degrees $0$ and $1$.

Define $\cF_{\bu}^{\flat *}$ to be the $\cO_{X'_{R}}$-linear dual of $\cF_{\bu}^{\flat}$, i.e., the cone of the dual morphism (in degrees $-1$ and $0$)
\[
 (\cF_1^* \oplus \ldots \oplus \cF_{r-1}^*)\rightarrow ((\cF_{1/2}^{\flat})^* \oplus \ldots \oplus (\cF_{r-1/2}^{\flat})^*).
\]
Note that for classical $R$, $\cF_{\bu}^{\flat *}$ is a coherent \emph{sheaf} on $X'_{R}$ concentrated in degree $0$ which may not be locally free.

We have a natural map of perfect complexes on $X'_{R}$
\begin{equation*}
\frp_{i}: \cF_{\bu}^{\flat}\to \cF_{i}, \quad i=0,1,\cdots, r
\end{equation*}
that is the composition of the projection to $\cF^{\flat}_{i-1/2}$ and the tautological map $\cF^{\flat}_{i-1/2} \rightarrow \cF_{i}$ when $i>0$, and the composition of the projection to $\cF^{\flat}_{i+1/2}$ and the tautological map $\cF^{\flat}_{i+1/2} \rightarrow \cF_{i}$ when $i<r$. Both constructions give the same map up to explicit chain homotopy when $0<i<r$. Dualizing $\frp_{i}$, we get a map of coherent sheaves on $X'_{R}$
\begin{equation}\label{eq: F bullet to F_i}
\frp^{*}_{i}: \cF_{i}^{*}\to \cF_{\bu}^{\flat*}.
\end{equation}

As $\cF_{\star}$ varies in $\Hk_{S}^{r}$, the construction $\cF_{\star}\mapsto \cF_{\bu}^{\flat *}$ gives a coherent sheaf $\cF_{\univ, \bu}^{\flat *}$ over $X'\times \Hk_{S}^{r}$. 


We now construct various spaces over $\Hk_S^r$. 
\begin{itemize} 
\item Define $\cU^{\flat}_{\Hk}:=\ul{\RHom(\cF_{\univ, \bu}^{\flat *}, \cE_1^*)}$ to be the perfect complex on $\Hk_S^r$ whose pullback to an $R$-point $(\cF_{\star}) \in \Hk_S^r(R)$ is naturally in $R$ isomorphic to $\RHom_{X_R'}(\cF^{\flat *}_{\bu}, \cE_1^*\ot R)$ regarded in $\Mod_R$. Let $\Hk_U^\flat := \Tot_{\Hk_S^r}(\cU^{\flat}_{\Hk})$ be the associated derived vector bundle over $\Hk_S^r$.
\item Define $\cV^{\flat}_{\Hk}:=\ul{\RHom(\cF^{\flat*}_{\univ, \bu},\wt{Q}_1)}$ to be the perfect complex on $\Hk_S^r$ whose pullback to an $R$-point $(\cF_{\star}) \in \Hk_S^r(R)$ is naturally in $R$ isomorphic to $\RHom_{X_R'}(\cF^{\flat *}_{\bu}, \wt{Q}_1\ot R)$ regarded in $\Mod_R$. Let $\Hk_V^\flat = \Tot_{\Hk_S^r}(\cV^{\flat}_{\Hk})$ be the associated derived vector bundle over $\Hk_S^r$. 
\item Define $\cW^{\flat}_{\Hk}:=\ul{\RHom(\cF^{\flat*}_{\univ, \bu}, \sigma^* \cE_2[1])}$ to be the perfect complex on $\Hk_S^r$ whose pullback to an $R$-point $(\cF_{\star}) \in \Hk_S^r(R)$ is naturally in $R$ isomorphic to $\RHom_{X_R'}(\cF^{\flat *}_{\bu}, \sigma^* \cE_2[1]\ot R)$ regarded in $\Mod_R$. Let $\Hk_W^\flat := \Tot_{\Hk_S^r}(\cW^{\flat}_{\Hk})$ be the associated derived vector bundle over $\Hk_S^r$. 
\end{itemize}
From \eqref{eq: descent tri}, we get an exact triangle of perfect complexes on $\Hk_S^r$,
\begin{equation}\label{eq: exact triangle on Hk_S E_1}
\cU^{\flat}_{\Hk}\to \cV^{\flat}_{\Hk}\to \cW^{\flat}_{\Hk}.
\end{equation}
At the level of total spaces, this induces maps $\Hk_U^\flat \xrightarrow{f} \Hk_V^\flat \xrightarrow{g} \Hk_W^\flat$ fitting into a derived Cartesian square of derived vector bundles over $\Hk_S^r$:
\begin{equation}\label{eq: Hk UVW Cart}
\begin{tikzcd}
\Hk_U^\flat \ar[r, "f"] \ar[d, "\pi"] & \Hk_V^\flat \ar[d, "g"] \\
\Hk_S^r \ar[r, "z"] & \Hk_W^\flat
\end{tikzcd}
\end{equation}

\begin{remark}\label{rem: derived hecke stack}
By \cite[\S 5.7]{FYZ2}, $\Hk_U^\flat$ is isomorphic to the derived pullback of the derived Hecke stack from \cite[\S 5]{FYZ2}, ``$\Hk_{\sM_{\GL(m)', U(n)}}^r \rightarrow \Bun_{\GL(m)'} \times \Hk_{U(n)}^r$'' 
along the map $\{\cE_1\} \times \Hk_S^r \rightarrow \Bun_{\GL(m)'} \times \Hk_{U(n)}^r$. 
\end{remark}

\subsubsection{Geometric properties} 
For each $i = 0, 1, \ldots, r$, let 
\begin{itemize}
\item $\wt{U}_i \rightarrow \Hk_S^r$ be the pullback of $U \rightarrow S$ along $h_i \co \Hk_S^r \rightarrow S$.
\item $\wt{V}_i \rightarrow \Hk_S^r$ be the pullback of $V \rightarrow S$ along $h_i \co \Hk_S^r \rightarrow S$.
\item $\wt{W}_i \rightarrow \Hk_S^r$ be the pullback of $W \rightarrow S$ along $h_i \co \Hk_S^r \rightarrow S$.
\end{itemize}
Let 
\begin{equation}\label{eq: hiUVW}
h_{i}^{U}: \wt U_{i}\to U, \quad h_{i}^{V}: \wt V_{i}\to V, \quad h_{i}^{W}: \wt W_{i}\to W,
\end{equation}
be the natural projections. When no ambiguity can result, we also write $U_i,V_i,W_i$ for copies of $U,V,W$ over the copy $S_i$ of $S$. 


The maps \eqref{eq: F bullet to F_i} induce natural maps of perfect complexes on $\Hk^{r}_{S}$
\begin{equation*}
\cU^{\flat}_{\Hk}\to h_{i}^{*}\cU, \quad \cV^{\flat}_{\Hk}\to h_{i}^{*}\cV,  \quad \cW^{\flat}_{\Hk}\to h_{i}^{*}\cW. 
\end{equation*}
At the level of total spaces, these induce maps of derived vector bundles over $\Hk^{r}_{S}$ for $i=0,\ldots, r$:
\begin{equation}\label{eq: tilde maps}
\wt a_{i}: \Hk^{\flat}_{U}\to \wt U_{i}, \quad \wt b_{i}: \Hk^{\flat}_{V}\to \wt V_{i},  \quad \wt c_{i}: \Hk^{\flat}_{W}\to \wt W_{i}. 
\end{equation}
Composing these maps with $h^{?}_{i}$, we get maps for $i=0,\cdots, r$
\begin{equation}\label{eq: no tilde maps}
a_{i}: \Hk^{\flat}_{U}\to U_{i}, \quad  b_{i}: \Hk^{\flat}_{V}\to V_{i},  \quad c_{i}: \Hk^{\flat}_{W}\to W_{i}. 
\end{equation}

\begin{lemma}\label{lem: hecke corr quasismooth tilde} We have the following properties of the morphisms in \eqref{eq: tilde maps}.
\begin{enumerate}
\item Each morphism $\wt{a}_i$ is a quasi-smooth closed embedding. 
\item Each morphism $\wt{b}_i$ is quasi-smooth and separated (in particular, representable in derived schemes). 
\item Each morphism $\wt{c}_i$ is a smooth vector bundle. 
\end{enumerate}
\end{lemma}

\begin{proof}
In all cases, the quasi-smoothness follows from the fact that the maps are induced by the dual of \eqref{eq: F bullet to F_i}. Indeed, for $(\cF_{\star})\in \Hk^{r}_{S}(R)$, write 
\begin{equation}\label{eq:T_i}
T_i := \Cofib(\cF_i^* \rightarrow \cF_{\bu}^{\flat *}).
\end{equation}
Note that along each fiber over $\Spec R$, the map $\cF_i^* \rightarrow \cF_{\bu}^{\flat *}$ is generically an isomorphism, hence $T_i$ is generically zero. Therefore, we can write $T_i = \iota_* (T_i^0)$ as a pushforward from some relative divisor $\iota \co Z \inj X'_R$, such that $Z$ is finite over $\Spec R$. The relative tangent complex of $\wt{a}_i$ at any $R$-point of $\Hk^{\flat}_{U}$ over $(\cF_{\star})$ is $\RHom(T_{i}, \cE_{1}^{*}\ot R)$. To analyze its tor-amplitude, write 
\[
\RHom_{X_R'}(T_{i}, \cE_{1}^{*}\ot R) = \RHom_Z(T_i^0, \iota^*\cE_{1}^{*}\ot R).
\]
Since $T_i^0$ is perfect of tor-amplitude $[-1, 0]$ on $Z$, $\RHom_Z(T_i^0,\iota^*\cE_1^*\ot R)$ is perfect of tor-amplitude $[0,1]$ over $Z$. Since $Z$ is finite over $\Spec R$, the direct image preserves tor-amplitude, hence $\RHom_{X_R'}(T_{i}, \cE_{1}^{*}\ot R)$ is perfect of tor-amplitude $[0, 1]$ over $R$. Let $T_{i,\univ}$ be the universal version of $T_{i}$ over $X'\times \Hk^{r}_{S}$, and form the perfect complex $\ul{\RHom(T_{i,\univ}, \cE_1^*)}$ on $\Hk^{r}_{S}$. Then the relative tangent complex of $\wt{a}_i$ is the pullback from $\Hk_S^r$ of $\ul{\RHom(T_{i,\univ}, \cE_1^*)}$, which is perfect of tor-amplitude $[0, 1]$. This shows that $\wt{a}_i$ is quasi-smooth.

Since being a closed embedding is a condition on classical truncations, the assertion that $\wt{a}_i$ is a closed embedding can be checked after base changing to the classical truncation of $S$. Then $T_i$ is an honest torsion sheaf (concentrated in degree $0$), so $\ul{\RHom(T_{i,\univ}, \cE_1^*)}$ is represented by a locally free coherent sheaf in degree $1$. By Lemma \ref{lem:closed-embedding}(3), this implies that $\wt{a}_i$ is a closed embedding.

The analysis of $\wt{b}_i$ is similar, except that its tangent complex is the pullback of $\ul{\RHom(T_{i,\univ}, \wt Q_1)}$, which is locally represented by a complex of locally free coherent sheaves in degrees $[0,1]$. 

The analysis of $\wt{c}_i$ is similar, except that its tangent complex is the pullback of $\ul{\RHom(T_{i,\univ}, \s^{*}\cE_{2}[1])}$, which is locally represented by a locally free coherent sheaf, which implies that $\wt{c}_i$ is a smooth vector bundle.

\end{proof}

By the argument of \cite[Lemma 6.9]{FYZ}, each map $h_i \co \Hk_S^r \rightarrow S$ is smooth. Since the maps $h_i^?$ in \eqref{eq: hiUVW} are all base changed from $h_i$, we get the following Corollary. 

\begin{cor}\label{cor: hecke corr quasismooth} We have the following properties of the morphisms in \eqref{eq: no tilde maps}.
\begin{enumerate}
\item Each morphism $a_i$ is quasi-smooth.
\item Each morphism $b_i$ is quasi-smooth and representable in derived schemes. 
\item Each morphism $c_i$ is smooth. 
\end{enumerate}
\end{cor}

\begin{cor}\label{cor: hecke corr base change assumption}  The diagram
\begin{equation*}
\xymatrix{ U_{0}\ar[ddd]^{\pi_{0}}\ar[dr]^{f_{0}} && \ar[ll]_-{a_{0}}\ar[rr]^-{a_{r}}\Hk^{\flat}_{U}\ar[ddd]^{\pi}\ar[dr]^{f} && U_{r}\ar[ddd]^{\pi_{r}}\ar[dr]^{f_{r}} \\
& V_{0}\ar[ddd]^{g_{0}} && \ar[ll]_-{b_{0}}\ar[rr]^-{b_{r}} \Hk^{\flat}_{V}\ar[ddd]^{g} && V_{r}\ar[ddd]^{g_{r}}\\
\\
S_{0}\ar[dr]^{z_{0}} && \ar[ll]^-{h_{0}}\ar[rr]_-{h_{r}} \Hk^{r}_{S}\ar[dr]^{z} && S_{r}\ar[dr]^{z_{r}}\\
& W_{0} && \ar[ll]_-{c_{0}}\ar[rr]^-{c_{r}} \Hk^{\flat}_{W} && W_{r} 
}
\end{equation*}
satisfies the conditions in \S\ref{sssec: descent assumptions}. Here, $z_{i}$ and $z$ are the inclusions of zero sections, and $\pi_{i}$ and $\pi$ are the natural projections.
\end{cor}
\begin{proof} We first check all maps in the above diagram are representable in derived schemes (this is abbreviated ``schematic'' below) and separated. 
\begin{itemize}
\item The horizontal maps. The maps $h_{i}$ are separated and schematic. By Lemma \ref{lem: hecke corr quasismooth tilde}, the maps $\wt a_{i}, \wt b_{i}$ and $\wt c_{i}$ are separated and schematic. Therefore the same is true for $a_{i}=h_{i}^{U}\c\wt a_{i}, b_{i}=h_{i}^{V}\c \wt b_{i}$ and $c_{i}=h_{i}^{W}\c\wt c_{i}$.
\item The maps $f_{i}, f$ and $z_{i}, z$. Since $W_{i}$ is a classical vector bundle stack over $S_{i}$, $z_{i}$ is separated and schematic. The same is true for $z$. By the Cartesian diagrams \eqref{eq: derived cartesian E_1 S} and \eqref{eq: Hk UVW Cart}, we see that $f_{i}$ and $f$ are also separated and schematic.
\item The maps $\pi_{i}, \pi$ and $g_{i}, g$. The map $g_{i}$ is given by the linear map  $\cV\to \cW$ of perfect complexes on $S$ whose cone is $\cU[1]$. To show $g_{i}$ and $\pi_{i}$ are separated and schematic, it suffices by Lemma \ref{lem:closed-embedding} to observe that $\cU$ has tor-amplitude $[0, \infty)$, which is clear from the definition. Similarly, the fact that $\cU^{\flat}_{\Hk}$ has tor-amplitude in $[0, \infty)$ implies that $g$ and $\pi$ are schematic and separated.
\end{itemize}

The vertical squares $(U_{0},V_{0}, S_{0}, W_{0})$, $(\Hk^{\flat}_{U},\Hk^{\flat}_{V},\Hk^{r}_{S},\Hk^{\flat}_{W})$ and $(U_{r}, V_{r}, S_{r}, W_{r})$ are derived Cartesian by \eqref{eq: derived cartesian E_1 S} and \eqref{eq: Hk UVW Cart}. It remains to check the pushability and pullability of various squares, which will all be seen to be cases of Example \ref{ex: pushable} or Example \ref{ex: pullable}. 
\begin{enumerate}
\item The square $(\Hk^{\flat}_{U}, U_{0}, \Hk^{\flat}_{V},V_{0})$ is pushable. For this it suffices to base change all relevant spaces to $\Hk^{r}_{S}$ and show instead that
\begin{equation*}
\xymatrix{ \wt U_{0}\ar[d]^{\wt f_{0}} & \Hk^{\flat}_{U} \ar[l]_{\wt a_{0}}\ar[d]^{f}\\
\wt V_{0} & \Hk^{\flat}_{V}\ar[l]_{\wt b_{0}}
}
\end{equation*}
is pushable. This follows from the fact that $\wt a_{0}$ is proper (Lemma \ref{lem: hecke corr quasismooth tilde}(1)) and $\wt b_{0}$ is separated.
\item The square $(\Hk^{r}_{S}, S_{0}, \Hk^{\flat}_{W},W_{0})$ is pushable. After base change to $\Hk^{r}_{S}$, it suffices to show that
\begin{equation*}
\xymatrix{ \Hk^{r}_{S}\ar[d]^{\wt z_{0}} & \Hk^{r}_{S} \ar@{=}[l]\ar[d]^{z}\\
\wt W_{0} & \Hk^{\flat}_{W}\ar[l]_{\wt c_{0}}
}
\end{equation*}
is pushable. This follows from the fact that $\Id_{\Hk^{r}_{S}}$ is proper and $\wt c_{0}$ is separated.
\item The square $(\Hk^{\flat}_{U}, U_{r}, \Hk^{r}_{S}, S_{r})$ is pullable. It suffices to check that
\begin{equation*}
\xymatrix{\Hk^{\flat}_{U} \ar[d]^{\pi}\ar[r]^{\wt a_{r}} & \wt U_{r}\ar[d]^{\wt \pi_{r}}\\
\Hk^{r}_{S} \ar@{=}[r]&\Hk^{r}_{S}}
\end{equation*}
is pullable. This follows from the fact that $\wt a_{r}$ is quasi-smooth (Lemma \ref{lem: hecke corr quasismooth tilde}(1)) and $\Id$ is smooth.  
\item The square $(\Hk^{\flat}_{V}, V_{r}, \Hk^{\flat}_{W}, W_{r})$ is pullable. It suffices to check that
\begin{equation*}
\xymatrix{\Hk^{\flat}_{V} \ar[d]^{g}\ar[r]^{\wt b_{r}} & \wt V_{r}\ar[d]^{\wt g_{r}}\\
\Hk^{\flat}_{W} \ar[r]^{\wt c_{r}} &\wt W_{r}
}
\end{equation*}
is pullable. This follows from the fact that $\wt c_{r}$ is smooth (Lemma \ref{lem: hecke corr quasismooth tilde}(3)) and $\wt b_{r}$ is quasi-smooth (Lemma \ref{lem: hecke corr quasismooth tilde}(2)).
\end{enumerate}
\end{proof}

\begin{lemma}\label{lem: dimension add}
For $i=0,\ldots, r$, we have $d(\wt{a}_i) + d(\wt{c}_i) = d(\wt{b}_i)$.
\end{lemma}
\begin{proof}	We have
\begin{align*}
d(\wt{a}_i) & = \rank(\cU^{\flat}_{\Hk}) - \rank(\cU) \\
d(\wt{b}_i)  &= \rank(\cV^{\flat}_{\Hk}) - \rank(\cV) \\
d(\wt{c}_i) &= \rank(\cW^{\flat}_{\Hk}) - \rank(\cW). 
\end{align*}
By \eqref{eq: exact triangle on S E_1} and \eqref{eq: exact triangle on Hk_S E_1}, we get the result. 
\end{proof}

\subsection{The stacks $U^{\perp}, \wh V, W^{\perp}$ and their Hecke correspondences}
\label{ssec: UVW dual} 

Dualizing the sequence \eqref{eq: descent ses 1} and applying $\sigma^*$, we get a short exact sequence
\begin{equation}\label{eq: descent ses 2}
0\to \s^{*}\cE_{1}\to \cE_{2}^{*}\to \wt Q_{2}\to 0
\end{equation}
where $\wt Q_{2}=\s^{*}\wt Q_{1}^{*}:=\s^{*}\cExt^{1}(\wt Q_{1},\cO_{X'})$.

\subsubsection{Definition of $U^{\perp}, \wh V$ and $W^{\perp}$}
We apply the same constructions in \S \ref{ssec: UVW} to get derived vector bundles $U^\perp, V'$ and $W^\perp$ over $S=\Bun^{\le \mu}_{U(n)}$ that fit into a derived Cartesian square
\begin{equation}\label{eq: U'V'W' SES}
\begin{tikzcd}
U^\perp \ar[r] \ar[d]  & V' \ar[d] \\
S \ar[r, "z_{W^\perp}"] & W^\perp
\end{tikzcd}
\end{equation}

To spell out the details, for any animated $\F_{q}$-algebra $R$,
\begin{itemize}
\item Define $\cU^{\perp}:=\ul{\RHom(\cF^*_{\univ}, \cE_2^*)}$ to be the perfect complex on $S$ whose pullback to an $R$-point $\cF \in S(R)$ is naturally in $R$ isomorphic to  $\RHom_{X_R'}(\cF^*, \cE_2^*\ot R)$ regarded in $\Mod_R$. Its associated vector bundle is denoted $U^{\perp} \rightarrow S$.

\item Define $\cV':=\ul{\RHom(\cF^*_{\univ},\wt{Q}_2)}$ to be the perfect complex on $S$ whose pullback to an $R$-point $\cF \in S(R)$ is naturally in $R$ isomorphic to $\RHom_{X_R'}(\cF^*, \wt{Q}_2\ot R)$ regarded in $\Mod_R$. Since $\Ext^1_{X_R'}(\cF^*,\wt{Q}_2\ot R) =0$ as $\wt{Q}_2$ is torsion, $\cV'$ is in fact a locally free coherent sheaf on $S$, corresponding to a vector bundle in the classical sense. Its associated vector bundle is denoted $V' \rightarrow S$.

\item Define $\cW^{\perp}:=\ul{\RHom(\cF^*_{\univ}, \sigma^* \cE_1[1])}$ to be the perfect complex on $S$ whose pullback to an $R$-point $\cF \in S(R)$ is naturally in $R$ isomorphic to $\RHom_{X_R'}(\cF^*, \sigma^* \cE_1[1]\ot R)$ regarded in $\Mod_R$. Its associated vector bundle is denoted $W^{\perp}\rightarrow S$.
\end{itemize}

\begin{lemma}\label{lem: dual vector bundles} 
As vector bundles over $S$, Serre duality identifies $V'$ with the dual vector bundle $\wh V$ of $V$. Under this identification, $U^{\perp}$ is identified with $\wh{W}$, and $W^{\perp}$ is identified with $\wh{U}$, and the derived fiber sequence \eqref{eq: U'V'W' SES} is identified with the dual fiber square to \eqref{eq: derived cartesian E_1 S}
\begin{equation*}
\begin{tikzcd}
\wh W \ar[r] \ar[d]  & \wh V \ar[d] \\
S \ar[r, "z_{\wh U}"] & \wh U
\end{tikzcd}
\end{equation*}
\end{lemma}

\begin{proof}
First we produce the identification $V'\cong \wh{V}$. For any $\cF \in S(R)$,  we have $\RHom_{X_R'}(\cF^*, \wt Q_{1}\ot R) \cong \RHom_{X_R'}(\wt Q^*_{1}[-1]\ot R, \cF)$ as perfect complexes over $R$. By relative Serre duality, the latter is $R$-dual to $\RHom_{X_R'}(\cF, \wt Q_{1}^* \otimes_{\cO_{X'}} \omega_{X'}\ot R)$. Using the Hermitian form $h_{\cF} \co \cF \xrightarrow{\sim} \sigma^* \cF^{\vee}$, we have 
\begin{eqnarray*}
\RHom_{X_R'}(\cF, \wt Q^*_{1} \otimes_{\cO_{X'}} \omega_{X'}\ot R) \cong \RHom_{X_R'}(\sigma^* \cF^{\vee}, \wt Q^*_{1} \otimes_{\cO_{X'}} \omega_{X'}\ot R) \\
\cong \RHom_{X_R'}(\cF^*, \sigma^* \wt Q^*_{1}\ot R)=\RHom_{X_R'}(\cF^*, \wt Q
_{2}\ot R).
\end{eqnarray*}
(For classical $R$, note that all $\RHom$'s appearing above are in fact concentrated in degree 0, hence identified with their respective $\Hom$'s.) This shows that $\RHom_{X_R'}(\cF^*, \wt Q_{1}\ot R)$ is $R$-dual to $ \RHom_{X_R'}(\cF^*, \wt Q_{2}\ot R)$.

Under the above identifications, the exact triangle 
\[
\RHom(\cF^*_{\univ}, \cE_1^*)  \rightarrow \RHom(\cF^*_{\univ}, \wt{Q}_1) \rightarrow \RHom(\cF^*_{\univ}, \sigma^* \cE_2[1])
\]
is dual to the exact triangle
\[
 \RHom(\cF^*_{\univ}, \cE_2^*)  \rightarrow \RHom(\cF^*_{\univ}, \wt{Q}_2) \rightarrow \RHom(\cF^*_{\univ}, \sigma^* \cE_1[1]).
\]
On total spaces, this says that $U \to V \to W$ is dual to $U^\perp \to V'\to W^\perp$.

\end{proof}

From now on, we will identify $V'$ with $\wh V$ using Lemma \ref{lem: dual vector bundles}. 

\subsubsection{More Hecke stacks} 

We let $\Hk_{U^{\perp}}^\flat \rightarrow \Hk_S^r$,  and $\Hk_{W^{\perp}}^\flat \rightarrow \Hk_S^r$, be the derived vector bundles on $\Hk_S^r$ defined similarly to $\Hk_U^\flat \rightarrow \Hk_S^r$, and $\Hk_W^\flat \rightarrow \Hk_S^r$, respectively, but interchanging $\cE_1$ and $\cE_2$. We let $\Hk_{\wh{V}}^\flat \rightarrow \Hk_S^r$ be the vector bundle defined similarly to $\Hk_V^\flat \rightarrow \Hk_S^r$, but replacing $\wt Q_1$ with $\wt Q_2$. 

Given a diagram $(\cF_{\star}) \in \Hk_{S}^r(R)$, define
\[
\cF^{\sh}_{\bu}:=\Cofib\left(\cF_{1/2}^{\flat} \oplus \ldots \oplus \cF_{r-1/2}^{\flat} \rightarrow \cF_0 \oplus \ldots \oplus \cF_{r}\right),
\]
where the map sends $(s_{1/2},\cdots, s_{r-1/2})$ to $(-s_{1/2}, s_{1/2}-s_{3/2}, \cdots, s_{r-3/2}-s_{r-1/2}, s_{r-1/2})$ (using the solid arrows in \eqref{eq: flat hecke diagram} to map $\cF^{\flat}_{i-1/2}$ to $\cF_{i-1}$ and $\cF_{i}$). For classical $R$, $\cF_{\bu}^{\sh}$ is a coherent sheaf on $X'_{R}$ concentrated in degree $0$. Let $\cF_{\bu}^{\sh *}$ be the 
$\cO_{X'_{R}}$-linear dual of $\cF_{\bu}^{\sh}$, i.e., $\cF_{\bu}^{\sh *}$ is the perfect complex on $X'_{R}$ in degrees $0$ and $1$,
\[
 (\cF_0^* \oplus \ldots \oplus \cF_{r}^*)\rightarrow ((\cF_{1/2}^{\flat})^* \oplus \ldots \oplus (\cF_{r-1/2}^{\flat})^*).
\] 
Note that even for classical $R$, the cohomology sheaves of $\cF_{\bu}^{\sh *}$ may be nontrivial in both degrees $0$ and $1$.

Comparing with the definition of $\cF_{\bu}^{\flat}$ given in \S\ref{sss: Hk UVW}, we have an exact triangle in $\Perf(X'\times \Hk_S^r)$:
\begin{equation}\label{eq: sharp exact triangle}
\cF_{\bu}^{\flat} \rightarrow \left( \cF_0 \oplus \cF_r \right) \rightarrow \cF_{\bu}^{\sharp}\to \cF_{\bu}^{\flat}[1]
\end{equation}

We now construct the $\sh$-version of the Hecke stacks over $\Hk_S^r$. 
\begin{itemize}
\item Define a perfect complexes $\cU^{\sh}_{\Hk}, \cV^{\sh}_{\Hk}$ and $\cW^{\sh}_{\Hk}$ on $\Hk_S^r$ similarly to $\cU^{\flat}_{\Hk}, \cV^{\flat}_{\Hk}$ and $\cW^{\flat}_{\Hk}$ respectively, replacing $\cF^{\flat*}_{\bu}$ with $\cF^{\sh*}_{\bu}$. Let $\Hk_U^{\sharp}:= \Tot_{\Hk_S^r}(\cU^{\sh}_{\Hk}), \Hk_{V}^{\sh}:=\Tot_{\Hk_S^r}(\cV^{\sh}_{\Hk})$ and $\Hk^{\sh}_{W}=\Tot_{\Hk_S^r}(\cW^{\sh}_{\Hk})$ be the associated derived vector bundles over $\Hk_S^r$.

\item Switching the roles of $\cE_{1}$ and $\cE_{2}$ and replacing $\wt Q_{1}$ by $\wt Q_{2}$, we define analogously the derived vector bundles $\Hk^{\sharp}_{U^{\perp}}, \Hk^{\sharp}_{\wh V}$ and $\Hk^{\sharp}_{W^{\perp}}$ over $\Hk_S^r$.
 
\end{itemize}

\subsubsection{} The exact triangle \eqref{eq: sharp exact triangle} induces three exact triangles in $\Perf(\Hk_S^r)$:
\begin{eqnarray*}
\cU^{\flat}_{\Hk}\to h_{0}^{*}\cU\op h_{r}^{*}\cU\to \cU^{\sh}_{\Hk}\to \cU^{\flat}_{\Hk}[1]\\
\cV^{\flat}_{\Hk}\to h_{0}^{*}\cV\op h_{r}^{*}\cV\to \cV^{\sh}_{\Hk}\to \cV^{\flat}_{\Hk}[1]\\
\cW^{\flat}_{\Hk}\to h_{0}^{*}\cW\op h_{r}^{*}\cW\to \cW^{\sh}_{\Hk}\to \cW^{\flat}_{\Hk}[1],
\end{eqnarray*}
which induce three (derived) Cartesian squares of derived vector bundles over $\Hk_S^r$:\begin{equation}\label{eq: three squares 1}
\begin{gathered}
\begin{adjustbox}{max width=\textwidth}
$\displaystyle
\begin{tikzcd}[ampersand replacement=\&]
\& \Hk_U^{\flat} \ar[dl, "\wt{a}_0"'] \ar[dr, "\wt{a}_r"] \\ 
\wt{U}_0 \ar[dr, "\wt{a}_r'"'] \& \& \wt{U}_r \ar[dl, "\wt{a}_0'"]  \\
\& \Hk_U^{\sharp} 
\end{tikzcd} \hspace{1cm}
\begin{tikzcd}[ampersand replacement=\&]
\& \Hk_V^{\flat} \ar[dl, "\wt{b}_0"'] \ar[dr, "\wt{b}_r"] \\ 
\wt{V}_0 \ar[dr, "\wt{b}_r'"'] \& \& \wt{V}_r \ar[dl, "\wt{b}_0'"]  \\
\& \Hk_V^{\sharp} 
\end{tikzcd}\hspace{1cm}
\begin{tikzcd}[ampersand replacement=\&]
\& \Hk_W^{\flat} \ar[dl, "\wt{c}_0"'] \ar[dr, "\wt{c}_r"] \\ 
\wt{W}_0 \ar[dr, "\wt{c}_r'"'] \& \& \wt{W}_r \ar[dl, "\wt{c}_0'"]  \\
\& \Hk_W^{\sharp} 
\end{tikzcd}
$
\end{adjustbox}
\end{gathered}
\end{equation}

Analogously, switching the roles of $\cE_{1}$ and $\cE_{2}$, and using the $\wt Q_{2}$ instead of $\wt Q_{1}$, the exact triangle \eqref{eq: sharp exact triangle} induces three (derived) Cartesian squares of derived vector bundles over $\Hk_S^r$:
\begin{equation}\label{eq: three squares 2}
\begin{gathered}
\begin{adjustbox}{max width=\textwidth}
$\displaystyle
\begin{tikzcd}[ampersand replacement=\&]
\& \Hk_{U^{\perp}}^{\flat} \ar[dl, "\wt{a}_0^{\perp}"'] \ar[dr, "\wt{a}_r^{\perp}"] \\ 
\wt{U}_0^{\perp} \ar[dr, "(\wt{a}_r')^{\perp}"'] \& \& \wt{U}_r^{\perp} \ar[dl, "(\wt{a}_0')^{\perp}"]  \\
\& \Hk_{U^{\perp}}^{\sharp} 
\end{tikzcd} \hspace{1cm}
\begin{tikzcd}[ampersand replacement=\&]
\& \Hk_{\wh V}^{\flat} \ar[dl, "\wt\b_{0}"'] \ar[dr, "\wt\b_{r}"] \\ 
\wh{\wt V}_0 \ar[dr, "\wt\b'_{r}"'] \& \& \wh{\wt V}_r \ar[dl, "\wt\b'_{0}"]  \\
\& \Hk_{\wh V}^{\sharp} 
\end{tikzcd}\hspace{1cm}
\begin{tikzcd}[ampersand replacement=\&]
\& \Hk_{W^{\perp}}^{\flat} \ar[dl, "\wt{c}_0^{\perp}"'] \ar[dr, "\wt{c}_r^{\perp}"] \\ 
\wt{W}_0^{\perp} \ar[dr, "(\wt{c}_r')^{\perp}"'] \& \& \wt{W}_r^{\perp} \ar[dl, "(\wt{c}_0')^{\perp}"]  \\
\& \Hk_{W^{\perp}}^{\sharp} 
\end{tikzcd}
$
\end{adjustbox}
\end{gathered}
\end{equation}

\begin{lemma}\label{lem: Hk dual vector bundles} 
Let $\wh{\Hk_U^{\flat}}, \wh{\Hk_V^{\flat}}$ and $\wh{\Hk_{W}^{\flat}}$ be the dual derived bundles to $\Hk_U^\flat, \Hk_{V}^{\flat}$ and $\Hk_{W}^{\flat}$ over $\Hk_S^r$.\footnote{Note $\wh{\Hk_V^{\flat}}$ has a different meaning from $\Hk_{\wh V}^{\flat}$.} Then we have identifications
\begin{enumerate}
\item $\Hk_{U^{\perp}}^{\sharp} \cong \wh{\Hk_W^{\flat}}$ and $\Hk_{U^{\perp}}^{\flat} \cong \wh{\Hk_W^{\sharp}}$.
\item $\Hk_{\wh V}^{\sharp} \cong \wh{\Hk_V^{\flat}}$ and $\Hk_{\wh V}^{\flat} \cong \wh{\Hk_V^{\sharp}}$.
\item $\Hk_{W^{\perp}}^{\sharp} \cong \wh{\Hk_U^{\flat}}$ and $\Hk_{W^{\perp}}^{\flat} \cong \wh{\Hk_U^{\sharp}}$. 
\end{enumerate}
Moreover, under these identifications, the first (resp. second, resp. third) derived Cartesian square in \eqref{eq: three squares 2} is the dual to the third (resp. second, resp. first) derived Cartesian square in \eqref{eq: three squares 1}. In particular, we have
\begin{equation*}
\wt\b_{0}=\wh{\wt b'_{r}}, \quad \wt\b_{r}=\wh{\wt b'_{0}}, \quad \wt\b'_{r}=\wh{\wt b_{0}}, \quad \wt\b'_{0}=\wh{\wt b_{r}}.
\end{equation*}
\end{lemma}

\begin{proof}
Similar to Lemma \ref{lem: dual vector bundles}.
\end{proof}

\subsubsection{Summary}\label{sssec: summary} Starting from \eqref{eq: descent ses 1}, we defined a collection of spaces and maps as in the diagram below: 
\begin{equation}\label{eq: big diagram for E_1}
\begin{tikzcd}
& & \Hk_U^{\flat} \ar[dl, "\wt{a}_0"', color=blue] \ar[dr, "\wt{a}_r"] \ar[ddd, bend left, "f"', color=orange] \\
& \ar[dl, "h_0^U"'] \wt{U}_0 \ar[ddd, "\wt{f}_0"', color=red] \ar[dr, "\wt{a}_r'"'] && \wt{U}_r \ar[dl, "\wt{a}_0'"] \ar[ddd, "\wt{f}_r"] \ar[dr, "h_r^U"] \\
U_0 \ar[ddd, "f_0"]  & & \Hk_U^{\sharp}  \ar[ddd, bend right, "f^{\sharp}"]   &  &  U_r  \ar[ddd, "f_r"] \\
& & \Hk_V^{\flat} \ar[dl, "\wt{b}_0"', color=green] \ar[dr, "\wt{b}_r"] \ar[ddd, bend left, "g"']   \\
& \wt{V}_0 \ar[ddd, "\wt{g}_0"']  \ar[dr, "\wt{b}_r'"']  \ar[dl, "h_0^V"'] && \wt{V}_r \ar[dl, "\wt{b}_0'"]  \ar[ddd, "\wt{g}_r"]  \ar[dr, "h_r^V"] \\
V_0 \ar[ddd, "g_0"] & & \Hk_V^{\sharp} \ar[ddd, bend right, "g^{\sharp}"]  & & V_r  \ar[ddd, "g_r"] \\
& & \Hk_W^{\flat}  \ar[dl, "\wt{c}_0"'] \ar[dr, "\wt{c}_r"] \\
& \wt{W}_0 \ar[dr, "\wt{c}_r'"'] \ar[dl, "h_0^W"']  & &  \wt{W}_r \ar[dl, "\wt{c}_0'"]  \ar[dr, "h_r^W"]  \\
W_0 & & \Hk_W^{\sharp}  & & W_r 
\end{tikzcd} 
\end{equation}
Here: 
\begin{itemize}
\item The maps in the columns come from exact triangles of perfect complexes. 
\item The three diamonds in the middle are derived Cartesian. 
\item The four parallelograms on the left and right sides are derived  Cartesian. 
\end{itemize}

 Starting from \eqref{eq: descent ses 2}, we defined a collection of spaces and maps as in the diagram below:
\begin{equation}\label{eq: big diagram for E_2}
\begin{tikzcd}
& & \Hk_{U^{\perp}}^{\flat} \ar[dl, "\wt{a}_0^{\perp}"'] \ar[dr, "\wt{a}_r^{\perp}"] \ar[ddd, bend left, "f^{\perp}"'] \\
& \ar[dl, "h_0^{U^{\perp}}"'] \wt{U}_0^{\perp} \ar[ddd, "\wt{f}_0^{\perp}"'] \ar[dr, "(\wt{a}_r')^{\perp}"'] && \wt{U}_r^{\perp} \ar[dl, "(\wt{a}_0')^{\perp}"] \ar[ddd, "\wt f_r^{\perp}"] \ar[dr, "h_r^{U^{\perp}}"]  \\
U_0^{\perp} \ar[ddd, "f_0^{\perp}"] & & \Hk_{U^{\perp}}^{\sharp}  \ar[ddd, bend right, "(f^{\sharp})^{\perp}"]   & & U_r^{\perp} \ar[ddd, "f_r^{\perp}"] \\
& & \Hk_{\wh{V}}^{\flat} \ar[dl, "\wt{\b}_0"'] \ar[dr, "\wt{\b}_r"] \ar[ddd, bend left, "g^{\perp}"']   \\
& \ar[dl, "h_0^{\wh V}"']  \wh{\wt{V}}_0 \ar[ddd, "\wt{g}_0^{\perp}"', color=red]  \ar[dr, "\wt\b'_{r}"', color=green]  && \wh{\wt{V}}_r \ar[dl, "\wt\b'_{0}"]  \ar[ddd, "\wt{g}_r^{\perp}"]  \ar[dr, "h_r^{\wh V}"]  \\
\wh{V}_0 \ar[ddd, "g_0^{\perp}"]  & & \Hk_{\wh{V}}^{\sharp} \ar[ddd, bend right, "(g^{\sharp})^{\perp}", color=orange] & & \wh{V}_r \ar[ddd, "g_r^{\perp}"]  \\	
& & \Hk_{W^{\perp}}^{\flat}  \ar[dl, "\wt{c}_0^{\perp}"'] \ar[dr, "\wt{c}_r^{\perp}"] \\
& \ar[dl, "h_0^{W^{\perp}}"']  \wt{W}_0^{\perp} \ar[dr, "(\wt{c}_r')^{\perp}"', color=blue]  & &  \wt{W}_r^{\perp} \ar[dl, "(\wt{c}_0')^{\perp}"]  \ar[dr, "h_r^{W^{\perp}}"]  \\
W_0^{\perp} & & \Hk_{W^{\perp}}^{\sharp} & & W_r^{\perp} 
\end{tikzcd}
\end{equation}
Again:
\begin{itemize}
\item The maps in the columns come from exact triangles of perfect complexes. 
\item The three diamonds in the middle are derived Cartesian. 
\item The four parallelograms on the left and right sides are derived Cartesian. 
\end{itemize}

Furthermore, \eqref{eq: big diagram for E_2} is the dual to the diagram \eqref{eq: big diagram for E_1}. The duality exchanges $U$ with $W^{\perp}$, $V$ with $\wh{V}$, and $W$ with $U^{\perp}$, and exchanges $\flat$ and $\sharp$ superscripts. Sample examples of dual morphisms in \eqref{eq: big diagram for E_1} and \eqref{eq: big diagram for E_2} are colored with the same color.

\subsection{Global presentability}\label{ssec: safety verification}
We will want to apply the derived Fourier theory of \S \ref{sec: FT} and \S \ref{sec: FT cc} to the ensemble of spaces and maps in \S \ref{sssec: summary}. In order to justify this we need to check that all the derived vector bundles are globally presented, and all the maps are globally presented. The reason for Harder-Narasimhan truncation (in this section) is to guarantee these properties.

The following observations are useful to perform this check. 

\begin{enumerate}
\item If $\cT$ is a torsion coherent sheaf on $X'$, then the perfect complex $\ul{\RHom}(\cF_{\univ}^*, \cT)$ on $\Bun_{U(n)}^{\leq \mu}$ is locally free.  
\item Given $\mu$ and any coherent sheaf $\cE$ on $X'$, for any effective divisor $D$ on $X'$ with $\deg D$ sufficiently large (depending on $\mu$ and $\cE$), $\ul{\RHom}(\cF_{\univ}^*, \cE(D))$ is locally free. 
\end{enumerate}
Now, from the exact triangle in $\Perf(X')$,
\[
\cE \rightarrow \cE(D) \rightarrow \cE |_D(D) 
\]
we get an exact triangle of complexes in $\Perf(\Bun_{U(n)}^{\leq \mu})$,
\begin{equation}\label{eq: 2-term complex}
\ul{\RHom}(\cF_{\univ}^*, \cE) \rightarrow \ul{\RHom}(\cF_{\univ}^*, \cE(D)) \rightarrow \ul{\RHom}(\cF_{\univ}^*, \cE |_D(D)).
\end{equation}
By the observations above, the second and third terms in \eqref{eq: 2-term complex} are vector bundles, so this presents $\ul{\RHom}(\cF_{\univ}^*, \cE) $ as a 2-term complex of vector bundles. Since the derived vector bundles $U_i, W_i, V_i, \wt{U}_i, \wt{W}_i, \wt{V}_i$ appearing in \S \ref{ssec: UVW} are all instances of this construction, they are all globally presented. 

Furthermore, if $\cE \rightarrow \cE'$ is a map of coherent sheaves on $X'$, then for any effective divisor $D$ on $X'$ with $\deg D$ sufficiently large (depending on $\mu, \cE, \cE'$) the diagram 
\[
\begin{tikzcd}
\ul{\RHom}(\cF_{\univ}^*, \cE) \ar[r] \ar[d] &  \ul{\RHom}(\cF_{\univ}^*, \cE(D)) \ar[r] \ar[d]  &  \ul{\RHom}(\cF_{\univ}^*, \cE |_D(D)) \ar[d] \\
\ul{\RHom}(\cF_{\univ}^*, \cE') \ar[r] &  \ul{\RHom}(\cF_{\univ}^*, \cE'(D)) \ar[r] &  \ul{\RHom}(\cF_{\univ}^*, \cE'|_D(D))  
\end{tikzcd}
\]
gives a global presentation for $\ul{\RHom}(\cF_{\univ}^*, \cE)  \rightarrow \ul{\RHom}(\cF_{\univ}^*, \cE')$ as a map of complexes of vector bundles. A similar trick gives a global presentation for all the commutative quadrilaterals involving the $U_i, V_i, W_i$ and $S_i$.

Since the universal perfect complex $\cF_{\bu}^{\flat}$ used to construct $\cU_{\Hk}^\flat, \cV_{\Hk}^\flat, \cW_{\Hk}^\flat$ is assembled out of the pullbacks of the $\cF_{\univ}$ from the various projections to $\Bun_{U(n)}^{\leq \mu}$, all the maps between $\Hk_U^\flat, \Hk_V^\flat, \Hk_W^\flat$ and $\Hk_S^r$ are also globally presented (this implicitly includes the statement that the individual derived vector bundles are globally presented). The same applies to all the variants of these spaces and maps considered in \S \ref{ssec: UVW}. 

\emph{We have now verified that the diagrams in \eqref{eq: big diagram for E_1} and \eqref{eq: big diagram for E_2} are globally presented, so that we may apply the results of \S \ref{sec: FT cc} to them.}

\subsection{Comparison of cohomological correspondences}\label{ssec: comparison of cc} We refer to the diagram \eqref{eq: big diagram for E_1} and its Fourier dual, diagram \eqref{eq: big diagram for E_2}.

By Corollary \ref{cor: hecke corr quasismooth}, the relative fundamental class of the quasi-smooth map $a_{r}$ defines a cohomological correspondence
\begin{equation*}
\cc_{U}=[a_{r}]\in \Corr_{\Hk^{\flat}_{U}}(\Qll{U_0}, \Qll{U_{r}}\tw{-d(a_r)}).
\end{equation*}
Similarly, the relative fundamental class of $a^{\perp}_{r}$ defines a cohomological correspondence
\begin{equation*}
\cc_{U^{\perp}}=[a^{\perp}_{r}]\in \Corr_{\Hk^{\flat}_{U^{\perp}}}(\Qll{U^{\perp}_0}, \Qll{U^{\perp}_{r}}\tw{-d(a^{\perp}_r)}).
\end{equation*}
By Corollary \ref{cor: hecke corr base change assumption}, the pushforward of cohomological correspondences along the morphism of correspondences $f: \Hk^{\flat}_{U}\to \Hk^{\flat}_{V}$ is defined, giving
\begin{equation*}
f_{!}(\cc_{U})\in \Corr_{\Hk^{\flat}_{V}}(f_{0!}\Qll{U_{0}}, f_{r!}\Qll{U_{r}}\tw{-d(a_{r})}).
\end{equation*}
Similarly, 
\begin{equation*}
f^{\perp}_{!}(\cc_{U^{\perp}})\in \Corr_{\Hk^{\flat}_{\wh V}}(f^{\perp}_{0!}\Qll{U^{\perp}_{0}}, f^{\perp}_{r!}\Qll{U^{\perp}_{r}}\tw{-d(a^{\perp}_{r})})
\end{equation*}
is defined.

Recall the notion of Fourier transform of cohomological correspondences from \S \ref{ssec: FT of cc}. We have
\begin{equation}\label{eq: ft f_! cc_U}
\FT_{\Hk^{\flat}_{V}}(f_{!}(\cc_{U}))\in \Corr_{\Hk^{\flat}_{\wh V}}(\FT_{V_{0}}(f_{0!}\Qll{U_{0}}), \FT_{V_{r}}(f_{r!}\Qll{U_{r}})\sm{[d(\wt b_{0})+d(\wt b_{r})](d(\wt b_{0}))}\tw{-d(a_{r})}).
\end{equation}
Here we use that $\wh{\Hk^{\sh}_{V}}\cong \Hk^{\flat}_{\wh V}$ (see Lemma \ref{lem: Hk dual vector bundles}). Note that $d(\wt b_{0})=d(\wt b_{r})$, therefore
\begin{equation*}
[d(\wt b_{0})+d(\wt b_{r})](d(\wt b_{0}))\tw{-d(a_{r})}=\tw{d(\wt b_{r})-d(a_{r})}.
\end{equation*}
Hence \eqref{eq: ft f_! cc_U} simplifies to
\begin{equation}
\FT_{\Hk^{\flat}_{V}}(f_{!}(\cc_{U}))\in \Corr_{\Hk^{\flat}_{\wh V}}(\FT_{V_{0}}(f_{0!}\Qll{U_{0}}), \FT_{V_{r}}(f_{r!}\Qll{U_{r}})\tw{d(\wt b_{r})-d(a_{r})}).
\end{equation}

Since $U^{\perp}$ is the orthogonal complement of $U$ relative to $V$ (in the derived sense), by \S\ref{sssec: FT functoriality} and Example \ref{ex: delta-constant duality}, we have canonical isomorphisms for $i=0,r$:
\begin{eqnarray*}
\FT_{V_{i}}(f_{i!}\Qll{U_{i}})\cong f^{\perp}_{i!}\Qll{U^{\perp}_{i}}\sm{[\rank(\cV)]}\tw{-\rank(\cU)}.
\end{eqnarray*}
Note the shift and twist on the right side is the same for $i=0$ and $i=r$.  Therefore $\FT_{\Hk^{\flat}_{V}}(f_{!}(\cc_{U}))$ can also be viewed as an element in 
\begin{equation*}
\Corr_{\Hk^{\flat}_{\wh V}}(f^{\perp}_{0!}\Qll{U^{\perp}_{0}}, f^{\perp}_{r!}\Qll{U^{\perp}_{r}}\tw{d(\wt b_{r})-d(a_{r})}).
\end{equation*}
 
On the other hand, by Lemma \ref{lem: Hk dual vector bundles}, $\wt c_{r}$ is dual to $(\wt a'_{0})^{\perp}$, which  has the same relative dimension as $\wt a^{\perp}_{0}$. Therefore, by Lemma \ref{lem: dimension add}, we have
\begin{equation*}
d(\wt b_{r})-d(\wt a_{r})=d(\wt c_{r})=-d(\wt a_{0}^{\perp})=-d(\wt a_{r}^{\perp}).
\end{equation*}
Therefore
\begin{equation*}
d(\wt b_{r})-d(a_{r})=d(\wt b_{r})-d(\wt a_{r})-d(h_{r})=-d(\wt a_{r}^{\perp})-d(h_{r})=-d(a_{r}^{\perp}).
\end{equation*}
We can therefore view $\FT_{\Hk^{\flat}_{V}}(f_{!}(\cc_{U}))$ as an element in
\[
\Corr_{\Hk^{\flat}_{\wh V}}(f^{\perp}_{0!}\Qll{U^{\perp}_{0}}, f^{\perp}_{r!}\Qll{U^{\perp}_{r}}\tw{-d(a_{r}^{\perp})}).
\]

The main result of this section is the following theorem.

\begin{thm}\label{thm: FT of pushforward correspondence} 
Under the above notations, we have
\[
\TT_{[d(f_0) + d(\pi_0)](d(\pi_0))}\FT_{\Hk_V^{\flat}}(f_{!} (\cc_U))   = f^{\perp}_{!} (\cc_{U^\perp})
\]
as elements of $\Corr_{\Hk_{\wh V}^\flat}(f_{0!}^\perp \Qll{U^{\perp}_0}, f_{r!}^\perp \Qll{U^{\perp}_r} \tw{-d(a_r^\perp)})$.
\end{thm}
\begin{proof}

Let $\frs\in \Corr_{\Hk^{r}_{S}}(\Qll{S}, \Qll{S}\tw{-d(h_{r})})$ be given by the relative fundamental class of $h_{r}$
\begin{equation*}
\frs=[h_{r}]: h^{*}_{0}\Qll{S}\to h^{!}_{r}\Qll{S}\tw{-d(h_{r})}.
\end{equation*}
Recall the maps of correspondences
\begin{equation*}
\begin{gathered}
\begin{aligned}
\pi: \Hk^{\flat}_{U}\to \Hk_{S}^{r},\quad
\pi^{\perp}: \Hk^{\flat}_{U^{\perp}}\to \Hk_{S}^{r},\\
z^{\perp}: \Hk^{r}_{S}\to \Hk^{\flat}_{W^\perp},\quad
g^{\perp}: \Hk_{\wh{V}}^{\flat}\to \Hk_{W^\perp}^{\flat}.
\end{aligned}
\end{gathered}
\end{equation*}

The theorem follows from a sequence of equalities of cohomological correspondences
\begin{align*}
\TT_{[d(f_0) + d(\pi_0)](d(\pi_0))}\FT(f_{!}\cc_{U}) &\stackrel{(1)}{=} 
\TT_{[d(f_0) + d(\pi_0)](d(\pi_0))} \FT(f_{!}\pi^{*}\frs)\\
&\stackrel{(2)}{=} (g^\perp)^{*}z^\perp_{!}\FT(\frs) \stackrel{(3)}{=}(g^\perp)^{*} z^{\perp}_! \frs\stackrel{(4)}{=}f_{!}^\perp(\pi^\perp)^{*}\frs\stackrel{(5)}{=}f_{!}^\perp \cc_{U^\perp}.
\end{align*}
We explain the reason for each equality:
\begin{enumerate}
\item[(1),(5)] follow from the equalities
\begin{equation*}
\pi^{*}\frs=\cc_{U}, \quad (\pi^{\perp})^{*}\frs=\cc_{U^{\perp}}
\end{equation*}
which will be proved in Lemma \ref{l:pi pullback s}.

\item[(2)] involves two applications of Proposition \ref{prop: cc push/pull functoriality}, namely
\begin{equation*}
\TT_{[d(f_0)]}\FT\c f_{!}=(g^\perp)^*\c\FT, \quad \TT_{[d(\pi_0)](d(\pi_0)) } \FT\c \pi^{*}=z^\perp_{!}\c\FT.
\end{equation*}
We used here that $g^{\perp\sh}:\Hk^{\sh}_{\wh V}\to \Hk^{\sh}_{W^{\perp}}$ is the dual of the map of correspondences $f: \Hk^{\flat}_{U}\to \Hk^{\flat}_{V}$, as summarized in the diagrams \eqref{eq: big diagram for E_1} and \eqref{eq: big diagram for E_2}. Under the identifications of Lemma \ref{lem: Hk dual vector bundles}, this is the map denoted $g^\perp$ in \eqref{eq: big diagram for E_2}.  Similarly, $\pi$ is dual to $z^{\perp\sh}$; in fact, the diagram
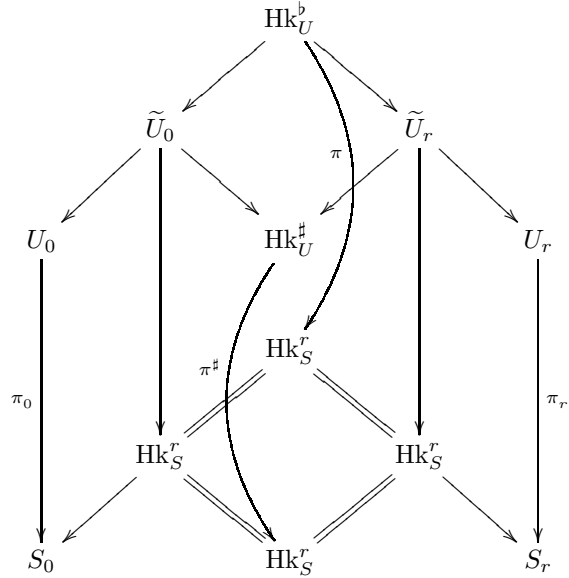
\begin{equation*}
\xymatrix{& & \Hk^{\flat}_{U}\ar[dl]\ar[dr]\ar@/^2pc/[ddd]_(.4){\pi}\\
& \wt U_{0}\ar[dl]\ar[dr]\ar[ddd] & & \wt U_{r}\ar[dl]\ar[dr]\ar[ddd]\\
U_{0}\ar[ddd]_{\pi^{}_{0}} && \Hk^{\sh}_{U}\ar@/_2pc/[ddd]_(.4){\pi^{\sh}} & & U_{r}\ar[ddd]^{\pi^{}_{r}}\\
& & \Hk^{r}_{S}\ar@{=}[dr]\ar@{=}[dl]\\
& \Hk^{r}_{S}\ar[dl]\ar@{=}[dr] & & \Hk^{r}_{S}\ar@{=}[dl]\ar[dr]\\
S_{0} && \Hk^{r}_{S} & & S_{r}}
\end{equation*}
is dual to the diagram
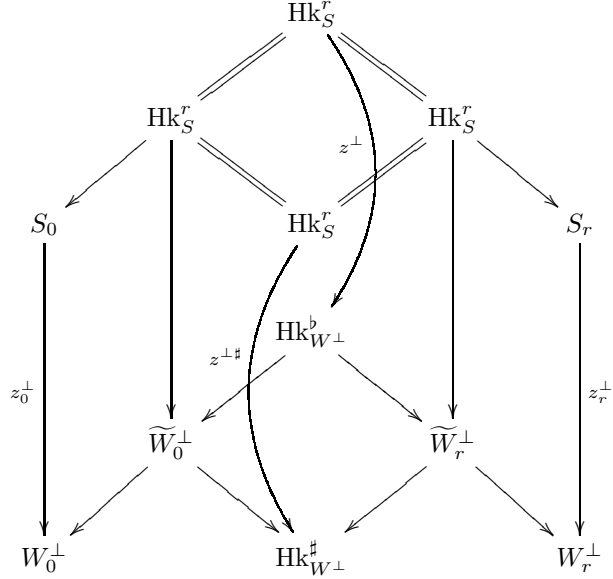
\begin{equation*}
\xymatrix{& & \Hk^{r}_{S}\ar@{=}[dr]\ar@{=}[dl]\ar@/^2pc/[ddd]_(.4){z^\perp}\\
& \Hk^{r}_{S}\ar[dl]\ar@{=}[dr]\ar[ddd] & & \Hk^{r}_{S}\ar@{=}[dl]\ar[dr]\ar[ddd]\\
S_{0}\ar[ddd]_{z_{0}^\perp}  && \Hk^{r}_{S}\ar@/_2pc/[ddd]_(.4){z^{\perp\sh}} & & S_{r}\ar[ddd]^{z_{r}^\perp}\\
& & \Hk^{\flat}_{W^\perp}\ar[dl]\ar[dr]\\
& \wt W^\perp_{0}\ar[dl]\ar[dr] & & \wt W^\perp_r\ar[dl]\ar[dr]\\
W^\perp_{0}&& \Hk^{\sh}_{W^\perp} & & W_{r}^\perp
}
\end{equation*}

\item[(3)] is the trivial equality $\frs=\FT_{\Hk_{S}^{r}}(\frs)$.

\item[(4)] follows from Theorem \ref{thm: descent of pushforward correspondence} and Lemma \ref{l:pi pullback s}. Note that we have verified in Corollary \ref{cor: hecke corr base change assumption} that the hypotheses of Theorem \ref{thm: descent of pushforward correspondence} hold in this situation. 
\end{enumerate}

\end{proof}

\begin{lemma}\label{l:pi pullback s}
With $\frs=[h_r]$ as above, we have
\[
\pi^{*}\frs=\cc_{U}\in \Corr_{\Hk^{\flat}_{U}}(\Qll{U_{0}}, \Qll{U_{r}}\tw{-d(a_{r})}).
\]
\end{lemma}
\begin{proof}
Unravelling the definition, we need to show that the composition 
\begin{equation*}
\Qll{\Hk^{\flat}_{U}}\xr{\pi^{*}[h_{r}]}\pi^{*}h_{r}^{!}\Qll{S_{r}}\tw{-d(h_{r})}\xr{\tu} a_{r}^{!}\pi^{*}_{r}\Qll{S_{r}}\tw{-d(a_{r})}
\end{equation*}
is equal to the relative fundamental class $[a_{r}]$. Here $\tu$ is the pull-pull base change map attached to the pullable (outer) square
\begin{equation}\label{HkUr factor}
\xymatrix{\Hk^{\flat}_{U}\ar@/^{1.5pc}/[rr]^{a_{r}}\ar[r]^{\wt a_{r}}\ar[dr]_{\pi} & \wt U_{r}\ar[r]^{h^{U}_{r}}\ar[d]^{\wt \pi_{r}} & U_{r}\ar[d]^{\pi_{r}}\\
& \Hk^{r}_{S}\ar[r]^{h_{r}} & S_{r} 
}
\end{equation}
By construction, $\tu$ is the composition of two steps
\begin{eqnarray*}
\pi^{*}h_{r}^{!}\Qll{S_{r}}\tw{-d(h_{r})}=\wt a_{r}^{*}\wt\pi^{*}_{r}h_{r}^{!}\Qll{S_{r}}\tw{-d(h_{r})}\xr{\di} \wt a_{r}^{*}(h_{r}^{U})^{!}\pi_{r}^{*}\Qll{S_{r}}\tw{-d(h_{r})}\\
\xr{[\wt a_{r}]}\wt a_{r}^{!}(h_{r}^{U})^{!}\pi^{*}_{r}\Qll{S_{r}}\tw{-d(\wt a_{r})-d(h_{r})}=a_{r}^{!}\pi^{*}_{r}\Qll{S_{r}}\tw{-d(a_{r})}.
\end{eqnarray*}
As explained in \S \ref{ssec: relative fundamental nt}, on the level of global sections the map $\di$ (adjoint to proper base change for the derived Cartesian square in \eqref{HkUr factor}) sends the relative fundamental class of $h_{r}$ to the relative fundamental class of $h^{U}_{r}$. Post-composing with $[\wt a_{r}]$, we get $[a_{r}]$.
\end{proof}

\section{Proof of the generic modularity theorem}\label{sec: AFT and modularity}
In this section we carry out Steps (6) and (7) of the proof outline from \S \ref{ssec: outline}, thus completing the proof of Theorem \ref{thm: transverse modularity}.

\subsection{Transverse Lagrangians}\label{ssec: AFT setup}
Let $\cG \in \Bun_{GU^{-}(2m)}(k)$. For notational simplicity, we assume the similitude line bundle of $\cG$ is trivial, hence the skew-Hermitian form takes the form $h_{\cG}: \cG\isom \s^{*}\cG^{*}$. Let $\wt{\cE}_1,  \wt{\cE}_2 \inj \cG$ be two transverse Lagrangian sub-bundles. The Modularity Conjecture asserts the equality of elements in $\Ch_{(n-m)r}(\Sht^{r}_{U(n)})_{\Qlbar}$
\begin{equation}\label{eq:modularity-eq}
\wt{Z}_m^r (\cG, \wt{\cE}_1) \stackrel{?}{=} \wt{Z}_m^r(\cG, \wt{\cE}_2).
\end{equation}
Our goal is to show that given any Harder-Narasimhan polygon $\mu$ for $\Bun_{U(n)}$, both sides of \eqref{eq:modularity-eq} have the same image in $\mBM_{2(n-m)r}(\Sht^{r,\le \mu}_{U(n)}\times_{(X')^{r}}(X'-T)^r;  \Qlbar)$ for a finite set of closed points $T$ depending only on $\wt{\cE}_{1}, \wt{\cE}_{2}$, and $\mu$.  This will establish Theorem \ref{thm: transverse modularity}.

We reproduce the diagram \eqref{eq: Q diagram} with $\cE_{i}$ replaced by $\wt{\cE}_{i}$:
\begin{equation}\label{eq: Q diagram for saturation}
\begin{tikzcd}
\wt{\cE}_1 \ar[rr, "b_{12}"] \ar[dr, hook] & & \sigma^* \wt{\cE}_2^* \ar[rr]\ar[dr, hook, "i_{1}"] & &  Q_{2}  \ar[dr, "\io_{2}"] \\
& \cG \ar[ur] \ar[rr] \ar[dr]  & &  \cG^{\sh}\ar[rr] && Q\\
\wt{\cE}_2 \ar[rr, "b_{21}"] \ar[ur, hook] & & \sigma^* \wt{\cE}_1^* \ar[rr] \ar[ur, hook, "i_{2}"'] & & Q_{1}\ar[ur, "\io_{1}"'] 
\end{tikzcd}
\end{equation}
Here $b_{12}$ is the composition
\[
\wt{\cE}_1 \rightarrow \cG\isom \s^{*}\cG^{*} \rightarrow \sigma^*  \wt{\cE}_2^*. 
\]
The map $b_{21}$ is defined similarly, and the three rows are short exact.  In particular, the maps $\io_{1}$ and $\io_{2}$ are isomorphisms of torsion sheaves on $X'$. As in \S\ref{sssec: self-duality of Q}, the duality between $Q_{1}$ and $Q_{2}$ equips $Q$ with two Hermitian structures $h_{12} \co Q \xrightarrow{\sim} \sigma^* Q^*$ and $h_{21} \co  Q \xrightarrow{\sim} \sigma^* Q^*$, related by $h_{12} = -h_{21}$.

For each $i \in \{1,2\}$ let $\cE_i \inj \wt{\cE}_i$ be a sub-sheaf with cokernel a torsion coherent sheaf $\cT_i$ on $X'$. Let $\cT^{*}_{i}=\RHom(\cT_{i}, \cO_{X'})[1]$ be its dual torsion sheaf on $X'$. Therefore, $\wt{\cE}_i$ is the saturation of $\cE_i$ in $\cG$, and in the diagram below 
\begin{equation}\label{eq: TL exact sequence}
\begin{tikzcd}
\cE_1  \ar[r, hook, "\cT_1"] &   \wt{\cE}_1 \ar[r, hook] \ar[rr, hook, bend right, "Q"'] & \cG \ar[r, twoheadrightarrow] & \s^{*}\wt{\cE}_{2}^* \ar[r, hook, "\sigma^* \cT_2^*"] &  \s^{*}\cE_2^*
\end{tikzcd}
\end{equation}
the arrows are labeled by their cokernels. 

\subsubsection{Assumptions on $\cT_{1}$ and $\cT_{2}$}\label{sss: assumptions Qi}
Fix a Harder-Narasimhan polygon $\mu$ for $\Bun_{U(n)}$, and write $S=\Bun_{U(n)}^{\le \mu}$ for the corresponding open substack of $\Bun_{U(n)}$.

We make the following assumptions:
 \begin{enumerate}
\item The supports $|Q|, |\cT_1|, |\cT_2|$ are disjoint after mapping to $X$.
\item For all $\cF \in S(\ov k)=\Bun_{U(n)}^{\leq \mu}(\ol{k})$ we have for $i=1,2$
\begin{equation}\label{eq: vanishing assumption 1}
\Ext^1_{X'_{\ol{k}}}(\cF^*, \cE_i^*) = 0.
\end{equation}
\end{enumerate}

Note that by the dualities in \eqref{eq: comparison of 5term}, \eqref{eq: vanishing assumption 1} is equivalent to
\begin{equation}\label{eq: vanishing assumption 2}
\Hom_{X'_{\ol{k}}}(\cF^*, \sigma^* \cE_i) = 0
\end{equation}
for all $\cF \in \Bun_{U(n)}^{\leq \mu}(\ol{k})$ and $i=1,2$.

\begin{remark} Given any $\mu$, the conditions \eqref{eq: vanishing assumption 1}, \eqref{eq: vanishing assumption 2} can be arranged by taking $\cE_i = \wt{\cE}_i(-D_i)$ for sufficiently large effective divisors $D_1, D_2$ whose supports over $X$ are disjoint from each other and from $|Q|$. By taking $D_i$ to be a sufficiently large multiple of a single closed point, we may even arrange that the support $|\cT_{i}|$ is a single closed point of $X'$.
\end{remark}

Let
\begin{eqnarray}
\wt Q_{1}:=Q^{*} \oplus \cT_1^* \oplus \s^{*}\cT_2,\\
\wt Q_{2}:= \sigma^* Q \oplus \s^{*}\cT_1 \oplus \cT_2^{*}.
\end{eqnarray}
From \eqref{eq: TL exact sequence} and the disjointness assumption in \S\ref{sss: assumptions Qi}, we have short exact sequences
\begin{equation}\label{eq: TL SES 1} 
0 \rightarrow  \s^{*}\cE_2 \rightarrow \cE_1^* \rightarrow \wt Q_{1} \rightarrow 0,
\end{equation}
\begin{equation}\label{eq: TL SES 2} 
0 \rightarrow \s^{*}\cE_1 \rightarrow \cE_2^* \rightarrow \wt Q_{2}\rightarrow 0
\end{equation}
  
\begin{remark}
The torsion coherent sheaves $\cT_1, \cT_2$ are auxiliary objects introduced for purely technical reasons (to ensure that certain spaces are smooth, and that certain maps are closed embeddings). They do not appear in \S \ref{ssec: r=0 transverse lagrangians}, but are necessary when $r>0$ in the argument as currently construed. 
\end{remark}

\subsubsection{Moduli spaces}

For $i \in \{0,r\}$ we define $U_i, \wt{U}_i, V_i, \wt{V}_i$, $W_i, \wt{W}_i$, etc. as in \S \ref{ssec: modularity for cc}. 

The vanishing assumption \eqref{eq: vanishing assumption 1} implies that $U,V$ and $W$ are all classical vector bundles over $S$, and we have a short exact sequence of classical vector bundles over $S$
\begin{equation}\label{eq: UVW SES}
0\to U \to V\to W\to 0.
\end{equation}
Similarly, we have a short exact sequence of classical vector bundles over $S$
\begin{equation}
0\to U^{\perp} \to \wh V\to W^{\perp}\to 0.
\end{equation}

\subsubsection{Open locus}
We denote
\begin{align}
X^{\c}&=X-\nu(|Q|\cup |\cT_{1}|\cup |\cT_{2}|)=X-\nu(|\wt Q_{1}|)=X-\nu(|\wt Q_{2}|);\\
X'^{\c}&=\nu^{-1}(X^{\c}).
\end{align}
For a stack $\cY$ over $X^r$, we denote
\begin{equation}
\cY^{\c}:=\cY\times_{X^{r}}(X^{\c})^{r}.
\end{equation}
In particular,   $\Hk_S^{r,\c} \subset \Hk_S^r$ denotes the open substack where the legs are all disjoint from $|\wt{Q}_1|\cup|\wt Q_{2}|$.

On $(X^{\c})^{r}$, the structure of $\Hk_V^r$ is simpler.

\begin{lemma}\label{l: HkV Xc}
For any $0\leq i \leq r$, the restriction $\wt b_i^{\c} \co \Hk_{V}^{\flat,\c} \rightarrow \wt{V}^{\c}_i$ of $\wt b_{i}$ and the restriction $\wt b'^{\c}_{i}: \wt{V}^{\c}_i\to \Hk_{V}^{\sh,\c} $ of $\wt b'_{i}$ are isomorphisms. 
\end{lemma}

\begin{proof} Let $(\cF_{\star})\in \Hk^{r,\c}_{S}(R)$. 
The projection $\wt b_i \co \Hk_V^{\flat} \rightarrow \wt{V}_i$, when base changed over $\cF_{\star}: \Spec R\to \Hk^{r}_{S}$, is induced by the map $\frp^{*}_{i}: \cF^{*}_i\to \cF_{\bu}^{\flat*}$ in \eqref{eq: F bullet to F_i}.
Recall from \eqref{eq:T_i} that $T_i:= \Cofib(\frp^{*}_{i}: \cF^{*}_i\to \cF_{\bu}^{\flat*})$ is a torsion sheaf (pushed forward from a relative divisor $\iota \co Z \inj X_R'$ finite over $\Spec R$) supported set-theoretically on the union of the legs of $\cF_{\star}$, which by assumption are disjoint from $\wt{Q}_1$. Therefore $\frp^{*}_{i}$ restricts to an isomorphism in an open neighborhood of $|\wt{Q}_1|\times \Spec R\subset X'_{R}$, and hence induces an isomorphism in $\Mod_R$
\[
\RHom_{X_R'}(\cF_{\bu}^{\flat *}, \wt{Q}_1\ot R) \xrightarrow{\sim} 
\RHom_{X_R'}(\cF_i^*, \wt{Q}_1\ot R) 
\]
This being true for any $R$-point $\cF_{\star}$ of $\Hk^{r,\c}_{S}$, we conclude that $\wt b^{\c}_i \co \Hk_{V}^{\flat,\c} \rightarrow \wt{V}^{\c}_i$ is an isomorphism. 

The argument for $\wt b'^{\c}_{i}$ is similar.
\end{proof}

We define $\Sht_U^r, \Sht_V^r, \Sht_W^r$ to be the derived fibered products
\[
\begin{tikzcd}
\Sht_U^r \ar[r] \ar[d] & \Hk_U^\flat \ar[d, "{(a_0, a_r)}"] \\
U \ar[r, "{(\Id, \Frob)}"] & U \times U
\end{tikzcd}
\quad 
\begin{tikzcd}
\Sht_V^r \ar[r] \ar[d] & \Hk_V^\flat \ar[d, "{(b_0, b_r)}"] \\
V \ar[r, "{(\Id, \Frob)}"] & V \times V
\end{tikzcd}
\quad 
\begin{tikzcd}
\Sht_W^r \ar[r] \ar[d] & \Hk_W^\flat \ar[d, "{(c_0, c_r)}"] \\
W \ar[r, "{(\Id, \Frob)}"] & W \times W
\end{tikzcd}
\]
Following the notation above, we define $\Sht^{r,\c}_{U},\Sht^{r,\c}_{V}, \Sht^{r,\c}_{W}$ to be the respective base changes to $(X^\c)^r$. We define $\Sht_{U^\perp}^r$, $\Sht_{\wh V}^r$, and $\Sht_{W^\perp}^r$ analogously. 

\begin{cor}\label{c: rel Fq vs} The projection map $\Sht^{r,\c}_{V}\to \Sht^{r,\c}_{S}=\Sht^{r,\le \mu,\c}_{U(n)}$ is a relative finite-dimensional $\F_{q}$-vector space. Hence the theory of the arithmetic Fourier transform (\S \ref{sec: arithmetic FT}) applies to it. 
\end{cor}

\subsection{Calculation of traces}\label{ssec: calculation of traces}

Recall from Remark \ref{rem: derived hitchin stack} that the spaces $U_i$ from \S \ref{sec: descent} can be viewed as the derived fiber of the derived Hitchin stack $\sM_{H_1, H_2}$ from \cite[\S 5]{FYZ2} over $\{ \cE_1 \} \times S   \rightarrow \Bun_{\GL(m)'} \times \Bun_{U(n)}$, where $H_1 = \GL(m)'$ and $H_2 = U(n)$. Similarly, we explained in Remark \ref{rem: derived hecke stack} that $\Hk_U^\flat$ was an open substack of the derived fiber of the derived Hecke stack $\Hk_{\sM_{H_1, H_2}}^r$ from \cite[\S 5]{FYZ2} over $\{ \cE_1 \} \times S   \rightarrow \Bun_{\GL(m)'} \times \Bun_{U(n)}$. Therefore, $\Sht_U^r$ is equipped with an open embedding in $\Sht_{\sM_{H_1, H_2}}^r$, and in particular is of virtual dimension $d(a_r)$. We then have two natural cycles in $\mBM_{2d(a_r)}(\Sht_U^r )$:
\begin{enumerate}
\item The intrinsic derived fundamental class $[\Sht_U^r] \in \mBM_{2d(a_r)}(\Sht_U^r )$, which agrees with the restriction of the intrinsic derived fundamental class $[\Sht_{\sM_{H_1, H_2}}^r] \in \mBM_{2d(a_r)}(\Sht_{\sM_{H_1, H_2}}^r)$ along the \'{e}tale map $\Sht_U^r \rightarrow \Sht_{\sM_{H_1, H_2}}^r$. 
\item The trace of the cohomological correspondence $\cc_U$, denoted $\Tr^{\Sht}(\cc_U) \in \mBM_{2d(a_r)}(\Sht_U^r )$, see \S\ref{sssec: fix vs sht}. Here, $\Tr^{\Sht}$ is calculated using the Weil structure on $\cc_U$ coming from the canonical identification $\Frob^* \Qll{U} = \Qll{U}$. 
\end{enumerate}

We assemble the earlier results to calculate the trace of our cohomological correspondences. For any space $?$ over $k$, we equip $\Qll{?}$ with the natural Weil structure $\Frob^* \Qll{?} = \Qll{?}$. This equips all of our cohomological correspondences with a Weil structure. 

The assumptions \eqref{eq: vanishing assumption 1} imply that the maps $U_i \rightarrow S$ and $U_i^\perp \rightarrow S$ are smooth, hence $U_i$ and $U_i^\perp$ are smooth. Then by Proposition \ref{prop: smooth derived local terms}  we have 
\begin{equation}\label{eq: trace of cc_U}
\Tr^{\Sht}(\cc_U) = [\Sht_U^r] \in \mBM_{2d(a_r)}(\Sht_U^r).
\end{equation}
In particular, $\Sht_U^r$ is an open substack of $\Sht_{\cM_{\cE_1}}^r$, so it is quasi-smooth and $[\Sht_U^r]$ is the restriction of $[\Sht_{\cM_{\cE_1}}^r]$, which was called $[\cZ_{\cE_1}^r]$ in \cite{FYZ}. 

Similarly, we have
\begin{equation}\label{eq: trace of cc_U perp}
\Tr^{\Sht}(\cc_{U^\perp}) = [\Sht_{U^\perp}^r] \in \mBM_{2d(a_r^\perp)}(\Sht_{U^{\perp}}^r),
\end{equation}
where $\Sht_{U^\perp}^r$ is defined by the derived Cartesian square 
\[
\begin{tikzcd}
\Sht_{U^\perp}^r \ar[r] \ar[d] & \Hk_{U^\perp}^\flat \ar[d, "{(a_0^\perp, a_r^\perp)}"]\\
U_0^\perp \ar[r, "{(\Id, \Frob)}"] & U_0^\perp \times U_r^\perp
\end{tikzcd}
\]

Next, the assumptions \eqref{eq: vanishing assumption 2} imply that the maps $f_i \co U_i \rightarrow V_i$, $f \co \Hk_U^\flat \rightarrow \Hk_V^\flat$, $f_i^{\perp} \co U_i^{\perp} \rightarrow \wh{V}_i$, and $f^\perp \co \Hk_{U^{\perp}}^\flat \rightarrow \Hk_{\wh{V}}^\flat$ are all closed embeddings, in particular proper. Therefore, by Proposition \ref{prop: trace commutes with proper push}, we have
\begin{equation}\label{eq: trace f_!}
\Tr^{\Sht}(f_! \cc_U) = \Sht(f)_! \Tr^{\Sht}(\cc_U) \stackrel{\eqref{eq: trace of cc_U}}= \Sht(f)_! [\Sht_{U}^r] \in \mBM_{2d(a_r)} (\Sht_V^r),
\end{equation}
where we write $\Sht(f) := \Fix(f^{(1)}) \co \Sht_U^r \rightarrow \Sht_V^r$ for the map induced by taking fixed points of the Frobenius-twisted map of geometric correspondences $f^{(1)}$, and similarly for other maps of correspondences. We similarly have 
\begin{equation}\label{eq: trace f_! perp}
\Tr^{\Sht}(f_!^\perp \cc_{U^\perp}) = \Sht(f^\perp )_! \Tr^{\Sht}(\cc_{U^\perp}) \stackrel{\eqref{eq: trace of cc_U perp}}= \Sht(f^\perp )_! [\Sht_{U^\perp}^r] \in \mBM_{2d(a_r^\perp)} (\Sht_{\wh V}^r). 
\end{equation}

Recall from Corollary \ref{c: rel Fq vs} that $\Sht^{r,\c}_{V}\to \Sht^{r,\c}_{S}$ is a relative $\F_{q}$-vector space. We now relate the cycle classes \eqref{eq: trace f_!} and \eqref{eq: trace f_! perp} under the arithmetic Fourier transform on $\Sht^{r,\c}_{V}$ as defined in \S \ref{sec: arithmetic FT}. 

\begin{thm}\label{thm: AFT of distribution}  
We have 
\[
\begin{gathered}
\begin{aligned}
\FT^{\arith}_{\Sht_V^{r,\c}}(\Sht(f)^{\c}_! [\Sht_U^{r,\c}])
&=  (-1)^{d(U/S)+d(f_0)} q^{d(U/S)}
\cdot  \Sht(f^{\perp})^{\c}_! [\Sht_{U^{\perp}}^{r,\c}]\\
&\in \mBM_{2d(a_r)}(\Sht_{\wh{V}}^{r,\c}).
\end{aligned}
\end{gathered}
\]
Here $\Sht(f)^{\c}: \Sht_U^{r,\c}\to \Sht^{r,\c}_{V}$ is the restriction of $\Sht(f)$, and similarly for $\Sht(f^{\perp})^{\c}$.
\end{thm}

\begin{remark}
A priori, $\Sht(f^{\perp})^{\c}_! [\Sht_{U^{\perp}}^{r,\c}] $ lies in $\mBM_{2d(a_r^\perp)}(\Sht_{\wh V}^{r,\c})$. We note that $d(a_r) = d(a_r^\perp) = (n-m)r$, so the statement makes sense.  
\end{remark}

\begin{proof} We apply Theorem \ref{thm: trace compatible with FT} with $E = V$, $C_E = \Hk_V^{\flat,\c}$ and $\cc =( f_! \cc_U)|_{\Hk_V^{\flat,\c}}$.  Then Theorem \ref{thm: trace compatible with FT} tells us that 
\begin{equation}\label{eq: aft trace 1}
\FT^{\arith}_{\Sht_V^{r,\c}}\left( \Tr^{\Sht}(f_! \cc_U)|_{\Sht^{r,\c}_{V}} \right)= 
\left(\Tr^{\Sht} \FT_{\Hk_V^\flat}(f_! \cc_U)\right)|
_{\Sht^{r,\c}_{\wh V}} \in \mBM_{2d(a_r)}(\Sht_{\wh{V}}^{r,\c}).
\end{equation}
By Theorem \ref{thm: FT of pushforward correspondence} we have
\[
\FT_{\Hk_V^{\flat}}(f_{!} \cc_U)   = \TT_{[-d(U/S)-d(f_0)](-d(U/S))} (f^{\perp})_{!} \cc_{U^\perp} 
\]
Putting this into \eqref{eq: aft trace 1} and then taking the trace, using \eqref{Tr Sht sh tw}, \eqref{eq: trace f_!} and \eqref{eq: trace f_! perp}, yields the result. 
\end{proof}

\subsubsection{Test functions}\label{sssec: test functions} We introduce some notation for functions on $\Sht_{V}$ and $\Sht_{\wh{V}}$. The decompositions $\wt{Q}_1 := Q^* \oplus \cT_1^* \oplus \sigma^* \cT_2$ and $\wt{Q}_2 := \sigma^* Q \oplus \sigma^* \cT_1 \oplus \cT_2^*$ induce the following. 
\begin{itemize}
\item A decomposition 
\[
V \cong V^{(0)} \times_S V^{(1)} \times_S V^{(2)},
\]
where $V^{(0)} \rightarrow S$ is the vector bundle associated to $\ul{\RHom(\cF^*_{\univ}, Q^*)}$, $V^{(1)} \rightarrow S$ is the vector bundle associated to $\ul{\RHom(\cF^*_{\univ}, \cT_1^*)}$, and $V^{(2)} \rightarrow S$ is the vector bundle associated to $\ul{\RHom(\cF^*_{\univ}, \sigma^* \cT_2)}$. 
\item A decomposition 
\[
\wh{V} \cong  \wh{V}^{(0)} \times_S \wh{V}^{(1)} \times_S  \wh{V}^{(2)},
\]
where $\wh{V}^{(0)} \rightarrow S$ is the vector bundle associated to $\ul{\RHom(\cF^*_{\univ}, \sigma^* Q)}$, $\wh{V}^{(1)} \rightarrow S$ is the vector bundle associated to $\ul{\RHom(\cF^*_{\univ},  \sigma^* \cT_1)}$, and $\wh{V}^{(2)} \rightarrow S$ is the vector bundle associated to $\ul{\RHom(\cF^*_{\univ},  \cT_2^*)}$. As the notation suggests,  $\wh{V}^{(i)} \rightarrow S$ is the dual bundle to $V^{(i)} \rightarrow S$. As in \S \ref{sssec: self-duality of Q}, the Hermitian structures $h_{12}$ and $h_{21}$ on $Q$ induce two isomorphisms $V^{(0)} \cong \wh{V}^{(0)}$, which are the negatives of each other. 
\item A decomposition 
\[
\Hk_V^\flat \cong  \Hk_{V^{(0)}}^\flat \times_{\Hk_S^r} \Hk_{V^{(1)}}^\flat \times_{\Hk_S^r}  \Hk_{V^{(2)}}^\flat
\]
where $\Hk_{V^{(0)}}^\flat $ is the vector bundle associated to $\ul{\RHom(\cF_{\univ, \bu}^{\flat *}, Q^*)}$, etc. 
\item A decomposition 
\[
\Hk_{\wh{V}}^\flat \cong  \Hk_{\wh{V}^{(0)}}^\flat \times_{\Hk_S^r} \Hk_{\wh{V}^{(1)}}^\flat \times_{\Hk_S^r} \Hk_{\wh{V}^{(2)}}^\flat,
\]
where $\Hk_{\wh{V}^{(0)}}^\flat $ is the vector bundle associated to $\ul{\RHom(\cF_{\univ, \bu}^{\flat *}, \sigma^* Q)}$, etc. The Hermitian structures $h_{12}$ and $h_{21}$ on $Q$ induce two isomorphisms $\Hk_{V^{(0)}}^\flat \cong \Hk_{\wh{V}^{(0)}}^\flat $, which are the negatives of each other. 

By Lemma \ref{l: HkV Xc}, we see that $\Hk_{\wh{V}^{(i)}}^{\flat,\c} \rightarrow \Hk^{r,\c}_{S}$ is the dual bundle to $\Hk_{V^{(i)}}^{\flat,\c} \rightarrow \Hk^{r,\c}_{S}$. 

\item A decomposition 
\[
\Sht_V^r = \Sht_{V^{(0)}}^r \times_{\Sht_S^r} \Sht_{V^{(1)}}^r \times_{\Sht_S^r} \Sht_{V^{(2)}}^r
\]
where $\Sht_{V^{(i)}}^r $ is defined by the derived fibered product 
\[
\begin{tikzcd}
\Sht_{V^{(i)}}^r  \ar[r] \ar[d] & \Hk_{V^{(i)}}^\flat \ar[d, "{(b_0^{(i)}, b_r^{(i)})}"] \\
V^{(i)} \ar[r, "{(\Id, \Frob)}"] & V^{(i)}  \times V^{(i)} 
\end{tikzcd}
\]
\item A decomposition 
\[
\Sht_{\wh{V}}^r = \Sht_{\wh{V}^{(0)}}^r \times_{\Sht_S^r} \Sht_{\wh{V}^{(1)}}^r \times_{\Sht_S^r}\Sht_{\wh{V}^{(2)}}^r
\]
where $\Sht_{\wh{V}^{(i)}}^r $ is defined by the derived fibered product 
\[
\begin{tikzcd}
\Sht_{\wh{V}^{(i)}}^r  \ar[r] \ar[d] & \Hk_{\wh{V}^{(i)}}^\flat \ar[d, "{(\wt\b_0^{(i)}, \wt\b_r^{(i)})}"] \\
\wh{V}^{(i)} \ar[r, "{(\Id, \Frob)}"] & \wh{V}^{(i)}  \times \wh{V}^{(i)} 
\end{tikzcd}
\]
We note that $\Sht_{\wh{V}^{(i)}}^{r,\c}$ is dual to $\Sht_{V^{(i)}}^{r,\c}$ as $\F_{q}$-vector spaces over $\Sht^{r,\c}_{S}$ in the sense of \S \ref{sssec: AFT setup}. The Hermitian structures $h_{12}$ and $h_{21}$ on $Q$ induce two isomorphisms $\Sht_{V^{(0)}}^r \cong \Sht_{\wh{V}^{(0)}}^r$, which are the negatives of each other. 

\end{itemize}

We denote $\mf{q}_{12} \co \Sht^{r}_{V^{(0)}} \rightarrow \ul{\F_q}$ and $\mf{q}_{21} \co \Sht^{r}_{V^{(0)}} \rightarrow \ul{\F_q}$ the two quadratic forms induced by $h_{12}$ and $h_{21}$, respectively. Namely, $\mf{q}_{12}$ is the composition
\begin{equation}
\mf{q}_{12}: \Sht^{r}_{V^{(0)}}\xr{(\Id, h_{12})}\Sht^{r}_{V^{(0)}}\times_{\Sht^{r}_{S}}\Sht^{r}_{\wh V^{(0)}}\to \ul{\F_{q}},
\end{equation}
and similarly for $\mf{q}_{21}$. They satisfy $\mf{q}_{12} = - \mf{q}_{21}$. 
\begin{itemize}
\item We let $\mf{q}_{12}^* \psi$ be the pullback of $\psi$ to $\Sht_{V^{(0)}}^r$ via $\mf{q}_{12}$, and similarly for $\mf{q}_{21}$. Abusing notation, we will also use the same notation $\mf{q}^*_{12} \psi$ to denote its pullback to $\Sht_V^r$ and to $\Sht^{r,\c}_{\wh V}$. The meaning will be clear from context. 
\item We let $\delta_{\Sht_{V^{(i)}}^{r,\c}}$ be the indicator function of the zero-section of the relative $\F_{q}$-vector space $\Sht_{V^{(i)}}^{r,\c} \rightarrow \Sht_{S}^{r,\c}$. Abusing notation, we will also use this same notation to denote its pullback to $\Sht_V^{r,\c}$. 
\item We let $\bbm{1}_{\Sht_{V^{(i)}}^{r,\c}}$ be the constant function on $\Sht_{V^{(i)}}^{r,\c}$ with value $1$. Abusing notation, we will also use this same notation to denote its pullback to $\Sht_V^{r,\c}$. 
\item We let $(\mf{q}_{21}^* \psi \cdot \delta_{\Sht_{V^{(1)}}^{r,\c}} \cdot \bbm{1}_{\Sht_{V^{(2)}}^{r,\c}})$ denote the product of the above functions, viewed as a locally constant function on $\Sht_V^{r,\c}$.
\item We use similar notation on $\Sht_{\wh{V}}^{r,\c}$. 
\end{itemize}

\begin{lemma}
Let $d^{(i)}$ be the rank of $\Sht_{V^{(i)}}^{r,\c}$ as an $\F_q$-vector space over $\Sht_{S}^{r,\c}$. Then we have the following identities.

\begin{equation}\label{eq: AFT delta}
\FT_{\Sht_{V^{(i)}}^{r,\c}}^{\arith}(\delta_{\Sht_{V^{(i)}}^{r,\c}}) = (-1)^{d^{(i)}} \bbm{1}_{\Sht_{\wh{V}^{(i)}}^{r,\c}}
\end{equation}

\begin{equation}\label{eq: AFT constant}
\FT_{\Sht_{V^{(i)}}^{r,\c}}^{\arith}(\bbm{1}_{\Sht_{V^{(i)}}^{r,\c}}) = (-1)^{d^{(i)}} q^{d^{(i)}} \delta_{\Sht_{\wh{V}^{(i)}}^{r,\c}}
\end{equation}

\begin{equation}\label{eq: AFT Gaussian}
\begin{gathered}
\begin{aligned}
\FT_{\Sht_{V^{(0)}}^{r,\c}}^{\arith}(\mf{q}_{21}^* \psi)
&= (-1)^{d^{(0)}}[1/2]^* \mf{q}_{21}^* [-1]^* \psi\\
&\quad \cdot q^{d^{(0)}/2} \eta_{F'/F}(D_Q)^n
\end{aligned}
\end{gathered}
\end{equation}
where we recall that $D_Q\in \Div(X)$ denotes the divisor of $Q$.
\end{lemma}

\begin{proof}
As in the proof of Lemma \ref{lem: involutivity}, we can reduce to the case where $\Sht_{V^{(0)}}^{r,\c} \rightarrow \Sht_{S}^{r,\c}$ is split, and then to the usual finite Fourier transform, which we handled in \S \ref{sssec: finite FT}. 
\end{proof}

\begin{cor}\label{cor: FT of test functions}
Let $d = d^{(0)} +  d^{(1)} + d^{(2)}$ be the rank of $ \Sht_V^{r,\c}$ as an $\F_q$-vector space over $\Sht_S^{r,\c}$. Then we have 
\[
\begin{gathered}
\begin{aligned}
\FT^{\arith} (\mf{q}_{21}^*  \psi \cdot \delta_{\Sht^{r,\c}_{V^{(1)}}} \cdot \bbm{1}_{\Sht^{r,\c}_{V^{(2)}}})
&= (-1)^{d} q^{d^{(2)} + \frac{1}{2}d^{(0)}} \eta_{F'/F}(D_Q)^n\\
&\quad \cdot ([1/2]^* \mf{q}_{21}^* [-1]^* \psi
\cdot \bbm{1}_{\Sht_{\wh{V}^{(1)}}^{r,\c}} \cdot \delta_{\Sht_{\wh{V}^{(2)}}^{r,\c}})
\end{aligned}
\end{gathered}
\]
as functions on $\Sht_{\wh{V}}^{r,\c}$.  
\end{cor}

\begin{proof}Multiply the equations \eqref{eq: AFT Gaussian}, \eqref{eq: AFT constant}, \eqref{eq: AFT delta} together. 
\end{proof}

\subsubsection{Higher theta series} Here we relate $[\Sht^{r,\c}_U]$ and $[\Sht^{r,\c}_{U^\perp}]$ to the special cycles $[\cZ^{r,\c}_{\wt{\cE}_1}]$ and $[\cZ^{r,\c}_{\wt{\cE}_2}]$, and then to the higher theta series associated to $\wt{\cE}_1$ and $\wt{\cE}_2$.

Recall that for $i \in \{1,2\}$ we had exact sequences of coherent sheaves on $X'$, 
\begin{equation}\label{eq: saturate SES i}
\cE_i  \rightarrow  \wt{\cE}_i \rightarrow  \cT_i.
\end{equation}
This induces an exact triangle in $\Perf(R)$ for any $\cF \in \Bun_{U(n)}(R)$, 
\begin{equation}\label{eq: saturate exact triangle i}
\RHom_{X_R'}(\wt{\cE}_i, \cF)  \rightarrow \RHom_{X_R'}(\cE_i, \cF)  \rightarrow \RHom_{X_R'}(\cT_i, \cF[1]).
\end{equation}
By linear duality we have $\RHom_{X_R'}(\cT_i, \cF[1]) \cong \RHom_{X_R'} (\cF^*, \cT_i^*)$. Since $\cT_i^*$ is a torsion sheaf (concentrated in degree $0$ on the classical scheme $X'$), $\RHom_{X_R'} (\cF^*, \cT_i^*)$ is equivalent to a locally free coherent sheaf on $X_R'$.  Let $\cN_{\cT_i}$ be the total space of $\ul{\RHom(\cT_i, \cF_{\univ}[1])} \cong \ul{\RHom(\cF_{\univ}^*, \cT_i^*)}$, a vector bundle over $\Bun_{U(n)}$.

Below we refer to the notation of \cite[\S 5]{FYZ2}: we will use the pair $(H_1, H_2) = (\GL(m)', U(n))$. For $i \in \{1,2\}$, we let 
\begin{itemize}
\item $\cM_{\cE_i}$ be the derived fiber of the Hitchin stack $\cM_{H_1, H_2}$ from \cite[Definition 5.14]{FYZ2} over $\cE_i \in \Bun_{H_1}(k)$, and
\item $\Hk_{\cM_{\cE_i}}^r$ be the derived fiber of the derived Hitchin stack $\Hk_{\cM_{H_1, H_2}}^r$ over $\cE_i \in \Bun_{H_1}(k)$.
\end{itemize}
We define $\cM_{\wt{\cE}_i}$ and $\Hk_{\cM_{\wt{\cE}_i}}^r$ similarly. 

Then the exact triangle \eqref{eq: saturate exact triangle i} corresponds at the level of total spaces to a derived Cartesian square 
\begin{equation}\label{eq: cartesian square 1}
\begin{tikzcd}
\cM_{\wt{\cE}_i}  \ar[r] \ar[d] & \cM_{\cE_i} \ar[d] \\
\Bun_{U(n)} \ar[r, "z"] & \cN_{\cT_i}
\end{tikzcd}
\end{equation}

Let $\Hk_{\cN_{\cT_i}}^r$ be the total space of $\ul{\RHom(\cT_i, \cF_{\univ, \bu}^\flat[1])}$ on $\Hk_{U(n)}^r$. The exact triangle 
\[
\ul{\RHom(\wt{\cE}_i, \cF_{\univ, \bu}^\flat)}  \rightarrow \ul{\RHom(\cE_i, \cF_{\univ, \bu}^\flat)}  \rightarrow \ul{\RHom(\cT_i, \cF_{\univ, \bu}^\flat[1])}
\]
corresponds at the level of total spaces to a derived Cartesian square 
\begin{equation}\label{eq: cartesian square 2}
\begin{tikzcd}
\Hk_{\cM_{\wt{\cE}_i}}^r  \ar[r] \ar[d]  & \Hk_{\cM_{\cE_i}}^r \ar[d] \\
\Hk_{U(n)}^r  \ar[r, "z"] & \Hk_{\cN_{\cT_i}}^r
\end{tikzcd}
\end{equation}
Consider the commutative diagram 
\begin{equation}\label{eq: cartesian square 3}
\begin{tikzcd}[column sep = huge]
\Hk_{\cM_{\cE_i}}^r  \ar[r, "{(h_0, h_r)}"]  \ar[d] & \cM_{\cE_i} \times \cM_{\cE_i}    \ar[d] &   \ar[l, "{(\Id, \Frob)}"']  \cM_{\cE_i}  \ar[d] \\
 \Hk_{\cN_{\cT_i}}^r  \ar[r, "{(h_0, h_r)}"]  & \cN_{\cT_i} \times \cN_{\cT_i} & \cN_{\cT_i}  \ar[l, "{(\Id, \Frob)}"']  \\
 \Hk_{U(n)}^r   \ar[u]  \ar[r, "{(h_0, h_r)}"] & \Bun_{U(n)} \times  \Bun_{U(n)} \ar[u] &  \Bun_{U(n)}  \ar[u] \ar[l, "{(\Id, \Frob)}"'] 
\end{tikzcd}
\end{equation}
By \eqref{eq: cartesian square 1} and \eqref{eq: cartesian square 2}, the derived fibered product of the columns of \eqref{eq: cartesian square 3} is 
\begin{equation}\label{eq: column fibered product}
\begin{tikzcd}[column sep = huge] 
\Hk_{\cM_{\wt{\cE}_i}}^r  \ar[r, "{(h_0, h_r)}"]   & \cM_{\wt{\cE}_i}    \times  \cM_{\wt{\cE}_i}& \ar[l]  \cM_{\wt{\cE}_i}   \ar[l, "{(\Id, \Frob)}"']
\end{tikzcd}
\end{equation}
By definition, the derived fibered product of the rows of \eqref{eq: cartesian square 3} is
\begin{equation}\label{eq: row fibered product}
\begin{tikzcd}
\Sht_{\cM_{\cE_i}}^r=\cZ^{r}_{\cE_{i}} \ar[d] \\
\Sht_{\cN_{\cT_i}}^r \\
\Sht_{U(n)}^r \ar[u] 
\end{tikzcd}
\end{equation}
In turn, the derived fibered products of \eqref{eq: column fibered product} and \eqref{eq: row fibered product} are canonically identified by the same proof as for \cite[Lemma A.9]{YZ}. This shows:

\begin{cor}\label{cor: cartesian square 4}
The commutative square 
\[
\begin{tikzcd}
\cZ^{r}_{\wt{\cE}_{i}}  \ar[r] \ar[d] & \cZ_{\cE_i}^r  \ar[d]  \\
\Sht_{U(n)}^r \ar[r, "0"] & \Sht_{\cN_{\cT_i}}^r
\end{tikzcd}
\]
is derived Cartesian. 
\end{cor}

To compare with earlier objects: 
\begin{itemize}
\item Restricting $\Sht_{\cN_{\cT_1}}^r \rightarrow \Sht_{U(n)}^r$ along the open embedding $\Sht_S^r \inj \Sht_{U(n)}^r$ recovers $\Sht_{V^{(1)}}^r \rightarrow \Sht_S^r$, i.e., we have a derived Cartesian square 
\begin{equation}\label{eq: cartesian square 5}
\begin{tikzcd}
\Sht_{V^{(1)}}^r \ar[r, hook, "\text{open}"] \ar[d] &   \Sht_{\cN_{\cT_1}}^r  \ar[d] \\
\Sht_S^r  \ar[r, hook, "\text{open}"] & \Sht^r_{U(n)}
\end{tikzcd}
\end{equation}
Restricting $\Sht_{\cN_{\cT_2}}^r \rightarrow \Sht_{U(n)}^r$ along the open embedding $\Sht_S^r \inj \Sht_{U(n)}^r$ recovers $\Sht_{\wh{V}^{(2)}}^r \rightarrow \Sht_S^r$. 

\item Restricting $ \cZ^r_{\cE_{1}} \rightarrow \Sht_{U(n)}^r$ along the open embedding $\Sht_S^r \inj \Sht_{U(n)}^r$ recovers $\Sht_{U}^r \rightarrow \Sht_S^r$, and restricting $\cZ_{\cE_{2}}^{r} \rightarrow \Sht_{U(n)}^r$ along the open embedding $\Sht_S^r \inj \Sht_{U(n)}^r$ recovers $\Sht_{U^\perp}^r \rightarrow \Sht_S^r$. In other words, we have derived Cartesian squares
\begin{equation}\label{eq: cartesian square 6}
\begin{tikzcd}
\Sht_{U}^r \ar[r, hook, "\text{open}"] \ar[d] &  \cZ_{\cE_{1}}^r  \ar[d] \\
\Sht_S^r  \ar[r, hook, "\text{open}"] & \Sht^r_{U(n)} 
\end{tikzcd}
\hspace{1cm}
\begin{tikzcd}
\Sht_{U^\perp}^r \ar[r, hook, "\text{open}"] \ar[d] &  \cZ_{\cE_{2}}^r  \ar[d] \\
\Sht_S^r  \ar[r, hook, "\text{open}"] & \Sht^r_{U(n)} 
\end{tikzcd}
\end{equation}
\end{itemize}

Abbreviate $\cZ_{\wt{\cE}_{i}}^{r, \leq \mu}$ for $\cZ^r_{\wt{\cE}_{i}} \times_{\Sht_{U(n)}^r} \Sht_S^r$.
We have derived Cartesian squares
\begin{equation}\label{eq: cartesian square 8}
\begin{tikzcd}
\cZ_{\wt{\cE}_1}^{r, \leq \mu} \ar[r, hook] \ar[d] &  \Sht_{U}^r  \ar[d] \\
\Sht_S^r  \ar[r, hook, "0"] & \Sht^r_{V^{(1)}} 
\end{tikzcd}
\hspace{1cm}
\begin{tikzcd}
\cZ_{\wt{\cE}_2}^{r, \leq \mu} \ar[r, hook] \ar[d] &  \Sht_{U^\perp}^r  \ar[d] \\
\Sht_S^r  \ar[r, hook, "0"] & \Sht^r_{\wh{V}^{(2)}}
\end{tikzcd}
\end{equation}

By Corollary \ref{c: rel Fq vs}, $\Sht_{V^{(i)}}^{r,\c}\rightarrow \Sht_{S}^{r,\c} $ is finite \'{e}tale, so the zero section $\Sht_{S}^{r,\c} \xrightarrow{0} \Sht_{V^{(i)}}^{r,\c}$ is open-closed, which by \eqref{eq: cartesian square 8} implies that $\cZ_{\wt{\cE}_1}^{r, \leq \mu,\c} =\cZ_{\wt{\cE}_1}^{r}|_{\Sht^{r,\c}_{S}}\inj   \Sht_{U}^{r,\c}$ is also open-closed; similarly for $\cZ_{\wt{\cE}_2}^{r, \leq \mu,\c} \inj   \Sht_{U^{\perp}}^{r,\c}$.

\begin{cor}\label{cor: cut down fundamental classes}  (1)  Viewing $\delta_{\Sht_{V^{(1)}}^{r,\c}} \in \rH^0(\Sht_V^{r,\c}; \Ql)$ and $\Sht(f)^{\c}_![\Sht_U^{r,\c}] \in \mBM_{2d(a_r)}(\Sht_V^{r,\c})$, we have 
\[
\Sht(f)^{\c}_! [\Sht_U^{r,\c}]  \cdot \delta_{\Sht_{V^{(1)}}^{r,\c}}  = [\cZ_{\wt{\cE}_{1}}^{r, \leq \mu,\c}] \in \mBM_{2d(a_r)} (\Sht_V^{r,\c}).
\]

(2) Viewing $\delta_{\Sht_{\wh{V}^{(2)}}^{r,\c}} \in \rH^0(\Sht_{\wh{V}}^{r,\c}; \Ql)$ and $\Sht(f^{\perp})^{\c}_![\Sht_{U^{\perp}}^{r,\c}] \in \mBM_{2d(a_r^\perp)}(\Sht_{\wh V}^{r,\c})$, we have 
\[
\Sht(f^{\perp})^{\c}_![\Sht_{U^{\perp}}^{r,\c}] \cdot \delta_{\Sht_{\wh{V}^{(2)}}^{r,\c}}  = [\cZ_{\wt{\cE}_{2}}^{r,\le\mu,\c}] \in \mBM_{2d(a_r)} (\Sht_{\wh{V}}^{r,\c}).
\]
\end{cor}

\begin{lemma}\label{lem: theta series as product} For $i=1,2$, let $\cA_{\wt{\cE}_i}$ be the Hitchin base as in \cite[\S 3.3]{FYZ2}. For $a\in \cA_{\wt{\cE}_{i}}(k)$, let $\cZ^{r, \leq \mu,\c}_{\wt{\cE}_i}(a)=\cZ^{r}_{\wt{\cE}_{i}}(a)\times_{\Sht^{r}_{U(n)}}\Sht^{r,\c}_{S}$.

(1) We have an equality in $\mBM_{2d(a_r)} (\Sht_V^{r,\c})$:
\[
\begin{gathered}
\begin{aligned}
\Sht(f)^{\c}_! [\Sht_{U}^{r,\c}]
\cdot  (\mf{q}_{21}^* \psi \cdot  \delta_{\Sht_{V^{(1)}}^{r,\c}} \cdot \bbm{1}_{\Sht_{V^{(2)}}^{r,\c}})
&=  \sum_{a \in \cA_{\wt{\cE}_1}(k)}
\psi(\langle e_{\cG,\wt{\cE}_1}, a \rangle )\\
&\quad \cdot \Sht(f)^{\c}_! [\cZ^{r, \leq \mu,\c}_{\wt{\cE}_{1}}(a)].
\end{aligned}
\end{gathered}
\]

(2) We have an equality in $\mBM_{2d(a_r)} (\Sht_{\wh{V}}^{r,\c})$:
\[
\begin{gathered}
\begin{aligned}
\Sht(f^{\perp})^{\c}_! [\Sht_{U^{\perp}}^{r, \c}]
\cdot (\mf{q}_{12}^* \psi  \cdot \bbm{1}_{\Sht_{\wh{V}^{(1)}}^{r,\c}}
\cdot \delta_{\Sht_{\wh{V}^{(2)}}^{r,\c}})
&=   \sum_{a \in \cA_{\wt{\cE}_2}(k)}
\psi(\langle e_{\cG,\wt{\cE}_2}, a \rangle )\\
&\quad \cdot \Sht(f^{\perp})^{\c}_! [\cZ^{r, \leq \mu,\c}_{\wt{\cE}_{2}}(a)].
\end{aligned}
\end{gathered}
\]

\end{lemma}

\begin{proof} 
(1) By Corollary \ref{cor: cut down fundamental classes}, we have
\begin{equation}\label{eq:10-2-7-eq1}
\Sht(f)^{\c}_! [\Sht_{U}^{r,\c}] \cdot  (\bbm{1}_{\Sht_{V^{(0)}}^{r,\c}} \cdot  \delta_{\Sht_{V^{(1)}}^{r,\c}} \cdot \bbm{1}_{\Sht_{V^{(2)}}^{r,\c}}) =  \sum_{a \in \cA_{\wt{\cE}_1}(k)}  \Sht(f)^{\c}_!  [\cZ^{r, \leq \mu,\c}_{\wt{\cE}_1}(a)].
\end{equation}
In the displayed equation of (1), the LHS is obtained by multiplying the LHS of \eqref{eq:10-2-7-eq1} by $\mf{q}_{21}^* \psi$. Here \eqref{eGE1 q12} and \eqref{q1221} are applied to the saturated transverse pair $(\wt{\cE}_1,\wt{\cE}_2)$ and the quotient $Q$ in \eqref{eq: Q diagram for saturation}. Observe that by \eqref{eGE1 q12} and \eqref{q1221}, the function $\mf{q}_{21}^* \psi$ on $\Sht_V^{r}$ coincides with the one sending $(\cF_\star, t) \in \Sht_V^r(R)$ to $\psi(\langle e_{\cG,\wt{\cE}_1}, a(t) \rangle)$. Thus $\mf{q}_{21}^* \psi$ factors over the map $\Sht_V^{r,\c} \rightarrow \cA_{\wt{\cE}_1}(k)$, as the map $a \mapsto \psi(\langle e_{\cG,\wt{\cE}_1}, a \rangle)$. Multiplying the RHS of \eqref{eq:10-2-7-eq1} by this function gives the result.

(2) The argument is similar, applying \eqref{eGE2 q21} to the same saturated pair. 
\end{proof}

\subsection[Conclusion of the proof of transverse modularity]{Conclusion of the proof of Theorem \ref{thm: transverse modularity}}\label{ssec: conclusion}

Let $d, d^{(i)}$, for $i \in \{0,1,2\}$, be as in \S \ref{sssec: test functions}. Note that $d^{(i)}$ is also the rank of $V^{(i)}$ as a vector bundle over $S$. 

By Lemma \ref{lem: theta series as product} we have
\begin{equation}\label{eq: conclusion -1}
\begin{gathered}
\begin{aligned}
\wt{Z}_m^r(\cG, \wt{\cE}_1) |_{\Sht^{r,\c}_{S}}
&= \chi(\det \wt{\cE}_1) q^{n (\deg \wt{\cE}_1  - \deg \omega_X)/2}\\
&\quad \cdot \langle \Sht(f)^{\c}_! [\Sht_U^{r,\c}],
\mf{q}_{21}^* \psi\cdot \delta_{\Sht_{V^{(1)}}^{r,\c}}
\cdot \bbm{1}_{\Sht_{V^{(2)}}^{r,\c}} \rangle
\end{aligned}
\end{gathered}
\end{equation}
and
\begin{equation}\label{eq: conclusion 0}
\begin{gathered}
\begin{aligned}
\wt{Z}_m^r(\cG, \wt{\cE}_2)|_{\Sht^{r,\c}_{S}}
&= \chi(\det \wt{\cE}_2) q^{n (\deg \wt{\cE}_2  - \deg \omega_X)/2}\\
&\quad \cdot \langle \Sht(f^\perp)^{\c}_! [\Sht_{U^\perp}^{r,\c}],
\mf{q}_{12}^*  \psi \cdot \bbm{1}_{\Sht_{\wh{V}^{(1)}}^{r,\c}}
\cdot \delta_{\Sht_{\wh{V}^{(2)}}^{r,\c}} \rangle.
\end{aligned}
\end{gathered}
\end{equation}

By the Plancherel formula of Lemma \ref{lem: Plancherel for cycles} and the near-involutivity of $\FT^{\arith}$ of Lemma \ref{lem: involutivity}, we have 
\begin{eqnarray}\label{eq: conclusion 1}
&& \langle \Sht(f)^{\c}_! [\Sht_U^{r,\c}], \mf{q}_{21}^*  \psi \cdot \delta_{\Sht_{V^{(1)}}^{r,\c}} \cdot \bbm{1}_{\Sht_{V^{(2)}}^{r,\c}} \rangle \\
\notag &=& \frac{1}{q^d} \langle \FT_{\Sht_V^{r,\c}}^{\arith} ( \Sht(f)^{\c}_! [\Sht_U^{r,\c}]) , \FT_{\Sht_V^{r,\c}}^{\arith}  \bigl( [-1]^*(\mf{q}_{21}^*  \psi \cdot \delta_{\Sht_{V^{(1)}}^{r,\c}} \cdot \bbm{1}_{\Sht_{V^{(2)}}^{r,\c}}) \bigr)\rangle.
\end{eqnarray}
Here $[-1]^*$ acts trivially on the second factor because $\mf q_{21}$ is quadratic, the zero section is fixed by $[-1]$, and $\bbm{1}_{\Sht_{V^{(2)}}^{r,\c}}$ is constant.
Then we use Theorem \ref{thm: AFT of distribution} and Corollary \ref{cor: FT of test functions} to rewrite the right side of  \eqref{eq: conclusion 1} as 
\begin{eqnarray}\label{eq: conclusion 2}
   \frac{1}{q^d} q^{d(U/S)} (-1)^{d(U/S) + d(f_0)} (-1)^d q^{d^{(2)}+  \frac{1}{2}d^{(0)}} \eta_{F'/F}(D_Q)^n \\
\hspace{1cm} \cdot  \notag \langle \Sht(f^\perp)^{\c}_! [\Sht_{U^\perp}^{r,\c}] , [1/2]^*\mf{q}_{21}^* [-1]^* \psi \cdot \bbm{1}_{\Sht_{\wh{V}^{(1)}}^{r,\c}} \cdot \delta_{\Sht_{\wh{V}^{(2)}}^{r,\c}} \rangle.
\end{eqnarray}
Since $\mf{q}_{21} =- \mf{q}_{12}$, we have $\mf{q}_{21}^* [-1]^* = \mf{q}_{12}^*$. Eliminating the $[1/2]^*$ does not affect the expression since $\Sht(f^\perp)^{\c}_! [\Sht_{U^\perp}^{r,\c}]$ is invariant under the scaling $\F_q^\times$ action on $\Sht_{\wh V}^{r,\c}$. Then clearly \eqref{eq: conclusion 2} agrees with \eqref{eq: conclusion 0} up to sign and integral power of $q$. Therefore, it remains to check the sign and the exponent of $q$. The exponent of $q$ in \eqref{eq: conclusion 2} is
\begin{equation}\label{eq: conclusion 3}
-d + d(U/S) + d^{(2)} + \frac{1}{2} d^{(0)} = -d^{(1)} - d^{(0)} + d(U/S) +\frac{1}{2} d^{(0)}.
\end{equation}
Recall that $d^{(1)} := \rank(\ul{\RHom(\cF^*, \cT_1^*)})$. By the exact triangle of perfect complexes over $S$
\[
\ul{\RHom(\cF^*_{\univ}, \wt{\cE}_1^*)} \rightarrow 
\ul{\RHom(\cF^*_{\univ}, \cE_1^*)} \rightarrow \ul{\RHom(\cF^*_{\univ}, \cT_1^*)}
\]
we have $-d^{(1)} = - \rank (\ul{\RHom(\cF^*_{\univ}, \cE_1^*)}) + \rank (\ul{\RHom(\cF^*_{\univ}, \wt{\cE}_1^*)} )$ and $\rank (\ul{\RHom(\cF^*_{\univ}, \cE_1^*)})   = d(U/S)$. Putting all this into \eqref{eq: conclusion 3} simplifies the exponent of $q$ to 
\begin{equation}\label{eq: conclusion 4}
 \rank (\ul{\RHom(\cF^*_{\univ}, \wt{\cE}_1^*)} ) -d^{(0)}  + \frac{1}{2} d^{(0)}.
\end{equation}
By the exact triangle 
\[
\ul{\RHom(\cF^*_{\univ}, \sigma^* \wt{\cE}_2)}  \rightarrow \ul{\RHom(\cF^*_{\univ}, \wt{\cE}_1^*)} \rightarrow \ul{\RHom(\cF^*_{\univ}, Q^*)}
\]
we have 
\begin{equation}\label{eq: conclusion 5}
\rank(\ul{\RHom(\cF^*_{\univ}, \wt{\cE}_1^*)} ) = \rank(\ul{\RHom(\cF^*_{\univ}, \sigma^* \wt{\cE}_2)} ) + \rank(\ul{\RHom(\cF^*_{\univ}, Q^*)}).
\end{equation}
Putting \eqref{eq: conclusion 5} into \eqref{eq: conclusion 4}, noting that $\rank(\ul{\RHom(\cF^*_{\univ}, Q^*)})= d^{(0)}$, the exponent of $q$ in \eqref{eq: conclusion 2} simplifies to 
\begin{equation}\label{eq: conclusion 6}
\rank(\ul{\RHom(\cF^*_{\univ}, \sigma^* \wt{\cE}_2)} )  + \frac{1}{2} d^{(0)}. 
\end{equation}
In \S \ref{sssec: conclusion r=0} we exactly showed that 
\[
n (\deg \wt{\cE}_1  - \deg \omega_X)/2 + \rank(\ul{\RHom(\cF^*_{\univ}, \sigma^* \wt{\cE}_2)} )  + \frac{1}{2} d^{(0)} =  n (\deg \wt{\cE}_2  - \deg \omega_X)/2 ,
\]
which after comparing \eqref{eq: conclusion -1}, \eqref{eq: conclusion 0} shows that the exponents match! Hence we have established that
\[
\wt{Z}_m^r(\cG, \wt{\cE}_1)|_{\Sht^{r,\c}_{S}}  = \pm \wt{Z}_m^r(\cG, \wt{\cE}_2)|_{\Sht^{r,\c}_{S}}.
\]

Finally, we match the signs. Let
\[
e:=n(\deg\wt{\cE}_1-\deg\omega_X)/2-d+d(U/S)+d^{(2)}+\frac{1}{2}d^{(0)}
=n(\deg\wt{\cE}_2-\deg\omega_X)/2,
\]
where the equality is the exponent comparison proved above. We found above that 
\begin{eqnarray}\label{eq: conclusion 7}
\wt{Z}_m^r(\cG, \wt{\cE}_1)|_{\Sht^{r,\c}_{S}} &=& \chi(\det \wt{\cE}_1) \eta_{F'/F}(D_Q)^n  (-1)^{d(U/S) + d(f_0)} (-1)^d q^e \\
\notag&&\langle \Sht(f^\perp)^{\c}_! [\Sht_{U^\perp}^{r,\c}], \mf{q}_{12}^*  \psi \cdot \bbm{1}_{\Sht_{\wh{V}^{(1)}}^{r,\c}} \cdot \delta_{\Sht_{\wh{V}^{(2)}}^{r,\c}}\rangle
\end{eqnarray}
while
\begin{equation}\label{eq: conclusion 8}
\begin{aligned}
\wt{Z}_m^r(\cG, \wt{\cE}_2) |_{\Sht^{r,\c}_{S}}
&= \chi(\det \wt{\cE}_2) q^e\\
&\quad \cdot \langle \Sht(f^\perp)^{\c}_! [\Sht_{U^\perp}^{r,\c}], \mf{q}_{12}^*  \psi\cdot \bbm{1}_{\Sht_{\wh{V}^{(1)}}^{r,\c}} \cdot \delta_{\Sht_{\wh{V}^{(2)}}^{r,\c}} \rangle,
\end{aligned}
\end{equation}
Next, we note that $d(U/S) - d(f_0) = d(V/S)=d$, hence $d(U/S) + d(f_0)+d\equiv0 \pmod{2}$. Furthermore, in \S \ref{sssec: conclusion r=0} we calculated that $\chi(\det \wt{\cE}_1) \eta_{F'/F}(D_Q)^n =  \chi(\det \wt{\cE}_2)$. Therefore, \eqref{eq: conclusion 7} and \eqref{eq: conclusion 8} are equal. \qed

\appendix 
\part{Appendices}

\section{Derived Fourier analysis: proofs}\label{app: A}

This appendix justifies assertions made in \S \ref{sec: FT} about the derived Fourier transform.

\textbf{We remind the reader of the notational conventions introduced in \S \ref{sec: FT} regarding derived vector bundles:} for perfect complexes denoted with calligraphic letters such as $\cE, \cE'$, etc., the corresponding Roman letters such as $E, E'$, etc. denote their associated total spaces.

\subsection{Preliminaries}
We record here some general facts about functors on sheaf categories induced by kernel sheaves. In particular, we introduce a configuration that we call a ``butterfly'', which induces natural transformations between such functors. 

\subsubsection{Correspondences}\label{sssec: correspondences} We begin with some preliminaries on correspondences. Given a correspondence
\[
\begin{tikzcd}
& C \ar[dl, "c_0"'] \ar[dr, "c_1"] \\
S & & T
\end{tikzcd}
\]
we associate the functor 
$$\Phi_{C}:=c_{1!} c_0^* \co D(S) \rightarrow D(T).$$
Consider the fibered product of correspondences 
\[
\begin{tikzcd}
& & C \times_T D \ar[dl] \ar[dr] \\ 
& C \ar[dl, "c_0"'] \ar[dr, "c_1"] & & D \ar[dl, "d_0"'] \ar[dr, "d_1"]  \\
S & & T & & U
\end{tikzcd}
\]
Then by proper base change, we have a natural isomorphism of functors 
\begin{equation}
\Phi_{C\times _{T}D}\cong \Phi_{D}\c \Phi_{C}: D(S)\to D(U).
\end{equation}

More generally, for any $\cK\in D(C)$, we have a functor
\begin{equation}
\Phi_{C, \cK}: c_{1!}(c_{0}^{*}(-)\ot \cK)\co D(S) \rightarrow D(T)
\end{equation}
using $\cK$ as the kernel sheaf. The functor considered in the previous paragraph corresponds to the special case $\cK=\Qll{C}$. By the projection formula, we have a natural isomorphism of functors $D(S)\to D(T)$
\begin{equation}
\Phi_{C, \cK}\cong \Phi_{S\times T, c_{!}\cK}
\end{equation}
where $c=(c_{0},c_{1}): C\to S\times T$.

What is a general mechanism to construct natural transformations between such functors?  For this we introduce the notion of a ``butterfly'', which is essentially a correspondence between correspondences.

\subsubsection{About butterflies}\label{sssec: butterflies} Suppose we have a commutative diagram of derived stacks
\begin{equation}\label{eq: butterfly nt}
\begin{tikzcd}
& & F \ar[dl, "u"'] \ar[dr, "v"] \\
& E  \ar[ddl, "\pr_S^E"'] \ar[ddrrr, "\pr_T^E"']  & & H \ar[ddr, "\pr_T^H"] \ar[ddlll, "\pr_S^H"]  \\
\\
S & & & & T
\end{tikzcd}
\end{equation}
where $u$ is proper, $v$ is quasi-smooth, and all morphisms are separated. We will refer to such a diagram as a \emph{butterfly} (we codify this notion because many butterflies will come up). We also let $\pr_S^F \co F \rightarrow S$ and $\pr_T^F \co F \rightarrow T$ be the obvious maps. We will construct a natural transformation
\begin{equation}\label{eq: nt by butterfly}
\Phi_{E}=\pr_{T!}^E \pr_S^{E*} \xrightarrow{\bigstar} \Phi_{H}\tw{-d(v)}=\pr_{T!}^H \pr_S^{H*} \tw{-d(v)}
\end{equation}
of functors $D(S) \rightarrow D(T)$, as the composition of the following maps. 
\begin{enumerate}
\item The unit for $(u^*, u_*)$ induces  
\[
\pr_{T!}^E \pr_S^{E*} \rightarrow \pr_{T!}^E u_* u^* \pr_S^{E*}
\]
\item Since $u_! = u_*$ because $u$ is proper, we have identifications 
\[
 \pr_{T!}^E u_* u^* \pr_S^{E*} =  \pr_{T!}^E u_! u^* \pr_S^{E*}  = \pr_{T!}^F \pr_{S}^{F*}  = \pr_{T!}^H v_! v^* \pr_S^{H*}.
 \]
 \item Since $v$ is quasi-smooth, its relative fundamental class induces (as explained in \S \ref{ssec: relative fundamental nt}) a natural transformation $[v] \co v^* \rightarrow v^!\tw{-d(v)}$, which gives a natural transformation
 \[
 \pr_{T!}^H v_! v^* \pr_S^{H*} \rightarrow \pr_{T!}^H v_! v^! \tw{-d(v)} \pr_S^{H*} .
 \]
 \item The counit for $(v_!, v^!)$ gives a natural transformation
 \[
 \pr_{T!}^H v_! v^! \tw{-d(v)} \pr_S^{H*}  \rightarrow \pr_{T!}^H \pr_S^{H*} \tw{-d(v)}.
 \]
\end{enumerate}

\begin{remark}\label{r: butterfly and cc} Let $\pr^{E}=(\pr^{E}_{S}, \pr^{E}_{T}): E\to S\times T$, and similarly define $\pr^{H}$ and $\pr^{F}$.  A butterfly of the form \eqref{eq: butterfly nt} can be viewed as a correspondence $F$ between $E$ and $H$ over $S\times T$
\begin{equation}
\xymatrix{ & F\ar[dl]_{u}\ar[dr]^{v}\\
E\ar[dr]_{\pr^{E}} && H\ar[dl]^{\pr^{H}}\\
& S\times T
}
\end{equation}
Assume that the morphisms in the butterfly are representable in derived schemes, in addition to being separated, so that we may use the pushforward formalism for cohomological correspondences from \S \ref{ssec: pushforward functoriality for CC}. The fundamental class of $v$ gives a cohomological correspondence
\begin{equation}
\cc_{v}\in \Corr_{F}(\Qll{E}, \Qll{H}\tw{-d(v)}).
\end{equation}
The assumption that $u$ is proper implies that the map of correspondences $(\pr^{E}, \pr^{F}, \pr^{H}): (E\leftarrow{}F\rightarrow{H} )\to (S\times T=S\times T=S\times T)$ is left pushable. Therefore
\begin{align*}
\pr^{F}_{!}:  \Corr_{F}(\Qll{E}, \Qll{H}\tw{-d(v)})&\to  \Corr_{S\times T}(\pr^{E}_{!}\Qll{E}, \pr^{H}_{!}\Qll{H}\tw{-d(v)})\\
&=\Hom_{S\times T}(\pr^{E}_{!}\Qll{E}, \pr^{H}_{!}\Qll{H}\tw{-d(v)})
\end{align*}
is defined. In particular we get a map
\begin{equation}
\pr^{F}_{!}\cc_{v}: \pr^{E}_{!}\Qll{E}\to \pr^{H}_{!}\Qll{H}\tw{-d(v)}
\end{equation}
which induces a natural transformation
\begin{equation}
\Phi_{E}=\Phi_{S\times T, \pr^{E}_{!}\Qll{E}}\to \Phi_{S\times T, \pr^{H}_{!}\Qll{H}\tw{-d(v)}}=\Phi_{H}\tw{-d(v)}.
\end{equation}
Unwinding the definitions we see that this construction recovers the map $\bigstar$ in \eqref{eq: nt by butterfly}.
\end{remark}


Next we establish some compatibilities between natural transformations induced by different butterflies.

\subsubsection{Butterfly and pushforward}\label{sssec: butterfly and pushforward} Consider the butterfly \eqref{eq: butterfly nt}. Suppose we are given a map $f \co S'\rightarrow S$. Then the construction in \S\ref{sssec: butterflies} gives a natural transformation
\begin{equation}\label{eq: app push butterfly 1}
\Phi_{E}f_{!}\xrightarrow{\bigstar f_!} \Phi_{H}f_{!}\tw{-d(v)}
\end{equation}
of functors $D(S') \rightarrow D(T)$. 

On the other hand, we have another natural transformation of the same form coming from another butterfly. Consider the diagram where the superscript $'$ means pullback along $f$. 
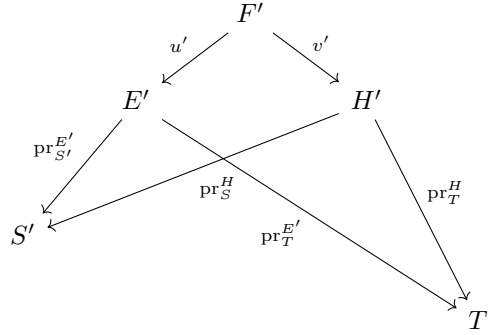
\[
\begin{tikzcd}
& & F'  \ar[dl, "u'"'] \ar[dr, "v'"]  \\
& E'  \ar[ddl, "\pr_{S'}^{E'}"'] \ar[dddrrr, "\pr_T^{E'}"']  &  & H' \ar[dddr, "\pr_T^{H'}"] \ar[ddlll, "\pr_S^{H'}"]  \\
&  & &   \\
S'\\	
 & & & & T
\end{tikzcd}
\]
This is also a butterfly ($u'$ is proper and $v'$ is quasi-smooth) and therefore gives another natural transformation 
\begin{equation}\label{eq: app push butterfly 2}
\Phi_{E'}\xrightarrow{\bigstar'} \Phi_{H'}\tw{-d(v')}=\Phi_{H'}\tw{-d(v)}
\end{equation}
of functors $D(S') \rightarrow D(T)$. 

The Cartesian squares  
\[
\begin{tikzcd}
S' \ar[d, "f"'] & E' \ar[d] \ar[l, "\pr_{S'}^{E'}"']  \ar[d, "f^E"] & S' \ar[d, "f"'] & H' \ar[d] \ar[l, "\pr_{S'}^{H'}"'] \ar[d, "f^H"]\\
S & E \ar[l, "\pr_S^E"]&  S & H \ar[l, "\pr_S^H"] 	
\end{tikzcd}
\]
induce proper base change identifications of  functors $D(S') \rightarrow D(T)$
\begin{eqnarray}\label{eq: app pbc 1}
\Phi_{E'}\cong \Phi_{E}f_{!},\\
\label{eq: app pbc 2}
\Phi_{H'}\cong \Phi_{H}f_{!}.
\end{eqnarray}

\begin{lemma}\label{lemma: butterfly pushforward}
With respect to the identifications \eqref{eq: app pbc 1} and \eqref{eq: app pbc 2}, the two natural transformations \eqref{eq: app push butterfly 1} and \eqref{eq: app push butterfly 2} agree.
\end{lemma}

\begin{proof}


Consider the diagram below, where all parallelograms are Cartesian. 
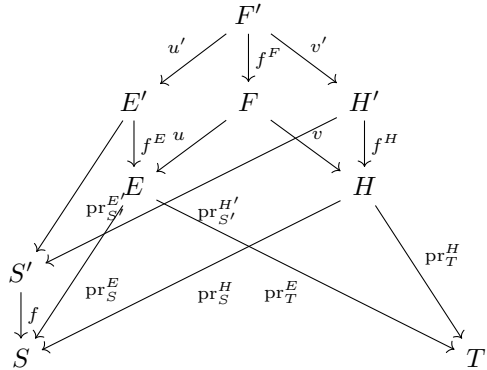
\[
\begin{tikzcd}
& & F' \ar[dl, "u'"'] \ar[dr, "v'"]  \ar[d, "f^F"] \\
& E' \ar[ddl, "\pr_{S'}^{E'}"]   \ar[d, "f^E"]  & F \ar[dl, "u"'] \ar[dr, "v"] & H' \ar[ddlll, "\pr_{S'}^{H'}"]  \ar[d, "f^H"] \\
& E  \ar[ddl, "\pr_S^E"] \ar[ddrrr, "\pr_T^E"']  & & H \ar[ddr, "\pr_T^H"] \ar[ddlll, "\pr_S^H"]  \\
S' \ar[d, "f"] \\	
S & & & & T
\end{tikzcd}
\]
This yields a diagram of natural transformations
\[
\begin{gathered}
\begin{adjustbox}{max width=\textwidth}
$\displaystyle
\begin{tikzcd}[ampersand replacement=\&]
\pr_{T!}^E \pr_S^{E*} f_{!} \ar[dd, "\text{unit}(u)"]  \ar[r, "\di_1"] \& \pr_{T!}^E f^E_! \pr_{S'}^{E'*} \ar[dd, "\text{unit}(u)"] \ar[r, equals] \& \pr_{T!}^{E'} \pr_{S'}^{E'*}  \ar[d, "\text{unit}(u')"] \\
\&\& \pr_{T!}^{E'} (u')_{*} (u')^*  \pr_{S'}^{E'*}\ar[d, equals]\\
 \pr_{T!}^E u_* u^*  \pr_S^{E*} f_! \ar[d, equals] \ar[r, "\di_1"] \& \pr_{T!}^E u_* u^* f^E_! \pr_{S'}^{E'*} \ar[r, "\di_2"] \& \pr_{T!}^E u_* f^F_! (u')^*  \pr_{S'}^{E'*} \ar[d, equals] \\
  \pr^F_{T!} \pr_S^{F*} f_! \ar[d, equals]\ar[rr, "\di_0"]  \& \& \pr_{T!}^{F} f^F_!  \pr_{S'}^{F'*} \ar[d, equals] \\
   \pr_{T!}^H v_! v^* \pr_S^{H*} f_! \ar[d, "{[v]}"]\ar[r, "\di_4"]  \&\pr_{T!}^H v_! v^* f^H_! \pr_{S'}^{H'*} \ar[d, "{[v]}"]  \ar[r, "\di_3"] \&  \pr_{T!}^{H} v_!  f^F_! (v')^* \pr_{S'}^{H' *} \ar[d, "{[v']}"] \\
     \pr_{T!}^H v_! v^! \pr_S^{H*} f_! \tw{-d(v)} \ar[dd, "\text{counit}(v)"]\ar[r, "\di_4"] \& \pr_{T!}^H v_! v^! f^H_! \pr_{S'}^{H'*} \tw{-d(v)} \ar[dd, "\text{counit}(v)"]   \&  \pr_{T!}^H  v_! f^{F}_{!} (v')^! \pr_{S'}^{H'*}\tw{-d(v)} \ar[d, equals] \ar[l]  \\
    \& \& \pr_{T!}^H  f^H_! v'_! (v')^! \pr_{S'}^{H'*}\tw{-d(v)} \ar[d, "\text{counit}(v')"]\\
      \pr_{T!}^H \pr_S^{H*} f_! \tw{-d(v)} \ar[r, "\di_4"]  \& \pr_{T!}^H f^H_!  \pr_{S'}^{H'*} \tw{-d(v)}  \& \ar[l, equals]  \pr_{T!}^{H'} \pr_{S'}^{H'*}\tw{-d(v)}  
\end{tikzcd}
$
\end{adjustbox}
\end{gathered}
\]

The natural transformation \eqref{eq: app push butterfly 1} is the composition along the left column, while the natural transformation \eqref{eq: app push butterfly 2} is the composition along the right column. All the arrows labeled by $\di_?$ are isomorphisms, 
\begin{enumerate}
\item $\di_0$ is induced by the base change natural isomorphism for the Cartesian square 
\[
\begin{tikzcd}
S' \ar[d] & F' \ar[l] \ar[d] \\
S & F \ar[l]
\end{tikzcd}
\]

\item $\di_1$ is induced by the base change natural isomorphism for the Cartesian square 
\[
\begin{tikzcd}
S' \ar[d] & E' \ar[l] \ar[d] \\
S & E \ar[l]
\end{tikzcd}
\]
\item $\di_2$ is induced by the base change natural isomorphism for the Cartesian square 
\[
\begin{tikzcd}
E' \ar[d] & F' \ar[l] \ar[d] \\ 
E & F \ar[l]
\end{tikzcd}
\]
\item $\di_3$ is induced by the base change natural isomorphism for the Cartesian square 
\[
\begin{tikzcd}
H' \ar[d] & F' \ar[l] \ar[d]   \\
H & F \ar[l]
\end{tikzcd}
\]
\item $\di_4$ is induced by the base change natural isomorphism for the Cartesian square 
\[
\begin{tikzcd}
S' \ar[d] & H' \ar[l] \ar[d]\\
S & H \ar[l]
\end{tikzcd}
\]
\end{enumerate}
It remains to show that the above diagram is commutative.  The only non-obvious squares to check are the two wide rectangles in the middle. Their commutativity follows from the observation that the two diagrams 
\[
\begin{tikzcd}
S' \ar[d] & E' \ar[l] \ar[d]  & F' \ar[d] \ar[l] \\
S & E \ar[l] & F \ar[l]
\end{tikzcd}
\]
and
\[
\begin{tikzcd}
S' \ar[d] & H' \ar[l] \ar[d] &  F' \ar[d] \ar[l] \\
S & H \ar[l] & F \ar[l]
\end{tikzcd}
\]
give two decompositions of the same Cartesian square
\[
\begin{tikzcd}
S' \ar[d] & F' \ar[d] \ar[l] \\
S & F \ar[l]
\end{tikzcd}
\]
so that the composition of their base change natural isomorphisms agree. 
\end{proof}

\subsubsection{Butterfly and pullback} Next suppose we are given $g \co T' \rightarrow T$. Then the  construction in \S\ref{sssec: butterflies} gives a natural transformation
\begin{equation}\label{eq: app pull butterfly 1}
g^{*}\Phi_{E}=g^* \pr_{T!}^E \pr_S^{E*} \xrightarrow{g^* \bigstar} g^{*}\Phi_{H}\tw{-d(v)}=g^* \pr_{T!}^H \pr_S^{H*} \tw{-d(v)}
\end{equation}
of functors $D(S) \rightarrow D(T')$. 

On the other hand, we have another natural transformation of the same form coming from another butterfly. Consider the diagram where superscript $'$ means pullback along $g$,
\[
\begin{tikzcd}
& & F'  \ar[dl, "u'"'] \ar[dr, "v'"]  \\
& E'  \ar[dddl, "\pr_{S}^{E'}"'] \ar[ddrrr, "\pr_{T'}^{E'}"']  &  & H' \ar[ddr, "\pr_{T'}^{H'}"] \ar[dddlll, "\pr_S^{H'}"]  \\
&  & &   \\
& & & & T'\\	
S & & & & 
\end{tikzcd}
\]
This is also a butterfly ($u'$ is proper and $v'$ is quasi-smooth) and therefore gives another natural transformation 
\begin{equation}\label{eq: app pull butterfly 2}
\Phi_{E'}\xrightarrow{\bigstar'}\Phi_{H'}\tw{-d(v')}=\Phi_{H'}\tw{-d(v)}
\end{equation}
of functors $D(S) \rightarrow D(T')$. 

The Cartesian squares
\[
\begin{tikzcd}
E' \ar[d, "g^E"] \ar[r] & T' \ar[d, "g"] & H' \ar[d, "g^H"] \ar[r] & T' \ar[d, "g"] \\
E \ar[r] & T  & H \ar[r]  & T
\end{tikzcd}
\]
induce proper base change identifications of functors $D(S) \rightarrow D(T')$
\begin{eqnarray}
\label{eq: app pbc pull 1}
\Phi_{E'}\cong g^{*}\Phi_{E},\\
\label{eq: app pbc pull 2}
\Phi_{H'}\cong g^{*}\Phi_{H}.
\end{eqnarray}


\begin{lemma}\label{lemma: pullback butterfly}
With respect to the identifications \eqref{eq: app pbc pull 1} and \eqref{eq: app pbc pull 2}, the two natural transformations \eqref{eq: app pull butterfly 1} and \eqref{eq: app pull butterfly 2} agree.
\end{lemma}

\begin{proof}The proof is analogous to that for Lemma \ref{lemma: butterfly pushforward}, considering instead the diagram below, all of whose parallelograms are Cartesian. 
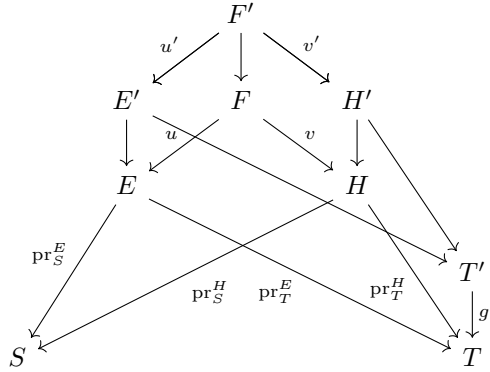
\[
\begin{tikzcd}
& & F' \ar[dl, "u'"'] \ar[dr, "v'"]\ar[dl] \ar[dr]  \ar[d] \\
& E'   \ar[d]  \ar[ddrrr] & F \ar[dl, "u"'] \ar[dr, "v"] & H'  \ar[ddr]   \ar[d] \\
& E  \ar[ddl, "\pr_S^E"'] \ar[ddrrr, "\pr_T^E"']  & & H \ar[ddr, "\pr_T^H"'] \ar[ddlll, "\pr_S^H"]  \\
& & & & T' \ar[d, "g"] \\	
S & & & & T
\end{tikzcd}
\]
\end{proof}

\subsection{Functorialities}\label{ssec: initial functorialities}
 The running notational conventions are now restored: $E,E'$ are derived vector bundles over the base $S$. Let $f \co E' \rightarrow E$ be a map of derived vector bundles over a derived stack $S$. Let $r$ be the rank of $E$ and $r'$ be the rank of $E'$. We begin by establishing some of the easier functorialities for the derived Fourier transform.

\subsubsection{$\FT f_!$ versus $\wh{f}^* \FT$} 

\begin{lemma}\label{lem: FT f_!} We have a natural isomorphism of functors $D(E') \rightarrow D(\wh{E})$
\[
  \FT_{E} \circ f_! [r'-r] \cong \wh{f}^* \circ \FT_{E'}.
\]
\end{lemma}

\begin{proof} The proof is the same as for classical vector bundles, cf. \cite[Th\'eor\`eme 1.2.2.4]{Lau87}. (In fact the same proof works if we replace $ \cL_{\psi}$ by any $\cL\in D(\A^1)$.)
\end{proof}

\subsubsection{$\FT f^*$ versus $\wh{f}_! \FT$: special cases}\label{sssec: FT f^*}
 Assume that $f \co E' \rightarrow E$ is either a closed embedding or is smooth. We explicate that since $f$ is a map of derived vector bundles, $f$ is a closed embedding if and only if $\bT_f$ has tor-amplitude in $[1, \infty)$, and $f$ is smooth if and only if $\bT_f$ has tor-amplitude in $(-\infty, 0]$. In particular, $f$ is a closed embedding if and only if $\wh{f}$ is smooth. 	

Therefore, if $f$ is a closed embedding then the diagram 
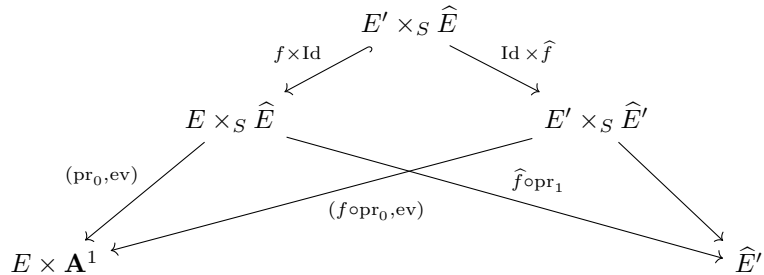
\[
\begin{tikzcd}
& & E' \times_S \wh{E} \ar[dl, hook, "f\times\Id"'] \ar[dr, "\Id\times \wh{f}"] \\
& E \times_S  \wh{E} \ar[ddrrr, "\wh{f} \circ \pr_1"] \ar[ddl, "{(\pr_0, \ev)}"'] && E' \times_S \wh{E}' \ar[ddr] \ar[ddlll, "{(f \circ \pr_0, \ev)}"] \\ \\
E \times \A^1 & & & & \wh{E}'
\end{tikzcd}
\]
is a butterfly. Composing its natural transformation $\bigstar$ with the functor $D(E) \rightarrow D(E \times \A^1)$ given by $\cK \mapsto \cK \boxtimes \cL_{\psi}$, it induces a natural transformation 
\begin{equation}\label{eq: app FT f^* functoriality left}
\FT_{E'} \circ f^* \leftarrow \wh{f}_! \circ \FT_E [r-r'] (r-r').
\end{equation}

\begin{lemma}\label{l: cl emb FT f^*}
If $f$ is a closed embedding, then the natural transformation \eqref{eq: app FT f^* functoriality left} is an isomorphism.
\end{lemma}

\begin{proof}Let $u :=f\times\Id \co E' \times_S \wh{E} \rightarrow E \times_S  \wh{E}$. Examining the definition of \eqref{eq: app FT f^* functoriality left}, and using that $\wh{f}$ is smooth, it suffices to check that the natural map 
\[
\wh{f}_! \pr_{1!} (\pr_0^*\cK\ot\ev^{*}\cL_{\psi}) \rightarrow \wh{f}_! \pr_{1!} u_* u^* (\pr_0^*\cK\ot\ev^{*}\cL_{\psi})
\]
is an isomorphism for all $\cK\in D(E)$.

Let $U := E \setminus E'$ and $j \co U \times_S \wh{E} \inj E \times_S \wh{E} $ be the open complement of $u$. Then it suffices to check that $\wh{f}_! \pr_{1!} j_! j^* (\pr_0^*\cK\ot\ev^*\cL_{\psi})=0$.

For this we can localize over stalks of the base $S$ and thus assume that $S$ is a geometric point. Let $C$ be the cone of $f$; then $\wh{C}$ is the (derived) fiber of $\wh{f}$. The fiber of $\wh{f}_! \pr_{1!}$ over $y \in \wh{E}'$ is $E \times (\wh{C}+y)$, so by proper base change the stalk of $\wh{f}_! \pr_{1!} j_! j^* (\pr_0^* \cK\ot\ev^*\cL_{\psi})$ at $y$ has cohomology groups $H_c^*(U \times (\wh{C}+y), \pr_0^* \cK \otimes \ev^* \cL_{\psi})$, which we want to show vanish. By the projection formula, it suffices to show that for the first projection map $\pr_0 \co U \times (\wh{C}+y) \rightarrow U$, we have $\pr_{0!} \ev^* \cL_{\psi} = 0$. Indeed, by definition any geometric point $u \in U$ has non-zero image in $C$, so $\ev^* \cL_{\psi}|_{\pr_0^{-1}(u)}$ is a nontrivial character sheaf, so its cohomology vanishes. 
\end{proof}

If $f$ is smooth, then $\wh{f} \co \wh{E} \rightarrow \wh{E}'$ is a closed embedding and the diagram 
\[
\begin{tikzcd}
& & E' \times_S \wh{E} \ar[dl, hook, "\Id\times\wh{f}"'] \ar[dr, "f\times\Id"] \\
&E' \times_S \wh{E}' \ar[ddrrr, "\pr_{1}"] \ar[ddl, "{(f\c\pr_{0}, \ev)}"'] &&   E \times_S  \wh{E}   \ar[ddr, "\wh{f}\c\pr_{1}"] \ar[ddlll, "{(\pr_{0}, \ev)}"] \\ \\
E \times \A^1 & & & & \wh{E}'
\end{tikzcd}
\]
is a butterfly, so it induces a natural transformation 
\begin{equation}\label{eq: app FT f^* functoriality right}
\FT_{E'} \circ f^* \rightarrow \wh{f}_! \circ \FT_E [r-r'] (r-r').
\end{equation}

A similar argument as for Lemma \ref{l: cl emb FT f^*} shows that

\begin{lemma}
If $f$ is smooth, then the natural transformation \eqref{eq: app FT f^* functoriality right} is an isomorphism.
\end{lemma}

\subsubsection{The right adjoint of $\FT$}\label{sssec: FT right adjoint}

Let $'\FT^\psi_E \co D(E) \rightarrow D(\wh{E})$ be the functor $\cK \mapsto \pr_{1*} (\ev^*\cL_{\psi} \otimes \pr_0^! (\cK))[-r]$, where maps are as in the diagram 
\[
\begin{tikzcd}
& E \times_S \wh{E} \ar[dl, "\pr_0"'] \ar[dr, "\pr_1"] \ar[r, "\ev"] & \A^1 \\
E & & \wh{E}
\end{tikzcd}
\]
By the compatibility of right adjoints with composition, $'\FT^{-\psi}_{\wh{E}} \co D(\wh{E}) \rightarrow D(E)$ is the right adjoint of $\FT^{\psi}_{E}$. When $\psi$ is understood, we abbreviate $'\FT_{\wh{E}} := '\FT^{-\psi}_{\wh{E}}$. 

Taking the right adjoint of Lemma \ref{lem: FT f_!} (using the convention $'\FT_{\wh E}:='\FT_{\wh E}^{-\psi}$) gives a natural isomorphism 
\begin{equation}
f^!   \circ ( '\FT_{\wh{E}} )[r-r']  \cong ('\FT_{\wh{E}'} ) \c \wh{f}_*. 
\end{equation}
This can also be proved directly by a similar argument as for Lemma \ref{lem: FT f_!}. 

We are especially interested in the case where $\wh{f}$ is a closed embedding (and dually $f$ is smooth), in which case we can replace $\wh{f}_*$ with $\wh{f}_!$ and $f^!$ with $f^*\tw{d(f)}$ above. 

Taking the right adjoint of \eqref{eq: app FT f^* functoriality left} and \eqref{eq: app FT f^* functoriality right} gives a natural isomorphism 
\begin{equation}
f_*    \circ ( '\FT_{\wh{E}'} )[r-r'](r-r')  \cong ('\FT_{\wh{E}} )\c \wh{f}^!
\end{equation}
if $f$ is smooth or a closed embedding. Alternatively, this isomorphism can be proved directly using analogous constructions to those in \S \ref{sssec: butterflies} and \S \ref{sssec: FT f^*}. 

We are especially interested in the case where $\wh{f}$ is smooth (and dually $f$ is a closed embedding), in which case we can replace $\wh{f}^!$ with $\wh{f}^*\tw{d(\wh f)}$ and $f_*$ with $f_!$ above.

\subsubsection{Involutivity}  The key to the functoriality properties of the derived Fourier transform is the (near) involutivity property. This appears to be significantly more subtle to establish than in the situation for classical vector bundles. It will be done later in \S \ref{ssec: involutivity}.

\begin{lemma}\label{lem: involutivity kernel}
Let $E$ be a derived vector bundle over $S$ of rank $r$. Recall $\d_{E}$ denotes $z_{!}\Qll{S}$ where $z: S\to E$ is the zero section. Then any isomorphism
\begin{equation}
\a_{E}: \d_{E}\cong \FT_{\wh E}(\Qll{\wh E})[r](r)
\end{equation}
determines an isomorphism of functors
\begin{equation}\label{FT inv}
\FT_{\wh E}\c \FT_{E}\cong [-1]^{*}(-r).
\end{equation}
\end{lemma}
\begin{proof}
Unravelling the definition of $\FT_{\wh E}\c\FT_{E}$ and using proper base change, we see that it is the endo-functor of $D(E)$ given by the kernel sheaf $\pr_{02!}(\ev_{01}+\ev_{12})^{*}\cL_{\psi}[2r]\in D(E\times_S E)$, where 
\begin{itemize}
\item $\pr_{02}: E\times_{S}\wh E\times_{S }E\to E\times_{S} E$ is the projection to the first and third factors, and
\item $\ev_{01}+\ev_{12}: E\times_{S}\wh E\times_{S }E\to \A^{1}$ is given by $(x,y,z)\mapsto \j{x+z, y}$.
\end{itemize}
On the other hand, $[-1]^{*}:D(E)\to D(E)$ is given by the kernel sheaf $\D^{-}_{!}\Qll{E}$,  where $\D^{-}: E\to E\times_{S}E$ is the anti-diagonal. To give an isomorphism \eqref{FT inv}, it therefore suffices to give an isomorphism of kernel sheaves
\begin{equation}\label{target kernel comp}
\pr_{02!}(\ev_{01}+\ev_{12})^{*}\cL_{\psi}[2r]\cong \D^{-}_{!}\Qll{E}(-r).
\end{equation}

If we let $a: E\times_{S} E\to E$ be the addition map, then proper base change for the Cartesian square 
\[
\begin{tikzcd}[column sep = huge] 
E \times_S \wh{E} \times_S E \ar[r, "{(\pr_0+\pr_2, \pr_1)}"] \ar[d, "\pr_{02}"] & E \times_S \wh{E} \ar[d, "\pr_E"] \\
E \times_S E \ar[r, "a"] & E
\end{tikzcd}
\]
supplies an isomorphism
\begin{equation}\label{inv kernel pullback a}
\pr_{02!}(\ev_{01}+\ev_{12})^{*}\cL_{\psi}[2r]\cong a^{*}\pr_{E!}\ev^{*}\cL_{\psi}[2r]
\end{equation}
Also, proper base change for the Cartesian square
\begin{equation}
\xymatrix{ E\ar[r]\ar[d]^{\D^{-}}& S\ar[d]^{z}\\
E\times_{S}E\ar[r]^{a} & E
}
\end{equation}
supplies an isomorphism
\begin{equation}\label{antidiag}
\D^{-}_{!}\Qll{E}\cong a^{*}\d_{E}.
\end{equation}
Now, an isomorphism $\a_{E} \co \d_{E}(-r)\cong \pr_{E!}\ev^{*}\cL_{\psi}[2r]$ induces an isomorphism $a^{*}\d_{E}(-r)\cong a^{*}\pr_{E!}\ev^{*}\cL_{\psi}[2r]$. Combining this with \eqref{inv kernel pullback a} and \eqref{antidiag}, we get an isomorphism \eqref{target kernel comp}. Since the canonicity of \eqref{target kernel comp} will be a recurring theme, we note that it only depended on $\a_{E}$.
\end{proof}



\begin{defn}\label{def: involutivity datum}
Let $E \rightarrow S$ be a derived vector bundle of virtual rank $r$. An \emph{involutivity datum} on $E$ is a pair of natural isomorphisms
\[
\eta_E \co \FT_{\wh{E}} \circ  \FT_E \xrightarrow{\sim} [-1]^*(-r) \quad \text{and} \quad  \eta_{\wh{E}} \co \FT_E \circ \FT_{\wh{E}}  \xrightarrow{\sim} [-1]^* (-r)
\]
such that the composition
\[
\FT_E  \xrightarrow{\FT_E \circ \eta_E^{-1}} \FT_E \circ \FT_{\wh{E}}  \circ \FT_E [-1]^*(r) \xrightarrow{\eta_{\wh{E}} \circ \FT_E} \FT_E 
\]
is multiplication by $(-1)^{r}$. 

\end{defn}


\begin{remark}
We will later equip every derived vector bundle with a canonical involutivity datum, but the construction is rather circuitous; for example, it will involve first constructing involutivity data that (a priori seem to) depend on auxiliary choices. 
\end{remark}

\subsubsection{Self-adjointness}

In \S \ref{sssec: FT right adjoint} we defined $'\FT_E$, and we explained that $'\FT_{\wh{E}}^{-\psi}$ was the right adjoint to $\FT_E^{\psi}$. 


\begin{lemma}\label{lem: self-adjointness}
Let $E \rightarrow S$ be a derived vector bundle of virtual rank $r$. Any involutivity datum $(\eta_E, \eta_{\wh{E}})$ on $E$ induces natural isomorphisms 
\[
\FT_{E} \xrightarrow{\sim} ('\FT_{E})(-r)
\]
and
\[
\FT_{\wh{E}} \xrightarrow{\sim} ('\FT_{\wh{E}})(-r).
\]
\end{lemma}

\begin{proof}
The given data of $\eta_E$ and $\eta_{\wh{E}}$ show that $\FT_E^\psi(r)$ and $\FT_{\wh{E}}^{-\psi}$ are inverses, and in particular adjoints. Since $'\FT_{E}^\psi$ was by definition the right adjoint of $\FT_{\wh{E}}^{-\psi}$, this induces the natural isomorphism $
\FT_{E}^\psi \xrightarrow{\sim} ('\FT_{E}^\psi)(-r)$; the other isomorphism is obtained similarly.
\end{proof}

\begin{example}[Compatibility with Verdier duality]
Let $E \rightarrow S$ be a derived vector bundle of virtual rank $r$ and $\DD_E$ (resp. $\DD_{\wh{E}}$) denote the Verdier duality functor on $E$ (resp. $\wh{E}$). We have
\[
\DD_{\wh{E}} \circ \FT_E^{\psi}   \cong ('\FT_E^{-\psi}) \circ \DD_E.
\]
By Lemma \ref{lem: self-adjointness}, an involutivity datum for $E$ equips $\FT_E$ and $\FT_{\wh{E}}$ with natural isomorphisms
\[
\DD_{\wh{E}} \circ  \FT_E^{\psi}  \cong \FT_E^{-\psi} \circ \DD_E (r).
\]
\end{example}

\subsubsection{More functoriality}\label{sssec: more functoriality} Let $f \co E' \rightarrow E$ be a linear map of derived vector bundles over $S$. Assume that $E,E'$ are equipped with involutivity data. We may then produce the remaining natural isomorphisms claimed in \S \ref{sssec: FT functoriality}. 

Recall that we always had the natural isomorphism
\begin{equation}\label{eq: init functoriality}
 \FT_{E} \circ f_! [r'-r] \cong \wh{f}^* \circ \FT_{E'}
\end{equation}
without any assumptions. Taking right adjoints in \eqref{eq: init functoriality}, applying Lemma \ref{lem: self-adjointness},  and relabeling terms gives 
\[
\wh{f}^! \circ \FT_{E'} \cong  \FT_E \circ f_*[r-r'](r-r').
\]
Pre-composing \eqref{eq: init functoriality} with $\FT_{\wh{E}'}$ and post-composing with $\FT_{\wh{E}}$, then applying 
\[
(-1)^{\rank E} \eta_{E}^{-1} \co [-1]^* (-r) \xrightarrow{\sim} \FT_{\wh{E}} \circ \FT_E
\]
 and 
 \[
 \eta_{\wh{E}'} \co \FT_{E'} \circ \FT_{\wh{E}'} \xrightarrow{\sim}  [-1]^* (-r'),
 \]
 and re-labeling terms, gives the isomorphism
\begin{equation}\label{eq: functoriality 3}
\FT_{E'} \circ f^* \cong \wh{f}_! \circ \FT_E [r-r'] (r-r').
\end{equation}
Taking right adjoints in \eqref{eq: functoriality 3}, applying Lemma \ref{lem: self-adjointness}, and relabeling terms gives
\[
\FT_{E'} \circ f^! \cong \wh{f}_* \circ \FT_E[r'-r].
\]

\subsection{Involutivity}\label{ssec: involutivity}
We will construct a canonical involutivity datum for any derived vector bundle $E \rightarrow S$. According to Lemma \ref{lem: involutivity kernel}, it suffices to produce isomorphisms 
\[
\alpha_E \co \delta_E \cong \FT_{\wh{E}}(\Qll{\wh{E}})[r](r), \quad \text{and} \quad \alpha_{\wh{E}} \co \delta_{\wh{E}} \cong \FT_{E}(\Qll{E})[r](r).
\]

\subsubsection{The case where $S$ is a point}\label{sssec: point}

Suppose $S$ is a point (i.e., the spectrum of a field). We will produce canonical $\alpha_{\wh{E}}$; the dual argument produces $\alpha_{E}$. This analysis will be used later to prove involutivity in general, and may be illuminating in any case.

Since $S$ is a point, $\cE$ is quasi-isomorphic to a formal complex $\bigoplus_i \cE^i [-i]$, with vanishing differentials. By the compatibility of $\FT$ with exterior tensor products (namely, that Fourier transform on a product of bundles takes an exterior tensor product of complexes to the exterior tensor product of the Fourier transforms of each factor), it suffices to treat the case where $\cE = \cE^i[-i]$, so we assume this to be the case.

If $i = 0$, then $\alpha_{\wh{E}}$ is the classical one of Laumon implicit in the proof of \cite[Th\'{e}or\`{e}me 1.2.2.1]{Lau87}. We will therefore focus on the case $i \neq 0$. 
\begin{itemize}
\item If $i>0$, then $E$ has classical truncation $S$. Therefore, $\pi_E^*$ induces an equivalence $D(S) \xrightarrow{\sim} D(E)$, whose inverse is its right adjoint $\pi_{E*}$. We then have $\pi_E^* = \pi_E^! \co D(S) \xrightarrow{\sim} D(E)$ and $\pi_{E*} = \pi_{E!} \co D(E) \xrightarrow{\sim} D(S)$. Since $\pi_E \circ z_E = \Id_S$, we deduce $z_{E*} = \pi_E^*$, $z_E^* = \pi_{E*}$, $z_{E!} = z_{E*}$, and $z_E^! = z_E^*$. 

\item If $i<0$, then $E$ is the $-i$-fold iterated classifying stack of a product of copies of $\G_a$. Since $\G_a$ is a connected unipotent group scheme, $\pi_E^*$ induces an equivalence $D(S) \xrightarrow{\sim} D(E)$, whose inverse is therefore $\pi_{E*}$. Since $\pi_E$ is smooth of rank $r$, we have $\pi^!_E = \pi^*_E \tw{r}$ and $\pi_{E!} \cong \pi_{E*}  \tw{-r}$. Since $\pi_E \circ z_E  = \Id_S$, we deduce that $z_{E*} \cong \pi_E^*$, $z_E^* \cong \pi_{E*}$, $z_E^! \cong z_E^*\tw{-r}$, and $z_{E!} \cong z_{E*}\tw{r}$. 
\end{itemize}
Below we will always use $\pi_E^*$ to identify $D(S) \cong D(E)$ and $\pi_{\wh{E}}^*$ to identify $D(S) \cong D(\wh{E})$. We will therefore view $\FT_E$ as an endofunctor of $D(S)$. 

\begin{itemize}
\item If $i>0$, then under the above identifications $\FT_E$ is simply the endofunctor $[r]$ of $D(S)$, so $\FT_E(\Qll{E})$ identifies with $\Qll{S}[r] \in D(S)$. On the other hand, $\delta_{\wh{E}} = z_{\wh{E}!} \Qll{S}$ identifies with $\Qll{S} \tw{r} \in D(S)$. Thus we obtain $\alpha_{\wh{E}}$.

\item If $i<0$, then under the above identifications $\FT_E$ is the endofunctor $[-r](-r)$ of $D(S)$, so $\FT_E(\Qll{E})$ identifies with $\Qll{S}[-r](-r) \in D(S)$. On the other hand, $\delta_{\wh{E}} = z_{\wh{E}!} \Qll{S}$ identifies with $\Qll{S}$. We take $\alpha_{\wh{E}}$ to be $(-1)^r$ times the obvious isomorphism. 
\end{itemize}

\subsubsection{Bootstrapping from vector bundles}\label{sssec: embed into classical} 
Suppose $\cE \in \Perf(S)$ admits a global presentation by a finite complex of vector bundles on $S$, say of the form 
\[
\ldots \rightarrow \cE^{-1} \rightarrow \cE^0 \rightarrow \cE^1 \rightarrow \ldots
\]

We will use the following device to bootstrap from the theory for classical vector bundles, where the theory was established by Laumon \cite{Lau87}. Consider the ``stupid truncations''
\[
\cE^{\geq 0} = (\cE^0 \rightarrow \cE^1 \rightarrow \ldots )
\]
and
\[
\cE^{\leq 0} = (\ldots \rightarrow \cE^{-1} \rightarrow \cE^0)
\]
and let $E^{\geq 0}, E^{\leq 0}$ be the associated derived vector bundles. Then we 
have a pullback (and pushout) square in $\Perf(S)$
\[
\begin{tikzcd}
\cE^{\geq 0} \ar[d]  \ar[r] & \cE^0 \ar[d] \\ 
\cE  \ar[r] & \cE^{\leq 0}
\end{tikzcd}
\]
which induces a derived Cartesian square of total spaces
\begin{equation}\label{eq: comparison square}
\begin{tikzcd}
E^{\geq 0} \ar[d, "p'"] \ar[r, "i'", hook] & E^0 \ar[d, "p"] \\ 
E  \ar[r, "i", hook] & E^{\leq 0} 
\end{tikzcd}
\end{equation}
Note that $E^0 \rightarrow S$ is a classical vector bundle, while $E^{\leq 0} \rightarrow S$ is smooth (and represented by classical stacks), and $E^{\geq 0} \rightarrow S$ is separated (and representable in derived schemes). Here $i$ is a closed embedding, so $i_! = i_* \co D(E) \inj D(E^{\leq 0})$ is fully faithful and similarly for $i'$. Also, $p$ is smooth with acyclic geometric fibers, so $p^* \co D(E^{\leq 0}) \inj D(E^0)$ is fully faithful and similarly for $p'$. Hence we have a fully faithful embedding 
\begin{equation}\label{eq: embed into classical}
p^* i_! \co D(E) \inj D(E^0).
\end{equation}

\subsubsection{Global presentations}\label{ssec: global presentations} 
We now establish involutivity for derived vector bundles $E$ that admit a global presentation. Note that, by the definition of a perfect complex, any derived vector bundle admits a global presentation locally on the base $S$. 

\begin{lemma}\label{lem: involutivity given global presentation}
Suppose $E$ is a derived vector bundle. 

(1) Any global presentation $E^{\bu}$ of $E$ induces an involutivity datum on $E$.

(2) If $S$ is a point, then the involutivity datum from (1) is naturally isomorphic to that from \S \ref{sssec: point}.
\end{lemma}

\begin{proof}
(1) As in \eqref{eq: embed into classical} we have a fully faithful embedding $p^* i_! \co D(E) \inj D(E^0)$. We produce a natural isomorphism 
\[
p^* i_! \FT_{\wh{E}} \circ \FT_E \cong p^* i_! [-1]^*(-r).
\]
By construction the maps $p, \wh{i}$ are smooth and the maps $i, \wh{p}$ are closed embeddings. Set $r_{\leq 0}:=\rank(E^{\leq 0})$ and $r_0:=\rank(E^0)$. We may apply the functorialities in \S \ref{ssec: initial functorialities} to produce a sequence of natural isomorphisms: 
\[
\begin{aligned}
p^*i_!\FT_{\wh E}\FT_E
&\cong p^*\FT_{\wh{E^{\leq 0}}}\wh i^*\FT_E
   [r_{\leq 0}-r](r_{\leq 0}-r) \\
&\cong \FT_{\wh{E^0}}\wh p_!\wh i^*\FT_E
   [2r_{\leq 0}-r-r_0](r_{\leq 0}-r) \\
&\cong \FT_{\wh{E^0}}\wh p_!\FT_{E^{\leq 0}}i_!
   [r_{\leq 0}-r_0](r_{\leq 0}-r) \\
&\cong \FT_{\wh{E^0}}\FT_{E^0}p^*i_!(r_0-r).
\end{aligned}
\]
Applying the classical involutivity datum $\eta_{E^0}\co \FT_{\wh{E^0}}\FT_{E^0}\xrightarrow{\sim}[-1]^*(-r_0)$ and using the equivariance of $p^*i_!$ with respect to multiplication by $-1$ gives
\[
p^*i_!\FT_{\wh E}\FT_E \cong p^*i_![-1]^*(-r).
\]
By the full faithfulness of $p^*i_!$, this reflects to an isomorphism $\eta_E\co \FT_{\wh E}\FT_E\xrightarrow{\sim}[-1]^*(-r)$. The natural isomorphism $\eta_{\wh E}$ is obtained by the same argument on the dual bundle. 

(2) By construction, the compatibility reduces to Lemma \ref{lem: ftft push/pull} below.
\end{proof}

\begin{lemma}\label{lem: ftft push/pull} Let $f \co E' \rightarrow E$ be a map of derived vector bundles.

(1) Suppose $f$ is a closed embedding. Then the diagram 
\begin{equation}\label{eq: push FTFT}
\begin{tikzcd}
f_! \FT_{\wh{E}'} \FT_{E'} \ar[r, "\sim"', "f_!\eta_{E'}"] \ar[d, "\sim"] &  f_! [-1]^* (-r')  \ar[dd, equals] \\
\FT_{\wh{E}} \wh{f}^* \FT_{E'} [r-r'](r-r') \ar[d, "\sim"]  &  \\
\FT_{\wh{E}} \FT_E f_!  (r-r') \ar[r, "\sim", "\eta_E"']  & f_! [-1]^* (-r')
\end{tikzcd}
\end{equation}
commutes. 

(2) Suppose $f$ is smooth. Then the diagram 
\begin{equation}\label{eq: pull FTFT} 
\begin{tikzcd}
f^* \FT_{\wh{E}} \FT_E \ar[r, "\sim"', "f^*\eta_E"] \ar[d, "\sim"] &  f^* [-1]^* (-r)  \ar[dd, equals] \\
\FT_{\wh{E}'} \wh{f}_! \FT_E [r-r'] \ar[d, "\sim"]  &  \\
\FT_{\wh{E}'} \FT_{E'} f^*  (r'-r) \ar[r, "\sim", "\eta_{E'}"']  & f^* [-1]^* (-r)
\end{tikzcd}
\end{equation}
commutes.

\end{lemma}

\begin{proof}
The two situations are similar so we just prove (1). The two paths are each given by an isomorphism of the respective kernel sheaves with the constant sheaf on the graph of $f$ in $E' \times_S E$. Hence they differ by a scalar in $H^0(S, \Qll{S})$. We can compute this scalar locally, thus reducing to $S = \pt$ by base change. Then $E$ and $E'$ factor as a product of derived vector bundles placed in degree $i$, so we may reduce to the case where they are each in degree $i$. If $i \neq 0$, then the Lemma follows from tracing through the explicit descriptions in \S \ref{sssec: point}. We henceforth focus on the case where $i=0$, so $E$ and $E'$ are classical vector bundles.
Observe that we can compute the scalar in question on any non-zero object; we will take the object $\delta_{E'}$. This reduces to the case where $E' = 0$. 

The isomorphism $\eta_{E} \co \FT_{\wh{E}} \FT_E \xrightarrow{\sim} [-1]^*(-r)$ is the map induced by the butterfly 
\begin{equation}\label{eq: ftft butterfly}
\begin{tikzcd}
& & E \times \wh{E} \ar[dl, hook'] \ar[dr, twoheadrightarrow] \\
&  E \times \wh{E} \times E  \ar[dl] \ar[drrr]   & & E \ar[dr] \ar[dlll]  \\
 E \times \A^1 & & & & E 
\end{tikzcd}
\end{equation}
using the Artin-Schreier sheaf on $\A^1$. According to \S \ref{sssec: butterfly and pushforward}, the bottom row of \eqref{eq: push FTFT} is then given by the butterfly obtained by pulling back along $f \co 0 \inj E$: 
\begin{equation}\label{eq: ftft push butterfly 1}
\begin{tikzcd}
& &  0  \times \wh{E} \ar[dl, hook'] \ar[dr, twoheadrightarrow] \\
&  0 \times \wh{E} \times E  \ar[dl] \ar[drrr]   & & 0 \ar[dr] \ar[dlll]  \\
0 \times  \A^1  & & & & E 
\end{tikzcd}
\end{equation}
Tracing through \S \ref{sssec: FT f^*}, we find that the isomorphism $\FT_{\wh{E}} \FT_E f_! \xrightarrow{\sim} f_! [-1]^*(-r)$ given by going around the left and top of \eqref{eq: push FTFT} is also given by a butterfly: 
\begin{equation}
\begin{tikzcd}
& &  \wh{E} \times 0  \ar[dl, hook'] \ar[dr, twoheadrightarrow] \\
&  \wh{E} \times E  \ar[dl] \ar[drrr]   & & \wh{0} \times 0 \ar[dr] \ar[dlll]  \\
 0 \times \A^1  & & & & E 
\end{tikzcd}
\end{equation}
By inspection, this agrees with the butterfly \eqref{eq: ftft push butterfly 1}. 
\end{proof}

\subsubsection{Base change}\label{sssec: FT bc proof} Consider the setup of \S\ref{sssec: FT base change} where we perform a base change $h:\wt S\to S$. We will prove the isomorphisms \eqref{eq: FT bc h^*}, \eqref{eq: FT bc h_!}, \eqref{eq: FT bc h_*} and \eqref{eq: FT bc h^!} \emph{without any assumptions}. (At an intermediate stage the proof invokes Lemma \ref{lem: involutivity given global presentation}, and the proof is used subsequently to construct the canonical involutivity datum for general derived vector bundles, which explains why it appears here.)

The isomorphisms \eqref{eq: FT bc h^*} and \eqref{eq: FT bc h_!} follow directly from proper base change.

Now consider \eqref{eq: FT bc h_*}. Consider the commutative diagram
\begin{equation}\label{FT bc proof}
\xymatrix{ &\wt E\times_{\wt S}\wh{\wt E}\ar[dl]_{\wt\pr_{0}}\ar[dr]^{\wt\pr_{1}}\ar[d]^{\wt h}\\
\wt E\ar[d]_{h^{E}} & E\times_{S}\wh E \ar[dl]_{\pr_{0}}\ar[dr]^{\pr_{1}}& \wh{\wt E}\ar[d]^{h^{\wh E}}\\
E & & \wh E 
}
\end{equation}
We have a natural transformation
\begin{equation}\label{eq: FT bc h_* map}
\FT_{E}\c h^{E}_{*}\to h^{\wh E}_{*}\c \FT_{\wt E}
\end{equation}
as the composition (for any $\cK\in D(\wt E)$)
\begin{align*}
\pr_{1!}(\pr_{0}^{*}h^{E}_{*}\cK\ot\ev^{*}\cL_{\psi}) & \xr{(1)}\pr_{1!}(\wt h_{*}\wt\pr_{0}^{*}\cK\ot\ev^{*}\cL_{\psi})\cong \pr_{1!}\wt h_{*}(\wt\pr_{0}^{*}\cK\ot\wt\ev^{*}\cL_{\psi}) \\
& \xr{(2)}h^{\wh E}_{*}\wt\pr_{1!}(\wt\pr_{0}^{*}\cK\ot\wt\ev^{*}\cL_{\psi}).
\end{align*}
Here $(1)$ is the natural transformation $\pr_{0}^{*}h^{E}_{*}\to \wt h_{*}\wt\pr_{0}^{*}$ coming from the left parallelogram in  \eqref{FT bc proof}, and $(2)$ is the natural transformation $\pr_{1!}\wt h_{*}\to h^{\wh E}_{*}\wt\pr_{1!}$ obtained by adjunction from the proper base change isomorphism attached to the right parallelogram in  \eqref{FT bc proof}. To check \eqref{eq: FT bc h_* map} is an isomorphism, one can work Zariski locally over $S$, hence reducing to the case where $E$ admits a global presentation. In this case, using Lemma \ref{lem: self-adjointness}, we may replace $\FT$ with $'\FT$ in \eqref{eq: FT bc h_* map} (up to a twist), which then is visibly an isomorphism by proper base change.

Now consider \eqref{eq: FT bc h^!}.  We have a natural transformation
\begin{equation}\label{eq: FT bc h^! map}
\FT_{\wt E}\c (h^{E})^{!}\to (h^{\wh E})^{!}\c \FT_{E}
\end{equation}
as the composition (for any $\cK\in D(E)$)
\begin{align*}
\wt\pr_{1!}(\wt\pr_{0}^{*}(h^{E})^{!}\cK\ot\wt\ev^{*}\cL_{\psi}) & \xr{(1)}\wt\pr_{1!}(\wt h^{!}\pr_{0}^{*}\cK\ot\wt\ev^{*}\cL_{\psi})\cong \wt\pr_{1!}\wt h^{!}(\pr_{0}^{*}\cK\ot\ev^{*}\cL_{\psi}) \\
& \xr{(2)}(h^{\wh E})^{!}\pr_{1!}(\pr_{0}^{*}\cK\ot\ev^{*}\cL_{\psi}).
\end{align*}
Here $(1)$ is the natural transformation $\wt\pr_{0}^{*}(h^{E})^{!}\to \wt h^{!}\pr_{0}^{*}$ obtained by adjunction from the proper base change isomorphism attached to the left parallelogram in  \eqref{FT bc proof}, and $(2)$ is the  natural transformation $\wt\pr_{1!}\wt h^{!}\to (h^{\wh E})^{!}\pr_{1!}$ coming from the right parallelogram in  \eqref{FT bc proof}. The proof that \eqref{eq: FT bc h^! map} is an isomorphism is similar to that of \eqref{eq: FT bc h_* map}.

\subsubsection{Involutivity in general} Finally we can produce the promised canonical involutivity datum for a general derived vector bundle.

Let $E \rightarrow S$ be a derived vector bundle of virtual rank $r$. Since $\delta_E := z_! \Qll{S}$, producing a map 
\begin{equation}\label{eq: a_E}
\a_{E}: \d_{E}\to \FT_{\wh E}(\Qll{\wh E})[r](r).
\end{equation}
is equivalent to producing the adjoint map 
\begin{equation}\label{eq: a_E'}
\a'_{E}: \Qll{S} \rightarrow z^{!}\FT_{\wh E}(\Qll{\wh E})[r](r).
\end{equation}

\begin{lemma}\label{lem: adjoint map isom}
If the map \eqref{eq: a_E} is an isomorphism, then the adjoint map \eqref{eq: a_E'} is an isomorphism. 
\end{lemma}

\begin{proof}
The map $\a_E'$ is the composition 
\[
\Qll{S} \rightarrow z^! z_! \Qll{S} \xrightarrow{z^! \a_E} z^! \FT_{\wh E}(\Qll{\wh E})[r](r).
\]
Since $\a_E$ is an isomorphism by assumption, it suffices to see that the unit map $\Id \rightarrow z^! z_!$ is an isomorphism. Equipping $E$ with the natural $\G_m$-action (by scaling), we will apply the ``stacky contraction principle'' Theorem \ref{thm:contraction-principle}. For $\pi \co E \rightarrow S$, it says thats on the subcategory of monodromic sheaves, $\pi_!$ is right adjoint to $z_!$. This applies here because $z_!\cK$ is monodromic for every $\cK \in D(S)$: its support is the fixed locus and the scaling action is trivial there. Thus, for every $\cK\in D(S)$, the unit map $\cK \to z^!z_!\cK$ is identified with the unit map $\cK \to \pi_!z_!\cK \cong \cK$, hence is an isomorphism. 
\end{proof}

\begin{lemma}
Let $E$ be a derived vector bundle over $S$ of virtual rank $r$. Then there is a canonical (i.e., independent of auxiliary choices) isomorphism
\begin{equation}
\a_{E}: \d_{E}\cong \FT_{\wh E}(\Qll{\wh E})[r](r).
\end{equation}
\end{lemma}

\begin{proof}
We will first construct a canonical isomorphism
\begin{equation}
\a'_{E}: \Qll{S}\cong z^{!}\FT_{\wh E}(\Qll{\wh E})[r](r).
\end{equation}

Let $\cK\in D(S)$ denote the right side. We claim that $\cK$ is a \emph{sheaf} concentrated in degree $0$, and Zariski-locally isomorphic to the constant sheaf $\Qll{S}$. The claim can be checked Zariski-locally, so to prove it we may assume that $E$ has a global presentation. Then Lemma \ref{lem: involutivity given global presentation} equips $E$ with an involutivity datum (which a priori depends on the choice of presentation), so \S \ref{sssec: more functoriality} applies and then Example \ref{ex: delta-constant duality} gives an isomorphism
\[
\a_E \co \d_{E}\xrightarrow{\sim} \FT_{\wh E}(\Qll{\wh E})[r](r).
\]
Then by Lemma \ref{lem: adjoint map isom}, the adjoint map $\Qll{S} \xrightarrow{\sim} z^{!}\FT_{\wh E}(\Qll{\wh E})[r](r)$ is an isomorphism. This establishes the claim. 

We now return to the case of general $E$. In the proof of the claim, we saw that Zariski-locally on $S$ we have global presentations of $E$, which induce local trivializations of $\cK$. By Lemma \ref{lem: involutivity given global presentation}(2), these local trivializations of $\cK$ agree on the stalk at any point of $S$. Since the claim implies in particular that $\cK$ is a local system on $S$, the local trivializations therefore glue to give an isomorphism $\Qll{S}\cong \cK$. (Note that we are gluing in the abelian category of local systems on $S$; it was crucial to first establish that $\cK$ lies in this subcategory in order to be able to glue.) Furthermore, Lemma \ref{lem: involutivity given global presentation}(2) shows that the resulting $\a_E'$ is independent of any and all choices of local trivializations. 

The canonical isomorphism $\a'_{E}$ then induces by adjunction a canonical map
\begin{equation}
\a_{E}: \d_{E}\to \FT_{\wh E}(\Qll{\wh E})[r](r).
\end{equation}
We claim that $\a_{E}$ is an isomorphism. This claim can be checked after base changing to a point of $S$. Hence we may assume that $E$ admits a global presentation, and then by its construction, the map $\a_E$ agrees with the map induced by the map $\d_{E}\to \FT_{\wh E}(\Qll{\wh E})[r](r)$ produced by the involutivity datum of Lemma \ref{lem: involutivity given global presentation} and Example \ref{ex: delta-constant duality}, which is an isomorphism. This finishes the proof.
\end{proof}

The canonical isomorphisms $\a_E$ and $\a_{\wh E}$ give the two involutivity isomorphisms by Lemma \ref{lem: involutivity kernel}. The sign condition in Definition \ref{def: involutivity datum} is local on $S$ and can be checked after base change to geometric points using the base-change compatibilities of \S \ref{sssec: FT bc proof}. After choosing local global presentations, the canonical maps agree with the presentation-dependent maps above by Lemma \ref{lem: involutivity given global presentation}(2). Thus all derived vector bundles have canonical involutivity data, independent of global presentations. Henceforth we implicitly equip all derived vector bundles with these canonical involutivity data.

\subsection{Fourier transform and proper base change} 
In this subsection we prove Proposition \ref{prop: ft pbc}. The notational conventions of the subsection are now reset. We maintain the notation and setup of \S \ref{sssec: pbc}.

\begin{lemma}\label{lem: factorize smooth/proper}
Let $f \co E'  \rightarrow E$ be a globally presented map of derived vector bundles. Then $f$ admits a factorization 
\[
f \co E' \xrightarrow{i} \wt{E} \xrightarrow{p} E
\]
such that $i$ is a closed embedding and $p$ is smooth.
\end{lemma}

\begin{proof}
Let $K := \Tot_S(\cK)$ be the derived kernel of $f$. By the hypothesis, $\cK$ admits a global presentation 
\[
\cK^{\bu} = \cdots \rightarrow \cK^{-1} \rightarrow \cK^0 \rightarrow  \cK^1 \rightarrow \cdots
\]
(for example, we can take for $\cK^{\bu}$ the usual cone construction of $f$ shifted by $1$). Then we have an exact triangle for naive truncations 
\[
\cK^{\geq 1} \rightarrow \cK \rightarrow \cK^{\leq 0}.
\]
We take $\wt{\cE} := \cone(\cK^{\geq 1} \rightarrow \cE')$, which then has a factorization 
\begin{equation}\label{eq: factorize closed embedding and smooth}
\cE' \rightarrow \wt{\cE} \rightarrow \cone(\cK \rightarrow \cE') \cong \cE.
\end{equation}
Note that the last isomorphism is in $\Perf(S)$ and we do not claim that it is represented by a map of the given presentations. Regardless, \eqref{eq: factorize closed embedding and smooth} induces a diagram of total spaces
\[
E' \xrightarrow{i} \wt{E} \xrightarrow{p} E
\]
whose composition is $f$, whose first map has derived kernel $\Tot_S(\cK^{\geq 1})$ hence is a closed embedding, and whose second map has derived kernel $\Tot_S(\cK^{\leq 0})$ hence is smooth. 
\end{proof}

\subsubsection{Formulation of the main statement}We expand the formulation of the compatibility \eqref{eq: FT BC compatibility}. Consider a Cartesian square of derived vector bundles: 
\begin{equation}\label{eq: ABCD}
\begin{tikzcd}
& B \ar[dl, "g'"'] \ar[dr, "f'"] \\
A  \ar[dr, "f"'] & & D  \ar[dl, "g"] \\
& C
\end{tikzcd}
\end{equation}
where all maps are linear. Set $d:=d(f)$, $\delta := d(g)$. Using \S \ref{sssec: FT functoriality} and proper base change, we have a hexagon of functors $D(A) \rightarrow D(\wh{D})$: 
\begin{equation}\label{eq: FT hexagon}
\begin{tikzcd}
\FT_D  g^* f_!  [d+\delta](\delta)  & \ar[l, "\sim"']   \FT_D  f'_! (g')^* [d+\delta](\delta)  \ar[r, "\sim"]  & (\wh{f}')^* \FT_B  \circ (g')^* [\delta](\delta)\\
\wh{g}_! \FT_C f_![d] \ar[u, "\sim"] &  \ar[l, "\sim"'] \wh{g}_! \wh{f}^* \FT_A  \ar[r, "\sim"]  & (\wh{f}')^* \wh{g}'_! \FT_A   \ar[u, "\sim"] 
\end{tikzcd}
\end{equation}

\begin{lemma}\label{lem: smooth/proper case}
Assume that each of $f,g$ is either smooth or a closed embedding. Then diagram \eqref{eq: FT hexagon} commutes. 
\end{lemma}

\begin{prop}\label{prop: FT hexagon}
Assume that $f$ and $g$ are globally presented. Then diagram \eqref{eq: FT hexagon} commutes. 
\end{prop}

\begin{remark}\label{rem: lemma vs prop}
We emphasize that Lemma \ref{lem: smooth/proper case} has a weaker hypothesis than Proposition \ref{prop: FT hexagon}; in particular, it does \emph{not} require maps $f,g$ to be globally presented. This is important because when reducing Proposition \ref{prop: FT hexagon} to Lemma \ref{lem: smooth/proper case}, we may not be able to guarantee that intermediate diagrams are globally presented compatibly with the original global presentation of \eqref{eq: ABCD}. Fortunately, because of the weaker hypothesis, there is no need to arrange such compatibility.
\end{remark}

\begin{proof}[Proof of Proposition \ref{prop: FT hexagon}, assuming Lemma \ref{lem: smooth/proper case}]

  We first prove an intermediate reduction: the desired hexagon commutes whenever at least one of $f$ and $g$ is smooth or a closed embedding. Suppose without loss of generality that $g$ is smooth or a closed embedding.

By Lemma \ref{lem: factorize smooth/proper}, we may factor $f$ as $A \xrightarrow{i} \wt{A} \xrightarrow{p} C$ where $i$ is a closed embedding and $p$ is smooth. This induces a factorization of the diagram \eqref{eq: ABCD} into a sequence of two Cartesian squares
\begin{equation}\label{eq: A decompose square}
\begin{tikzcd}
& B \ar[dl, "g'"'] \ar[ddrr, bend left, "f'"] \ar[dl]  \ar[dr, "i'"']\\
A  \ar[dr, "i"] \ar[ddrr, "f"', bend right] & & \wt{D}  \ar[dl] \ar[dr, "p'"'] \\
& \wt{A}  \ar[dr, "p"]   & & D  \ar[dl, "g"] \\
& &  C
\end{tikzcd}
\end{equation}

Now both inner Cartesian squares have the property that their sides are all smooth or closed embeddings, so we may now apply Lemma \ref{lem: smooth/proper case} to them. Lemma \ref{lem: smooth/proper case} is compatible with compositions of Cartesian squares, so the commutativity of the analogous hexagons to \eqref{eq: FT hexagon} for the two inner squares of \eqref{eq: A decompose square} implies the commutativity of \eqref{eq: FT hexagon}.

It remains to remove the auxiliary assumption that one of $f$ and $g$ is already smooth or a closed embedding. As above, we may factor $f$ as $A \xrightarrow{i} \wt{A} \xrightarrow{p} C$ where $i$ is a closed embedding and $p$ is smooth. This produces a commutative diagram \eqref{eq: A decompose square} both of whose inner Cartesian squares have a set of parallel sides which are smooth or closed embeddings, and we may repeat the argument of the previous paragraph. 
\end{proof}

It remains to prove Lemma \ref{lem: smooth/proper case}. We make some initial reductions. We claim that the result is immediate if $f$ is a closed embedding, or if $g$ is smooth. Indeed:
\begin{itemize}
\item If $f$ is a closed embedding then $\wh{f}$ is smooth and we may replace $f_!$ by $f_*$ everywhere in \eqref{eq: FT hexagon}. Since $f'$ is also a closed embedding, the relevant Fourier functoriality isomorphisms may be written in the normalized form
\[
\FT_C f_*[d]\xrightarrow{\sim}\wh f^*\FT_A
\qquad\text{and}\qquad
\FT_D f'_*[d]\xrightarrow{\sim}(\wh f')^*\FT_B.
\]
After applying the adjunctions $(g')^*\dashv g'_*$ and $\wh g'_!\dashv \wh g'^!$, the commutativity of \eqref{eq: FT hexagon} is therefore equivalent to the commutativity of
\[
\begin{tikzcd}
\FT_C f_*[d]   \ar[r, "\textup{unit}(g')"] & \FT_C f_* g'_* (g')^* [d]  \\
\wh{f}^* \FT_A  \ar[u, "\sim"] \ar[r, "\textup{unit}(\wh g')"] & \wh{f}^* \wh{g}'^! \wh{g}'_! \FT_A  \ar[u, "\sim"]
\end{tikzcd}
\]
This last square is precisely the compatibility of the functoriality isomorphism $\FT_C f_*[d]\simeq \wh f^*\FT_A$ with the unit of adjunction, which is incorporated into its construction in \S \ref{sssec: FT functoriality}.
\item If $g$ is smooth then $\wh{g}$ is a closed embedding and we may replace $\wh{g}_!$ by $\wh{g}_*$ everywhere in \eqref{eq: FT hexagon}, and a similar argument applies. 
\end{itemize}
\emph{We may and do henceforth assume that $f$ is smooth and $g$ is a closed embedding}.

\subsubsection{More general statement} In fact, we can formulate a more general statement. The transform $\FT$ is defined with respect to $\cL_{\psi}$ on $\A^1$ via pull-push in the diagram 
\[
\begin{tikzcd}
A \times \A^1  & A \times_S \wh{A} \ar[l] \ar[r] & \wh{A}
\end{tikzcd}
\]

More generally, define the functor $\Conv \co D(A \times \A^1) \rightarrow D(\wh{A})$ as $\pr_{1!} (\pr_0, \ev)^*$ (see \S\ref{sssec: correspondences}) for the diagram 
\[
\begin{tikzcd}
A \times \A^1 & A \times_S \wh{A} \ar[l, "{(\pr_0, \ev)}"'] \ar[r, "\pr_1"] & \wh{A}.
\end{tikzcd}
\]
Composing $\Conv$ with the functor $D(A) \rightarrow D(A \times \A^1)$ given by external product with the Artin-Schreier sheaf $\cL_{\psi}$ recovers $\FT^{\psi}$, up to shift. 

We will construct a hexagon
\begin{equation}\label{eq: conv hexagon}
\begin{tikzcd}
\Conv g^* f_! \tw{\delta}   & \ar[l, "\sim"']   \Conv f'_! (g')^* \tw{\delta}\ar[r, "\sim"]  & (\wh{f}')^* \Conv (g')^* \tw{\delta}\\
\wh{g}_! \Conv f_! \ar[u] &  \ar[l, "\sim"'] \wh{g}_! \wh{f}^* \Conv  \ar[r, "\sim"]  & (\wh{f}')^* \wh{g}'_! \Conv \ar[u] 
\end{tikzcd}
\end{equation}
that recovers \eqref{eq: FT hexagon} in the above manner; note however that in general the vertical arrows are not isomorphisms. The top-left and bottom-right horizontal arrows are the natural isomorphisms induced by proper base change. The other arrows are explained in \S \ref{sssec: Conv f_!} and \S \ref{sssec: Conv g^*} below.

We are still assuming that $g, g'$ are closed embeddings (so $\wh{g}, \wh{g}'$ are smooth) and $f,f'$ are smooth (so $\wh{f}, \wh{f}'$ are closed embeddings).

\subsubsection{}\label{sssec: Conv f_!} For any $f \co A \rightarrow C$, we will define a natural isomorphism $\wh{f}^* \Conv \cong \Conv f_!$ of functors $D(A \times \A^1) \rightarrow D(\wh{C})$. (We do not need the smoothness of $f$ here.) \emph{This supplies the bottom left and top right horizontal arrows of \eqref{eq: conv hexagon}.}
\begin{itemize}
\item The functor $\wh{f}^* \Conv$ is given by the correspondence 
\[
\begin{tikzcd} 
& & A \times_S \wh{C} \ar[dl] \ar[dr] \\
& A \times_S \wh{A} \ar[dr, "\pr_1"] \ar[dl, "{(\pr_0, \ev)}"']  & & \wh{C} \ar[dl, "\wh{f}"']  \ar[dr, equals] \\
A \times \A^1  & & \wh{A}  & & \wh{C}
\end{tikzcd}
\]
\item The functor $\Conv f_!$ is given by the correspondence
\[
\begin{tikzcd} 
& & A \times_S \wh{C} \ar[dl, "{(\pr_0,\ev \circ (f\times \Id))}"'] \ar[dr] \\
& A \times \A^1  \ar[dr, "f \times \Id"] \ar[dl, equals]  & & C \times_S \wh{C} \ar[dl, "{(\pr_0, \ev)}"'] \ar[dr, "\pr_1"] \\
A \times \A^1  & & C \times \A^1  & & \wh{C}
\end{tikzcd}
\]
\end{itemize}
These correspondences between $A\times\A^{1}$ and $\wh C$ visibly agree. 

\subsubsection{}\label{sssec: Conv g^*} For any closed embedding $g' \co B \rightarrow A$, we will define a natural transformation $\wh{g}'_! \Conv  \rightarrow \Conv (g')^*\tw{-d(\wh g')}$ of functors $D(A \times \A^1) \rightarrow D(\wh{B})$. (We \emph{will} use the smoothness of $\wh{g}'$ here.) \emph{This supplies the vertical arrows of \eqref{eq: conv hexagon}.}
\begin{itemize}
\item The functor $\Conv (g')^*$ is given by the correspondence
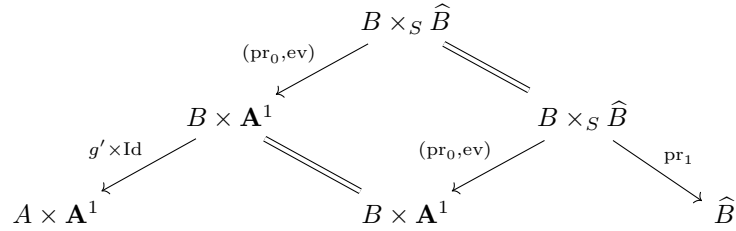
\[
\begin{tikzcd}
& & B \times_S \wh{B} \ar[dl, "{(\pr_0, \ev)}"'] \ar[dr, equals] \\
& B\times \A^1 \ar[dl, "g' \times \Id"'] \ar[dr, equals]  & & B \times_S \wh{B} \ar[dl, "{(\pr_0, \ev)}"'] \ar[dr, "\pr_1"] \\
A \times \A^1 & & B \times \A^1 & & \wh{B} 
\end{tikzcd}
\]
\item The functor $\wh{g}'_! \Conv$ is given by the diagram 
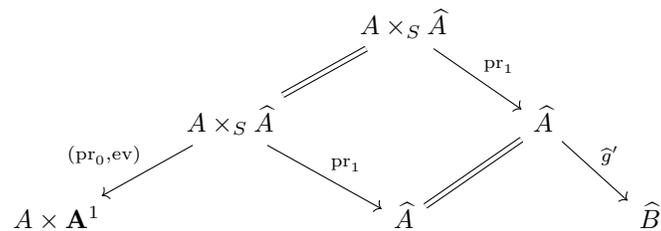
\[
\begin{tikzcd}
& & A \times_S \wh{A} \ar[dl, equals] \ar[dr, "\pr_1"] \\
& A \times_S \wh{A} \ar[dl, "{(\pr_0, \ev)}"'] \ar[dr, "\pr_1"] & & \wh{A} \ar[dl, equals] \ar[dr, "\wh{g}'"]  \\
A  \times \A^1& & \wh{A} & & \wh{B}
\end{tikzcd}
\]
\end{itemize}
To compare these, we consider the diagram
\[
\begin{tikzcd}
& & B \times_S \wh{A} \ar[dr , "\Id\times\wh{g}'"] \ar[dl, "g'\times\Id"']  \\
&A \times_S \wh{A}  \ar[ddl] \ar[ddrrr] & & B \times_S \wh{B}  \ar[ddr] \ar[ddlll] \\
\\
A  \times \A^1 & & & & \wh{B} 
\end{tikzcd}
\]
where the maps are as above. This is a butterfly thanks to the assumption that $g'$ is a closed embedding (so $\wh{g}'$ is smooth), and therefore gives a natural transformation 
\[
\wh{g}'_! \Conv  \xrightarrow{\bigstar} \Conv (g')^* \tw{-d(\wh g')}=\Conv (g')^* \tw{\delta}.
\]

\subsubsection{}  
Returning to the hexagon \eqref{eq: conv hexagon}, the top row of natural isomorphisms all come from the correspondence
\[
\begin{tikzcd} 
& B  \times_S \wh{D}  \ar[dl, "{(g', \ev \circ (\Id\times\wh{f}'))}"']  \ar[dr, "\pr_1"] \\
A \times \A^1 & &  \wh{D}  
\end{tikzcd}
\]
via \S \ref{sssec: correspondences}, as seen in the diagrams below: 
\begin{itemize}
\item $(\wh{f}')^* \Conv (g')^*$ is induced via \S \ref{sssec: correspondences} by the fibered product of correspondences 
\[
\begin{gathered}
\begin{adjustbox}{max width=\textwidth}
$\displaystyle
\begin{tikzcd}[ampersand replacement=\&]
\& \& \& B \times_S \wh{D} \ar[dl] \ar[dr, equals] \\
 \& \& B \times_S \wh{B} \ar[dl] \ar[dr, equals] \& \& B \times_S \wh{D} \ar[dl] \ar[dr]  \\
 \& B \times \A^1  \ar[dl, "g'\times\Id"']  \ar[dr, equals] \& \& B \times_S \wh{B} \ar[dl] \ar[dr] \& \& \wh{D} \ar[dl, "\wh{f}'"] \ar[dr, equals] \& \\
 A \times \A^1 \& \& B \times \A^1 \& \& \wh{B} \& \& \wh{D} 
\end{tikzcd}
$
\end{adjustbox}
\end{gathered}
\]
\item $\Conv f'_! (g')^*$ is induced via \S \ref{sssec: correspondences} by the fibered product of correspondences
\[
\begin{tikzcd}
& &  B \times_S \wh{D} \ar[dl] \ar[dr] \\
& B  \times \A^1 \ar[dl, "g'\times\Id"'] \ar[dr, "f'\times\Id"] & & D \times_S \wh{D} \ar[dl] \ar[dr] \\
A \times \A^1 & & D\times \A^1   & & \wh{D} 
\end{tikzcd}
\] 

\item $\Conv g^* f_!$ is induced via \S \ref{sssec: correspondences} by the fibered product of correspondences
\[
\begin{gathered}
\begin{adjustbox}{max width=\textwidth}
$\displaystyle
\begin{tikzcd}[ampersand replacement=\&]
\& \& \& B \times_S \wh{D} \ar[dl] \ar[dr] \\
\& \& B \times \A^1  \ar[dl, "g'\times\Id"'] \ar[dr, "f'\times\Id"]  \& \& D \times_S \wh{D} \ar[dl] \ar[ddrr] \\
 \& A \times \A^1  \ar[dl, equals] \ar[dr, "f\times\Id"] \&  \&  D \times \A^1  \ar[dl, "g\times\Id"']  \\
 A \times \A^1 \& \& C\times \A^1  \& \& \& \&  \wh{D}
\end{tikzcd}
$
\end{adjustbox}
\end{gathered}
\]
\end{itemize}

On the other hand, the bottom row of natural isomorphisms in the hexagon all come from the correspondence 
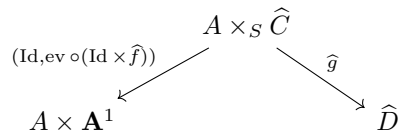
\[
\begin{tikzcd}
 & A \times_S \wh{C} \ar[dl, "{(\Id, \ev \circ (\Id\times\wh{f}))}"'] \ar[dr, "\wh{g}"] &\\
 A \times \A^1  & &  \wh{D} 
\end{tikzcd}
\]
via \S \ref{sssec: correspondences}, as seen in the diagrams below: 
\begin{itemize}
\item $(\wh{f}')^*(\wh{g}')_! \Conv$ is induced via \S \ref{sssec: correspondences} by the fibered product of correspondences
\[
\begin{gathered}
\begin{adjustbox}{max width=\textwidth}
$\displaystyle
\begin{tikzcd}[ampersand replacement=\&]
\& \& \& A \times_S \wh{C} \ar[dl] \ar[dr] \\
\& \& A \times_S \wh{A} \ar[dl, equals] \ar[dr] \& \& \wh{C} \ar[dl] \ar[dr] \\
\& A \times_S \wh{A} \ar[dl] \ar[dr] \& \& \wh{A} \ar[dl, equals] \ar[dr, "\wh{g}'"] \& \& \wh{D} \ar[dl, "\wh{f}'	"'] \ar[dr, equals] \\
A \times \A^1  \& \& \wh{A} \& \& \wh{B} \& \& \wh{D} 
\end{tikzcd}
$
\end{adjustbox}
\end{gathered}
\]
\item $\wh{g}_! (\wh{f})^* \Conv $ is induced via \S \ref{sssec: correspondences} by the fibered product of correspondences
\[
\begin{gathered}
\begin{adjustbox}{max width=\textwidth}
$\displaystyle
\begin{tikzcd}[ampersand replacement=\&]
\& \& \& A \times_S \wh{C} \ar[dl, equals] \ar[dr] \\
\& \& A \times_S \wh{C} \ar[dl] \ar[dr]  \&\& \wh{C} \ar[dl, equals] \ar[dr, equals]\\
\& A \times_S \wh{A} \ar[dl] \ar[dr] \& \& \wh{C} \ar[dl, "\wh{f}"]  \ar[dr, equals]  \& \& \wh{C} \ar[dl, equals] \ar[dr, "\wh{g}"] \\
A  \times \A^1\& \& \wh{A} \& \& \wh{C} \& \& \wh{D}
\end{tikzcd}
$
\end{adjustbox}
\end{gathered}
\]

\item $\wh{g}_! \Conv f_!$ is induced via \S \ref{sssec: correspondences} by the fibered product of correspondences
\[
\begin{gathered}
\begin{adjustbox}{max width=\textwidth}
$\displaystyle
\begin{tikzcd}[ampersand replacement=\&]
\& \& \& A \times_S \wh{C} \ar[dl] \ar[dr] \\
\& \& A \times_S \wh{C} \ar[dl] \ar[dr] \& \& C \times_S \wh{C}\ar[dl] \ar[dr] \\
\& A  \times \A^1\ar[dl, equals] \ar[dr, "f\times\Id"] \& \& C \times_S \wh{C} \ar[dl] \ar[dr] \& \& \wh{C} \ar[dl, equals] \ar[dr, "\wh{g}"] \\
A \times \A^1 \& \& C  \times \A^1 \& \& \wh{C} \& \& \wh{D} 
\end{tikzcd}
$
\end{adjustbox}
\end{gathered}
\]
\end{itemize}

\begin{prop}
Hexagon \eqref{eq: conv hexagon} commutes. 
\end{prop}

\begin{proof} By the preceding discussion, the natural transformation 
\[
\wh{g}_! \Conv f_! \rightarrow \Conv g^* f_! \tw{\delta}
\]
is the pre-composition of the natural transformation $\wh{g}_! \Conv \xrightarrow{\bigstar} \Conv g^* \tw{d(g)}$ coming from the butterfly
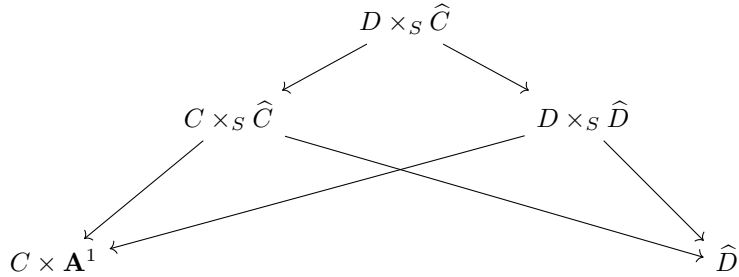
\[
\begin{tikzcd}
& & D \times_S \wh{C} \ar[dl] \ar[dr] \\
& C \times_S \wh{C} \ar[ddl] \ar[ddrrr] & & D \times_S \wh{D} \ar[ddr]  \ar[ddlll]\\ 
\\
C  \times \A^1  & & & & \wh{D}
\end{tikzcd}
\]
with $f_!$. Computing the pullback of this butterfly along $f$ gives the diagram
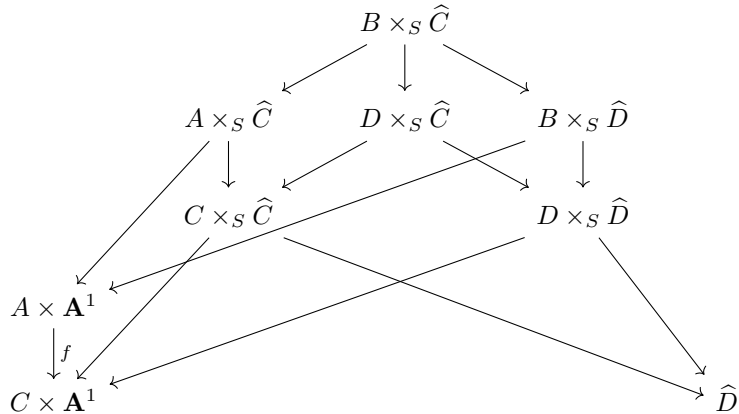
\[
\begin{tikzcd}
& & B \times_S \wh{C} \ar[dl] \ar[dr]  \ar[d] \\
& A \times_S \wh{C} \ar[ddl]   \ar[d]  & D \times_S \wh{C} \ar[dl] \ar[dr] & B \times_S  \wh{D} \ar[ddlll]  \ar[d] \\
& C \times_S \wh{C}  \ar[ddl] \ar[ddrrr] & & D \times_S \wh{D} \ar[ddr] \ar[ddlll]  \\
A \times \A^1  \ar[d, "f\times\Id"]  \\	
C  \times \A^1 & & & & \wh{D}
\end{tikzcd}
\] 
and Lemma \ref{lemma: butterfly pushforward} identifies the resulting natural transformation with the one from the upper butterfly 
\begin{equation}\label{eq: butterfly 1}
\begin{tikzcd}
& & B \times_S \wh{C} \ar[dl] \ar[dr] \\
& A \times_S \wh{C} \ar[ddl] \ar[ddrrr] & & B \times_S  \wh{D} \ar[ddr]  \ar[ddlll] \\ 
\\
A  \times \A^1  & & & & \wh{D}
\end{tikzcd}
\end{equation}

On the other hand, the natural transformation 
\[
(\wh{f'})^* \wh{g'}_! \Conv \rightarrow (\wh{f'})^* \Conv (g')^* \tw{\delta} 
\]
is the composition of $(\wh{f'})^*$ with the natural transformation $\wh{g'}_! \Conv \xrightarrow{\bigstar} \Conv (g')^* \tw{\delta} $ coming from the butterfly 
\[
\begin{tikzcd}
& & B \times_S \wh{A} \ar[dl, "g'\times\Id"'] \ar[dr, "\Id\times\wh{g'}"] \\
& A \times_S \wh{A}  \ar[ddl] \ar[ddrrr] & & B \times_S \wh{B} \ar[ddr] \ar[ddlll] \\ 
\\
A  \times \A^1  & & & & \wh{B}
\end{tikzcd}
\]
Computing the pullback of this butterfly along $\wh{f'}$ gives the diagram
\[
\begin{tikzcd}
& & B \times_S \wh{C} \ar[dl] \ar[dr]  \ar[d] \\
& A \times_S \wh{C}   \ar[d]  \ar[ddrrr] & B \times_S \wh{A} \ar[dl] \ar[dr] & B \times_S \wh{D}  \ar[ddr]   \ar[d] \\
& A \times_S \wh{A}  \ar[ddl] \ar[ddrrr]  & & B \times_S \wh{B} \ar[ddr] \ar[ddlll]  \\
& & & & \wh{D} \ar[d, "\wh{f'}"] \\	
A  \times \A^1  & & & & \wh{B}
\end{tikzcd}
\]
and Lemma \ref{lemma: pullback butterfly} identifies the resulting natural transformation with the one from the upper butterfly 
\[
\begin{tikzcd}
& & B \times_S \wh{C} \ar[dl] \ar[dr] \\
& A \times_S \wh{C} \ar[ddl]  \ar[ddrrr] & & B \times_S  \wh{D} \ar[ddr]  \ar[ddlll] \\ 
\\
A \times \A^1 & & & & \wh{D}
\end{tikzcd}
\]
which is visibly the same butterfly as \eqref{eq: butterfly 1}.
\end{proof}

\section{The stacky contraction principle}
At the behest of the referee, we will reprove the ``stacky contraction principle'' of Drinfeld--Gaitsgory  from \cite[Appendix C]{DG15} in our current context. 

Throughout this subsection, let $\cY$ be a locally finite type derived Artin stack over a field $k$, with an action of $\A^1$ viewed as a multiplicative monoid. (The derived structure will play no role in this discussion.) 

\subsection{The contracting stack} We may view $\{0,1\}$ as a monoid under multiplication. By restriction along $\{0,1\} \subset \A^1$, $\cY$ acquires an action of $\{0,1\}$. 

Let $M$ be a finite module for the monoid $\{0,1\}$. We denote by $\cY^M$ the stack of $\{0,1\}$-equivariant maps $M \rightarrow \cY$. Concretely, for a $k$-algebra $R$, $\cY(R)$ carries an action of the monoid $\{0,1\}$ and $\cY^M(R)$ is the groupoid of $\{0,1\}$-equivariant functors $M \rightarrow \cY(R)$.

\begin{example}\label{ex:regular-monoid-isom}
Viewing $\{0,1\}$ as a module over itself in the natural way, the natural forgetful map $\cY^{\{0,1\}} \rightarrow \cY$ is an isomorphism. 
\end{example}

\begin{defn}Viewing $\{0\}$ as a module for the monoid $\{0,1\}$ under multiplication, the \emph{contracting stack} of $\cY$ is $\cY^0$. Restricting along the $\{0,1\}$-equivariant maps $\{0,1\} \rightarrow \{0\} \inj \{0,1\}$, and invoking Example \ref{ex:regular-monoid-isom}, we obtain maps
\[
\cY \xleftarrow{i} \cY^0 \xleftarrow{\pi} \cY
\]
such that $\pi \circ i = \Id_{\cY^0}$.
\end{defn}

If $\cY$ is represented by a derived scheme, then $i$ is the closed embedding of the ``contracting locus''. In general, however, $i$ need not be a closed embedding, so we refer to $\cY^0$ as the ``contracting stack''. 

\begin{example}
Let $\cY = E$ be a derived vector bundle over a base $S$. Then $\A^1$ acts on $E$ by scaling, and $\cY^0$ identifies with the zero-section of $E$, which is a closed embedding if and only if $\cE$ comes from a perfect complex $\cE \in \Perf(S)$ with tor-amplitude in $[0, \infty)$.
\end{example} 

\subsection{Formulation of the contraction principle}\label{sssec:contraction-principle}

By restriction along $\G_m \subset \A^1$, we have a $\G_m$-action on $\cY$. For $q \co \cY \rightarrow \cY/\G_m$, define the \emph{monodromic derived category} 
\[
D(\cY  \rightdash \G_m) \subset D(\cY)
\]
to be the full subcategory of $D(\cY)$ generated by objects pulled back from  $D(\cY/\G_m)$ via $q^*$. Note that since the map $q$ is smooth, the definition is equivalent if we use $q^!$ instead of $q^*$. 

\begin{example}
If $\G_m$ acts trivially on $\cY$, then $D(\cY \rightdash \G_m) = D(\cY)$. In particular, the induced action of $\A^1$ on a contracting stack $\cY^0$ is always trivial, so we always have $D(\cY^0 \rightdash \G_m) = D(\cY^0)$. 
\end{example}

We have natural transformations in $\cK_0 \in D(\cY)$, $\cK_1 \in D(\cY^0)$,
\begin{align}\label{eq:contracting-adjunctions}
\Hom_{\cY^0}(i^* \cK_0, \cK_1) \cong \Hom_{\cY}(\cK_0, i_* \cK_1) & \rightarrow \Hom_{\cY^0}(\pi_* \cK_0, \pi_* i_* \cK_1) \nonumber \\
& \cong \Hom_{\cY^0}(\pi_* \cK_0, \cK_1).
\end{align}
From the composite, we obtain a natural transformation $\pi_* \rightarrow i^* $. 

\begin{thm}\label{thm:contraction-principle} With the assumptions and definitions above, the natural transformations are isomorphisms:
\[
\pi_* \cong i^*  \co D(\cY \rightdash \G_m) \rightarrow D(\cY^0)
\]
and
\[
i^! \cong \pi_!  \co D(\cY \rightdash \G_m) \rightarrow D(\cY^0)
\]
\end{thm}

In the D-module setting, this is a result of Drinfeld--Gaitsgory \cite[Theorem C.5.3]{DG15}. The proof is essentially a close copy of Drinfeld--Gaitsgory, but for completeness we spell it out (after prodding from the referee). We will focus on proving the first statement; the argument for the second is of a dual nature. 

\subsection{A special case}
We will reduce Theorem \ref{thm:contraction-principle} to the special case where the $\A^1$ monoid action factors through an action of the \emph{monoidal stack} $\A^1/\G_m$. 

Let $\ol \cY$ be a derived Artin stack equipped with an action of the monoidal stack $\A^1/\G_m$. Thanks to the monoid homomorphism $\A^1 \rightarrow \A^1/\G_m$, this equips $\ol \cY$ with a natural action of $\A^1$, so that the contracting stack $\ol \cY^0$ is defined and we have natural maps 
\[
\ol\cY \xleftarrow{i} \ol\cY^0 \xleftarrow{\pi} \ol\cY
\]
such that $\pi \circ i = \Id_{\ol \cY^0}$. As in \S \ref{sssec:contraction-principle}, we obtain natural transformations $\pi_* \rightarrow i^*$ and $i^! \rightarrow \pi_!$. Then Theorem \ref{thm:contraction-principle} specializes to the following assertion. 

\begin{thm}\label{thm:different-contraction-principle}
Let $\ol \cY$ be a derived Artin stack equipped with an action of the monoidal stack $\A^1/\G_m$. Then the natural transformations above are isomorphisms
\[
\pi_* \cong i^* \co D(\ol \cY) \rightarrow D(\ol \cY^0)
\]
and
\[
i^! \cong \pi_! \co D(\ol \cY) \rightarrow D(\ol \cY^0)
\]
\end{thm}

\begin{example}\label{ex:B-base-case}
We verify Theorem \ref{thm:different-contraction-principle} for $\ol \cY = \A^1/\G_m$ with the tautological $\A^1/\G_m$-action, which will be used later as a ``base case'' for the proof in general. In this case, $\ol \cY^0 = 0/\G_m$ and we want to show that the natural transformation $\pi_* \rightarrow i^*$ is an isomorphism. Let $j \co \G_m/\G_m \rightarrow \A^1/\G_m$ be the open embedding of the complement of $\ol \cY^0$. Writing any $\cK \in D(\ol \cY)$ as an extension
\[
j_! j^* \cK \rightarrow \cK \rightarrow i_* i^* \cK
\]
we reduce to showing the statement for $\cK$ of the form $i_* \cK'$ or $\cK = j_! \cK'$. In the first case, the assertion is clear from the fact that $\pi \circ i = \Id$. In the second case, $i^* j_! \cK'=0$ so it suffices to show that $\pi_* j_! \cK' = 0$. To prove this vanishing for an arbitrary $\cK' \in D(\G_m/\G_m)=D(\pt)$, pull back along the smooth atlas $\pt \rightarrow 0/\G_m$. The base change is $R\Gamma(\A^1, j_!\Q_\ell)\otimes \cK'$, by the projection formula. Thus it suffices to show $R\Gamma(\A^1,j_!\Q_\ell)=0$, which follows from the localization triangle and the fact that restriction to the origin induces an isomorphism $R\Gamma(\A^1,\Q_\ell)\xrightarrow{\sim} R\Gamma(\pt,\Q_\ell)$.
\end{example}

\subsection[Reductions of the contraction principle]{Reductions of Theorem \ref{thm:contraction-principle} to Theorem \ref{thm:different-contraction-principle}}\label{ssec:B-reduction}

We will prove Theorem \ref{thm:contraction-principle} assuming Theorem \ref{thm:different-contraction-principle}. Examining \eqref{eq:contracting-adjunctions}, it suffices to show that the natural transformation 
\begin{equation}\label{eq:contracting-upstairs}
\Hom_{\cY}(\cK_0, i_* \cK_1) \rightarrow \Hom_{\cY^0}(\pi_* \cK_0, \pi_* i_* \cK_1) \cong \Hom_{\cY^0}(\pi_* \cK_0, \cK_1)
\end{equation}
is an isomorphism for all $\cK_0 \in D(\cY \rightdash \G_m )$ and $\cK_1 \in D(\cY^0 \rightdash \G_m ) = D(\cY^0)$. Consider the Cartesian square
\begin{equation}\label{eq:Y-diagram}
\begin{tikzcd}
\cY^0 \ar[r, "i"] \ar[d, "q^0"] & \cY \ar[d, "q"] \ar[l, bend right, "\pi"'] \\
\cY^0/\G_m \ar[r, "\ol i"] & \cY/\G_m  \ar[l, bend right, "\ol\pi"']
\end{tikzcd}
\end{equation}
Let $\ol \cY := \cY/\G_m$; note that its contracting stack is $\ol \cY^0 \cong \cY^0/\G_m$. 

We will recast \eqref{eq:contracting-upstairs} in terms of the lower row of \eqref{eq:Y-diagram}. Since $D(\cY\rightdash \G_m)$ is generated by the $q^!$-pullbacks and both sides of \eqref{eq:contracting-upstairs} are exact in $\cK_0$, it suffices to check the claim on such generators. Thus we may assume $\cK_0=q^!\ol\cK_0$ with $\ol\cK_0\in D(\ol\cY)$. Since the action on $\cY^0$ is trivial, we may also write $\cK_1=(q^0)^!\ol\cK_1$ for some $\ol\cK_1\in D(\ol\cY^0)$.
Then we have
\begin{align*}
\Hom_{\cY}(\cK_0, i_* \cK_1) & \cong 
\Hom_{\cY}(q^! \ol \cK_0, i_* (q^0)^! \ol \cK_1)  \\ 
& \cong \Hom_{\cY}(q^! \ol \cK_0, q^! \ol i_*  \ol \cK_1) \\
& \cong \Hom_{\cY}(q^* \ol \cK_0 \tw{1}, q^! \ol i_*  \ol \cK_1)   \\
& \cong \Hom_{\ol \cY}(\ol \cK_0 , q_*q^! \ol i_*  \ol \cK_1 \tw{-1} )
\end{align*}
where we used that $q^! \cong q^* \tw{1}$ because $q$ is smooth (of relative dimension $1$). Let $p\co \pt \rightarrow \pt/\G_m$ be the universal $\G_m$-torsor and put $\sA=p_*p^!\Ql\tw{-1}$. If $s\co \ol\cY\rightarrow \pt/\G_m$ and $s_0\co \ol\cY^0\rightarrow \pt/\G_m$ are the structural maps, set $\cA_{\ol\cY}:=s^*\sA$ and $\cA_0:=s_0^* \sA$. The maps $\ol i$ and $\ol\pi$ are over $\pt/\G_m$, so $\cA_{\ol\cY}\cong\ol\pi^*\cA_0$ and $\ol i^*\cA_{\ol\cY}\cong\cA_0$. By the projection formula, $q_*q^!(-)\tw{-1}$ is tensoring with $\cA_{\ol\cY}$. Thus we may rewrite
\begin{equation}\label{eq:downstairs-adjunctions-1}
\Hom_{\cY}(\cK_0, i_* \cK_1)   = \Hom_{\ol \cY}(\ol \cK_0 , \ol\pi^*\cA_0 \otimes  \ol i_*  \ol \cK_1 )
\end{equation}
and similarly
\begin{equation}\label{eq:downstairs-adjunctions-2}
\begin{aligned}
\Hom_{\cY^0}(\pi_* \cK_0, \cK_1)
&= \Hom_{\ol \cY^0}(\ol \pi_* \ol \cK_0 , \cA_0 \otimes \ol \cK_1 ).
\end{aligned}
\end{equation}
Theorem \ref{thm:different-contraction-principle} asserts that $\ol \pi_* \cong \ol i^*$. In particular, the map $\ol \pi_* (\ol \pi^*(A) \otimes \ol\cK) \leftarrow A \otimes \ol \pi_* \ol \cK$ is an isomorphism for any $A \in D(\ol \cY^0)$. Hence for any $A \in D(\ol \cY^0)$, we have a priori a natural transformation
\begin{align}\label{eq:downstairs-hom-with-A}
\Hom_{\ol \cY}(\ol \cK_0, \ol\pi^*A \otimes \ol i_* \ol \cK_1) &\rightarrow \Hom_{\ol \cY^0}(\ol \pi_* \ol \cK_0 , \ol \pi_* (\ol\pi^*A \otimes \ol i_* \ol \cK_1)) \nonumber \\
& \cong \Hom_{\ol \cY^0}(\ol \pi_* \ol \cK_0 , A \otimes  \ol \cK_1),
\end{align}
which for $A=\cA_0$ agrees with the map \eqref{eq:contracting-upstairs} under our identifications \eqref{eq:downstairs-adjunctions-1} and \eqref{eq:downstairs-adjunctions-2}. Therefore, it suffices to show that \eqref{eq:downstairs-hom-with-A} is an isomorphism for $A=\cA_0$. Theorem \ref{thm:different-contraction-principle} tells us that it is an isomorphism for $A = \Ql$, but $\cA_0$ has a two-step filtration whose graded pieces are shifts and twists of $\Ql$, because $\sA$ computes the Borel--Moore homology of $\G_m$. Since both sides of \eqref{eq:downstairs-hom-with-A} are exact in $A$, the isomorphism for $A=\Ql$ implies the isomorphism for $A=\cA_0$. \qed

\subsection{Generalities on idempotent functors}\label{ssec:idempotent-functors}

Thanks to \S \ref{ssec:B-reduction}, it now suffices to prove Theorem \ref{thm:different-contraction-principle}. This will essentially follow from abstract nonsense. We recall some relevant definitions from \cite[\S C.1]{DG15}.

Given a monoidal category $(\cC, \otimes)$, there is a notion of \emph{algebra} $A \in \cC$. We say that $A$ is \emph{idempotent} if the multiplication map $A \otimes A \rightarrow A$ is an isomorphism. An idempotent algebra can be thought of equivalently as a monoidal functor $F  \co \{0,1\} \rightarrow \cC$: $A = F(0)$ is an idempotent algebra in the above sense, while $F(1) = 1_{\cC}$ by definition. 

\begin{defn}
Given any category $\cC$, we can form a monoidal category of endofunctors $\End(\cC)$. An \emph{idempotent endofunctor} of $\cC$ is an idempotent algebra in $\End(\cC)$.
\end{defn}

\begin{example}\label{ex:idempotent-endofunctor}
Let $\cD$ and $\cC$ be categories and 
\[
\cD \xrightarrow{i_* } \cC \xrightarrow{\pi_*} \cD
\]
a sequence of functors equipped with an isomorphism $\epsilon \co \pi_* \circ i_* \xrightarrow{\sim} \Id_{\cD}$. Then $\ul{0}_* := i_* \circ \pi_*$ is an idempotent endofunctor of $\cC$: $\epsilon$ supplies an isomorphism 
\[
\ul{0}_* \circ \ul{0}_* \cong (i_* \pi_* ) (i_* \pi_*) \cong i_* (\pi_* i_* ) \pi_* \xrightarrow{\epsilon} i_* \pi_* = \ul{0}_*.
\]
\end{example}

In fact, every idempotent endofunctor arises from the construction of Example \ref{ex:idempotent-endofunctor}. To see this, let $\cC$ be a category equipped with an idempotent endofunctor $\ul{0}_* \in \End(\cC)$. As discussed above, the data of $\ul{0}_*$ is equivalent to that of a monoidal functor $\{0,1\} \rightarrow \End(\cC)$, i.e., an action of the monoid $\{0,1\}$ on $\cC$. Let $\cC^0$ be the category of $\{0,1\}$-equivariant functors $\{0\} \rightarrow \cC$. Then restriction along the $\{0,1\}$-equivariant maps $\{0\} \rightarrow \{0,1\} \rightarrow \{0\}$ induces a diagram
\[
\cC^0 \xleftarrow{\pi_*} \cC \xleftarrow{i_*} \cC^0
\]
whose composition is $\Id_{\cC^0}$. 

\begin{lemma}\label{lem:idempotent-endofunctor-unit}
Let $\cC$ be a category and $\ul{0}_*$ be an idempotent endofunctor of $\cC$. Factorize $\ul{0}_* = i_* \pi_*$ with $\epsilon \co \pi_* i_* \xrightarrow{\sim} \Id_{\cC^0}$ as above. Regarding $\ul{0}_*$ as an algebra in $\End(\cC)$, promoting $\ul{0}_*$ to a unital algebra is equivalent to promoting $(\pi_*, i_*)$ to an adjunction with counit $\epsilon$. 
\end{lemma}

\begin{proof}
To promote $(\pi_*, i_*)$ to an adjunction with counit $\epsilon$ is to give a unit $\eta \co \Id_{\cC} \rightarrow i_* \pi_*$ and commutative triangles
\begin{equation}\label{eq:adjunction-triangle}
\begin{tikzcd}
\pi_* \ar[r, "\pi_* \eta"] \ar[dr, "\Id"'] & \pi_* (i_* \pi_*) = (\pi_* i_*) \pi_* \ar[d, "\epsilon \pi_*"] \\
& \pi_* 
\end{tikzcd} \quad \text{and} \quad 
\begin{tikzcd}
i_* \ar[r, "\eta i_* "] \ar[dr, "\Id"'] &  (i_* \pi_*)i_*  = i_* (\pi_* i_*)  \ar[d, "i_* \epsilon"] \\
& i_*
\end{tikzcd}
\end{equation}
The unit of the adjunction is taken to be the unit of $\ul{0}_*$. To check that it satisfies the axioms of a unit, we must produce commutative triangles
\begin{equation}\label{eq:unital-triangle}
\begin{tikzcd}
\ul{0}_* \ar[r, "\ul{0}_* \eta"]  \ar[dr, "\Id"'] & \ul{0}_* (\ul{0}_*)  \ar[d, "\mu"] \\
& \ul{0}_* 
\end{tikzcd}\quad \text{and} \quad 
\begin{tikzcd}
\ul{0}_* \ar[r, "\eta \ul{0}_*"]  \ar[dr, "\Id"] & (\ul{0}_*) \ul{0}_*  \ar[d, "\mu"] \\
& \ul{0}_* 
\end{tikzcd}
\end{equation}
For the first triangle, we compose the first triangle of \eqref{eq:adjunction-triangle} with $i_*$ to obtain  
\begin{equation}\label{eq:unit-triangle-left-expand}
\begin{tikzcd}
i_* \pi_* \ar[r, "i_* \pi_* \eta"] \ar[dr, "\Id"'] & i_* \pi_* (i_* \pi_*) = i_* (\pi_* i_*) \pi_* \ar[d, "i_* \epsilon \pi_* "] \\
& i_* \pi_* 
\end{tikzcd}
\end{equation}
as desired. The second triangle is handled similarly, using the second triangle of \eqref{eq:adjunction-triangle}. 

To go in the other direction, we start with commutative triangles \eqref{eq:unital-triangle}. Composing \eqref{eq:unit-triangle-left-expand} with $\pi_*$ and applying the counit $\epsilon$ (and using that it is an isomorphism) gives the left triangle of \eqref{eq:adjunction-triangle}, and a similar procedure works for the other commutative triangle. 
\end{proof}

In the setting of Theorem \ref{thm:different-contraction-principle}, let $\ul{0} \co \ol\cY \rightarrow \ol\cY$ be the composition
\[
\ol\cY \xrightarrow{\pi} \ol\cY^0 \xrightarrow{i} \ol\cY.
\]
This induces the setup of Example \ref{ex:B-base-case}, so $\ul{0}_*$ is an idempotent endofunctor of $D(\ol\cY)$. Since $i^*$ is left adjoint to $i_*$, we can interpret Theorem \ref{thm:different-contraction-principle} as saying that $\pi_*$ is left adjoint to $i_*$. By Lemma \ref{lem:idempotent-endofunctor-unit},  this assertion can in turn be reformulated as saying that $\ul{0}_*$ can be promoted to a unital algebra in $\End(D(\ol\cY))$. This proves the isomorphism $\pi_*\cong i^*$; the isomorphism $i^!\cong\pi_!$ follows by Verdier duality, since Verdier duality interchanges $*$- and $!$-functors and dualizes the natural transformation $\pi_*\rightarrow i^*$. In these terms, Example \ref{ex:B-base-case} gives such a promotion for the special case $\ol\cY = \A^1/\G_m$. We shall pass from the ``base case'' of 
Example \ref{ex:B-base-case} to the general case using a general mechanism. 

\subsection{Generalities on monoidal stacks}

Recall that a monoidal stack is a monoid in the category of stacks. In particular, a monoidal stack $\cY$ is equipped with a multiplication $m \co \cY \times \cY \rightarrow \cY$. 

\begin{example}
A monoid scheme (such as a group scheme) is a monoidal stack. Also $\A^1/\G_m$ is a monoidal stack; it parametrizes line bundles equipped with a cosection, and the monoid operation is tensor product. 
\end{example}

Given a monoidal derived Artin stack $\cY$, the derived category $D(\cY)$ acquires a monoidal structure $\star$ under convolution: 
\[
\cK_0 \star \cK_1 := m_* (\cK_0 \boxtimes \cK_1). 
\]
We denote by $1_{\cY}$ the unit of $(D(\cY), \star)$. If $f \co \cY \rightarrow \cY'$ is a monoidal homomorphism of such monoidal stacks, then $f_* \co D(\cY) \rightarrow D(\cY')$ promotes to a monoidal functor for $\star$. In particular, we have $f_*(1_{\cY}) = 1_{\cY'}$. If $f$ is only a semigroup homomorphism (i.e., compatible with the multiplication but not the unit), then $f_*(1_{\cY})$ is an idempotent algebra in $(D(\cY'), \star)$. 

Given $y \in \cY(k)$, let $\delta_y$ be the direct image of $\Q_\ell \in D(\pt)$ under the map $y \rightarrow \cY$. 

\begin{example}
For $y = 1$ the monoidal unit, the map $1 \rightarrow \cY$ is a monoidal homomorphism, so $\delta_1 = 1_{\cY}$ is the monoidal unit of $(D(\cY), \star)$. 
\end{example}

\begin{example}\label{ex:delta_0-algebra}
If $\cY$ is a monoidal stack and the map denoted $0\rightarrow\cY$ is a semigroup homomorphism from an idempotent one-object monoidal stack, equivalently if the multiplication composite $0\times 0\rightarrow \cY\times\cY\xrightarrow{m}\cY$ agrees with $0\rightarrow\cY$, then its direct image $\delta_0$ is an idempotent algebra in $(D(\cY),\star)$.
\end{example}

\subsection[Proof of the quotient contraction principle]{Proof of Theorem \ref{thm:different-contraction-principle}}
We can now complete the proof of Theorem \ref{thm:different-contraction-principle}.  With $\ol\cY$ as in the setup of Theorem \ref{thm:different-contraction-principle}, we have a monoidal action of $\A^1/\G_m$ on $\ol\cY$. This induces an action of $D(\A^1/\G_m)$ on $D(\ol\cY)$: letting 
\begin{equation}\label{eq:B-action-map}
a \co \A^1/\G_m \times \ol\cY \rightarrow \ol\cY
\end{equation}
denote the action map, the action functor
\[
D(\A^1/\G_m) \times D(\ol\cY) \rightarrow D(\ol\cY)
\]
is given by 
\[
\cK_0, \cK_1 \mapsto a_*(\cK_0 \boxtimes \cK_1).
\]

The above action of $D(\A^1/\G_m)$ on $D(\ol\cY)$
can be formulated equivalently as a monoidal functor 
\begin{equation}\label{eq:B-action-homomorphism}
D(\A^1/\G_m) \rightarrow \End(D(\ol\cY)).
\end{equation}
Recall that using Lemma \ref{lem:idempotent-endofunctor-unit}, we reformulated Theorem \ref{thm:different-contraction-principle} as saying that the idempotent algebra $\ul{0}_* \in \End(D(\ol\cY))$ can be promoted to a unital algebra. 

\begin{lemma}\label{lem:idempotent-algebra-as-delta}
The idempotent algebra $\ul{0}_* \in \End(D(\ol\cY))$ is the image of $\delta_0$ under \eqref{eq:B-action-homomorphism}. 
\end{lemma}

\begin{proof}
By definition, the image $\delta_0$ under \eqref{eq:B-action-homomorphism} is the endomorphism of $D(\ol\cY)$ induced by pushforward along endomorphism of $\ol\cY$, which is in turn obtained by pulling back the action map \eqref{eq:B-action-map} along $\{0\} \rightarrow \A^1/\G_m$. The resulting map 
\[
\{0\} \times \ol\cY \rightarrow \ol\cY
\]
is precisely the endomorphism $\ul{0} \co \ol\cY \xrightarrow{\pi} \ol\cY^0 \xrightarrow{i} \ol\cY$, as desired. 
\end{proof}

In the base case $\ol\cY = \A^1/\G_m$, we saw at the end of \S \ref{ssec:idempotent-functors} that Example \ref{ex:B-base-case} could be interpreted as promoting $\ul{0}_* \in \End(D(\A^1/\G_m))$ to a unital algebra. By Lemma \ref{lem:idempotent-algebra-as-delta}, we can interpret this as a unital algebra structure on $\delta_0$ regarded as an idempotent algebra in $(D(\A^1/\G_m), \star)$ via Example \ref{ex:delta_0-algebra}. This is then carried by the monoidal functor \eqref{eq:B-action-homomorphism} to a unital algebra structure on $\ul{0}_* \in \End(D(\ol \cY))$, completing the proof. \qed 

 \qed

\bibliographystyle{amsalpha}
\bibliography{Bibliography}

\end{document}